\documentclass{article}
\usepackage{graphicx} 
\usepackage[margin=1in]{geometry}
\usepackage{amsmath}
\usepackage{amsthm}
\usepackage{amssymb}
\usepackage{bbold}
\usepackage{mathtools}
\usepackage{enumitem}
\usepackage{caption}
\usepackage{subcaption}
\usepackage{comment}
\usepackage{centernot}
\usepackage{xcolor}
\usepackage{enumitem}
\usepackage{hyperref}
\usepackage{booktabs}
\usepackage[page]{appendix}
\newtheorem{theorem}{Theorem}[section]
\newtheorem{proposition}{Proposition}[section]
\newtheorem{lemma}{Lemma}[section]

\newtheorem{definition}{Definition}[section]

\newtheorem{corollary}{Corollary}[section]

\DeclareRobustCommand{\bigO}{%
  \text{\usefont{OMS}{cmsy}{m}{n}O}%
}
\DeclareMathOperator{\diag}{Diag}
\DeclareMathOperator{\diver}{div}
\DeclareMathOperator{\backwardscoloneqq}{=:}

\DeclareMathOperator{\tr}{tr}

\DeclareMathOperator*{\argmin}{arg\,min}
\DeclareMathOperator{\dif}{\mathrm{d}}

\title{Ordered Diffusion Kernels}

\begin{document}

\maketitle

\begin{center}
Jack H. Soulsby\textsuperscript{1*}, Andreas C.S. Jørgensen\textsuperscript{1,2}, Atiyo Ghosh\textsuperscript{1},
Vahid Shahrezaei\textsuperscript{1}\\

\bigskip
\textbf{$^1$} Department of Mathematics, Imperial College London, South Kensington, London SW7 2BX, UK\\
\textbf{$^2$} Laboratoire des Sciences du Climat et de l’Environnement, CEA CNRS UVSQ, Orme des Merisiers, 91191 Gif-sur-Yvette, France \\
$^*$ Correspondence to be addressed to jack.soulsby18@imperial.ac.uk\\

\end{center}

\begin{abstract}
    We introduce Ordered Diffusion Kernels (ODKs), a novel class of local kernels that can approximate the infinitesimal generator of an arbitrary Itô Stochastic Differential Equation (SDE).  ODKs are designed to be applied to data sampled from dynamical systems where little dynamical information is available a priori.  The Laplacian of classical diffusion kernels approximates the Laplace–Beltrami operator on the underlying manifold; adjusting the normalisation introduces an advection term that depends on the sampling density; recently, TMDmap generalised this normalisation to target an arbitrary measure, but at the cost of coupling advection to diffusion.  More general local kernels can learn arbitrary second-order elliptic operators but are formulated in terms of known velocity fields --- making the first step of any analysis a potentially ill-posed inference problem.  To formulate ODK, we first relax the problem of potential estimation to the more tractable task of inferring an ordering of the data, which we represent through an ordering function.  We prove ODK’s Laplacian converges to the infinitesimal generator of a gradient-flow SDE with state-dependent isotropic diffusion, without coupling advection and diffusion.  We provide various extensions of ODK to: arbitrary drifts via local ordering functions; anisotropic diffusions via a Strang splitting scheme; multiple ordering functions; and self-tuning bandwidths.  In addition, we introduce two loss functions which exploit the structure of ODKs to solve a non-parametric inference problem.  We validate this framework on synthetic data from deterministic and stochastic systems, demonstrating accurate recovery of operators, velocity fields, extrinsic curvature, and spatially dependent drift and diffusion coefficients. 

    \vspace{0.5em}
\end{abstract}

\section{Introduction}

Over the last 20 years, a large body of research has been focused on applying kernels to complex, high-dimensional data to extract geometric and dynamical features. This endeavour has been successful, leading to the development of many techniques now standard to a data scientist's tool chest, tackling problems such as dimensionality reduction \cite{scholkopf1998nonlinear, coifman2006diffusion, belkin2003laplacian, tenenbaum2000global}, dimensionality estimation \cite{hein2005intrinsic, granata2016accurate, hino2017ider, facco2017estimating, allegra2020data, di2026scale}, density estimation \cite{rosenblatt1956, sheather1991reliable, berry2016variable, berry2017density}, and clustering \cite{shi2000normalized, ng2001spectral, dhillon2004kernel, schiebinger2015geometry, couillet2016kernel, vankadara2020optimality, pillaud2023kernelized}, among others.

Much mathematical machinery, such as the tools of calculus, require the algebraic structure of a vector (or linear) space, in order to be constructed; however, many geometries of interest cannot be represented as linear spaces, so the theory of manifolds was developed as a generalisation of linear spaces to those that are only locally linear and of homogeneous dimension. The relaxation of global linearity to local linearity allows the tools of linear spaces to be brought over to the more complex class of manifolds, where they become the union of their local applications --- making sure this union is well-defined is the basis of differential geometry, and is encoded in the manifold's atlas (set of charts). The foundational assumption underlying most kernel methods is that the numerical data we observe clusters on or around an unobserved $d$-dimensional submanifold $\mathcal{M}\subset\mathbb{R}^n$ with $d\ll n$ and $n \gg 1$ being the dimension of the data --- this is often called the \textit{Manifold Hypothesis}. Kernels leverage the manifold structure of data analogously to how manifolds themselves are defined. Although the data we are provided is high-dimensional and globally complex, we can use fast-decaying kernels to define affinities between pairs of points which will encode some aspect of the local neighbourhood for every data point. These local affinities are easily interpretable and vary smoothly on the underlying manifold. The global geometric complexity of the data is then encoded through the aggregation of these locally simple descriptions --- in much the same way that a manifold is built by smoothly glueing together locally linear patches. The result is a smooth and interpretable representation of the data that respects its intrinsic geometry. The central question this raises is whether such a representation, constructed entirely from local relations, faithfully captures the global structure of the data: in short, is the union of our intuitions intuitive?

The pursuit of formalising this question --- and in particular, understanding the theoretical foundations of the Laplacian Eigenmaps algorithm \cite{belkin2003laplacian} --- led researchers to establish that the graph Laplacian of an appropriately normalised symmetric kernel asymptotically approaches the Laplace-Beltrami operator on the underlying manifold (with metric inherited from ambient space) \cite{belkin2008towards}. This first convergence result was pointwise, but subsequent work extended this to spectral and uniform convergence \cite{belkin2008towards, belkin2006convergence, gine2006empirical}.  Coifman \& Lafon extended these results to arbitrary compactly supported sampling distributions, $q(x)$, and through this introduced the class of $\alpha$-normalisations that add a density dependent drift term into the limiting operator, $L_{\varepsilon,\alpha}f(x)\xrightarrow{\varepsilon\to0}2(1-\alpha)\nabla (\log q(x))\cdot\nabla f(x)+\Delta f(x)$, where it becomes clear the $\alpha\in[0,1]$ controls the magnitude of this additional term.  Operators of this form can be viewed as the infinitesimal generators of It{\^o} Stochastic Differential Equations (SDEs): first-order terms become drifts (or advections), second-order terms become diffusions, and the sampling distribution becomes the transient distribution of the process.  If we make such a connection for the $\alpha$-normalised operators, we find the coupling between the drift and sampling density is highly restrictive, limiting us to only considering potential systems in steady state, with potential satisfying $U(x)=-\log q(x)$.  Alternative normalisations have been considered, for example by Banisch et al. \cite{banisch2020diffusion}, whose Target Measure Diffusion Maps (TMDmap) kernel is capable of approximating the drift term of any gradient flow; however, TMDmap is limited to only considering unit diffusions, as otherwise there will be an inextricable coupling between the drift and diffusion.  The most general result of this kind came from Berry \& Sauer \cite{berry2016local}, who define the very general class of local kernels which satisfy the conditions of being fast-decaying (local) but need not be symmetric. In their work they show that appropriately constructed kernels can approximate operators of the form $\mathcal{L}f(x)=\mu(x)\cdot\nabla f(x)+D^{ij}(x)\nabla_i\nabla_j f$ where $\mu\in \mathcal{T}\mathcal{M}$ and $D\in \mathcal{T}^2\mathcal{M}$ --- which is the general form of a symmetric second-order elliptic operator and is the infinitesimal generator of an arbitrary It{\^o} SDE of the form $dx_t=\mu(x_t)\mathrm{d}t+\sigma(x_t)dW_t$ where $D(x)=\frac{1}{2}\sigma(x)\sigma(x)^T$.

To apply local kernels as conceptualised in \cite{berry2016local}, one is required to know the drift vector field $\mu(x)$ and diffusion tensor field $D(x)$ beforehand. In practice, this limits the applicability of local kernels to complex, novel data where these quantities are not known a priori and are otherwise highly challenging to infer.  In response to these constraints, we introduce \textit{Ordered Diffusion Kernels} (ODKs) which, in their standard form, approximate operators of the form $\mathcal{L}f(x)=\frac{\tau(x)}{\kappa(x)}\nabla l(x)\cdot\nabla f(x)+\frac{1}{2}\rho(x)^2\Delta f(x)$, where $l(x)$, $\tau(x)$, $\kappa(x)$, and $\rho(x)^2$ are scalar fields that define the ordering, speed, order scale-regularisation, and diffusive strength, respectively.  In place of a potential, $U(x)$, we introduce what we call an \textit{Ordering Function}, $l(x)$, which in theory would capture the same level set structure as $U(x)$, but whose scale is not important.  Through the scale-regularisation term $\kappa(x)$, we can make ODKs agnostic to the scale of $l(x)$, thus, the quality of an ordering function is only in the directionality it implies --- and importantly, by relaxing the constraints on $l(x)$ we significantly reduce the difficulty of the initial inference problem. As a result, ODKs allow us to approximate operators with arbitrary gradient flows, without coupling the drift to the diffusion, and can be extended to capture arbitrary drifts.  Thus, ODKs are a novel class of local kernels that are much more practically applicable to data generated from an unknown dynamical process.

An existing class of datasets to which this method is particularly well suited is single-cell RNA-sequencing (scRNA-seq) data \cite{luecken2019current}. In scRNA-seq, the measured features are the counts of RNA molecules associated with each gene in a given cell, with each sample corresponding to a single cell that must be destroyed in order to make the measurement. The destructive nature of this measurement process means that the RNA expression of the same cell cannot be measured twice, yielding data that is a collection of independent expression snapshots containing little to no direct dynamical information about gene expression dynamics. This has motivated a large and rapidly growing body of computational methods aimed at recovering the underlying dynamics from such static observations \cite{saelens2019comparison, street2018slingshot, lavenant2024toward}. The most common approach is to construct a so-called \textit{pseudotime} --- a hallucinated ordering of cells that is meant to correlate with the underlying sampling time of the system, and typically is structured as a distance from a collection of root, or progenitor, cells \cite{trapnell2014dynamics}. Notably, one of the most widely used pseudotime methods, Diffusion PseudoTime \cite{haghverdi2016diffusion}, already exploits diffusion maps by aggregating diffusion distances over all time scales \cite{coifman2006diffusion}.  Pseudotime, as it is conceptualised within the scRNA-seq field, is a clear example of an ordering function, and thus such data fits naturally within the ODK framework.

A central goal of single-cell genomics is to go beyond static descriptions of cell states and recover the underlying gene regulatory networks (GRNs) that govern cell fate determination dynamics. Learning a GRN is fundamentally equivalent to learning the dynamical system and vector field that drives gene expression --- that is, inferring the drift and diffusion terms of the stochastic process underlying the data. This perspective has motivated a number of recent computational frameworks: Dynamo \cite{qiu2022mapping} attempts to reconstructs the full vector field of gene expression dynamics from metabolic labelling data and uses it to infer regulatory Jacobians and predict perturbation outcomes, while RegVelo \cite{wang2026regvelo} integrates the inference of splicing kinetics with GRN structure in a joint deep learning model, enabling in silico perturbation predictions across cell types. Despite the power of these approaches, neither provides a principled probabilistic framework for inferring spatially varying diffusion, which reflects cell-to-cell variability in transcriptional noise and is rarely modelled explicitly in this field. ODKs offer a complementary approach: by using pseudotime as an ordering function and fitting the speed and diffusion parameters either from an external inference method or directly via gradient descent, we reduce the inference problem from $\frac{1}{2}mn(n+3)$ parameters --- corresponding to a full vector field and a symmetric 2-tensor field --- to only $2m$ parameters corresponding to two scalar fields, where $m$ is the number of samples, and $n$ the number of ambient dimensions. This dramatic reduction in the size of the inference problem makes it far more likely to be well-posed and tractable on the characteristically noisy and high-dimensional count matrices that arise in scRNA-seq experiments. The application of ODKs to scRNA-seq data is not explored in this paper, but this discussion illustrates the natural fit between the ODK framework and this existing class of data, motivating the development of ODKs as a principled tool for recovering the dynamics of complex systems in the absence of dynamical information.

In this paper, we will make clear the assumptions on our data and introduce the idea of an ordering function as a relaxation of a potential in section \ref{sec:problem_setup}.  Following this, in section \ref{sec:main_results}, we will define ODKs and explain the notational choices we make to provide practitioners fine control over the first- and second-order terms of the normalised Laplacian; we will then describe the limiting behaviour of this kernel and prove that its Laplacian approximates the infinitesimal generator of a gradient-flow SDE with isotropic diffusion.  In section \ref{subsec:nongradflow_extension}, we will extend ODK to capture arbitrary drifts through both a local ordering function and a multiple ordering function formulation; in section \ref{subsec:nonisotropic_extension}, we will generalise to anisotropic diffusions, which introduces a coupling between the drift and diffusion which we can remove by composing two different kernels within a Strang splitting scheme.  In addition, in section \ref{subsec:selftuning_odk}, we extend ODK to self-tuning bandwidths of the form $\rho(x)\rho(y)$.  Finally, in section \ref{sec:numerical_results}, we will apply ODKs to a variety of numerical examples on data sampled from dynamical and stochastic systems; of particular importance, in section \ref{subsec:numerical_results_sde}, we will introduce two objectives for gradient-based optimisation schemes, which we use to infer the drift and diffusion parameters of an unknown SDE with a spatially dependent diffusion.

\section{Problem Setup}\label{sec:problem_setup}
In this section, we extend the discussion of the introduction and specify the sampling assumptions of the data we wish to study before providing a more complete discussion of the ordering function and why we distinguish it from the more usual potential function formulation.

\subsection{Sampling of the Data}\label{subsec:setup_sampling_data}
We consider a dataset $\mathcal{D}=\{x_i\}_{i=1}^m\subset\mathbb{R}^n$ that is composed of $m$ independent samples from a deterministic or stochastic dynamical system; we assume the distribution over initial states of the system is known and has bounded support, denoted $p_0(x)$ with $X_0\coloneqq\operatorname{supp}p_0$.  The sampling procedure is assumed to be as follows: for a given initial state $x_0\sim p_0$, we evolve a single trajectory of the dynamical system forward up to a random stopping time $t\sim r$, where we call $r(t)$ the \textit{time-sampling distribution} and which is compactly supported on the interval $[0,T]$, where $T<\infty$ is the upper time horizon.  For deterministic data, the sampling distribution satisfies:
\begin{equation}
    q(x)=\int_0^Tr(t)p_0(\Phi_t^{-1}(x))\left|\frac{d\Phi_t^{-1}}{dx}(x)\right|\mathrm{d}t=\int_0^Tr(t)p_0(\Phi_t^{-1}(x))\exp\left(-\int_0^t\nabla\cdot v(\Phi_s(\Phi_t^{-1}(x)))\mathrm{d}s\right)\mathrm{d}t,
\end{equation}
where $\Phi_t(x)$ is the flow map of the dynamical system $\frac{dx}{dt}=v(x)$ over the time interval $[0,t]$, and satisfies: $\Phi_t(x_0)=x_0+\int_0^tv(\Phi_s(x_0))\mathrm{d}s$ and $\Phi_0(x_0)=x_0$.  For data generated by a stochastic system, the sampling distribution satisfies:
\begin{equation}\label{eq:sampling_distribution_integral}
    q(x)=\int_0^Tr(t)\int_{X_0}p_t(X_t=x|X_0=x_0)p_0(x_0)\mathrm{d}x_0\mathrm{d}t.
\end{equation}

We assume that our data lies on or near a $d$-dimensional Riemannian submanifold $\mathcal{M}$ of $\mathbb{R}^n$, with metric inherited from the standard Euclidean metric on $\mathbb{R}^n$.\footnote{By ``near", we mean that there can be some additional measurement noise that slightly perturbs the observed samples away from a smooth underlying manifold.}  More precisely, we assume there exists an isometric embedding $\iota:\mathcal{M}\to\mathbb{R}^n$ such that our data clusters around the image $\iota(\mathcal{M})$ --- since $\iota$ is an isomorphism of Riemannian manifolds, when it is not ambiguous to do so, we abuse notation and let $\mathcal{M}$ denote $\iota(\mathcal{M})$.  As noted in the introduction, the assumption of the existence of $\mathcal{M}$ is often referred to as the \textit{Manifold Hypothesis} and is at the core of the theoretical development of kernel-based methods and manifold learning in general \cite{bengio2013representation, fefferman2016testing}.  The justification of this assumption is simple: nature is inherently efficient, and if our data $\mathcal{D}\subset\mathbb{R}^n$, with $n\gg 1$, has been generated by some physical system, then by Occam's Razor it is unlikely that there will be $n$ independent degrees of freedom, and thus it is wise to consider the \textit{intrinsic dimensionality} of the system as being much smaller than $n$, and which we always denote as $d$.  The validation of this statement might reasonably be found in the success of methods that presume it: in the field of deep learning, methods such as AutoEncoders operate by first compressing the data with an encoder $f:\mathbb{R}^n\to\mathbb{R}^d$, with $d\ll n$, before undoing this transformation with a decoder $g:\mathbb{R}^d\to\mathbb{R}^n$, the objective being optimised is then some variation of $\lVert g\circ f-\text{id}\rVert$.  The fact that these bottlenecked architectures have been shown to work across a wide variety of datasets empirically supports the conclusion that most high-dimensional data will, in reality, have far fewer degrees of freedom than their expressed dimensionality would indicate \cite{hinton2006reducing,vincent2010stacked,kingma2013auto}.  After these considerations, we are emboldened to assume the ubiquitous existence of an \textit{underlying manifold}.

\subsection{Ordering Functions}\label{subsec:ordering_functions}
As discussed in the introduction, if the velocities generating the underlying data are unknown, and one has no strong a priori knowledge of the underlying system, inferring them is a highly challenging task whose complexity scales linearly with the system's dimensionality.  A more tractable task may be the inference of a potential function, $U(x)$, that governs the deterministic behaviour of the system by favouring lower energy states; but a potential must simultaneously imply the correct directionality and rate of change of the system, which requires either a priori knowledge or some additional dynamical constraints, examples of such additional information include: time marginals that allow you to set up an Optimal Transport problem \cite{schiebinger2019optimal}, velocity estimates \cite{bergen2020generalizing}, or an assumption of stationarity, where the sampling distribution should correspond to the potential through the Boltzmann distribution \cite{jin2018scepath} --- this latter example would also correspond to using $\alpha=\frac{1}{2}$ in the $\alpha$-normalised Diffusion Maps transition matrix \cite{coifman2006diffusion}.  In summary, inferring a reasonable potential remains a difficult task in practice and often requires building or tailoring an inference machine to the specific data at hand.

To navigate this obstacle, we propose relaxing the problem from inferring a potential to inferring --- or simply constructing --- an \textit{ordering} of the data, which is only required to imply the correct directionality of the system and whose scale is irrelevant.  Implying the correct directionality is just about making statements of the form ``given I'm at some $x\in\mathcal{M}$, which of my neighbours are closer to where I'm ultimately going than I currently am" (or equivalently, further away from where I came), this can be phrased as binary relation that tells you whether another point $y\in\mathcal{M}$ is further along in its evolution that my current $x$.  The study of different notions of ordering a mathematical space is tackled in the field of order theory.  The most familiar generalisation of the canonical ordering on $\mathbb{R}$ to more complex spaces is the notion of a Linear Quasi-Order (LQO)\footnote{Terms vary in the literature, with quasi-orders often referred to as pre-orders, and the property of linearity as totality.} of a topological space, $X$, which is a binary relation, $\leq:X\times X\to\{0,1\}$, that is reflexive ($\forall x\in X,\ x\leq x$), transitive ($\forall x,y,z\in X,\ x\leq y\ \land\ y\leq z \implies x\leq z$), and linear ($\forall x,y\in X,\ x\leq y\ \lor\ y\leq x$).  This ordering is continuous if its upper and lower sets are topologically closed, i.e. the order is adapted to the underlying topology in a natural manner; see \cite{ward1954partially} or Appendix \ref{appBsubsec:ordering_functions} for more details.  It is true that any continuous LQO on a topological space can be generated by a continuous function $l:X\to\mathbb{R}$ in the following sense:
\begin{equation}\label{eq:lqo_impliedbyorder}
    x\leq_l y\iff l(x)\leq l(y).
\end{equation}
We call such a continuous function an \textit{ordering function} of the topological space $X$.  We note that an ordering function, as we have conceptualised it, has been studied before as a utility function or, in more modern terminology, a utility representation of Quasi-Ordered Topological Spaces $(X,\leq)$.  It was shown by Debreu \cite{debreu1983representation} that $(X,\leq)$ admits a utility function of the form \eqref{eq:lqo_impliedbyorder} if and only if $\leq$ is linear and continuous, and $X$ is second-countable.  When we specialise to topological manifolds, $\mathcal{M}$, we can imagine a suitable ordering function being such that the level sets of $l(x)$ slice the manifold into disjoint connected upper and lower sets with a certain degree of regularity; we could formulate this concisely as: the level sets of an ordering function should (almost) form a codimension-one foliation of $\mathcal{M}$.\footnote{A $p$-dimensional foliation of a $d$-dimensional manifold $\mathcal{M}$ is a decomposition of $\mathcal{M}$ into a union of disjoint and connected $p$-dimensional submanifolds called the leaves of the foliation. The codimension of a foliation is $q=d-p$.}  We say \textit{almost} a foliation to allow for singular leaves which allow, for example, an ordering function to have extremal elements; we refer the reader to Appendix \ref{appBsubsec:ordering_functions} for a more in-depth discussion.  Hence, when we refer to an \textit{ordering function}, it is functions of this sort that one should be imagining, and it is the ordering function that we propose to use in place of potential functions.

The reason we make the distinction between orderings and potentials is to greatly reduce the difficulty of the initial inference step, as the complex dynamical constraints that must be satisfied by a potential have been removed.  As an example, if we knew the distribution over initial states of our system, $p_0(x)$, we could construct an ordering function as $l(x)=\mathbb{E}_{x_0\sim p_0}[d_\mathcal{G}(x_0,x)]$ where $d_\mathcal{G}(x,y)$ denotes the geodesic distance on some symmetric connected graph $\mathcal{G}=(\mathcal{V},\mathcal{E})$ of the data (see Appendix \ref{appBsubsec:heuristic_orderingfunction} for more details).  A good ordering function needs only to capture the rough flow of the dynamics, and thus we have much greater freedom when inferring (or constructing) it, and can conduct our analysis much more freely --- albeit at a coarser level.  Naturally, one should use whatever dynamical information is available and, where appropriate, set up a domain-specific inference problem for the ordering function or potential --- in the latter case, we would then set $l(x)=-U(x)$.

In section \ref{subsec:nongradflow_extension}, we generalise the notion of an ordering function of the form $l(x)$ to that of a \textit{local ordering function} $l(x,y)$, that satisfies $l(x,x)=0$ for all $x\in\mathcal{M}$, and defines an ordering in a local neighbourhood, $U_x$, of any $x\in\mathcal{M}$ through: $x\leq_{l_x} y\iff l(x,y)\geq0$.  This generalisation is entirely compatible with the former conception of an ordering function since any global (linear) ordering function $l(x)$ can be written in the form of a local ordering function through $l(x,y)=l(y)-l(x)$, with $U_x=\mathcal{M}$ for all $x\in\mathcal{M}$.  The implications of this are twofold: first, if a quasi-order is not linear, it is not always possible to define a single global ordering function, but it is guaranteed that there exists a local ordering function capable of representing the order.  Furthermore, if we wished to consider a more general notion of ordering on a topological space --- such as d-spaces \cite{grandis2009directed}, or streams \cite{krishnan2009convenient}, both of which relax the assumption of transitivity --- that is more local in nature (i.e. we can only compare points that are \textit{close} to one another), or if we wished to study more complex order topologies (such as cyclic orders which require the condition of transitivity to be relaxed), then a single global ordering function is unlikely to be sufficient (or is guaranteed to be insufficient in the case of cyclic orders) in representing this new complexity.  Second, as we will see in the next section, using a global ordering function restricts us to only approximate the dynamics of gradient-flow SDEs, i.e. infinitesimal generators of the form $\mathcal{L}=-\nabla U\cdot\nabla+\frac{1}{2}\rho^2\Delta$.  As we will see in section \ref{subsec:nongradflow_extension}, generalising to local ordering functions allows us to approximate operators with arbitrary drifts.

\section{Main Results}\label{sec:main_results}

In Local Kernels (\cite{berry2016local}), Berry \& Sauer define a prototypical local kernel, with isotropic diffusion, to be:
\begin{equation}\label{eq:prototypical_local_kernel_isotropic}
    K_\varepsilon(x,y)=\exp\left(-\frac{\lVert y-x-\varepsilon b(x)\rVert^2}{2\varepsilon\rho(x)^2}\right),
\end{equation}
where $b(x)$ is a vector field, and $\rho(x)$ a scalar field bounded uniformly away from zero. As they note, it is essential for convergence as $\varepsilon\to 0$ that $b(x)$ and $\rho(x)^2$ scale linearly with $\varepsilon$.  This is analogous to the short-timescale approximation of an It\^o transition $P(X_{dt}=y|X_0=x)$ where both the drift and diffusion scale linearly with the time increment $dt$. Therefore, intuitively: $\varepsilon$ is a timescale, $b(x)$ a drift, and $\tfrac{1}{2}\rho(x)^2I_d$ an isotropic diffusion.  With this structure in mind, we define an Ordered Diffusion Kernel (ODK) as:
\begin{definition}[Ordered Diffusion Kernel]\label{def:isotropic_odk}
    Let $\mathcal{M}$ be a $d$-dimensional Riemannian manifold embedded in $\mathbb{R}^n$ and with metric inherited from the ambient space, and let $l:\mathcal{M}\to\mathbb{R}$ be a $C^2(\mathcal{M})$ ordering function.  An Ordered Diffusion Kernel (ODK) on $\mathcal{M}$ with respect to the ordering $l$, is defined as:
    \begin{equation}\label{eq:general_odk}
        K_\varepsilon(x,y)\coloneqq h\left(\frac{\lVert y-x\rVert^2}{\varepsilon\rho(x)^2}+\frac{(l(y)-l(x)-\varepsilon\tau(x)\kappa(x))^2}{\varepsilon\rho(x)^2\kappa(x)^2}-\frac{(l(y)-l(x))^2}{\varepsilon\rho(x)^2\kappa(x)^2}\right),
    \end{equation}
    where $h:\mathbb{R}_{\geq0}\to\mathbb{R}_{\geq 0}$ is a smooth function bounded above by an exponential (i.e. $\exists C,\sigma>0:h(u)\leq C\exp(-u/\sigma^2)$); $\tau(x)$ is a $C^1(\mathcal{M})$ scalar function; $\kappa(x)$ is a $C^1(\mathcal{M})$ positive scalar function bounded uniformly away from zero (i.e. $\exists c_\kappa$ such that $\forall x\in\mathcal{M},\kappa(x)\geq c_\kappa$); and $\rho(x)$ is a $C^2(\mathcal{M})$ positive scalar function bounded uniformly away from zero (i.e. $\exists c_\rho>0$ such that $\forall x\in\mathcal{M},\rho(x)\geq c_\rho$).
\end{definition}
Denote by $(\star)$ the expression inside of $h(\cdots)$ in equation \eqref{eq:general_odk}, to understand the reasoning behind definition \ref{def:isotropic_odk} we will break it down term-by-term: the first term of $(\star)$ drives the diffusion analogous to a diffusion maps kernel \cite{coifman2006diffusion} (or, more precisely, the left-formulation of a Variable Bandwidth Kernel \cite{berry2016variable}), while the second term of $(\star)$ introduces a bias towards higher-order samples analagous to a $1$-dimensional local kernel \cite{berry2016local}.  The third term of $(\star)$ acts as a form of dimensionality correction, which is required since the first two terms together effectively constitute a $(d+1)$-dimensional Euclidean metric on the embedded graph of $l(x)$, and is dimensionally inconsistent with the underlying manifold.  More precisely, we note that since $l(x)$ is an ordering function, the (almost) foliating property implies the existence of a vector field normal to the level sets of $l$, denoted $\nabla l$, suppose around some point $x$, one takes a normal coordinate system $(x_1,\dots,x_d)$ such that $\partial_d$ is parallel to $\nabla l$ with respect to the inherited metric, after setting $b(x)=\tau(x)\partial_d$ in the local kernel formulation, we find:
\begin{equation}\label{eq:heuristic_expansion_odk_intuition}
    \lVert y-x-\varepsilon b(x)\rVert^2=\sum_{i=1}^{d-1}x_i^2+(x_d-\varepsilon\tau(x))^2=\sum_{i=1}^dx_i^2+(x_d-\varepsilon\tau(x))^2-x_d^2= \lVert y-x\rVert^2+(x_d-\varepsilon\tau(x))^2-x_d^2
\end{equation}
which is the exact form of $(\star)$ in definition \ref{def:isotropic_odk}.  So, these three terms together can be seen to be implicitly extracting the $x_d$ coordinate parallel to $\nabla l$, adding a bias along this coordinate, which resembles a $(d+1)$-dimensional norm, and then subtracting the extra, unbiased, squared distance along $x_d$. Thus, we see that \eqref{eq:general_odk}, is approximately equivalent to:
\begin{equation}
    K_\varepsilon(x,y)\approx h\left(\frac{1}{\varepsilon\rho(x)^2}\left\lVert y-x-\varepsilon\frac{\tau(x)}{\kappa(x)}\nabla l(x)\right\rVert^2\right)
\end{equation}
with $\varepsilon$ controlling the timescale, $\tfrac{1}{2}\rho(x)^2$ the diffusion, and $\tfrac{\tau(x)}{\kappa(x)}\nabla l(x)$ the drift.  To justify the introduction of the $\kappa(x)$ term first note that $l(y)-l(x)\approx\langle y-x,\nabla l(x)\rangle=\lVert\nabla l(x)\rVert x_d$, so $(l(y)-l(x))^2$ depends on the scale of $l(x)$, setting $\kappa(x)=\lVert\nabla l(x)\rVert$ in equation \eqref{eq:general_odk}, allows us to regularise the drift to be agnostic to the scale of $l$.  We do not fix $\kappa(x)=\lVert\nabla l (x)\rVert$, since in cases where one's ordering is thought to contain dynamical information (such as when $l(x)=-U(x)$ is a negative potential) setting $\tau(x)=\kappa(x)\equiv 1$ will give the desired dynamics.  Additionally, if $l(x)$ was equal to a timing function $t(x)$,\footnote{If our underlying dynamics are deterministic, the timing function can be defined unambiguously; otherwise, if our underlying dynamics are stochastic, we must decide on a convention, such as the mean first passage time, to define $t(x)$.} defined as the time to go from an initial state to the state $x$, then setting $\tau(x)\equiv 1$ and $\kappa(x)=\lVert\nabla l(x)\rVert^2\approx\lVert\nabla U(x)\rVert^{-2}$ might be most appropriate, since $\nabla t(x)\cdot \nabla U(x)\equiv -1$ will always be true, and if $\lVert \nabla U(x)\rVert$ is constant on level sets of $U$, then $\nabla t(x)=-\tfrac{\nabla U(x)}{\lVert\nabla U(x)\rVert^2}$, and this setup will reconstruct the dynamics exactly.  As a result $\kappa(x)$ should be treated more like a hyperparameter of the kernel, needing to be set by the user depending on whether ones ordering function is \textit{potential-like}, where one should set $\kappa(x)\equiv 1$, \textit{time-like}, where one should set $\kappa(x)\propto\lVert\nabla l(x)\rVert^2$, or is purely an representation of an ordering, in which case one should set $\kappa(x)\propto\lVert\nabla l(x)\rVert$.  In practice it is unlikely that you will have an analytic expression for $l(x)$, so one must estimate $\lVert \nabla l(x)\rVert$ numerically from the data, which is non-trivial in high dimensions; in Appendix \ref{appBsubsec:gradest} we overview a method we use to compute these estimates in practice, and which we have found to be superior to other alternatives when applied to real data.

To understand the limiting behaviour of ODKs as $\varepsilon\to 0$, we first define a $\varepsilon>0$ family of (backward\footnote{Calling an operator the \textit{backward} or \textit{forward} operator refers to the Kolmogorov equations \cite{kolmogoroff1931analytischen}, but when working with discrete data it is equivalent to the Markov transition probability matrix acting on the left or the right of a given test function, respectively.}) integral operators, $\{G_\varepsilon\}_{\varepsilon>0}$, associated to a given ODK.  For some bounded $f\in C^3(\mathcal{M})$, we define $G_\varepsilon$ as:
\begin{equation}\label{eq:def_backward_integral_operator}
    G_\varepsilon f(x)\coloneqq\varepsilon^{-d/2}\int_\mathcal{M}K_\varepsilon(x,y)f(y)\mathrm{dy}.
\end{equation}
For each $\varepsilon>0$, we can define a normalised integral operator corresponding to $G_\varepsilon$ through:
\begin{equation}\label{eq:def_normalised_backward_integral_operator}
    P_\varepsilon f(x)\coloneqq(G_\varepsilon\mathbb{1}(x))^{-1}G_\varepsilon f(x),
\end{equation}
where $\mathbb{1}(x)=1,\forall x\in\mathcal{M}$.  The expansion of these operators in terms of $\varepsilon$ is given in the following theorem:
\begin{theorem}[Expansion of Ordered Diffusion Kernels]\label{thm:odk_backward_expansion}
    Suppose $\mathcal{M}$ is a $d$-dimensional Riemannian manifold embedded in $\mathbb{R}^n$ and with the metric inherited from the standard Euclidean metric in ambient space. Let $l\in C^2(\mathcal{M})$ be an ordering function, and let $K_\varepsilon(x,y)$ be the corresponding ODK as defined in \eqref{eq:general_odk}.  Then for any $f\in C^3(\mathcal{M})$ the associated family of backward integral operators $\{G_\varepsilon\}_{\varepsilon> 0}$ admits the following asymptotic expansion as $\varepsilon\to 0$ when applied to $f$:
    \begin{equation}\label{eq:odk_backward_operator_expansion}
        G_\varepsilon f(x) = m(x)f(x) + \varepsilon\left( \omega_l(x)f(x) + m(x)\mathcal{L}f(x) \right)  + \varepsilon^{3/2}\Omega_l(x)+ \bigO(\varepsilon^2),
    \end{equation}
    where $\omega_l(x)$ depends on the manifold and kernel $K_\varepsilon(x,\cdot)$ at $x$, and $\Omega_l(x)$ depends on the manifold, kernel $K_\varepsilon(x,\cdot)$, and the test function $f$ at $x$. The differential operator $\mathcal{L}$ is a second-order elliptic operator of the form:
    \begin{equation}\label{eq:backward_differential_operator}
        \mathcal{L}f(x) = \mu(x)\cdot\nabla f(x) + D(x)\Delta f(x)
    \end{equation}
    where $m(x)$, $\mu(x)$ and $D(x)$, satisfy:
    \begin{align}
        m(x)&=\int_{T_x\mathcal{M}}h\left(\tfrac{\lVert z\rVert^2}{\rho(x)^2}\right)\mathrm{d}z\\
        \mu(x)&=\frac{\tau(x)}{\kappa(x)}\nabla l(x)\label{eq:odk_backward_expansion_drift}\\
        D(x)&=\ \frac{1}{2}m(x)^{-1}\int_{T_x\mathcal{M}}z_1^2\ h\left(\tfrac{\lVert z\rVert^2}{\rho(x)^2}\right)\mathrm{d}z\backwardscoloneqq \tfrac{1}{2}\sigma_h(x)^2\rho(x)^2\label{eq:odk_backward_expansion_diffusion}\\
        &\ \ \quad\text{where}\quad\sigma_h(x)^2=\frac{1}{\int_{T_x\mathcal{M}}h(\lVert z\rVert^2)\dif z }\int_{T_x\mathcal{M}}z_1^2h(\lVert z\rVert^2)\dif z.\label{eq:sigma_h_def}
    \end{align}
    The normalised backward integral operator can thus be expanded in terms of $\varepsilon$ as:
    \begin{equation}
        P_\varepsilon f(x) = f(x)+\varepsilon\mathcal{L}f(x)+\bigO(\varepsilon^{3/2}).
    \end{equation}
\end{theorem}
\begin{proof}
    See Appendix \ref{appAthm:odk_backward_expansion}.
\end{proof}
If we use the prototypical ODK, i.e. if we set $h(u)=\exp(-u/2)$, then $D(x)=\tfrac{1}{2}\rho(x)^2$ since $\sigma_h(x)^2\equiv1$ for the standard Gaussian kernel, thus $\mathcal{L}$ would be the infinitesimal generator of the It{\^o} SDE $dx_t=\frac{\tau(x)}{\kappa(x)}\nabla l(x)\mathrm{d}t+\rho(x)dw_t$, which gives a very explicit connection between the parameters of the kernel and the underlying stochastic dynamics they imply.

Although our choice of normalisation for $P_\varepsilon$ is the most canonical (when applied to discrete data, this normalisation is equivalent to row-normalising the transition matrix as we will see in section \ref{subsec:discrete_continuous}), it should be noted that other normalisations exist and have the same effect of removing the unwanted $\omega_l(x)f(x)$ term from the expansion of $G_\varepsilon f(x)$ (see, for example, \cite{belkin2003laplacian,coifman2006diffusion,berry2016local}). The utility of these normalisations are that they preclude the need to compute the dimensionality of the manifold, since the normalisation constant $\varepsilon^{-d/2}$ cancels, allowing us to directly use unnormalised kernels that do not integrate to 1 on the tangent space, and which otherwise would require us to estimate the volume element over the manifold, increasing the computational complexity and introducing new sources of error.

We also note that we are using a different convention from what is standard in the literature by defining $\mathcal{L}f=\tfrac{1}{m}[\tilde{\mu}\cdot\nabla f+\tilde{D}\Delta f]$, where $\tilde{\mu}$ and $\tilde{D}$ are the first- and second-order moments of the unnormalised kernel (as given equation 1 of \cite{berry2016local}).  We are motivated to do this because the symmetry enforced by wrapping the expression $(\star)$ in $h(\cdot)$ allows for highly interpretable forms for $\mu(x)$ and $D(x)$, and we do not wish to obscure this notationally by using the more general form utilised in local kernels \cite{berry2016local}.

We denote by $G_\varepsilon^*$ the formal adjoint of the (backward) integral operator $G_\varepsilon$ defined in equation \eqref{eq:def_backward_integral_operator}, i.e. for any $f,g\in L^2(\mathcal{M})\cap C^3(\mathcal{M})$, $G_\varepsilon^*$ satisfies $\langle f,G_\varepsilon g\rangle_{L^2}=\langle G_\varepsilon^*f,g\rangle_{L^2}$, we call this operator the forward integral operator and note that it necessarily takes the following form:
\begin{equation}\label{eq:def_forward_integral_operator}
    G_\varepsilon^* f(x)\coloneqq\varepsilon^{-d/2}\int_\mathcal{N}K_\varepsilon(y,x)f(y)\mathrm{dy}.
\end{equation}
Furthermore, due to our choice of convention for the differential operator $\mathcal{L}$, we are afforded the following far more intuitive form for the normalised adjoint operator:
\begin{equation}
    P_\varepsilon^*f(x)\coloneqq G_\varepsilon^*((G_\varepsilon\mathbb{1})^{-1}f)(x).
\end{equation}
In a discrete setting, this normalisation is equivalent to taking the matrix transpose of the Markov transition probability matrix given in \eqref{eq:def_normalised_backward_integral_operator} --- such a convenient formulation is afforded to us through our choice of normalisation.  Analogous to the expansion of the backward operator, we can derive the following expansion for the forward operators in the weak sense:
\begin{theorem}\label{thm:odk_forward_expansion}
    Let $K_\varepsilon(x,y)$ be an Ordered Diffusion Kernel under the same assumptions as Theorem \ref{thm:odk_backward_expansion}; then, in the weak sense, the family of adjoint operators, $\{G_\varepsilon^*\}$ admit the following expansion as $\varepsilon\to 0$ when applied to a test function $f\in C^3(\mathcal{M})\cap L^2(\mathcal{M})$:
    \begin{equation}\label{eq:odk_adjoint_expansion}
        G_\varepsilon^*f(x)\coloneqq\varepsilon^{-d/2}\int_\mathcal{M}K_\varepsilon(y,x)f(y)\mathrm{d}y=m(x)f(x)+\varepsilon(\omega_l(x)f(x)+\mathcal{L}^*(fm)(x))+\bigO(\varepsilon^{3/2}),
    \end{equation}
    where $\mathcal{L}^*$ is the formal adjoint of $\mathcal{L}$ as given in equation \eqref{eq:backward_differential_operator}, this operator takes the form (\cite{evans2022partial}, p320):
    \begin{equation}\label{eq:odk_adjoint_laplacian_form}
        \mathcal{L}^*f(x)=-\text{div}\left(f(x)\mu(x)\right)+\Delta(D(x)f(x)),
    \end{equation}
    where $\mu(x)$ and $D(x)$ are the normalised first- and second-order moments of $K_\varepsilon$, as in defined in equations \eqref{eq:odk_backward_expansion_drift} and $\eqref{eq:odk_backward_expansion_diffusion}$, respectively. The normalised adjoint operator can thus be expanded in terms of $\varepsilon$ as:
    \begin{equation}
        P_\varepsilon^*f(x)=f(x)+\varepsilon\mathcal{L}^*f(x)+\bigO(\varepsilon^{3/2}).
    \end{equation}
\end{theorem}
\begin{proof}
    See Appendix \ref{appAthm:odk_forward_expansion}.
\end{proof}
The reason we can only take the expansion of $G_\varepsilon^* f(x)$ up to $\bigO(\varepsilon^{3/2})$ is due to $\Omega_l(x)$ depending non-trivially on the second test function required to formulate weak convergence, and thus one cannot simply recombine this into an expression for $G_\varepsilon^*f(x)$.  Given the expansions of the normalised backwards and forwards integral operators of theorems \ref{thm:odk_backward_expansion} and \ref{thm:odk_forward_expansion}, and defining their respective \textit{Laplacians} to be:
\begin{equation}
    L_\varepsilon\coloneqq \tfrac{1}{\varepsilon}(P_\varepsilon-I),\qquad\text{\&}\qquad L_\varepsilon^*\coloneqq\tfrac{1}{\varepsilon}(P_\varepsilon^*-I)
\end{equation}
It is then trivial to determine the limiting behaviour of ODKs:
\begin{corollary}\label{cor:laplacian_expansions}
    Let $K_\varepsilon(x,y)$ be an Ordered Diffusion Kernel under the same assumptions as Theorem \ref{thm:odk_backward_expansion}; then:
    \begin{align}
        L_\varepsilon f(x)&\coloneqq\frac{(G_\varepsilon\mathbb{1})^{-1}G_\varepsilon f(x)-f(x)}{\varepsilon}=\mathcal{L}f(x)+\bigO(\varepsilon^{1/2})\\
        L_\varepsilon^*f(x)&\coloneqq\frac{G_\varepsilon^*((G_\varepsilon \mathbb{1})^{-1}f)(x)-f(x)}{\varepsilon}=\mathcal{L}^*f(x)+\bigO(\varepsilon^{1/2}).
    \end{align}
    where $\mathcal{L}$ and $\mathcal{L}^*$ are as in equations \eqref{eq:backward_differential_operator} and \eqref{eq:odk_adjoint_laplacian_form}, respectively.
\end{corollary}
\begin{proof}[Proof of Corollary \ref{cor:laplacian_expansions}]
    Both results follow immediately from the expansions of $P_\varepsilon f(x)$ and $P_\varepsilon^*f(x)$, given in Theorem \ref{thm:odk_backward_expansion} and Theorem \ref{thm:odk_forward_expansion}, respectively.
\end{proof}

This last corollary completes the proof that the family of differential operators $\{L_\varepsilon\}_{\varepsilon>0}$ asymptotically approaches the infinitesimal generator $\mathcal{L}$ of an It{\^o} process with gradient-flow drift and isotropic diffusion.  It is important to note that the nature of the results given in this section are about a family of continuous operators (i.e. operators defined directly on the underlying manifold) asymptotically approaching another continuous operator $\mathcal{L}=\mu\cdot\nabla+D\Delta$, and we have not yet said anything about the convergence of the discrete operators we define on real data to the continuous operators we've discussed here.  The nature of this discrete to continuous convergence is of great importance and has been studied in detail in the literature (see, for example, \cite{belkin2006convergence,gine2006empirical,singer2006graph,hein2007graph,ting2011analysis,wang2015spectral,tang2018limit,inagaki2025spectral}), where almost all results are formulated as convergence in probability, i.e. the probability of the discrepancy between the discrete and continuous objects being greater than any small error tends to zero as $m\to\infty$.  Different authors are concerned with different mathematical objects, so if one is primarily focused on spectral properties, it might be most fitting to consider the convergence of discrete eigenvectors to continuous eigenfunctions; but if you are more concerned with the operators themselves, then (pointwise or uniform) convergence of the discrete operator to the continuous operator would be more fitting.  For our purposes, we are satisfied with pointwise convergence in probability of the Laplacian operator; in section \ref{subsec:discrete_continuous}, we will show that under very mild assumptions on the sampling distribution the discrete Laplacian, $L_{m,\varepsilon}$, formed by evaluating ODK on dataset of $m$ independent samples from $q(x)$, will converge to the continuous Laplacian, $L_\varepsilon$, as $m\to\infty$, i.e. for a fixed and bounded $f\in C^3(\mathcal{M})$, we will show $L_{m,\varepsilon}f\xrightarrow[m\to\infty]{P}L_\varepsilon f$.  It is quite reasonable to assume this would always be the case, as if one had arbitrarily many data points, then one could resolve any function arbitrarily well from its evaluations across all these points.  The complexity of such questions arises when one wants to make the convergence rate optimal, especially under further regularity assumptions on $\mathcal{M}$ (such as bounded Ricci curvature \cite{inagaki2025spectral}) and $f$ (such as bounded Sobolev norms); then one can significantly tighten one's results and find regimes where the convergence rate is much faster.  Results such as these are not the main focus of this paper, so we do not consider such extensions, nor do we attempt to prove a tight convergence rate; but, in section \ref{subsec:discrete_continuous}, we provide a more detailed discussion of the bridge between discrete and continuous data, and state the pointwise convergence in probability in its most generic form.

Up until now, we have ignored how the boundary of $\mathcal{M}$ affects our asymptotics.\footnote{Topologically, a boundary is $\partial X\coloneqq\bar X\setminus\mathring X $, where $\bar X$ is the closure of $X$ and $\mathring X$ its interior, i.e. all points of closure that are not strictly inside of $X$.  As an example, if we consider $D^3=\{x\in\mathbb{R}^3:||x||_2\leq1\}$, then $\partial D^3=S^2=\{x\in\mathbb{R}^3:||x||_2=1\}$, which is the 2-sphere, but we see that $\partial S^2=\varnothing$, since all points of closure are also inside $S^2$.  For topological manifolds, this is generalised by the addition of charts to the half space $\mathbb{R}^{d-1}\times\mathbb{R}_{\geq 0}$ at the boundary, it is known that if $\mathcal{M}$ is a $d$-dimensional manifold with boundary, then $\partial\mathcal{M}$ is a $(d-1)$-dimensional manifold.}  In the proof of Theorem \ref{thm:odk_backward_expansion}, we assume that we were a certain \textit{large enough} distance away from the boundary - more precisely, we assumed we were at least $\varepsilon^\gamma$ away from the boundary for some $\gamma\in(0,1/2)$ --- so that when $\varepsilon<1$ we can define an adequate region around a ball of radius $\varepsilon^{1/2}$ which does not intersect the boundary.  In practice, the data will only allow you to set $\varepsilon$ to be so small before you cannot adequately resolve the kernel on the data --- i.e. all the variation and detail of the kernel is happening in the space between observed samples; thus, the previous condition provides a limit on how close we can get to the boundary before the asymptotics break down.  This has been studied in the literature (see, for example, \cite{coifman2006diffusion,zhou2011behavior,belkin2012toward}) and it has been found that the leading-order errors of the asymptotics given in Theorems \ref{thm:odk_backward_expansion} and \ref{thm:odk_forward_expansion} become of order $\bigO(\varepsilon^{1/2})$ rather than $\bigO(\varepsilon)$ near the boundary; as an example, if $\varepsilon=10^{-4}$ then the error near the boundary would be of order $\sim 10^{-2}$ which is 100 times larger than in the interior.  Critically, most datasets sampled from dynamical systems will likely have boundaries, since if the initial conditions are fairly regular, the orbits (paths) will trace out (cluster around) something similar to a hypercylinder, which necessarily has a boundary.  The only way to correct for these boundary effects is to sample more data points near the boundary so you can resolve smaller $\varepsilon$'s, but this will just reduce the volume of the problematic boundary region, not remove it entirely.  Thus, it is important to note that ODK will likely experience these boundary effects, and in section \ref{subsec:numerical_results_ode} we show results in figures \ref{fig:vf_ode_velests} and \ref{fig:vf_ode_mcv} where the breakdown in resolution near the boundary is clearly visible.  Musingly, these boundary effects can be mitigated if the data is noisy, such as when we sample from an ODE in the presence of some measurement noise, or if the data is sampled from an SDE with non-negligible intrinsic noise.  In these cases, the noise effectively blurs the boundary region to the point where one cannot really tell exactly where the boundary might lie; thus, the kernel sees \textit{boundary points} more so as outliers in the data rather than as a fundamental topological limitation\footnote{We do not rigorously consider the effects of noise on convergence rates near the boundary, but this would be interesting future work.}.  This results in a less obvious breakdown in the kernel's resolution and can be observed in figures \ref{fig:ou_sde_example1} and \ref{fig:nl_results} of section \ref{subsec:numerical_results_sde}, where, to our eye, the data clearly has some sort of boundary but the kernel still infers reasonable values for $l(x)$ and $\rho(x)$.  Ultimately, if the data clearly resolves a boundary, one must expect higher errors in these regions, but if the boundary is blurry, the kernel will perform better than the current theory would suggest.

We do not discuss practical methods for parameterising ODK until section \ref{subsec:numerical_results_sde}, but we refer the reader to Appendix \ref{appBsubsec:heuristic_parametrisations} for a discussion of several heuristic parameterisations that one could utilise for $l(x)$, $\varepsilon$, $\tau(x)$ and $\rho(x)$.

\subsection{Generalising ODK to Arbitrary Drifts}\label{subsec:nongradflow_extension}

In Theorems \ref{thm:odk_backward_expansion} and $\ref{thm:odk_forward_expansion}$ of section \ref{sec:main_results}, we use an ordering function of the form $l(x)$; we do this because the philosophy of an ordering function is to make ODKs easily applicable to data with very little dynamical information.  This results in the first-order advection term of an ODK's Laplacian being proportional to $\nabla l(x)$; thus, we are restricted to studying the infinitesimal generators of gradient-flow SDEs.  In this section, we will extend ODK to capture arbitrary drifts (analogous to local kernels \cite{berry2016local}) by generalising $l(x)$ to a \textit{local ordering function}.

As discussed briefly in section \ref{subsec:ordering_functions} (and in more detail in Appendix \ref{appBsubsec:ordering_functions}), if one wished to consider more complex classes of orderings (such as cyclic orders), or one wished to work on domains with non-trivial homotopy structures (such as a torus), then one needs to relax the assumptions of linearity and global transitivity of $\leq_l$.  Particularly, the relaxation of transitivity to be a local property, leads us to consider \textit{local ordering functions}, $l:\bigsqcup_{x\in\mathcal{M}} U_x\to\mathbb{R}$, that satisfy $l(x,x)=0$ for all $x\in\mathcal{M}$ (which we will see implies reflexivity), and at every $x\in\mathcal{M}$ defines a locally linear order on $U_x$ through $x\leq_{l_x} y\iff l_x(y)\geq 0$ where $l_x: U_x\to\mathbb{R}$ denotes $l(x,\cdot)$ with the first argument fixed.  We define a locally Ordered Diffusion Kernel (lODK) to be:
\begin{definition}[Locally Ordered Diffusion Kernel]\label{def:nongradflow_odk}
    Let $\mathcal{M}$ be a $d$-dimensional Riemannian manifold embedded in $\mathbb{R}^n$ and with metric inherited from the ambient space; furthermore, assume that the $\exists c>0$ such that $\inf_{x\in\mathcal{M}}\text{inj}(x)=c>0$, where $\text{inj}(x)$ is the injectivity radius of the exponential map at $x$ --- we call this property bounded geometry.  Let $\{U_x\}_{x\in\mathcal{M}}$ be a family of open neighbourhoods of every point in $\mathcal{M}$ such that the diameters of the $ U_x$ are uniformly bounded away from zero.  Let $\text{Dom}(l)=\bigsqcup_{x\in\mathcal{M}} U_x$ and let $l:\text{Dom}(l)\to\mathbb{R}$ be a $C^2$ local ordering function of two inputs satisfying $l(x,x)=0$ for all $x\in\mathcal{M}$.  Under the same constraints on $h(u),\tau(x),\rho(x)$ and $\kappa(x)$ as in definition \ref{def:isotropic_odk}, we define a locally Ordered Diffusion Kernel (lODK) to be:
    \begin{equation}
        K_\varepsilon(x,y)=h\left(\frac{\lVert y-x\rVert^2}{\varepsilon\rho(x)^2}+\frac{(l(x,y)-\varepsilon\tau(x)\kappa(x))^2}{\varepsilon\rho(x)^2\kappa(x)^2}-\frac{l(x,y)^2}{\varepsilon\rho(x)^2\kappa(x)^2}\right)
    \end{equation}
\end{definition}
The purpose of assuming the injectivity radius of the exponential map is bounded away from zero on $\mathcal{M}$ (a property sometimes referred to as bounded geometry) is to ensure that $l(x,y)$ can always be well-defined.  It is true that for an arbitrary smooth vector field $b:\mathcal{M}\to T\mathcal{M}$, we can define $l(x,y)$ such that $\nabla l_x(x)=b(x)$ for all $x\in\mathcal{M}$.  We discuss this in more detail in Appendix \ref{appBsubsec:ordering_functions}, but the gist is to define each $l_x$ locally as $l_x(y)=y_ib_i(x)$ where $y_1,\dots,y_d$ are normal coordinates about $x$.  By the assumption of bounded geometry, we can always choose a family $\{U_x\}_{x\in\mathcal{M}}$ such that for every $x\in\mathcal{M}$, the exponential map will be a diffeomorphism over $U_x$ and $\text{diam}(U_x)$ will be bounded uniformly away from zero.  By the smoothness of $b$ and non-degeneracy of $\{U_x\}_{x\in\mathcal{M}}$, we can smoothly glue these locally defined functions into a globally defined function $l:\bigsqcup_{x\in\mathcal{M}}U_x\to\mathbb{R}$ with a partition of unity such that $\forall x,y\in\mathcal{M}$ we will have $l(x,y)=l_x(y)$ and $l(x,x)=0$.

Before stating the asymptotic expansions associated to lODKs, it is important to comment on the domain of both the backward operator and its adjoint: since the local ordering function is defined on the domain $\text{Dom}(l)=\bigsqcup_{x\in\mathcal{M}} U_x$ for a family of subsets $\{U_x\}_{x\in\mathcal{M}}$ with diameter bounded away from zero (always possible by the assumption of bounded geometry), then the domain of integration for the backwards operator at a point $x\in\mathcal{M}$ is clearly $U_x$, since $l(x,y)$ is not defined for $y\in\mathcal{M}\setminus U_x$.  Likewise, for a given $x\in\mathcal{M}$ the domain of the adjoint (forwards) operator applied to $f$ at $x$ would be $V_x\coloneqq \{y\in\mathcal{M}: x\in U_y\}$, and since we assume every $U_y$ is an open neighbourhood of $y$ with diameter bounded away from zero, and $x\in U_x$ for all $x\in\mathcal{M}$ (by assumption), then $V_x$ will be an open neighbourhood of $x$ with bounded diameter (since otherwise the diameter being bounded away from zero would be violated).  Furthermore, the intersection $U_x\cap V_x\backwardscoloneqq M_x$ will be a neighbourhood of $x$ with diameter bounded away from zero and, when viewing $l(x,y)$ through the order it implies, will be the largest subset of $\mathcal{M}$ containing $x$, such that all elements $y\in M_x$ are comparable to $x$ from the perspective of the local orders at both $x$ and $y$.  In order for the adjoint to be well-defined, we would wish for the domain of integration to be equal and thus we assume the backwards and forwards operators are defined as integrals of $K_\varepsilon(x,y) f(y)$ and $K_\varepsilon(y,x)f(y)$ over the domain $M_x$.  By referencing the proof of Theorem \ref{thm:odk_backward_expansion} given in Appendix \ref{appAthm:odk_backward_expansion}, for a given value of $\gamma\in(0,1/2)$ --- which is chosen to bound the operator away from the boundary and to define a region of integration, that being the ball of radius $\varepsilon^\gamma$ around some $x\in\mathcal{M}$, outside of which contributions are of order $\bigO(\varepsilon^2)$ --- we must have that $B(x,\varepsilon^\gamma)\subset M_x$, since $M_x$ is not arbitrarily small, there always exists a $\varepsilon_x>0$ below which this will always be satisfied.  We can view this as if we have taken some sensible extension of $l(x,y)$ to all of $\mathcal{M}\times\mathcal{M}$, since the integral of $K_\varepsilon(x,y)f(y)$ ($K_\varepsilon(y,x)f(y)$) outside of $B(x,\varepsilon^\gamma)$ contributes $\bigO(\varepsilon^2)$ to $G_\varepsilon f(x)$ ($G_\varepsilon^* f(x)$) the terms depending on the choice of extension are hidden away in the higher-order terms, and thus the limiting behaviour of $L_\varepsilon f(x)$ ($L_\varepsilon^* f(x)$) is independent of the choice of extension.  It is this latter statement that allows us to talk about the limiting behaviour of lODKs without ambiguity.

\begin{theorem}[Expansion of locally Ordered Diffusion Kernels]\label{thm:nongradflow_odk_expansion}
    Let $\mathcal{M}$ be a $d$-dimensional Riemannian manifold with bounded geometry embedded in $\mathbb{R}^n$, and with the metric inherited from the standard Euclidean metric.  Let $l(x,y)$ be a $C^2(\text{Dom}(l))$ local ordering function of two inputs and let $b(x)=\nabla l_x(x)$.  For the local ODK $K_\varepsilon(x,y)$, as defined in \ref{def:nongradflow_odk}, and for any $f\in C^3(\mathcal{M})$, we have the following asymptotic expansion of the associated family of backward integral operators $\{G_\varepsilon\}_{\varepsilon>0}$ as $\varepsilon\to 0$, applied to $f$:
    \begin{equation}
        G_\varepsilon f(x) = m(x)f(x) + \varepsilon\left( \omega_l(x)f(x) + m(x)\mathcal{L}_bf(x) \right)  + \varepsilon^{3/2}\Omega_l(x)+ \bigO(\varepsilon^2)
    \end{equation}
    where the second-order differential operator $\mathcal{L}_b$ takes the form:
    \begin{equation}
        \mathcal{L}_b f(x)=\frac{\tau(x)}{\kappa(x)}b(x)\cdot\nabla f(x)+D(x)\Delta f(x)
    \end{equation}
    In the above, $m(x)$ and $D(x)$ are as in Theorem \ref{thm:odk_backward_expansion}.  The normalised backward integral operator and its associated Laplacian (defined in the usual manner) can be expanded as:
    \begin{equation}
        P_\varepsilon f(x) =f(x)+\varepsilon\mathcal{L}_b f(x)+\bigO(\varepsilon^{3/2})\implies L_\varepsilon f(x)=\mathcal{L}_b f(x)+\bigO(\varepsilon^{1/2})
    \end{equation}
    Analogous to the previous theorems, the family of adjoint operators, $\{G_\varepsilon^*\}_{\varepsilon>0}$, satisfies the following asymptotic expansion when applied to $f\in C^3(\mathcal{M})\cap L^2(\mathcal{M})$, in the weak sense:
    \begin{equation}
        G_\varepsilon^*f(x)=m(x)f(x)+\varepsilon(\omega_l(x)f(x)+\mathcal{L}_b^*(fm)(x))+\bigO(\varepsilon^{3/2})
    \end{equation}
    where $\mathcal{L}_b^*$ is the formal adjoint of $\mathcal{L}_b$ and takes the form:
    \begin{equation}
        \mathcal{L}_b^* f(x)=-\diver\left(f(x)\frac{\tau(x)}{\kappa(x)}b(x)\right)+\Delta(D(x)f(x))
    \end{equation}
    The normalised adjoint integral operator $P_\varepsilon^*f=G_\varepsilon^*((G_\varepsilon \mathbb{1})^{-1} f)$ and its Laplacian can be expanded as:
    \begin{equation}
        P_\varepsilon^* f(x)=f(x)+\varepsilon\mathcal{L}_b^* f(x)+\bigO(\varepsilon^{3/2})\implies L_\varepsilon^* f(x)=\mathcal{L}_b^* f(x)+\bigO(\varepsilon^{1/2})
    \end{equation}
\end{theorem}
\begin{proof}
    See Appendix \ref{appAthm:nongradflow_odk_expansion}.
\end{proof}
When we were considering the global ordering function $l(x)$, we introduced $\kappa(x)$ as a means to cleanly make ODK's independent of the scale of $l(x)$, particularly, if we set $\kappa(x)=\lVert\nabla l(x)\rVert$ then we can make the ODK entirely agnostic to the scale of $l$; the same principle holds for lODKs, except one would now set $\kappa(x)=\lVert\nabla l_x(x)\rVert$.  We note that in Theorem \ref{thm:nongradflow_odk_expansion}, we assumed that $l(x,y)$ was $C^2$ in both arguments, we can actually relax this to only having $\tfrac{l(x,y)}{\kappa(x)}$ be $C^2$ in both arguments.\footnote{This condition trivially still implies $l(x,y)$ is $C^2$ in the second argument which is essential in the proof of Theorem \ref{thm:nongradflow_odk_expansion}.}  Indeed, suppose $a:\mathcal{M}\to\mathbb{R}_{>0}$, was any discontinuous strictly positive function, then setting $\tilde l(x,y)=a(x)l(x,y)$ and $\tilde\kappa(x)=a(x)$ for some $l(x,y)$ that is $C^2$ in both arguments, will result in $\tfrac{\tilde l(x,y)}{\tilde \kappa(x)}=l(x,y)$ still being $C^2$ in both arguments in spite of $\tilde l(x,y)$ being discontinuous in the first argument.  The takeaway from this is that if one is constructing a local ordering function, one can be less careful when it comes to ensuring $l(x_1,\cdot)$ and $l(x_2,\cdot)$ are commensurate when $x_1$ and $x_2$ are close to one another, so long as $\kappa(x)$ is set as to correct for any ad hoc differences in scale.

Essentially, if one wishes to study dynamics more general than gradient flows, or the topology of the data (such as domains with non-trivial homotopy groups, e.g. a torus) and/or the ordering (such as cyclic orderings, e.g. ordering the numbers on a clock) makes the application of a global ordering function unreasonable, then using a local ordering function $l(x,y)$ to encode directionality, and using the associated lODK to define transition probabilities and to approximate the infinitesimal generator $\mathcal{L}_b=b\cdot\nabla+D\Delta$, is the best course of action.  As an example of such a dataset, we could consider cell-cycle scRNA-seq data, where cells go through distinct and well-defined phases ($G1\to S\to G2\to M\to G1\to\dots$) as they grow and divide; if during each phase it has been shown experimentally that changes in expression of a particular subset of marker genes uniquely identifies how far along each phase one is, then one could very reasonably define a local ordering function $l(x,y)=\sum_{C\in\{G1,S,G2,M\}}\mathbb{1}(x\in C)(l_{C}(y)-l_C(x))$ where $\mathbb{1}(x\in C)$ is an indicator function that determines which phase we are in, and each $l_C$ is an ordering function constructed based off the heuristics we have for cells in each respective phase.  Furthermore, in this example, the cyclic nature of the cell-cycle makes it theoretically impossible to define a single global ordering function, $l(x)$, since global transitivity is violated in cyclic order topologies.  In section \ref{subsec:numerical_results_operators}, we demonstrate the convergence of the isotropic and anisotropic ODK as $\varepsilon\to 0$ on synthetic data sampled uniformly on a torus; due to the cyclic domain of the torus, a global ordering would produce discontinuities at the boundaries of the glueing and thus we construct a local ordering function $l(x,y)$ to match the known topology of the data.

As a final comment, if one is unsure whether to use a local ordering function or a global ordering function on one's specific data, before committing to one or the other, we recommend applying the method of Persistent Homology (from the field of Topological Data Analysis, see \cite{edelsbrunner2002topological,zomorodian2004computing,carlsson2009topology} for more details) to ascertain if non-trivial topologies are present in one's data, and if they are, at what scales do they emerge.  Since rotational dynamics, cyclic orderings, and non-trivial homotopies will almost certainly coincide given our assumptions on data generation outlined in section \ref{subsec:setup_sampling_data}, such a preliminary analysis will provide practitioners with the information they need to decide which formulation is most appropriate given their specific problem.

\subsubsection{Multi-Ordered Diffusion Kernels}\label{subsec:multi_odk}
In this short subsection, we will write ODK in a manner such that they are compatible with multiple ordering functions, denoted $\{l_i:\mathcal{M}\to\mathbb{R}\}_{i=1}^k$.  The practical motivation for this would be if each $l_i$ indicated a preference towards (or maybe (negative) proximity to) a set of distinguished states $\{s_i\}_{i=1}^k\subset \mathcal{M}$, such as during stem cell differentiation, where we know the expected gene expression of each differentiated cell type and thus could compute a distance measure (such as geodesic distance on a connected neighbourhood graph) for each differentiated state to form each $l_i$.  Another example, if one interprets each $l_i$ as being a negative potential, could be to model a process where there are multiple potential wells that combine non-trivially to affect the evolution of a system; if each potential can be computed independently but we were unsure of the relative importance of each field, then it would be useful to apply this directly rather than attempt to combine the different scales of each potential.  

More theoretically, in Appendix \ref{appBsubsec:ordering_functions}, we discuss how a global ordering function, $l(x)$, as we have conceptualised it, is a single utility function representation of a Linear Quasi Order, and the intuition of relaxing global transitivity to form a locally linear quasi-order was the motivation for defining the local ordering function $l(x,y)$.  A multi-utility representation of a --- not necessarily linear --- quasi-order, $\leq$, is a family of functions $\{f_i\}_{i\in I}$ that satisfies, for all $x,y\in X$: $x\leq y\iff f_i(x)\leq f_i(y)\quad \forall i\in I$; it is true that if $X$ is a locally compact Hausdorff space that is also $\sigma$-compact, then any continuous quasi-order will admit a continuous multi-utility representation.  Relaxing the linearity constraint leads us to consider a (finite) family of ordering functions, $\{l_i\}_{i=1}^k$, just as relaxing the transitivity constraint led us to consider local ordering functions.  We define a multi-Ordered Diffusion Kernel (mODK) to be:
\begin{definition}[Multi-Ordered Diffusion Kernel]\label{def:multi_ODK}
    Let $\mathcal{M}$ be a $d$-dimensional Riemannian manifold embedded in $\mathbb{R}^n$ and with metric inherited from the ambient space; and let $\{l_i:\mathcal{M}\to\mathbb{R}\}_{i=1}^k\subset C^2(\mathcal{M})$ be a finite set of ordering functions.  A multi-Ordered Diffusion Kernel (mODK) is defined as:
    \begin{equation}\label{eq:multi_odk}
        K_\varepsilon(x,y)\coloneqq h\left(\frac{\lVert y-x\rVert ^2}{\varepsilon\rho(x)^2}+\sum_{i=1}^k\frac{(l_i(y)-l_i(x)-\varepsilon\tau_i(x)\kappa_i(x))^2}{\varepsilon\rho(x)^2\kappa_i(x)^2}- \sum_{i=1}^k\frac{(l_i(y)-l_i(x))^2}{\varepsilon\rho(x)^2\kappa_i(x)^2}\right),
    \end{equation}
    where $h:\mathbb{R}_{\geq0}\to\mathbb{R}_{\geq 0}$ is a smooth function bounded above by an exponential ($\exists C,\sigma>0:h(u)\leq C\exp(-u/\sigma)$); $\{\tau_i(x)\}_{i=1}^k$ and $\{\kappa_i(x)\}_{i=1}^k$ are sets of k $C^1(\mathcal{M})$ scalar functions with each $\kappa_i(x)$ bounded uniformly away from zero; and $\rho(x)$ is a smooth positive scalar functions bounded uniformly away from zero.
\end{definition}
We do not provide the asymptotic expansions of all intermediary results and instead simply state the form of the limiting operator, which we denote $\mathcal{L}_m$, as:
\begin{equation}
    \mathcal{L}_mf(x)=\left(\sum_{i=1}^k\tau_i(x)\frac{\nabla l_i(x)}{\kappa_i(x)}\right)\cdot\nabla f(x)+D(x)\Delta f(x)
\end{equation}
where $D(x)$ is as in Theorem \ref{thm:odk_backward_expansion}, and from hereon we denote by $\mu_m(x)$ the advection term of $\mathcal{L}_m$.  Once again, we include the functionality of the $\kappa_i(x)$ terms, which allows us to make the kernel agnostic to the scale of each $l_i(x)$ by setting $\kappa_i(x)=\lVert\nabla l_i(x)\rVert$ if one so wishes.  The drift is thus a weighted combination of the directions of increasing order; the fact that each ordering function is weighted independently by a distinct, state-dependent $\tau_i(x)$ allows for complex combinations of the order directions that could be interpreted as encoding interactions between different potential fields.  

Interestingly, in theory, if one had chosen $\{l_i\}_{i=1}^k$ such that $\text{Span}(\{dl_i|_x:i=1,\dots,k\})=T_x^*\mathcal{M}$ for all $x\in\mathcal{M}$, then we could choose $\{\tau_i\}_{i=1}^k$ such that $\mu_m(x)$ is any arbitrary vector-field.  We can show this constructively: by Whitney's embedding theorem we know that all Hausdorff, second-countable $d$-dimensional smooth manifolds can be embedded into $\mathbb{R}^{2d}$, and the assumptions of our manifold ensure there is an embedding to $\mathbb{R}^n$, so if $k=\min(n,2d)$,\footnote{For most manifolds this value of $k$ will be much higher than is necessary.} and $\iota:\mathcal{M}\to\mathbb{R}^k$ was the associated embedding, then setting $l_i(x)=\iota_i(x)$, where $\iota_i(x)$ is the $i^\text{th}$ component of $\iota(x)$, is sufficient for $\{d\iota_i\}_{i=1}^k$ to span the cotangent bundle, and thus for any vector field $b\in\mathcal{T}\mathcal{M}$, there exist scalar functions $\{\tau_i\}_{i=1}^k$ such that $b(x)=\sum_{i=1}^k\tau_i(x)\nabla l_i(x)$ for all $x\in\mathcal{M}$.

It is worth noting that in the proof of Theorem \ref{thm:odk_backward_expansion}, for a given $x\in\mathcal{M}$ we expand $l(y)$ about $x$ and end up with something approximating $K_\varepsilon(x,y)=h\left((\varepsilon\rho(x)^2)^{-1}\lVert y-x-\varepsilon\mu(x)\rVert^2+\bigO(\varepsilon)\right)$, where the $\bigO(\varepsilon)$ error depends on the choice of ordering function and the geometry of $\mathcal{M}$.  From the point of view of the asymptotics, this is inconsequential, but on real data there is a fundamental limit to how small $\varepsilon$ can be made; thus, if one decides to use a staggering number of ordering functions in an mODK, then under certain conditions this $\bigO(\varepsilon)$ error can accumulate to the point of skewing the kernel and obscuring the underlying dynamics.

\subsection{Generalising ODK to Anisotropic Diffusions}\label{subsec:nonisotropic_extension}
In previous theorems, we restricted ourselves to using isotropic diffusions of the form $D(x)=\tfrac{1}{2}\rho(x)^2I_d$.  We can extend this to anisotropic diffusions in a natural way but at the cost of coupling the first- and second-order moments of the Laplacian.  We define anisotropic Ordered Diffusion Kernels (anODKs) to be:
\begin{definition}[Anisotropic Ordered Diffusion Kernel]\label{def:nonisotropic_odk}
    Under the same constraints as definition \ref{def:isotropic_odk}, the anisotropic Ordered Diffusion Kernel (anODK) is defined as:
    \begin{equation}\label{eq:non_isotropic_odk}
        K_\varepsilon(x,y)\coloneqq h\left(\frac{(y-x)^TA(x)^{-1}(y-x)}{\varepsilon\rho(x)^2}+\frac{(l(y)-l(x)-\varepsilon\tau(x)\kappa(x))^2}{\varepsilon\rho(x)^2\kappa(x)^2}-\frac{(l(y)-l(x))^2}{\varepsilon\rho(x)^2\kappa(x)^2}\right),
    \end{equation}
    where $A(x)$ is a symmetric and uniformly elliptic ($\exists c>0$ such that $\lambda_\text{min}(A(x))\geq c$ for all $x\in\mathcal{M}$) type $\binom{2}{0}$-tensor field.\footnote{We use the convention that a type $\binom{k}{l}$-tensor field, $T$, acts on $k$-covariant and $l$-contravariant vectors, i.e. $T:\bigsqcup_{x\in\mathcal{M}}(T_x^*\mathcal{M})^{k}\times(T_x\mathcal{M})^l\to\mathbb{R}$, the space of all such tensors is written $\mathcal{T}_l^k\mathcal{M}$; a $\binom{2}{0}$-tensor, $A$, is uniformly elliptic if all of its eigenvalues are bounded away from zero, in other words, there exists $c>0$ such that all $\lambda$ which are the solutions to the equation $g_{jk}A^{ki}v^j=A_j^iv^j=\lambda v^j$, satisfy $\lambda>0$ --- if $g_{ij}=\delta_{ij}$ this is equivalent to just taking the spectrum of the matrix $A(x)$.}
\end{definition}

In this definition, we make a separation between the diffusion driving terms $A(x)$ and $\rho(x)^2$, we do this for two reasons: firstly, for practical convenience, since the biased order terms are not \textit{aware} that the coordinates have been rescaled by $A$, this separation can help to minimise the coupling that $A$ imparts on the inferred drift --- this will become clear in Theorem \ref{thm:noniso_odk_expansion}; secondly, we do this so that we can make a distinction between $A(x)$ acting like a diffusion and like a generalised Riemannian metric.  In the former formulation $\rho(x)^2$ and $A(x)$ act in tandem and combine into the diffusion tensor $D_A(x)=\tfrac{1}{2}\rho(x)^2A(x)$; but in the latter formulation, since $A(x)$ only appears in the computation of the distance metric $(y-x)^TA(x)^{-1}(y-x)$ we can view this as working on the Riemannian manifold $(\mathcal{M},A^{-1})$ with non-trivial metric.  In the remainder of this section, we adopt the former interpretation, but in Appendix \ref{appAsec:noniso_metric} we discuss the latter interpretation in more detail.

\begin{theorem}[Expansion of Anisotropic Ordered Diffusion Kernels]\label{thm:noniso_odk_expansion}
    Under the same assumptions as theorem \ref{thm:odk_backward_expansion}, the anisotropic ODK has the expansion:
    \begin{equation}
        G_\varepsilon f(x)=m_A(x)f(x)+\varepsilon(\omega_l(x)f(x)+m_A(x)\mathcal{L}_Af(x))+\varepsilon^{3/2}\Omega_l(x)+\bigO(\varepsilon^2),
    \end{equation}
    where $m_A(x)=m(x)\sqrt{|A(x)|}$ where $m(x)$ is the zeroth-order moment from theorem \ref{thm:odk_backward_expansion}.  The second-order elliptic operator, $\mathcal{L}_A$, takes the form:
    \begin{equation}\label{eq:non_isotropic_laplacian}
        \mathcal{L}_Af(x)=\mu_A(x)\cdot\nabla f(x)+\tr_{D_A}\nabla^2f(x)=\mu_A^i(x)\partial_i f(x)+D_A^{ij}(x)\partial_i\partial_j f(x)
    \end{equation}
    where in the second equality we have written this operator in Euclidean normal coordinates (and we note that the trace of a $\binom{0}{2}$-tensor $T$ with respect to a $\binom{2}{0}$-tensor $B$ is defined to equal $\tr_B T=B^{ij}T_{ji}$).  The first- and second-order moments $\mu_A(x)$ and $D_A(x)$ satisfy:
    \begin{align}
        \mu_A(x)&=\frac{\tau(x)}{\kappa(x)}A(x)\nabla l(x)=A(x)\mu(x)\\
        D_A(x)&=\frac{1}{2}m_A(x)^{-1}\int_{T_x\mathcal{M}}zz^Th\left(\frac{z^TA(x)^{-1}z}{\rho(x)^2}\right)\dif z=\tfrac{1}{2}\sigma_h(x)^2\rho(x)^2A(x)=A(x)D(x)
    \end{align}
    where $\mu(x)$ and $D(x)$ are the same as in Theorem \ref{thm:odk_backward_expansion}.  The normalised backward integral operator and its Laplacian can be expanded in terms of $\varepsilon$ as:
    \begin{equation}
        P_\varepsilon f(x) = f(x)+\varepsilon\mathcal{L}_A f(x)+\bigO(\varepsilon^{3/2})\implies L_\varepsilon f(x)=\mathcal{L}_A f(x)+\bigO(\varepsilon^{1/2})
    \end{equation}
    Analogously to Theorem \ref{thm:odk_forward_expansion}, in the weak sense the adjoint operator satisfies:
    \begin{equation}
        G_\varepsilon^*f(x)=m_A(x)f(x)+\varepsilon(\omega_l(x)f(x)+\mathcal{L}_A^*(fm_A)(x))+\bigO(\varepsilon^{3/2}),
    \end{equation}
    where $\mathcal{L}_A^*$ is the formal adjoint of the operator $\mathcal{L}_A$ in theorem \ref{thm:noniso_odk_expansion}, satisfying $\langle f,\mathcal{L}_A^*g\rangle_{L^2}=\langle \mathcal{L}_Af,g\rangle_{L^2}$ for any $f,g\in L^2(\mathcal{M})$, this operator takes the form (\cite{evans2022partial}, p320):
    \begin{equation}
        \mathcal{L}_A^*f(x)=-\diver\left(f(x)\mu_A(x)\right)+\Delta(D_A(x)f(x))=-\partial_i(f(x)\mu_A^i(x))+\partial_i\partial_j(D_A^{ij}(x)f(x))
    \end{equation}
    and the normalised adjoint integral operator and its Laplacian can be expanded as:
    \begin{equation}
        P_\varepsilon^*f(x)=f(x)+\varepsilon\mathcal{L}_A^*f(x)+\bigO(\varepsilon^{3/2})\implies L_\varepsilon^* f(x)=\mathcal{L}_A^* f(x)+\bigO(\varepsilon^{1/2}).
    \end{equation}
\end{theorem}
\begin{proof}
    See Appendix \ref{appAthm:noniso_odk_backward_expansion}.
\end{proof}
From the form of $\mu_A$, we see that in Euclidean normal coordinates the drift is coupled with the anisotropic diffusion matrix in an analogous manner to how the drift is coupled with the diffusion in the isotropic TMDmap kernels.  Intuitively, this happens because the biased drift term isn't \textit{aware} of the ambient coordinate rescaling introduced by $A(x)$.  As it is, this may be useful in some scenarios, but for our purposes, unless $A(x)$ has been defined in such a way as to skew the direction of increasing order towards more dynamically relevant regimes, it seems like a loss of generality when compared to anisotropic Local Kernels with gradient flows.  We also note that identical results would hold if we had used a local or multi-ordering function, with $\mu_A=A\mu_b$ or $\mu_A=A\mu_m$ in each case, respectively.

We can correct for this coupling by composing two separate kernels together: first, we can construct the normalised integral operator $P_{\varepsilon,\alpha}^\mu$ corresponding to the isotropic ODK with $\tau(x)$, $\kappa(x)$ and $l(x)$ identical, but with the diffusion set to $\bar\rho(x)^2\coloneqq\alpha\rho(x)^2\lambda_\text{min}(A(x))$, where $\alpha\in(0,1)$ and $\lambda_\text{min}(A(x))>0$ gives the smallest eigenvalue of the matrix $A(x)$ --- we assume that $A(x)$ is such that $\lambda_\text{min}(A(x))$ is bounded uniformly away from zero.  We can then also construct the normalised integral operator $P_{\varepsilon,\alpha}^A$ corresponding to the anisotropic local kernel with no drift and covariance matrix $\rho(x)^2\bar A(x)\coloneqq \rho(x)^2(A(x)-\alpha\lambda_\text{min}(A(x))I_d)=\rho(x)^2A(x)-\bar\rho(x)^2I_d$.  Let the diffusions and Laplacians corresponding to these kernels be denoted by $D_\alpha^\mu=\tfrac{1}{2}\sigma_h^2\bar\rho^2I_d$, $\mathcal{L}_\alpha^\mu=\mu\cdot\nabla+D_\alpha^\mu\Delta$, and $D_\alpha^A=\tfrac{1}{2}\sigma_h^2\rho^2\bar A$ and $\mathcal{L}_\alpha^A=\tr_{D_\alpha^\mu}\nabla^2$, respectively.  We also note that indeed:  
\begin{equation}
    D_\alpha^\mu(x)+D_\alpha^A(x)=\tfrac{1}{2}\sigma_h(x)^2\bar\rho(x)^2I_d+\tfrac{1}{2}\sigma_h(x)^2\left[\rho(x)^2A(x)-\bar\rho(x)^2I_d\right]=\tfrac{1}{2}\sigma_h(x)^2\rho(x)^2A(x)=D_A(x).
\end{equation}
It is clear from the results of theorems \ref{thm:odk_backward_expansion} and \ref{thm:noniso_odk_expansion} that $\bar{\mathcal{L}}_A\coloneqq\mathcal{L}_\alpha^\mu+\mathcal{L}_\alpha^A=\mu(x)\cdot\nabla +\tr_{D_A}\nabla^2$, is the operator we desire to reconstruct, and which is the infinitesimal generator of the It{\^o} SDE $dx_t=\mu(x_t)\mathrm{d}t+\rho(x_t)A^{1/2}(x_t)\mathrm{d}w_t$, where $A^{1/2}(x)$ satisfies $A(x)=A^{1/2}(x)A^{1/2}(x)^T$ for all $x\in\mathcal{M}$.  However, it is unsatisfactory to simply sum the two approximations for of $\mathcal{L}_\alpha^\mu$ and $\mathcal{L}_\alpha^A$, since this breaks the chain of operators which we usually descend, i.e. we will not have an expression for the matrix $\bar{P}_{\varepsilon,\alpha}$ whose Laplacian $\varepsilon^{-1}(\bar{P}_{\varepsilon,\alpha}-I)\to\bar{\mathcal{L}}_A$ as $\varepsilon\to 0$; thus we seek to first construct $\bar{P}_{\varepsilon,\alpha}$ such that this property holds.  

The key insight is that if $\bar P_{\varepsilon,\alpha}$ is such an operator, then it should be true that $\bar P_{\varepsilon,\alpha}=\exp(\varepsilon\bar{\mathcal{L}}_A)$, where $\exp$ is the operator exponential defining in terms of the power-series of $e^x$.  Since $\bar{\mathcal{L}}_A=\mathcal{L}_\alpha^\mu+\mathcal{L}_\alpha^A$, we can construct $\bar P_{\varepsilon,\alpha}$ in terms of products of $P_{\varepsilon,\alpha}^\mu\approx\exp(\varepsilon\mathcal{L}_\alpha^\mu)$ and $P_{\varepsilon,\alpha}^A\approx\exp(\varepsilon\mathcal{L}_\alpha^A)$, but this construction will not be exact unless $\mathcal{L}_\alpha^\mu$ and $\mathcal{L}_\alpha^A$ commute, since in general $\exp(\varepsilon\mathcal{L}_\alpha^\mu+\varepsilon\mathcal{L}_\alpha^A)\neq\exp(\varepsilon\mathcal{L}_\alpha^\mu)\exp(\varepsilon\mathcal{L}_\alpha^A)$.  Approximating the exponential of an operator from products of exponentials of its components is a well studied problem in the field of numerical integration (see, for example, \cite{hairer2006numerical}), and amounts to generalising different integration schemes (such as the Euler method, the Verlet scheme, etc).  The generalisation of the Euler method gives the so-called Lie-Trotter approximation (see \cite{trotter1959product}), given by $\exp(\varepsilon\mathcal{L}_\alpha^\mu)\exp(\varepsilon\mathcal{L}_\alpha^A)\approx P_{\varepsilon,\alpha}^\mu P_{\varepsilon,\alpha}^A$, but we notice that the order of operations matter (i.e. $P_{\varepsilon,\alpha}^\mu P_{\varepsilon,\alpha}^A\neq P_{\varepsilon,\alpha}^A P_{\varepsilon,\alpha}^\mu$) and particularly, the order seems to switch upon taking the adjoint, i.e. if $\bar{P}_{\varepsilon}=P_{\varepsilon,\alpha}^\mu P_{\varepsilon,\alpha}^A$ (or $P_{\varepsilon,\alpha}^A P_{\varepsilon,\alpha}^\mu$) then $\bar{P}_{\varepsilon} f$ applies the anisotropic diffusion and then the drift (or applies the drift and then the anisotropic diffusion), but then $\bar{P}_{\varepsilon}^*f$ would apply the drift and then the anisotropic diffusion (or the anisotropic diffusion and then the drift).  To ameliorate this asymmetry, and to reduce the magnitude of the error arising from the choice of composition scheme (see Appendix \ref{appBsubsec:strang_splitting} for more details), we choose to use a higher-order method that generalises the Verlet scheme, and which is often referred to as the Strang Splitting scheme \cite{strang1968construction}:
\begin{equation}\label{eq:strang_operator_formula}
    \bar{P}_\varepsilon\approx\exp\left(\frac{\varepsilon}{2}\mathcal{L}_\alpha^A\right)\exp\left(\varepsilon\mathcal{L}_\alpha^\mu\right)\exp\left(\frac{\varepsilon}{2}\mathcal{L}_\alpha^A\right)\approx P_{\frac{\varepsilon}{2},\alpha}^AP_{\varepsilon,\alpha}^\mu P_{\frac{\varepsilon}{2},\alpha}^A.
\end{equation}
By applying the Baker-Campbell-Hausdorff (BCH) expansion (see \cite{hairer2006numerical}, Chapter 3.4) we can see that the Lie-Trotter scheme has leading order error $\bigO(\varepsilon^2)$ in approximating $\exp(\varepsilon\bar{\mathcal{L}}_A)$, but by repeated application of the BCH expansion, we see that the symmetric form of the Strang Splitting scheme causes all even order errors to cancel, and thus the leading order error is $\bigO(\varepsilon^3)$, thus providing a better approximation for $\exp(\varepsilon\bar{\mathcal{L}}_A)$.  Furthermore, we note that the observed asymmetry between the order of applications between the backward operator and its adjoint are no longer present; however the order of operations still matters, and we must choose whether the inner exponential is with respect to $P_{\varepsilon,\alpha}^\mu$ as we write it in equation \eqref{eq:strang_operator_formula}, or with respect to $P_{\varepsilon,\alpha}^A$.  In section \ref{subsec:numerical_results_operators}, we apply the Strang splitting scheme to approximate a known Laplacian, in this example we found that the form given in equation \eqref{eq:strang_operator_formula} was more stable to the choice of $\alpha$ and converged faster to the desired Laplacian $\bar{\mathcal{L}}_A$, thus we chose to use this order for our composition.  It is our belief that the best order of operations will most likely depend on which of $\mu(x)$ and $A(x)$ is most variable across the data, with the inner factor of the Strang splitting scheme being the more variable kernel.  In most scenarios, we expect the order of operations given in equation \eqref{eq:strang_operator_formula} to give the best approximation.  We introduce a parameter $\alpha\in(0,1)$ to determine how the anisotropic diffusion is split between the two operators $P_{\varepsilon/2,\alpha}^A$ and $P_{\varepsilon,\alpha}^\mu$.  In theory, this is a free parameter, but in practice there will be an optimal range of $\alpha$ where the Strang splitting scheme is most effective; this range will be determined by the manifold, $\mathcal{M}$, and the eigenstructure of $A(x)$ (see Appendix \ref{appBsubsec:choosing_alpha} for a more detailed discussion).  Ultimately, if $A(x)$ was isotropic, setting $\alpha=1$ would be optimal, hence the optimal splitting of the diffusion will likely amount to making $\alpha$ as large as possible before $P_{\varepsilon/2,\alpha}^A$ degenerates and is no longer able to resolve $\rho(x)^2\bar A(x)$ on the data --- thus, the process of choosing an optimal $\alpha$ is analogous to choosing an optimal bandwidth.

To summarise, we have extended ODKs to anisotropic diffusions at the cost of coupling the drift and diffusion, analogous to how TMDmap couples drift and diffusion in the isotropic case.  We then propose implementing a Strang splitting scheme to remove this coupling by splitting the anisotropic diffusion between a standard isotropic ODK, as defined in \ref{def:isotropic_odk}, and a symmetric anisotropic local kernel.  In Appendix \ref{appAsec:noniso_metric}, we interpret $A^{-1}$ as being a non-trivial Riemannian metric and show how three variations of ODK can be constructed to approximate the more complex infinitesimal generator in this case.

\subsubsection{Self-Tuning Ordered Diffusion Kernels}\label{subsec:selftuning_odk}
In what follows, we will demonstrate the limiting behaviour of an ODK after replacing the $\rho(x)^2$ term in equation \ref{eq:general_odk} with $\rho_1(x)\rho_2(y)$, where for $i=1,2$, $\rho_i(x):\mathcal{M}\to\mathbb{R}_{\geq c}$ is a $C^2$ scalar function bounded uniformly away from zero by $c>0$.  Such a diffusion doesn't have a standardised name: it has been called a \textit{self-tuning}, adaptive, and \textit{locally scaling} bandwidth (for example in \cite{zelnik2004self}); but when $\rho_1=\rho_2$, which is almost always the case, it is often called a symmetric bandwidth (for example in \cite{berry2016variable}); since we do not impose that $\rho_1=\rho_2$ we chose to call it the self-tuning bandwidth.\footnote{In this section alone we favour the term bandwidth over diffusion because in most applications a self-tuning bandwidth is more of an ad hoc construction than an intrinsic quality of the underlying process, thus we feel the term bandwidth creates some necessary separation between the two interpretations.}  A self-tuning Ordered Diffusion Kernel (stODK) is thus defined as:
\begin{definition}[Self-Tuning Ordered Diffusion Kernel]\label{def:selftuning_odk}
    Let $\mathcal{M}$ be a $d$-dimensional Riemannian manifold embedded in $\mathbb{R}^n$ and with metric inherited from the ambient space; and let $l:\mathcal{M}\to\mathbb{R}$ be an ordering function.  A self-tuning Ordered Diffusion Kernel (stODK) is defined as:
    \begin{equation}\label{eq:selftuning_odk}
        K_\varepsilon^{st}(x,y)\coloneqq h\left(\frac{\lVert y-x\rVert^2}{\varepsilon\rho_1(x)\rho_2(y)}+\frac{(l(y)-l(x)-\varepsilon\tau(x)\kappa(x))^2}{\varepsilon\rho_1(x)\rho_2(y)\kappa(x)^2}-\frac{(l(y)-l(x))^2}{\varepsilon\rho_1(x)\rho_2(y)\kappa(x)^2}\right),
    \end{equation}
    where $h:\mathbb{R}_{\geq0}\to\mathbb{R}_{\geq 0}$ is a smooth function bounded above by an exponential ($\exists C,\sigma>0:h(u)\leq C\exp(-u/\sigma)$); $\tau(x)$ is a smooth scalar function; $\kappa(x),\rho_1(x)$, and $\rho_2(x)$ are smooth positive scalar functions bounded away from zero.
\end{definition}
The motivation for using such kernels is usually to facilitate variable bandwidths (which is useful on non-uniformly sampled data as it improves both the fit and the convergence rate, \cite{berry2016variable}) when working with unbiased (zero-drift) kernels without making them asymmetric.  If we try to make the Diffusion Maps kernel $h(\varepsilon^{-1}\lVert x-y\rVert^2)$ into a variable bandwidth kernel in a manner consistent with the theory, then one might define $K_\varepsilon(x,y)=h((\varepsilon\rho(x)^2)^{-1}\lVert x-y\rVert^2)$, but it would then not be guaranteed that $K_\varepsilon(x,y)=K_\varepsilon(y,x)$, which although not a problem theoretically, will introduce problem computationally, since methods such as diffusion maps seek to take the eigendecomposition of their (normalised) kernel matrices.  Especially when these matrices are large, it is computationally very important --- for reasons such as numerical stability, computational speed, and to avoid a complex spectrum --- that the kernel matrix be symmetric.  Thus, the symmetric bandwidth $\rho(x)\rho(y)$ was introduced in place of $\rho(x)^2$ to solve this issue.  

Other than computational convenience, there are good reasons why one might want to use such a bandwidth: for example, if your data is very sparse or inhomogeneously sampled (such as scRNA-seq data), then if you are using your kernel to numerically integrate the Fokker-Planck equation from a known initial condition, then it may become very challenging to parametrise your transition probabilities so as to transition into these sparsely sampled regions; if one sets $\rho\propto q^{-1/d}$, which is standard in the literature (see, for example, \cite{berry2016consistent}) and is effectively equivalent to setting $\rho(x)=\lVert x-x^k\rVert$ where $x^k$ is the $k^\text{th}$ nearest neighbour to $x$ in $\mathcal{D}$, (and which follows from $q(x)\propto \lVert x-x^k\rVert^{-d}$) --- which is just stating that if a point's neighbours are further away then that region has been less densely samples --- then going from $\rho(x)^2$ to $\rho(x)\rho(y)$ will increase the bandwidth relative to more densely sampled states, thus increasing the probability of transitioning into sparser regimes.  For such ill-conditioned datasets, a self-tuning kernel is thus quite well justified.  The expansion of stODKs is given in the following Theorem:
\begin{theorem}[Expansion of Self-Tuning Ordered Diffusion Kernels]\label{thm:self_tuning_expansion}
    Suppose $\mathcal{M}$ is a $d$-dimensional Riemannian manifold embedded in $\mathbb{R}^n$ and with metric inherited from the standard Euclidean metric in ambient space. Let $l\in C^2(\mathcal{M})$ be an ordering function (or a local ordering function), and let $K_\varepsilon^{st}(x,y)$ be the corresponding stODK, as defined in \eqref{def:selftuning_odk}, with self-tuning bandwidth of the form $\rho_1(x)\rho_2(y)$ for smooth functions $\rho_1,\rho_2\in C^2(\mathcal{M})$ bounded uniformly away from zero. Then for any $f\in C^3(\mathcal{M})$ the associated family of backward integral operators $\{G_\varepsilon\}_{\varepsilon> 0}$ admits the following asymptotic expansion as $\varepsilon\to0$ when applied to $f$:
    \begin{equation}
        G_\varepsilon f(x) = m(x)f(x) + \varepsilon\left( \omega_l(x)f(x) + m(x)\mathcal{L}_{st}f(x) \right)  + \varepsilon^{3/2}\Omega_l(x)+ \bigO(\varepsilon^2),
    \end{equation}
    at a point $x\in\mathcal{M}$ a distance larger than $\varepsilon^\gamma$ from $\partial\mathcal{M}$, for a fixed $\gamma\in(0,1/2)$; and where $\omega_l(x)$ depends on the manifold and the kernel $K_\varepsilon(x,\cdot)$ at $x$, and $\Omega_l(x)$ depends on the manifold, the kernel $K_\varepsilon(x,\cdot)$ at $x$, and also the test function $f(x)$. The differential operator $\mathcal{L}$ is a second-order elliptic operator of the form:
    \begin{equation}\label{eq:full_selftuning_operator}
        \mathcal{L}_{st}f(x) = \left[\mu(x)+(d+2)D(x)\frac{\nabla \rho_2(x)}{\rho_2(x)}\right]\cdot\nabla f(x) + D_{st}(x)\Delta f(x)
    \end{equation}
    where $m(x)$, $\mu(x)$ and $D(x)$, satisfy:
    \begin{equation}\label{eq:selftuning_moments}
        m(x)=\int_{T_x\mathcal{M}}h\left(\tfrac{\lVert z\rVert^2}{\rho_1(x)\rho_2(x)}\right)\mathrm{d}z,\qquad\mu(x)=\frac{\tau(x)}{\kappa(x)}b(x),\qquad D_{st}(x)=\ \tfrac{1}{2}\sigma_h(x)^2\rho_1(x)\rho_2(x)\\
    \end{equation}
    The normalised backward integral operators and their Laplacians can thus be expanded in terms of $\varepsilon$ as:
    \begin{equation}
        P_\varepsilon f(x) = f(x)+\varepsilon\mathcal{L}_{st}f(x)+\bigO(\varepsilon^{3/2})\implies L_\varepsilon^{st}f(x)=\mathcal{L}_{st}f(x)+\bigO(\varepsilon^{1/2})
    \end{equation}
    Analogous to previous theorems, the family of adjoint operators $\{G_\varepsilon^*\}_{\varepsilon>0}$ satisfies the following asymptotic expansion, in the weak sense, when applied to $f\in C^3(\mathcal{M})\cap L^2(\mathcal{M})$:
    \begin{equation}
        G_\varepsilon^*f(x)=m(x)f(x)+\varepsilon(\omega_l(x)f(x)+\mathcal{L}_{st}^*(fm)(x))+\bigO(\varepsilon^{3/2})
    \end{equation}
    where $\mathcal{L}_{st}^*$ is the formal adjoint of $\mathcal{L}_{st}$ and takes the form:
    \begin{equation}
        \mathcal{L}_{st}^* f(x)=-\text{div}\left(f(x)\left[\mu(x)+(d+2)D_{st}(x)\frac{\nabla \rho_2(x)}{\rho_2(x)}\right]\right)+\Delta(D_{st}(x)f(x))
    \end{equation}
    The normalised adjoint integral operator $P_\varepsilon^*f=G_\varepsilon^*((G_\varepsilon \mathbb{1})^{-1} f)$ and its Laplacian can be expanded as:
    \begin{equation}
        P_\varepsilon^* f(x)=f(x)+\varepsilon\mathcal{L}_{st}^* f(x)+\bigO(\varepsilon^{3/2})\implies L_\varepsilon^* f(x)=\mathcal{L}_{st}^* f(x)+\bigO(\varepsilon^{1/2})
    \end{equation}
\end{theorem}
\begin{proof}
    See Appendix \ref{appAthm:self_tuning_expansion}.
\end{proof}
An immediate corollary of this result is that it strengthens the proof given by Berry \& Harlim in ``Variable Bandwidth Kernels" (VBKs, \cite{berry2016variable}) since in their Appendix A.3 they prove the expansion of the \textit{right-formulation} in the weak sense, whereas we do not, and by setting $\tau(x)\equiv 0$, $l(x)\equiv 0$ and $\kappa(x)\equiv 1$ we recover the form of a VBK.  The reason they use the weak formulation is that they seek to make a substitution of the form $\hat y=\tfrac{y-x}{\sqrt{\rho_2(y)}}+x$, but such a substitution need not be invertible and thus cannot be directly applied.  We do not need to do this to derive our expansion, but if we made the substitution $z\mapsto \tilde z/\sqrt{\rho_1(x)\rho_x(x)}$ in the formula for $m(x)$ as given in equation \eqref{eq:selftuning_moments}, then we would recover the necessary normalisation factors of $\rho_1^{d/2}\rho_2^{d/2}$ that constitute the difference in convention between ODKs and VBKs.

In equation \eqref{eq:full_selftuning_operator}, the drift term $\mu(x)$ may equal $\tfrac{\tau(x)}{\kappa(x)}\nabla l(x)$, $\tfrac{\tau(x)}{\kappa(x)}\nabla l_x(x)$, or $\sum_{i=1}^k\tfrac{\tau_i(x)}{\kappa_i(x)}\nabla l_i(x)$ depending on whether you use a global, local, or multi-ordering function, respectively.  We can simplify equation \ref{eq:full_selftuning_operator} by assuming $h(u)=\exp(-u/2)$ is the prototypical kernel (which implies $\sigma_h(x)^2=1$), and that $\rho_1=\rho_2$, we can then write the limiting operator as:
\begin{equation}
    \mathcal{L}_{st}=(\mu+\tfrac{1}{2}(d+2)\rho\nabla \rho)\cdot\nabla f+\tfrac{1}{2}\rho^2\Delta f=\left(\mu-\tfrac{1}{2}\tfrac{d+2}{d}\tfrac{\nabla q}{q}\right)\cdot\nabla f+\tfrac{1}{2}\rho^2\Delta f
\end{equation}
where in the second equality we consider the special case of $\rho\propto q^{-1/d}$.  It is important to note that the additional advection term coming from the $y$-dependent bandwidth, and which appears alongside the drift in equation \eqref{eq:full_selftuning_operator}, is not \textit{actually} a drift from the view of the transition operator $P_\varepsilon$ and is more of an \textit{effective drift}, the transition probabilities are still centred about $x+\varepsilon\mu(x)$ it is just the probabilities surrounding the mean (which we usually attribute to the noise of an SDE) are skewed in such a way as to to create a probabilistic flux in the direction proportional to $\nabla \rho_2$.

Let us assume for convenience that $\sigma_h(x)\equiv 1$, then we can introduce a $\beta$-regularisation that acts analogously to diffusion maps $\alpha$-regularisation, where $\beta\in[0,1]$, and we use a left-bandwidth of the form $\rho_1(x)^{2\beta}$ and a right-bandwidth of the form $\rho_2(y)^{2(1-\beta)}$.  The results are most interpretable when $\rho_1=\rho_2$, so assuming this, the Laplacian takes the following form:
\begin{equation}
    \mathcal{L}_{\beta}f=\mu\cdot\nabla f+(1-\beta)(d+2)\rho\nabla\rho\cdot\nabla f+\tfrac{1}{2}\rho^2\Delta f
\end{equation}
We see that when $\beta=1$ (i.e we have a bandwidth $\rho(x)^2$) the additional advection term disappears and we are back in the standard isotropic ODK case where the diffusion only depends on ones current state (think of this as the \textit{left-formulation} of VBK); but when $\beta=0$ ($\rho(y)^2$) we can think of this as being the \textit{right-formulation} of VBK where the diffusion is determined by where the state is about to transition, allowing the kernel to push probability towards higher-bandwidth states (where the meaning of this depends on how you've defined your bandwidth) - this is the least Markov of all values of $\beta$.  When $\beta=1/2$ ($\rho(x)\rho(y)$), we can think of this as the standard self-tuning bandwidth, or \textit{symmetric formulation} of VBK where the bandwidth is an equally weighted mixture of the bandwidth in my current state and in my future state.  Finally, we note that when you set $\beta=\tfrac{2d+3}{2d+4}=1-\tfrac{1}{2(d+2)}$, then our additional advection terms combines nicely with the diffusion term:
\begin{equation}
    (d+2)(1-\tfrac{2d+3}{2d+4})\rho\nabla\rho\cdot \nabla f+\tfrac{1}{2}\rho^2\Delta f=\tfrac{1}{2}\rho\nabla \rho\cdot\nabla f+\tfrac{1}{2}\rho^2\Delta f=\tfrac{1}{2}\rho\nabla\cdot(\rho\nabla f)
\end{equation}
which one might notice as the diffusion term for a Stratonovich SDE $dx_t=\mu(x_t)\mathrm{d}t+\rho(x_t)\circ dw_t$.  This is interesting, particularly since when one actually defines an SDE intrinsically on the tangent space of a manifold, the Stratonovich formulation of noise is much preferred \cite{elworthy1998stochastic}.  However, following our discussion of the $\tfrac{1}{2}\rho\nabla\rho$ only being an effective drift, the transition probabilities $P_\varepsilon$ would not be consistent with the Stratonovich SDE associated to the generator $L_{\varepsilon,1-(2(d+2))^{-1}}$, since the additional advection from the drift will be real in the Stratonovich sense.

\subsection{Bridging Discrete and Continuous Domains}\label{subsec:discrete_continuous}

The results we have given up until now are done in a continuous setting, i.e. all functions can be evaluated on the entirety of $\mathcal{M}$ and the integrals are exact; in other words, there are no gaps in our underlying space.  When we move to discrete data $\mathcal{D}=\{x_i\}_{i=1}^m$ composed of $m$ independent samples $x\sim q(x)$, where $\text{supp}(q)\subseteq\mathcal{M}$ is the underlying process's transient distribution --- or more broadly just the sampling distribution of the data.  Then we replace our continuous functions with vectors of their evaluations across $\mathcal{D}$, and our integrals become finite sums that only approximate the integral in the Monte-Carlo sense.  We will use the following square bracket notation to represent the discretisation of continuous functions onto $\mathcal{D}$: assume we fix an indexing of $\mathcal{D}$, then for any function $F:\mathcal{M}^k\to\mathbb{R}$, we define $[F]\in\mathbb{R}^{m^{\otimes k}}$, i.e. a $k$-dimensional array, that satisfies $[F]_{ij\dots k}=F(x_i,x_j,\dots,x_k)$.  We use this notation only in this subsection to make it absolutely clear whether an object is to be seen as discrete or continuous.

The essence of Monte-Carlo integration, is that if I want to approximate $\int_\mathcal{M} f(x)q(x)\mathrm{d}x$ for any $f\in L^2(\mathcal{M})$ and where $q(x)$ is a probability density function; then I can sample a \textit{large number} of points $x_i\sim q(x)$ from this sampling distribution and then compute the arithmetic mean of $f(x)$ evaluated on these points \cite{robert2004monte}.  Precisely, we will have:
\begin{equation}
    \lim_{m\to\infty} \frac{1}{m}\sum_{i=1}^mf(x_i)\overset{a.s.}{=}\int_{\mathcal{M}}f(x)q(x)\ \mathrm{d}x
\end{equation}
when $x_i\sim q(x)$.  The intuition behind this is that the arithmetic mean implicity uses a uniform measure of the data, but because the data has been sampled from $q(x)$, it will naturally cluster in high-density regions and then \textit{overcount} the integral in this region; as a result, you implicitly weight your integral by the distribution from which your data is sampled.  From our perspective, the integral operator $G_\varepsilon$ could be approximated as:
\begin{equation}
    \frac{1}{m}[K_\varepsilon]\cdot[f]\approx\left[\int_\mathcal{M}K_\varepsilon(x,y)f(y)q(y)\mathrm{d}y\right]=\varepsilon^{d/2}[G_\varepsilon(fq)]
\end{equation}
Depending on one's data, this implicit weighting by the sampling distribution could be seen as a problem; but, supposing we know $q(x)$, then we could define $K_\varepsilon^R(x,y)=K_\varepsilon(x,y)q(y)^{-1}$ as we have done before.  We call this the right-normalised kernel, since letting $D_R\coloneqq \text{Diag}([q])$ denote the diagonal matrix with diagonal entries equalling $q(x_i)$, then $[K_\varepsilon^R]=[K_\varepsilon]D_R^{-1}$ - i.e. the sampling distribution normalisation acts on the right of the discretised $K_\varepsilon$.  In this case we will have $\frac{1}{m}[K_\varepsilon^R]\cdot[f]\approx\varepsilon^{d/2}[G_\varepsilon f]$, and all that is left is to normalise for the mass of the kernel through $D_L\coloneqq\text{Diag}\left(\frac{1}{m}[K_\varepsilon^R]\cdot [\mathbb{1}]\right)\approx\varepsilon^{d/2}[G_\varepsilon \mathbb{1}]$; therefore:
\begin{equation}\label{eq:qregsec_backward}
    [P_\varepsilon f]=\text{Diag}([G_\varepsilon\mathbb{1}])^{-1}[G_\varepsilon f]=\text{Diag}(\varepsilon^{d/2}[G_\varepsilon\mathbb{1}])^{-1}\varepsilon^{d/2}[G_\varepsilon f]\approx D_L^{-1}[K_\varepsilon]D_R^{-1}\cdot[f]\backwardscoloneqq P_{m,\varepsilon}[f]
\end{equation}
where in the final equality we are defining $P_{m,\varepsilon}=\frac{1}{m}[P_\varepsilon]\in\mathbb{R}_{\geq 0}^{m\times m}$, and since $[P_\varepsilon f]\approx \frac{1}{m}[P_\varepsilon]\cdot[f]=P_{m,\varepsilon}[f]=D_L^{-1}[K_\varepsilon]D_R^{-1}\cdot[f]$ is satisfied for any $[f]$ we have that $P_{m,\varepsilon}=D_L^{-1}[K_\varepsilon]D_R^{-1}$ - or written elementwise we have $(P_{m,\varepsilon})_{ij}=\frac{K_\varepsilon(x_i,x_j)}{Z(x_i)q(x_j)}$ where $Z(x_i)=\sum_{j=1}^mK_\varepsilon(x_i,x_j)q(x_j)^{-1}$, which should be recognised as the standard form of a row-normalised transition matrix.  To explain the difference intuitively, as $\varepsilon$ gets smaller $P_\varepsilon(x,\cdot)$ gets sharper and more $\delta$-like as a distribution, and thus its maximum value will begin tending to infinity, but the discrete row-normalised transition probability matrix cannot have entries greater than $1$, since all entries must be positive and the row-sums must equal $1$; division by $m$ is what maps the discretisation of the continuous transition operator to the standard discrete transition matrix.  

If we wished to approximate the adjoint operator $[P_\varepsilon^* f]=[G_\varepsilon^*((G_\varepsilon\mathbb{1})^{-1} f)]$ it should now be easy to see that:
\begin{equation}
    [P_\varepsilon^* f]=\varepsilon^{-d/2}\left [\int_\mathcal{M}q(y)\frac{K_\varepsilon(y,x)}{q(y)}\frac{f(y)}{G_\varepsilon\mathbb{1}(y)}\ \mathrm{d} y\right]\approx \frac{\varepsilon^{-d/2}}{m}D_R^{-1}[K_\varepsilon]^T\cdot\left[\frac{f}{G_\varepsilon\mathbb{1}}\right]=D_R^{-1}[K_\varepsilon]^TD_L^{-1}\cdot[f]\backwardscoloneqq P_{m,\varepsilon}^*[f]
\end{equation}
As we hinted at in previous sections, we can now see explicitly that our choice of normalisation results in the convenient fact that $P_{m,\varepsilon}^*=P_{m,\varepsilon}^T$, which is exactly as we would expect for a Markov transition probability matrix.  As a final note, $L_\varepsilon\coloneqq\varepsilon^{-1}(P_\varepsilon-I)\approx\mathcal{L}$ and therefore:
\begin{equation}
    [\mathcal{L}f]\approx[L_\varepsilon f]=\varepsilon^{-1}[P_\varepsilon f-f]=\varepsilon^{-1}\left([P_\varepsilon f]-[f]\right)\approx \varepsilon^{-1}(P_{m,\varepsilon}-I_m)[f]\backwardscoloneqq L_{m,\varepsilon}[f]
\end{equation}
where $I_m$ is the identity matrix in $\mathbb{R}^{m\times m}$.  Since this is true for any $[f]$ we can write the familiar expression for the discrete Laplacian: $L_{m,\varepsilon}=\varepsilon^{-1}(P_{m,\varepsilon}-I_m)$; as a point of interest we can write this in the less familiar form of $L_{m,\varepsilon}=\tfrac{1}{m\varepsilon}([P_\varepsilon]-m I_m)$, where we notice that the second term $m I_m$ is approximating the Kronecker $\delta$-function as $m\to\infty$ since if $i\neq j$ ($x\neq y$) then $(mI_m)_{ij}=0$ ($\delta(x,y)=0$), but for any $i$ we have $(mI_m)_{ii}=m\to\infty$ ($\delta(x,x)=\infty$).  Otherwise stated, if $\delta_q(x,y)\coloneqq \delta(x,y)q(y)^{-1}$, then $[\delta_q]=m I_m$ and $[L_\varepsilon]=\varepsilon^{-1}[P_\varepsilon-\delta_q]$.  We also note again that $L_{m,\varepsilon}^*=L_{m,\varepsilon}^T$ since, assuming everything exists and is well-defined, the adjoint of the generator will be the generator of the adjoints - i.e. it follows directly from the fact $P_{m,\varepsilon}^*=P_{m,\varepsilon}^T$.\footnote{We have avoided using the notation $[L_\varepsilon]$ as conceptually we view $[L_\varepsilon]$ as being a differential operator and not an integral operator and thus talking about $L_\varepsilon(x,y)$ does not agree with this formulation.  It is not entirely wrong to write $[L_\varepsilon]$ since $L_\varepsilon$ is effectively an integral against the kernel $P_\varepsilon(x,y)-\delta_q(x,y)$.}

Since $P_{m,\varepsilon}$ is a discrete Markov matrix, we can approximate the derivative of some probability distribution evolving under $P_\varepsilon$ and denoted as $p_t$ at time $t$, as $[p_{t+\varepsilon}]=P_{m,\varepsilon}^T[p_t]$, written as a derivative this becomes $\frac{\partial p_t}{\partial t}\approx\frac{[p_{t+\varepsilon}]-[p_t]}{\varepsilon}=\varepsilon^{-1}(P_{m,\varepsilon}^T-I_m)[p_t]=L_{m,\varepsilon}^T[p_t]$, when approximated through an Euler scheme - this is the discrete Fokker-Planck equation (Kolmogorov forward equation) for a Markov process.  We can perform a similar approximation for the Feynman-Kac equation (Kolmogorov backward equation), yielding $\frac{\partial u_t}{\partial t}=-L_{m,\varepsilon}[u_t]$, where the negative sign comes from the convention of time going backwards in the Feynman-Kac formula.

In the literature, particularly in Diffusion Maps \cite{coifman2006diffusion}, they introduce an $\alpha\in[0,1]$ parameter to tune the strength of density regularisation, with $\alpha=0$ being no density regularisation, $\alpha=1/2$ being their \textit{Fokker-Planck} normalisation, and $\alpha=1$ being fully regularised.  This amounts to defining $K_{\varepsilon,\alpha}(x,y)\coloneqq K_\varepsilon(x,y)q(y)^{-\alpha}$ and only makes sense if we change the integrals for the backward and forward integral operators to be with respect to the underlying sampling measure $q(x)$.  In other words, we now consider the integral operator:
\begin{equation}
    G_{\varepsilon,\alpha}f(x)\coloneqq \varepsilon^{-d/2}\int_\mathcal{M}K_{\varepsilon,\alpha}(x,y)f(y)q(y)\ \mathrm{d}y=\varepsilon^{-d/2}\int_\mathcal{M} K_\varepsilon(x,y) f(y)q(y)^{1-\alpha}\ \mathrm{d}y=G_\varepsilon(fq^{1-\alpha})(x)
\end{equation}
The limiting operator of $\alpha$-regularised Laplacian takes the form (see Appendix \ref{appAsec:densreg}, for a derivation):
\begin{equation}
    \mathcal{L}_\alpha f=\mu\cdot\nabla f+2(1-\alpha)D\frac{\nabla q}{q}\cdot\nabla f+D\Delta f
\end{equation}
We see that if one doesn't fully regularise for the sampling density, then there will be another advection term proportional to $\nabla q$ that is weighted by a factor of $(1-\alpha)$.  If one progressively increases $\alpha$ from $0$ to $1$, then the contribution of this additional term decreases.  If one wants an exact reconstruction of the underlying dynamics then one should always consider the $\alpha=1$ case; but if one's data is highly noisy (many outliers) then pushing the effective drift towards the core concentration of the data may be wise depending on your ultimate goal.  The $\alpha$-normalised adjoint has an expansion analogous to the previous results given in Theorem \ref{thm:odk_forward_expansion}.

We can also put this together with the $\beta$-regularisation proposed in section \ref{subsec:selftuning_odk}, this would result in the infinitesimal generator:
\begin{align}
    \mathcal{L}_{\alpha,\beta}f&=\mu\cdot\nabla f+2(1-\beta)(d+2)D\tfrac{\nabla \rho}{\rho}\cdot\nabla f+2(1-\alpha)D\tfrac{\nabla q}{q}\cdot\nabla f+D\Delta f\\
    &=\mu\cdot\nabla f+2D\left[(1-\alpha)-(1-\beta)\tfrac{d+2}{d}\right]\tfrac{\nabla q}{q}\cdot\nabla f+D\Delta f\label{eq:alpha_beta_operator}
\end{align}
where the second equality comes from assuming $\rho\propto q^{-1/d}$, and we can see how the two \textit{effective drifts} (or skewed diffusions) act in opposite directions, i.e. setting $\alpha=1$ pushes you towards denser regions, and setting $\beta=1$ pushes you towards sparser regions.  Hence, ones choice of $\alpha$ and $\beta$ depends on the data itself and the analysis you wish to do: if your data is highly noisy, you should lower $\alpha$; if your data is very sparse, you should lower $\beta$; if your data is neither then unless essential for further analysis, setting $\alpha=\beta=1$ seems most reasonable; and if your data is both highly noisy and very sparse, then from equation \eqref{eq:alpha_beta_operator}, we can see that the seemingly sensible choice of setting $\alpha=\beta=0$ will actually be achieving very little when $d\gg 1$, thus further exploration of the data is required to determine a good parametrisation.

Such results are not novel and have been reproduced many times (for example in \cite{coifman2006diffusion,berry2016variable}), and it just shows how the sampling density - and the resulting normalisations - affect the differential operator your Laplacian approximates.  Recall, the sampling distribution was introduced into our integral as an artifact of how Monte-Carlo integration aggregates under non-uniform sampling distributions.  But note, the results of, for example, equation \eqref{eq:alphanormalised_TP} are still computed with the continuous operators $P_{\varepsilon,\alpha}:L^2(\mathcal{M},q)\cap C^3(\mathcal{M})\to L^2(\mathcal{M},q)\cap C^3(\mathcal{M})$, and the results of our theorems are still viewed continuously; to truly bridge the gap between the discrete operators and the continuous operators, we must show that in the limit of large data\footnote{Whenever we say ``limit of large data" we are implicity taking a limit of $m\to\infty$ where $m=|\mathcal{D}|$.} the action of increasingly large matrices, i.e. $L_{m,\varepsilon}$, on $[f]$ becomes arbitrarily close to the discretisation of the continuous operator, i.e. $[L_\varepsilon f]$, and this should be true independent of any particular sampling of $m$ data points from $q(x)$.  This question has been asked and answered before in the literature, but we state a similar result here for completeness:
\begin{theorem}[Discrete Operator Convergence]\label{thm:conv_in_prob}
    Let $\mathcal{M}$ be a $d$-dimensional Riemannian manifold embedded in $\mathbb{R}^n$ and with metric inherited from the standard Euclidean metric in ambient space.  Suppose $\mathcal{D}=\{x_i\}_{i=1}^m$ is a dataset of $m$ independent samples from the sampling distribution $q(x)$ which we assume satisfies $\inf_{x\in\mathcal{M}}q(x)\backwardscoloneqq c_q>0$, i.e. $\text{supp}(q)=\mathcal{M}$ and $q$ doesn't decay to zero anywhere on $\mathcal{M}$ - this condition necessitates that $\mathcal{M}$ be bounded.  Suppose $K_\varepsilon(x,y)$ is an Ordered Diffusion Kernel and $L_\varepsilon$ its associated Laplacian, with the sampling density fully regularised away ($\alpha=1$), and suppose, after fixing an indexing of $\mathcal{D}$, $L_{m,\varepsilon}\in\mathbb{R}^{m\times m}$ is the discretisation of $L_\varepsilon$ onto $\mathcal{D}$ (i.e. $[L_\varepsilon f]\approx L_{m,\varepsilon}[f]=\frac{1}{m}[L_\varepsilon]\cdot[f]$).  Then for any bounded $f\in L^2(\mathcal{M},q)\cap C^3(\mathcal{M})$, i.e. $\exists C_f>0$ such that $\sup_{x\in\mathcal{M}}|f(x)|=C_f$, and any $\delta>0$:
    \begin{equation}
        P\left(\left\lVert L_{m,\varepsilon}[f]-[L_\varepsilon f]\right\rVert_\infty>\delta\right)\leq2m\exp\left(-\frac{\tfrac{1}{2}m\varepsilon^{d/2+2}\delta^2}{C_f(C_K+\varepsilon^{d/2})(C_f+\tfrac{1}{3}\varepsilon\delta)}\right)\backwardscoloneqq (\dagger)
    \end{equation}
    where $C_K>0$ is a constant that depends on $K_\varepsilon$ and $q$ and is independent of $\mathcal{D}$.
\end{theorem}
\begin{proof}
    See Appendix \ref{appAthm:discrete_conv_proof}.
\end{proof}
The result of this is that for a fixed value of $\varepsilon>0$ and for any $\delta>0$, we will have $(\dagger)=2m\exp(-Cm)$, where $C>0$ groups all terms inside the exponential that do not depend on $m$.  Since $\lim_{m\to\infty} 2me^{-Cm}=0$ for any $C>0$ we have that $L_{m,\varepsilon}f$ converges (pointwise) in probability to $L_\varepsilon f$.  We can also ask how we may decrease $\varepsilon$ as $m$ increases, represented as a function $\varepsilon(m)$ so $L_{m,\varepsilon(m)}$ is converging in probability to the continuous operators $L_{\varepsilon(m)}$, and $L_{\varepsilon}$ is converging asymptotically to $\mathcal{L}$.  For small enough $\varepsilon>0$ we can group the terms of the denominator inside the exponential of $(\dagger)$, ignoring the $\varepsilon$ dependence, since (after expanding $\tfrac{1}{a+b\varepsilon^k}$) the $\varepsilon^{d/2+2}$ term will dominate.  The convergence criteria, for any fixed $\delta>0$, is that $m$ and $\varepsilon$ must satisfy:
\begin{equation}
    \lim_{m\to\infty}\exp(\log (m)-C'm\varepsilon^{d/2+2})=0\iff\lim_{m\to\infty}\log(m)\left(1-C'\frac{m}{\log(m)}\varepsilon^{d/2+2}\right)= -\infty
\end{equation}
so at the very least we must have $\varepsilon\geq (\log(m)/m)^{2/(d+4)}/C'$, but to make this convergence decisive (i.e. independent of $C'$ and $o(1)$) then we must have $m\varepsilon^{d/2+2}/\log(m)\to\infty$, but we also desire $\varepsilon(m)\to 0$ as $m\to\infty$.  Thus for any $g(m)$ such that $\lim_{m\to\infty} g(m)=\infty$, and $g(m)=o\left((m/\log(m))^{2/(d+4)}\right)$, we can write:
\begin{equation}
    \varepsilon(m)=\left(\frac{\log m}{m}\right)^{2/(d+4)}g(m)
\end{equation}
which clearly satisfies $m\varepsilon(m)^{d/2+2}/\log(m)=g(m)\to\infty$, and because $g(m)$ is dominated by $(m/\log(m))^{2/(d+4)}$, we will still have that $\varepsilon(m)\xrightarrow{m\to\infty} 0$.  Thus, for any bounded $f\in L^2(\mathcal{M},q)\cap C^3(\mathcal{M})$, in the limit of large data the discrete operator applied to $f$ will satisfy $L_{m,\varepsilon(m)}f\xrightarrow{m\to\infty} \mathcal{L}f=\mu\cdot\nabla f+D\Delta f$, since as $m\to\infty$ we have $L_{m,\varepsilon(m)}f\xrightarrow{P}L_{\varepsilon(m)}f$ and $\varepsilon(m)\to0$ as $m\to\infty$, implies we will also have $L_{\varepsilon(m)}f\to\mathcal{L}f$.  This completes the story for the convergence of ODK's, applied to discrete data, to the continuous second-order differential operator $\mathcal{L}$.

\section{Numerical Results}\label{sec:numerical_results}
In this section, we will take several idealised examples and demonstrate the limiting behaviour and utility of ODKs when applied to synthetic data. We will begin in subsection \ref{subsec:numerical_results_operators} by mirroring the example given in local kernel's paper \cite{berry2016local} and demonstrating the equivalent convergence for ODKs with both isotropic and non-isotropic diffusions; in subsection \ref{subsec:numerical_results_ode} we will demonstrate how to extract geometric and dynamical features from data sampled from an ODE; and finally, in subsection \ref{subsec:numerical_results_sde}, we will leverage the structure of ODKs to build two simple gradient based optimisation schemes for scenarios where you do, and do not, have access to sampling times; we use these optimisers to infer the ordering and diffusion for the underlying SDE that generated the data.

\subsection{Numerical Example: Laplacian Operator}\label{subsec:numerical_results_operators}
In our first numerical example, we will demonstrate the convergence of the discrete Laplacian to the continuous analytical Laplacians of Theorems \ref{thm:odk_backward_expansion} (isotropic diffusion) and \ref{thm:noniso_odk_expansion} (anisotropic diffusion). We will use the same process as was used in local kernels to demonstrate equivalent convergence \cite{berry2016local}.

\begin{figure}
    \centering
    \includegraphics[width=0.9\textwidth]{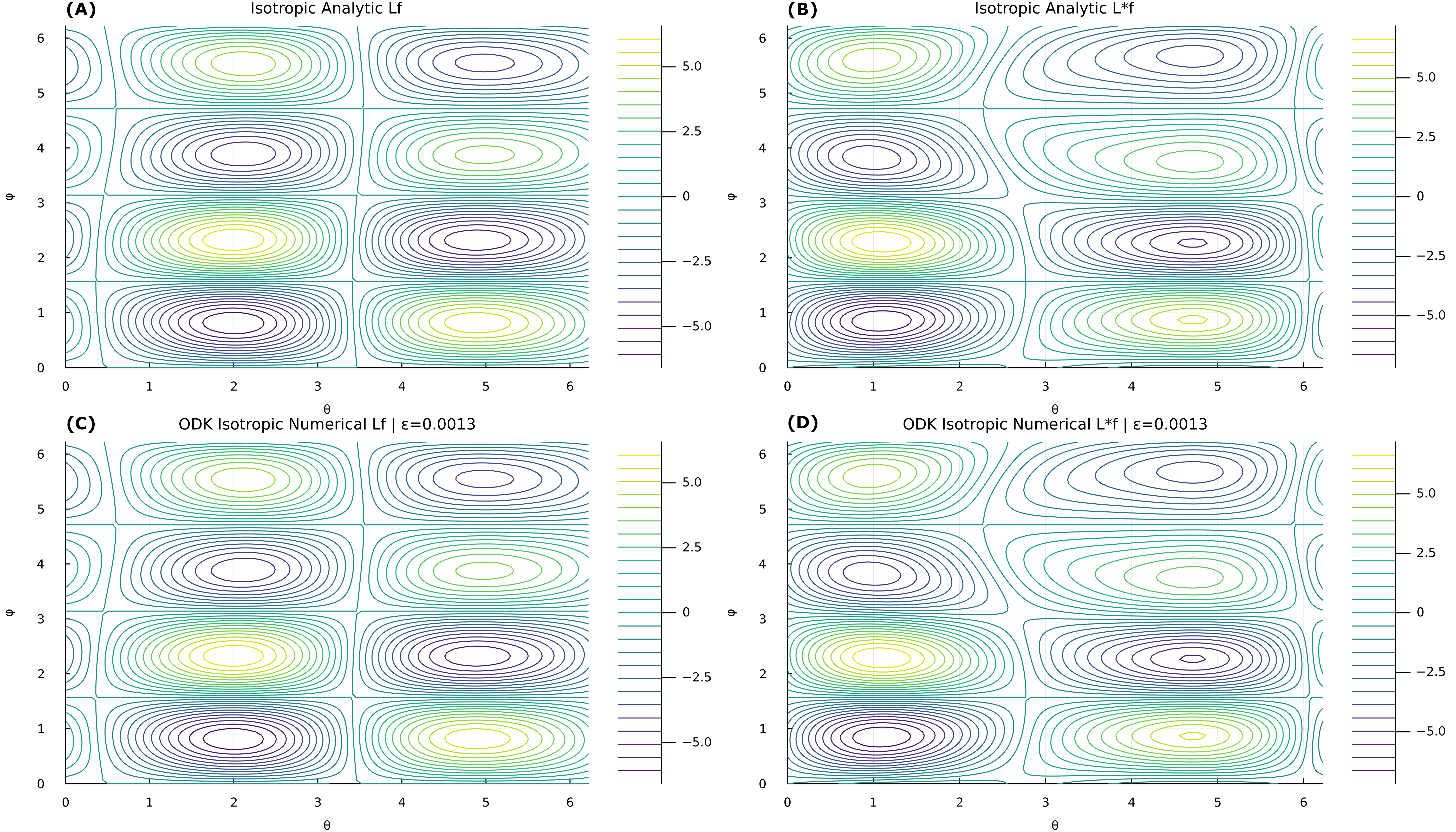}
    \caption{In this figure, we show the results of applying the isotropic ODK, as defined in equation \eqref{eq:numerical_iso_odk}, with $\varepsilon=0.0013$, to $10^4$ data points sampled uniformly on a torus embedded in $\mathbb{R}^4$.  We use a local ordering function $\tilde l(x,y)$, as given in equation \eqref{eq:numerical_local_ordering_function}, and define the isotropic diffusion to be $\tfrac{1}{2}\bar\rho(x)^2\coloneqq=\tfrac{1}{4}\tr(C(x))$, which is the best isotropic approximation to an anisotropic diffusion when viewed through the metric of KL-divergence.  In this figure we apply the Laplacians and their adjoints to the test function $f(\theta,\phi)=\sin(\theta)\sin(2\phi)$: in \textbf{Panel A}, we show a contour plot of $\mathcal{L}f$ as given in equation \eqref{eq:numerical_iso_analytic_laplacian}, and in \textbf{Panel B} we do the same for $\mathcal{L}^* f$ as given in equation \eqref{eq:numerical_iso_analytic_laplacian_adjoint}.  In \textbf{Panel C}, we show the numerical approximation of ODKs $L_\varepsilon f$ to $\mathcal{L}f$, and in \textbf{Panel D}, we do the same for $L_\varepsilon^* f$.  Broadly, we see very little deviation between the rows of this figure with the $L^2$-norms between the analytic and numerical approximations being $\lVert\mathcal{L}f-Lf\rVert_{L^2}=0.00791$ and $\lVert\mathcal{L}^*f-L^*f\rVert_{L^2}=0.01032$, respectively.}
    \label{fig:isotropic_laplacian_results}
\end{figure}

Our dataset will be constructed as a uniform mesh of a $2$-torus, denoted $T^2$, embedded in $\mathbb{R}^4$: to do this we first construct a uniform grid of 10,000 points $x_i=(\theta_i,\phi_i)\in[0,2\pi)^2\subset[0,2\pi]^2/\sim\backwardscoloneqq T^2$ --- this is our intrinsic manifold which we take to have the Euclidean metric;\footnote{$[0,2\pi]^2/\sim$ denotes the torus constructed by gluing opposite sides of a square in $\mathbb{R}^2$ together, we first represent this as a quotient space $[0,2\pi]^2/\sim$.} we then embed these points into $\mathbb{R}^4$ through the isometric embeddeding:
\begin{equation}\label{eq:operator_results_embedding}
 \iota(\theta,\phi)\mapsto(\sin(\theta),\cos(\theta),\sin(\phi),\cos(\phi)),\qquad   D\iota(\theta,\phi)=\begin{pmatrix}
        \cos(\theta) & -\sin(\theta) & 0 & 0\\
        0 & 0 & \cos(\phi) & -\sin(\phi)
    \end{pmatrix}^T.
\end{equation}
Since $\mathcal{M}\coloneqq\iota(T^2)$ is assumed to inherit the Euclidean metric on $\mathbb{R}^4$, isometry is equivalent to $D\iota(x)^TD\iota(x)=I_2$ for all $x\in[0,2\pi)^2$, this can be verified easily from equation \eqref{eq:operator_results_embedding}. We desire an isometric embedding so that the uniformity of our grid on $[0,2\pi)^2$ is preserved on the embedded torus in $\mathbb{R}^4$, this will simplify the analysis. In \cite{berry2016consistent} they define the following dynamics:
\begin{equation}\label{eq:numericalresults_laplacian_dynamics}
    \mu(\theta,\phi)=\begin{pmatrix}
        2+\sin(\theta)\\
        0
    \end{pmatrix},
    \qquad C(\theta,\phi)=\begin{pmatrix}
        3+\sin(\phi) & 1\\
        1 & 1
    \end{pmatrix}
\end{equation}
To replicate these dynamics using an ordering function, we first treat $l$ like a potential and define the global ordering function $l(\theta,\phi)=2\theta-\cos(\theta)$, which naturally satisfies $\nabla l(x)=\mu(x)$; however, since the torus has non-trivial homotopy, and the orbits under $\mu(\theta,\phi)$ form cycles, a single global ordering function will be insufficient due to the branch cut at $\theta=2\pi$.  Thus, we adapt the global $l(\theta,\phi)$ into a local ordering function $\tilde{l}(x_i,x_j)$ through:
\begin{equation}\label{eq:numerical_local_ordering_function}
    \tilde{l}(x_i,x_j)=\begin{cases}
        l(\theta_j,\phi_j)-l(\theta_i,\phi_i) & \text{if } |\theta_j-\theta_i|\leq\pi\\
        l(\theta_j-2\pi,\phi_j)-l(\theta_i,\phi_i) & \text{if } \theta_j>\theta_i+\pi\\
        l(\theta_j+2\pi,\phi_j)-l(\theta_i,\phi_i) & \text{if } \theta_j<\theta_i-\pi
    \end{cases}
\end{equation}
We will first consider the standard isotropic ODK with $\bar\rho^2\coloneqq\frac{1}{2}\tr(C(x))=\frac{1}{2}(4+\sin(\phi))$ --- this is the best possible isotropic approximation to $C(x)$ as measured by the KL-divergence.  Since we are treating $l(\theta,\phi)$ like a potential, we set $\tau(x)=\kappa(x)\equiv 1$, after taking $h(u)=\exp(-u/2)$ our ODK can be written as:
\begin{equation}\label{eq:numerical_iso_odk}
    K_\varepsilon(x_i,x_j)=\exp\left(-\frac{\lVert\iota(x_i)-\iota(x_j)\rVert^2}{2\varepsilon\bar\rho(x)^2}-\frac{(\tilde{l}(x_i,x_j)-\varepsilon)^2}{2\varepsilon\bar\rho(x)^2}+\frac{\tilde{l}(x_i,x_j)^2}{2\varepsilon\bar\rho(x)^2}\right)
\end{equation}
By the results of Theorems \ref{thm:odk_backward_expansion} and \ref{thm:odk_forward_expansion}, we would expect the Laplacian of these kernels to approach:
\begin{align}
    \mathcal{L}f(x)&=\mu(x)\cdot\nabla f(x)+\frac{1}{2}\bar\rho(x)^2\Delta f(x)=(2+\sin(\theta))\frac{\partial f}{\partial \theta}(x)+\frac{1}{4}(4+\sin(\phi))\left(\frac{\partial^2 f}{\partial\theta^2}(x)+\frac{\partial^2 f}{\partial\phi^2}(x)\right)\label{eq:numerical_iso_analytic_laplacian}\\
    \mathcal{L}^*f(x)&=-\diver(\mu(x)f(x))+\frac{1}{2}\Delta(\bar\rho^2f(x))=\frac{\partial}{\partial\theta}((2+\sin(\theta))f(x))+\frac{1}{4}\Delta((4+\sin(\phi))f(x))\label{eq:numerical_iso_analytic_laplacian_adjoint}
\end{align}
We take $f(\theta,\phi)=\sin(\theta)\sin(2\phi)$ to be our test function, and we found numerically that $\varepsilon\approx 0.0013$ was the optimal value of $\varepsilon$ (i.e. the smallest value of $\varepsilon$ before the kernel begins losing resolution on the data).  In Figure \ref{fig:isotropic_laplacian_results}, we show contour plots of $\mathcal{L}f$ (Panel A), $\mathcal{L}^*f$ (Panel B), $Lf$ (Panel C) and $L^*f$ (Panel D), where $L$ and $L^*=L^T$ denote ODKs discrete Laplacian and its adjoint, respectively.  We can immediately see that there is very little discernible difference at all between the analytic and the numerical results, with the $L^2$ norm between each being $\lVert\mathcal{L}f-Lf\rVert_{L^2}=0.00791$ and $\lVert\mathcal{L}^*f-L^*f\rVert_{L^2}=0.01032$, respectively.

For our second example, we will apply the anisotropic extension of ODK to the same data using $\rho(x)^2=|C(x)|$ and $A(x)=|C(x)|^{-1}C(x)$, where $C(x)$ is given in equation \eqref{eq:numericalresults_laplacian_dynamics}.  From section \ref{subsec:nonisotropic_extension} we know that to approximate the operators:
\begin{align}
    \mathcal{L}_Af(x)&=\mu(x)\cdot\nabla f(x)+\tfrac{1}{2}\rho(x)^2\sum_{i,j\in\{\theta,\phi\}}A^{ij}(x)\partial_i\partial_jf(x)\\
    \mathcal{L}_A^*f(x)&=\diver(\mu(x)f(x))+\tfrac{1}{2}\sum_{i,j\in\{\theta,\phi\}}\partial_i\partial_j(\rho(x)^2A^{ij}(x)f(x)).
\end{align}
We construct two normalised transition matrices $P_{\varepsilon,\alpha}^\mu$ and $P_{\varepsilon/2,\alpha}^A$: where $P_{\varepsilon,\alpha}^\mu$ is an isotropic ODK with drift $\mu(x)$ and diffusion $\tilde{\rho}(x)^2\coloneqq\alpha\rho(x)^2\lambda_\text{min}(A(x))$; and $P_{\varepsilon/2,\alpha}^A$ is a symmetric local kernel with diffusion $\rho^2(x)\tilde{A}(x)\coloneqq\rho(x)^2A(x)-\tilde{\rho}(x)^2I$.  We then apply our proposed Strang Splitting scheme to approximate the corresponding normalised transition probability matrix, $\bar P_\varepsilon^A$, using the following order of composition $\bar P_\varepsilon^A\coloneqq P_{\varepsilon/2,\alpha}^AP_{\varepsilon,\alpha}^\mu P_{\varepsilon/2,\alpha}^A$ --- we found numerically that this order of operations was more stable to the choice of $\alpha$ and converged to the true Laplacian $\mathcal{L}_Af$ (and its adjoint $\mathcal{L}_A^* f$) at a slightly faster rate.  We found numerically that the optimal value of $\alpha$ in this example was approximately $\tfrac{2}{3}$.  We also note that the optimal value of $\varepsilon$ increases in this example $\varepsilon\approx 0.0063$, which makes sense since the minimal isotropic diffusion of $P_{\varepsilon,2/3}^\mu$, has single eigenvalue of magnitude $\tfrac{2}{3}\lambda_\text{min}(C)$, whereas the isotropic diffusion of the ODK defined in equation \eqref{eq:numerical_iso_odk} has a single eigenvalue of magnitude $\tfrac{1}{2}(\lambda_\text{max}(C)+\lambda_\text{min}(C))$ and the mean ratio between the two takes the form:
\begin{equation}
    \mathbb{E}_{x\sim q}\left[\frac{\tfrac{1}{2}(\lambda_\text{max}(C)+\lambda_\text{min}(C))}{\tfrac{2}{3}\lambda_\text{min}(C)}\right]=\frac{3}{4}\mathbb{E}\left[\frac{\lambda_\text{max}(C)}{\lambda_\text{min}(C)}\right]-\frac{3}{4}\approx 5.36
\end{equation}
and $5.36\times 0.0013\approx0.00697$ is close to the optimal value of $\varepsilon$ found numerically in the anisotropic example.

\begin{figure}
    \centering
    \includegraphics[width=0.9\textwidth]{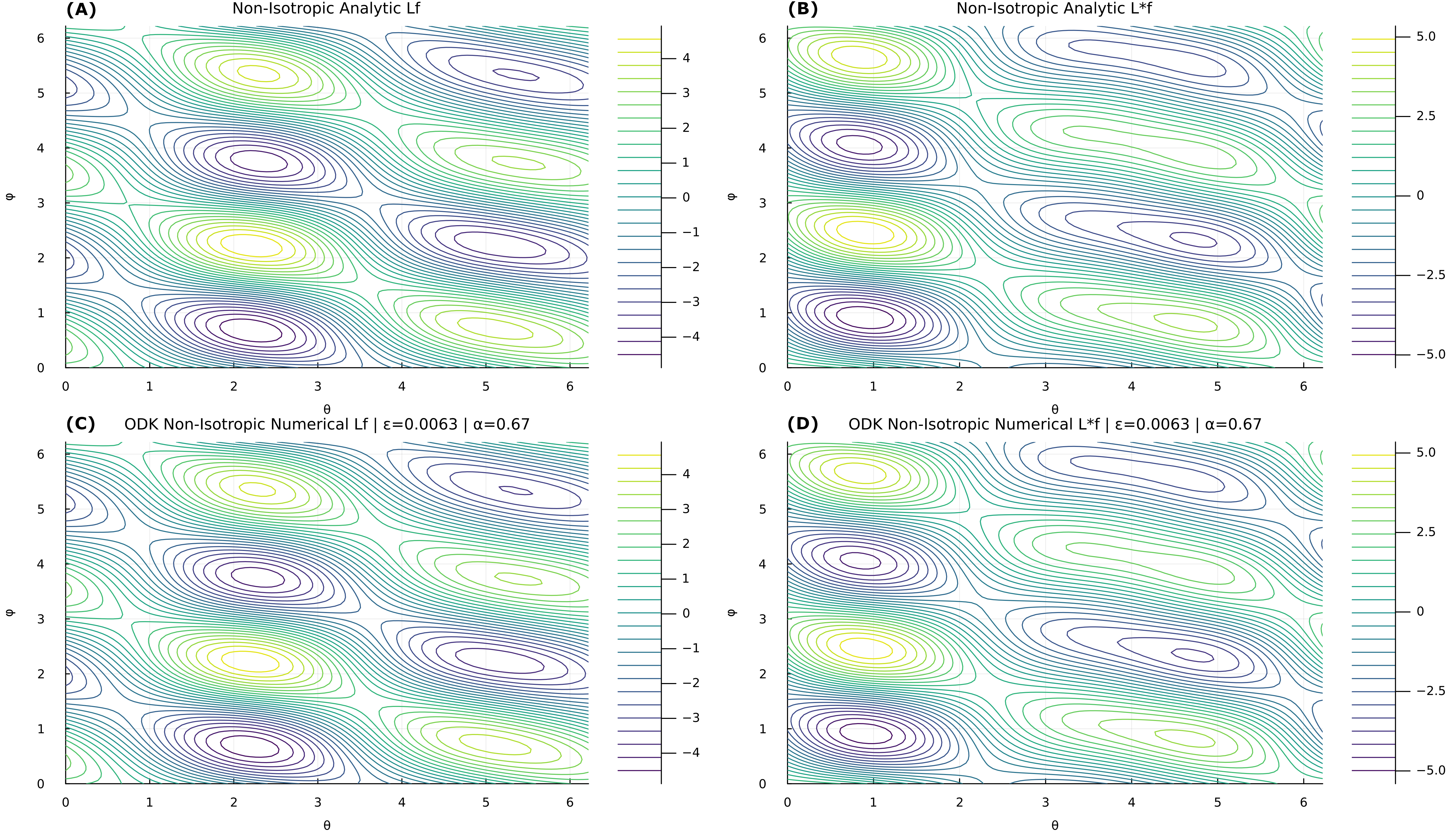}
    \caption{In this figure, we show the results of applying an anisotropic ODK, as discussed in section \ref{subsec:nonisotropic_extension}, with $\varepsilon=0.0067$, to $10^4$ data points sampled uniformly on a torus embedded in $\mathbb{R}^4$.  We use a local ordering function $\tilde l(x,y)$, as given in equation \eqref{eq:numerical_local_ordering_function}, and seek to approximate the diffusion $\tfrac{1}{2}C(x)$, where $C(x)$ is given in equation \eqref{eq:numericalresults_laplacian_dynamics}.  We set $\rho(x)^2=|C(x)|$, and $A(x)=|C(x)|^{-1}C(x)$, we then approximate the normalised anisotropic transition matrix through the Strang Splitting scheme with the following order of operations $\bar P_{\varepsilon,\alpha}^A=P_{\varepsilon/2,\alpha}^AP_{\varepsilon,\alpha}^\mu P_{\varepsilon/2,\alpha}^A$, where we take $\alpha=\tfrac{2}{3}$.  We then apply the associated Laplacians to the test function $f(\theta,\phi)=\sin(\theta)\sin(2\phi)$.  In \textbf{Panel A}, we show a contour plot of $\mathcal{L}_A f$, and in \textbf{Panel B}, we do the same for $\mathcal{L}^* f$.  In \textbf{Panel C}, we show the numerical approximation of ODKs $L_{\varepsilon,\alpha} f$ to $\mathcal{L}_A f$, and in \textbf{Panel D}, we do the same for $L_{\varepsilon,\alpha}^*$.  Broadly, we see very little deviation between the rows of this figure, with $L^2$-norms between the analytic and numerical approximations being $\lVert\mathcal{L}_Af-L_Af\rVert_{L^2}=0.0305$ and $\lVert\mathcal{L}_A^*f-L_A^*f\rVert_{L^2}=0.03518$.}
    \label{fig:nonisotropic_laplacian_results}
\end{figure}

In Figure \ref{fig:nonisotropic_laplacian_results}, we show contour plots of $\mathcal{L}_Af$ (Panel A), $\mathcal{L}_A^*f$ (Panel B), $L_Af$ (Panel C) and $L_A^*f$ (Panel D), where $L_A$ and $L_A^*=L_A^T$ denote the non-isotropic ODKs discrete Laplacian and adjoint, respectively.  We can again immediately see that there is very little discernible difference between the analytic and numerical results, with the $L^2$ norms being $\lVert\mathcal{L}_Af-L_Af\rVert_{L^2}=0.0305$ and $\lVert\mathcal{L}_A^*-L_A^*\rVert_{L^2}=0.03518$, where the slight increase in the non-isotropic $L^2$ differences from the isotropic $L^2$ differences can likely be attributed to an increased optimal $\varepsilon$.

\subsection{Numerical Example: ODEs}\label{subsec:numerical_results_ode}

We will now consider a more complex example and demonstrate how ODKs can be used to extract dynamic and geometric features from data sampled from an ODE.  Our intrinsic manifold will be the square $\mathcal{M}\coloneqq[1,2]\times[0,1]$ and our embedded manifold will be a 2-dimensional spiral in $\mathbb{R}^3$, formed from the embedding:
\begin{equation}\label{eq:vf_embedding}
    \iota(u_1,u_2)=\begin{pmatrix}
        u_1\cos(2\pi u_1)\\
        u_2\\
        u_1\sin(2\pi u_1)
    \end{pmatrix}
    ,\qquad D\iota(u_1,u_2)=\begin{pmatrix}
        \cos(2\pi u_1)-2\pi u_1\sin(2\pi u_1)&0\\
        0&1\\
        \sin(2\pi u_1)+2\pi u_1\cos(2\pi u_1)&0
    \end{pmatrix}
\end{equation}
where we will denote $\mathcal{N}\coloneqq\iota(\mathcal{M})$ for clarity.  When constructing our dataset $\mathcal{D}\subset\mathcal{N}$ we first sample a $100\times100$ grid in $[1,2]\times[0,1]$ in such a way that our data $\mathcal{D}$ is approximately a uniform mesh of the embedded spiral $\mathcal{M}$: we achieve this be sampling $\tilde{u}_1\in[2,4]$ uniformly, and then we define $u_1\coloneqq\sqrt{\tilde{u}_1}$ --- this is not an exact constraint of uniformity but rather a good enough approximation.  

Since $\mathcal{N}$ will inherit the Euclidean metric, we note that $\iota:(\mathcal{M},I)\to(\mathcal{N},I)$ is not isometric; so to compute anything intrinsically we must do so with respect to the metric $g_\mathcal{M}(u)\coloneqq \iota^*g_\mathcal{N}(u)=D\iota(u)^TD\iota(u)=\diag(1+4\pi^2u_1^2,1)$ --- i.e. the pullback of the metric on $\mathcal{N}$ --- which guarantees $\iota:(\mathcal{M},g_\mathcal{M})\to(\mathcal{N},I)$ is an isometry. Hence, the intrinsic manifold will be $(\mathcal{M},D\iota^TD\iota)$ and all distances and derivatives computed intrinsically must be done so with respect to this metric; for example, if $f:\mathcal{M}\to\mathbb{R}$ and $v(u)\in T\mathcal{M}$, then:
\begin{align}
    \lVert v(u)\rVert^2&=v(u)^TD\iota(u)^TD\iota(u)v(u)=(1+4\pi^2u_1^2)v_1(u)^2+v_2(u)^2,\\
    \nabla_\mathcal{M}f(u)&\coloneqq\text{grad}_\mathcal{M}(f)(u)=(D\iota(u)^TD\iota(u))^{-1}\nabla f=
    \frac{\partial_1f(u)}{1+4\pi^2u_1^2}e_1+\partial_2 f(u)e_2.
\end{align}
where $e_1,e_2$ are the standard Euclidean basis vectors on $\mathcal{M}\subset\mathbb{R}^2$.  We use $\nabla f$ to denote the derivative of $f$ with respect to the metric $g=I_2$ (i.e. the derivative of $f$ on $(\mathcal{M},I)$), and $\nabla_\mathcal{M}f\coloneqq\text{grad}_{g_\mathcal{M}}(f)=(D\iota^TD\iota)^{-1}\nabla f$ the derivative of $f$ with respect to the metric $g_\mathcal{M}=D\iota^TD\iota$ (i.e. the derivative of $f$ on $(\mathcal{M},g_\mathcal{M})$).  If we view $\tilde{f}=f\circ\iota^{-1}$ as being $f$ defined on $\mathcal{N}$ then the extrinsic derivative of $f$ restricted to $\mathcal{N}$, which we denote as $\nabla_\text{ext} \tilde{f}$ for clarity, will satisfy $\nabla_\text{ext}\tilde{f}=(D\iota^{-1})^T\nabla f=D\iota(D\iota^TD\iota)^{-1}\nabla f=D\iota\nabla_\mathcal{M} f$, where in the second equality we noted that since $D(\iota^{-1}\circ\iota)=D\iota^{-1}D\iota=I_2$ this implies $D\iota^{-1}=(D\iota)^+=(D\iota^TD\iota)^{-1}D\iota^T$ is a left pseudo-inverse of $D\iota$.  Furthermore, the magnitude of these two gradients with respect to each of their metrics will also be equal, i.e.:
\begin{equation}\label{eq:vf_metric_consistency}
    \lVert\nabla_\text{ext}\tilde f\rVert_2^2=\nabla_\mathcal{M}f^TD\iota^TD\iota\nabla_\mathcal{M}f=g_\mathcal{M}(\nabla_\mathcal{M}f,\nabla_\mathcal{M}f)=\lVert\nabla_\mathcal{M} f\rVert_{g_\mathcal{M}}^2,
\end{equation}
This is all to say that in practice we compute extrinsic derivatives on $(\mathcal{N},I)$, but these are identical (after pushing forward through $D\iota$) to their respective intrinsic derivatives on $(\mathcal{M},g_\mathcal{M})$. Let us define $\hat t_{||}(u)\propto D\iota(u)e_1$ to be the unit vector parallel to the spiral but perpendicular to the y-axis, and $\hat t_\perp(u)\propto D\iota(u)e_1\times D\iota(u)e_2$ to be the inwards-facing unit normal vector to $\mathcal{N}$.

We imagine that on the intrinsic manifold, $\mathcal{M}$, the data is generated by evolving the following ODE from the initial condition set $\mathcal{M}_0\coloneqq\{1\}\times[0,1]$:
\begin{equation}\label{eq:vf_intrinsic_ode}
    \frac{du}{dt}=\begin{pmatrix}
        u_1(2-u_1)\\
        0
    \end{pmatrix}\backwardscoloneqq v_\text{in}(u)
\end{equation}
We also define $s_\text{in}(u)=\lVert v_\text{in}(u)\rVert_2=u_1(2-u_1)$. This is equivalent to the spiral being generated in $\mathbb{R}^3$ by evolving the ODE:
\begin{equation}\label{eq:vf_embedded_ode}
    \frac{dx}{dt}=D\iota(\iota^{-1}(x))v_\text{in}(\iota^{-1}(x))=u_1(2-u_1)(1+4\pi^2u_1^2)^{1/2}\ \hat t_{||}(x)\backwardscoloneqq v_\text{ext}(x),
\end{equation}
where for notational convenience we use $u_1=(\iota^{-1}(x))_1$ and $\hat t_{||}(x)=\hat t_{||}(\iota^{-1}(x))$; we also define $s_\text{ext}(x)=\lVert v_\text{ext}(x)\rVert_2=u_1(2-u_1)\sqrt{1+4\pi^2u_1^2}$.  For a given point $u\in\mathcal{M}$, the timing function $t(u)$ satisfies $1+\int_0^{t(u)}\left[v_\text{in}(u(t'))\right]_1 dt'=u_1$ and can be solved by separation of variables to be:
\begin{equation}\label{eq:intrinsic_timing_function}
    t(u)=\frac{1}{2}\log\left(\frac{u_1}{2-u_1}\right).
\end{equation}
If we generated the data directly on the spiral, then solving the timing function with respect to the ODE given in equation \eqref{eq:vf_embedded_ode} would very naturally result in an identical timing function; thus, there is no ambiguity between whether we defined $t(u)$ intrinsically or extrinsically.  The gradient of the timing function is naturally related to the underlying dynamics through: 
\begin{equation}
    \nabla_\mathcal{M}t(u)=(D\iota^TD\iota)^{-1}\nabla t(u)=\left(\frac{1}{1+4\pi^2u_1^2}e_1e_1^T+e_2e_2^T\right)\frac{e_1}{s_\text{in}(x)}=\frac{1}{s_\text{in}(u)(1+4\pi^2u_1^2)}e_1
\end{equation}
therefore: $\nabla_\text{ext}t(x)=D\iota\nabla_\mathcal{M} t=\left[s_\text{in}(u)^2(1+4\pi^2u_1^2)\right]^{-1/2}\hat t_{||}(x)$, and $\lVert\nabla_\text{ext}t(x)\rVert_2^2=\left[s_\text{in}(u)^2(1+4\pi^2u_1^2)\right]^{-1}$, so in particular the gradient of the timing function is related to the underlying velocity through:
\begin{equation}
    \frac{\nabla_\text{ext}t(x)}{\lVert\nabla_\text{ext}t(x)\rVert^2}=\frac{s_\text{in}(u)^2(1+4\pi^2u_1^2)}{s_\text{in}(u)(1+4\pi^2u_1^2)^{1/2}}\hat t_{||}(x)=s_\text{in}(u)(1+4\pi^2u_1^2)^{1/2}\hat t_{||}(x)=v_\text{ext}(x)
\end{equation}
Such a convenient relationship will hold if and only if $v_\text{in}\lVert v_\text{in}\rVert^2$ is curl-free.  Thus, letting $\tau(x)\equiv 1$, $l(x)=t(x)$, and $\kappa(x)=\lVert\nabla_\text{ext} t(x)\rVert^2$, allows us to define the ODK:
\begin{equation}\label{eq:vf_time_odk}
    K_\varepsilon(u,v)=\exp\left(-\frac{\lVert\iota(u)-\iota(v)\rVert^2}{2\varepsilon\rho^2}-\frac{(t(v)-t(u)-\varepsilon\lVert\nabla_\text{ext} t(\iota(u))\rVert^2)^2}{2\varepsilon\rho^2\lVert\nabla_\text{ext} t(\iota(u))\rVert^4}+\frac{(t(v)-t(u))^2}{2\varepsilon\rho^2\lVert\nabla_\text{ext} t(\iota(u))\rVert^4}\right)
\end{equation}
by Theorem \ref{thm:odk_backward_expansion} we know the associated drift will equal $\mu(x)=\tfrac{\tau(x)}{\kappa(x)}\nabla_\text{ext} l(x)=v_\text{ext}(x)$.  We would now like to extract this tangent vector from the kernel: since we can think of the normalised backward operator, $P_\varepsilon$, as a transition probability matrix, it would seem natural to define a pointwise velocity estimate as the mean displacement under the measure $P_\varepsilon(u,\cdot)$, i.e.:
\begin{equation}\label{eq:vf_def_velocity_computation}
    \check{v}_\text{md}(u)\coloneqq\frac{1}{\varepsilon}\sum_{v\in\mathcal{D}}P_\varepsilon(u,v)(v-u)
\end{equation}
letting $\iota_i(u)$ denote the i-th element of the vector $\iota(u)\in\mathbb{R}^3$, we can rewrite equation \ref{eq:vf_def_velocity_computation} as:
\begin{equation}
    \check{v}_\text{md}(u)_i=\frac{1}{\varepsilon}\sum_{v\in\mathcal{D}}P_\varepsilon(u,v)(\iota_i(v)-\iota_i(u))=\frac{1}{\varepsilon}\sum_{v\in\mathcal{D}}P_\varepsilon(u,v)\iota_i(v)-\frac{1}{\varepsilon}\iota_i(u)=\frac{(P_\varepsilon-I)}{\varepsilon}\iota_i(u)=L_\varepsilon \iota_i(u).
\end{equation}
If for a multivariate function $F(x)=(f_1(x),\dots,f_n(x))^T$, we define $L_\varepsilon F(x)\coloneqq(L_\varepsilon f_1(x),\dots,L_\varepsilon f_n(x))^T$, then the formula for the mean displacement becomes:
\begin{equation}\label{eq:vf_meandisplacement_expansion}
    \check{v}_\text{md}(u)=L_\varepsilon \iota(u)=g_\mathcal{M}(\mu(u),\nabla_\mathcal{M}\iota(u))+\frac{1}{2}\rho^2\Delta_\mathcal{M}\iota(u)+\bigO(\varepsilon^{1/2})=D\iota(u)\mu(u)+\frac{1}{2}\rho^2H(u)+\bigO(\varepsilon^{1/2})
\end{equation}
where $H(u)$ is the Mean Curvature Vector (MCV) of the spiral, and encodes the curvature introduced through the embedding: this term appears since the Laplace-Beltrami operator applied coordinate wise to an embedding into a Euclidean space satisfies $\Delta_\mathcal{M}\iota(u)=\tr_{g_\mathcal{M}}\Pi(u)\backwardscoloneqq H(u)$ where $\Pi(u)$ is the second fundamental form (the proof of this is given on Pg. 41 of \cite{chen2011pseudo}) --- by definition of the second-fundamental form, $H(u)$ will be normal to the embedded manifold $\mathcal{M}$.  As we've already discussed, the first-order term in this expansion will equal $v_\text{ext}(\iota(u))$, and thus, in the limit, the mean displacement estimate of the velocity will equal:
\begin{equation}
    \lim_{\varepsilon\to0}L_\varepsilon\iota(u)=v_\text{ext}(\iota(u))+\frac{1}{2}\rho^2H(u)
\end{equation}
The first term is the embedded velocity field we want to estimate, and the second term is perpendicular to this velocity and proportional to the MCV.  To untangle the contribution of each component, we define the following two improved estimates of the embedded velocity field:
\begin{align}
    \check{v}_\text{sym}(u)&=\left(L_\varepsilon-L_\varepsilon^\text{sym}\right)\iota(u)\label{eq:symcorrected_velocity}\\
    \check{v}_\text{proj}(u)&=\left(I-\hat N(u)\hat N(u)^T\right)L_\varepsilon\iota(u)\label{eq:projcorrected_velocity}
\end{align}
where the first estimate subtracts a symmetric local kernel (which should use the same kernel wrapper $h(\cdot)$): we do this because $\lim_{\varepsilon\to 0}L_\varepsilon^\text{sym}\iota(u)=\frac{1}{2}\rho^2H(u)$ and thus $\lim_{\varepsilon\to 0}\check{v}_\text{sym}=v_\text{ext}(u)$. In the second estimate, the columns of the matrix $\hat N(u)$ form an orthonormal basis of $N_{\iota(u)}\mathcal{M}\cong\mathbb{R}^{n-d}$, if such a matrix can be easily computed given ones data --- such as for our $2$-dimensional spiral embedded in $\mathbb{R}^3$ in which $\hat N(u)=\hat t_\perp (u)$ --- then $\check{v}_\text{proj}(u)$ will remove the MCV without the need to evaluate another kernel; the issue is when $d$ is unknown and whatever it is, we will have $n\gg d$, this may result in us needing to estimate the intrinsic dimensionality of $\mathcal{M}$, which in general we do not wish to do, and the extraction of a large basis of the normal bundle.  Due to this latter consideration, we favour the use of the \text{symmetric corrected} velocity estimate given in \eqref{eq:symcorrected_velocity}, over the \text{projection corrected} velocity estimate.

If instead of the extrinsic velocity field, one wished to estimate the intrinsic velocity field, $v_\text{in}(u)$, then one would need to have knowledge of the embedding itself since $v_\text{ext}(\iota(u))=D\iota(u)v_\text{in}(u)$, which implies $v_\text{in}(u)=(D\iota^TD\iota)^{-1}D\iota(u) v_\text{ext}(\iota(u))$, again using $D\iota^{-1}=(D\iota)^+$, we could then esimate $v_\text{in}(u)$ through: 
\begin{equation}\label{eq:embeddingaware_intrinsic_vf_est}
    \check{v}_\text{in}(u)=(D\iota^TD\iota)^{-1}D\iota^T L_\varepsilon\iota(u)=(D\iota^TD\iota)^{-1}D\iota^TD\iota\ v_\text{in}(u)+\tfrac{1}{2}\rho^2(D\iota^TD\iota)^{-1}D\iota^TH(u)=v_\text{in}(u)
\end{equation}
where we use the fact that $H(u)$ is normal to $\mathcal{M}$, and thus $D\iota(u)^TH(u)=0$.  If we don't have knowledge of the embedding but we did have access to the intrinsic data $\iota^{-1}(\mathcal{D})$, then we could define a new Laplacian $L_\varepsilon^\text{in}$ on this data, and since we have assumed their exist global coordinates on $\mathcal{M}$, the metric will be the trivial Euclidean metric and thus:
\begin{equation}
    \check{v}_\text{in}(u)_i=L_\varepsilon^\text{in} u_i(u)=v_\text{in}(u)\cdot\nabla u_i(u))+\tfrac{1}{2}\rho^2\Delta u_i(u)=v_\text{in}(u)\cdot e_i=v_\text{in}(u)_i.
\end{equation}

\begin{figure}
    \centering
    \includegraphics[width=0.98\textwidth]{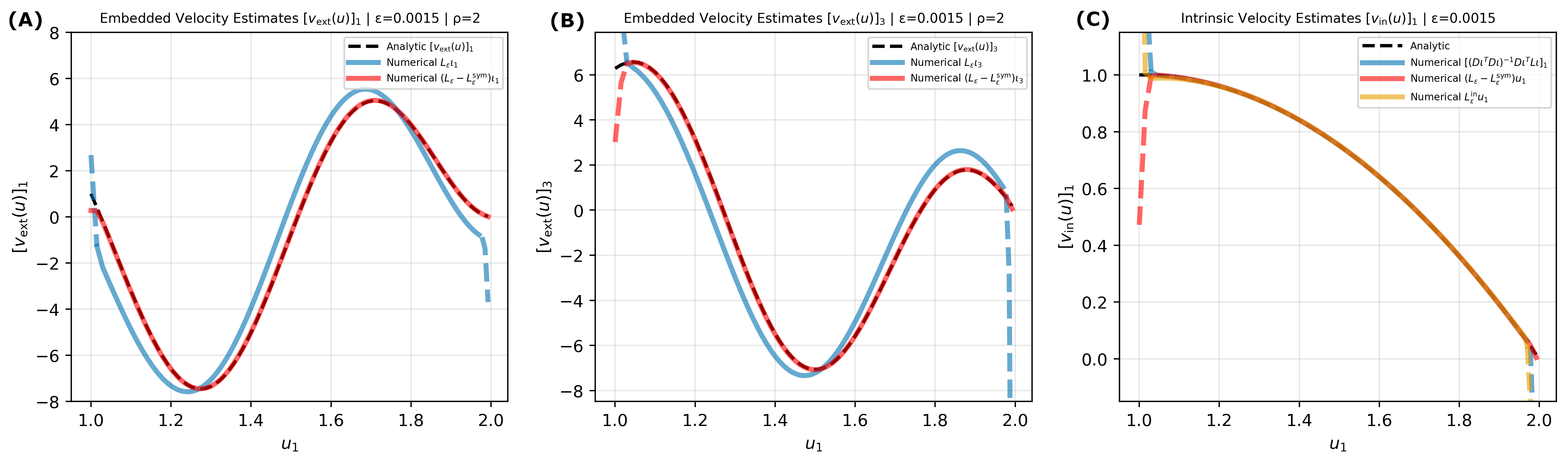}
    
    \caption{In this figure, we show the results of applying an ODK to a lattice of $10^4$ points (sampled uniformly with respect to the metric $D\iota^TD\iota$) non-isometrically embedded in $\mathbb{R}^3$ through $\iota:\mathcal{M}\to\mathbb{R}^3$ as given in equation \eqref{eq:vf_embedding}, to form a $2$-dimensional spiral.  We supposed this data was generated by evolving the ODE $\tfrac{dx}{dt}=v_\text{ext}(x)$ from the initial condition set $\iota(\{0\}\times[0,1])\subset\mathbb{R}^3$, where $v_\text{ext}(x)$ is given in equation \eqref{eq:vf_embedded_ode}.  We parametrised our ODK, which is given in equation \eqref{eq:vf_time_odk}, through: $l(u)=t(u)$, where $t(u)$ is the timing function of the ODE, $\kappa(u)=g_\mathcal{M}(u)(\nabla_\mathcal{M} t(u),\nabla_\mathcal{M} t(u))$ and $\tau(u)\equiv 1$, which amounts to the ODKs drift being equal to $v_\text{ext}(\iota(u))$; and $\rho(x)\equiv2$.  In \textbf{Panel A} and \textbf{Panel B}, we show two estimates of the extrinsic velocity (in the first and third dimensions) obtained from our ODK: the dashed black line shows the analytic components of the velocity; the solid blue line shows the mean displacement estimate, $\check{v}_\text{md}(u)$; and the solid red line shows the symmetric corrected estimate $\check{v}_\text{sym}(u)=(L_\varepsilon-L_\varepsilon^\text{sym})\iota(u)$, where $L_\varepsilon^\text{sym}$ is the Laplacian corresponding to a symmetric local kernel with diffusion $\tfrac{1}{2}\rho^2=2$.  Across both panels, we see that the addition of the normal Mean Curvature Vector (MCV) to $v_\text{ext}(u)$ displaces $\check{v}_\text{md}(u)$ away from the analytic values, whereas $\check{v}_\text{sym}(u)$ corrects for this and is almost entirely aligned with the analytic.  The only discrepancy occurs at the boundaries, where the kernel can no longer resolve the diffusion, but we see that this effect is mitigated in $\check{v}_\text{sym}(u)$ because boundary errors cancel when subtracting the two Laplacians.  In \textbf{Panel C}, we show three estimates of the first component of the intrinsic velocity field $\mu_\text{in}(u)$, as given in equation \eqref{eq:vf_intrinsic_ode}, again the black dashed line shows the analytic expression, the solid blue line shows the embedding aware estimate $D\iota^{-1}(\iota(u))L_\varepsilon\iota(u)$, the orange line shows the mean displacement estimate applied directly to the intrinsic data (using the intrinsic Laplacian $L_\varepsilon^\text{in}$ with $\rho_\text{in}\equiv\frac{2}{l_\text{spiral}}$ where $l_\text{spiral}\approx9.48$ is the arclength of the spiral), and the solid red line shows the symmetric corrected estimate, $(L_\varepsilon-L_\varepsilon^\text{sym})u$, applied directly to the intrinsic data.  All three estimates are highly aligned with the true values, but only the symmetric-corrected estimate can mitigate boundary effects.}
    \label{fig:vf_ode_velests}
\end{figure}

In Figure \ref{fig:vf_ode_velests}, we show the results of applying the ODK given in equation \eqref{eq:vf_time_odk} with $\rho(u)\equiv2$ and $\varepsilon=0.0015$.  In Panels A and B, we show the analytic expression for $\left[v_\text{ext}(u)\right]_i$, for $i=1,3$, respectively, as the dashed black line; the solid blue line shows the mean displacement estimate $\check{v}_\text{md}(u)=L_\varepsilon\iota(u)\approx v_\text{ext}(u)+2H(u)$ which can be seen to capture the rough outline of the embedded velocity field but is clearly offset in both the first and third dimensions by the addition of the normal MCV; the solid red line is our symmetric kernel corrected velocity estimate, $\check{v}_\text{sym}(u)$, as given in equation \eqref{eq:symcorrected_velocity}, and we can see clearly that there is a very high degree of alignement between this estimate and the analytic solution, with the only deviations being at the left boundary (i.e. $u_1\approx 1$).  This latter observation is a realisation of the boundary effects discussed in section \ref{sec:main_results}; as we noted, these effects are mostly unavoidable when the boundary is properly resolved by the data.  However, when comparing $\check{v}_\text{sym}(u)=(L_\varepsilon-L_\varepsilon^\text{sym})\iota(u)$ to $\check{v}_\text{md}(u)=L_\varepsilon\iota(u)$, in Panels A and B, we notice that the observed erros near the boundary for $\check{v}_\text{sym}$ are of a much smaller magnitude that for $\check{v}_\text{md}$: this is because the boundary error is realised in $L_\varepsilon$ and $L_\varepsilon^\text{sym}$ in roughly equal proportions, and thus when we subtract the two Laplacians the two sources of error cancel; the only error that we still observe is when $u_1\approx1$, and in this case $v_\text{ext}(u)$ is large and thus ODKs diffusion resolves better (less error) than the symmetric local kernel because $u+\varepsilon v_\text{ext}(u)$ is further away from the boundary than $u$.  We do not show the estimates of $\check{v}_\text{proj}(u)$ as it is mostly identical to $\check{v}_\text{sym}(u)$, with the only difference being that near the boundary $\check{v}_\text{proj}(u)$ behaved similarly to $\check{v}_\text{md}(u)$, since we do not have the convenient error cancellations as we do with $\check{v}_\text{sym}(u)$.

In Panel C of Figure \ref{fig:vf_ode_velests}, we show 3 estimates of the intrinsic velocity field $v_\text{in}(u)$, where again the dashed black line gives the true values, the solid blue line gives the embedding aware estimate as given in equation \eqref{eq:embeddingaware_intrinsic_vf_est}, the solid red line gives the velocity estimate after applying the symmetric correction directly to the intrinsic data (i.e. $(L_\varepsilon-L_\varepsilon^\text{sym})u_1(u)$), and the orange line shows the result of recomputing the Laplacian directly on the intrinsic data to form $L_\varepsilon^\text{in}$ (we set $\rho_\text{in}=\tfrac{2}{l_\text{spiral}}$ in this construction, where $l_\text{spiral}\approx9.48$ is the arclength of the spiral), and applying it to the coordinate functions.  We see that all three estimates are highly aligned with the analytic solution, with the only divergences occurring at the boundaries: the symmetric estimate is the least affected for the same reason as before.

\begin{figure}
    \centering
    \includegraphics[width=0.98\textwidth]{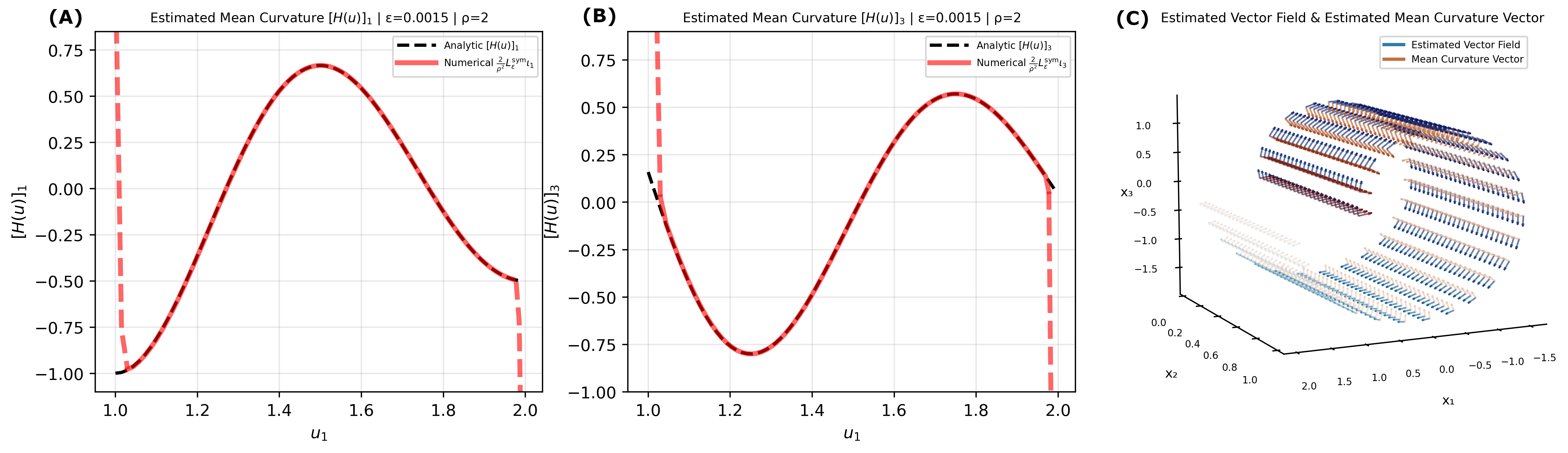}
    \caption{In this figure, we demonstrate how a symmetric local kernel can be used to estimate the Mean Curvature Vector (MCV) of a manifold, $\mathcal{M}$, under an unknown embedding in $\iota:\mathcal{M}\to\mathbb{R}^n$.  We apply a symmetric local kernel $K_\varepsilon^\text{sym}=\exp\left(-(\varepsilon\rho^2)^{-1}\lVert x-y\rVert^2\right)$, with constant diffusion $\tfrac{1}{2}\rho^2\equiv 1$, to a lattice of $10^4$ points sampled in $\mathcal{M}\coloneqq[1,2]\times[0,1]$, non-isometrically embedded in $\mathbb{R}^3$ through $\iota:\mathcal{M}\to\mathbb{R}^3$ as given in equation \eqref{eq:vf_embedding}, we sampled our lattice on $\mathcal{M}$ with respect to the pullback of the Euclidean metric through $\iota$ (i.e. with respect to $D\iota^TD\iota$).  In \textbf{Panel A} and \textbf{Panel B}, the dashed black line shows the analytic MCV in the first and third dimensions, as given in equation \eqref{eq:analytic_mcv}; the solid red line shows the estimated MCV, computed as in equation \eqref{eq:estimated_mcv}.  We see there is very little misalignment between the analytic and estimated MCV, except at the boundaries where the kernel can no longer be resolved sufficiently and the estimate degrades (we show this degradation with dashed red lines).  In \textbf{Panel C} we show a quiver plot evaluated on a $25\times25$ subgrid of $\mathcal{D}$, where arrows with a blue hue show the velocities estimated through $\check{v}_\text{sym}(u)$, as in equation \eqref{eq:symcorrected_velocity}, and arrows with a red hue show the estimated MCVs, $\check{H}(u)$, as in equation \eqref{eq:estimated_mcv}.  The strength of the hue of the arrows is determined by their magnitude, and we note that we have rescaled the arrows by a multiplicative factor of $0.025$ (velocities) and $0.2$ (MCVs) to make their lengths comparable.}
    \label{fig:vf_ode_mcv}
\end{figure}

In Panels A and B of figure \ref{fig:vf_ode_mcv}, the dashed black line shows the analytic MCV in the first and third dimensions, respectively, where the analytic expression of the MCV is computed from the embedding and is given below:
\begin{equation}\label{eq:analytic_mcv}
    H(u)=\frac{4\pi(1+2\pi^2u_1^2)}{(1+4\pi^2u_1^2)^2}
    \begin{pmatrix}
        -\sin(2\pi u_1)-2\pi u_1\cos(2\pi u_1)\\
        0\\
        \cos(2\pi u_1)-2\pi u_1\sin(2\pi u_1)
    \end{pmatrix}
    =\frac{4\pi(1+2\pi^2u_1^2)}{(1+4\pi^2u_1^2)^{3/2}}\hat t_\perp(u)
\end{equation}
The solid red line in both of these figures shows the estimated MCV (in the first and third dimensions) obtained from computing:
\begin{equation}\label{eq:estimated_mcv}
    \check H(u)=\frac{2}{\rho^2}L_\varepsilon^\text{sym}\iota(u)
\end{equation}
From both of these figures, we see that the estimated and analytic forms only diverge from one another near the boundaries, and are otherwise almost perfectly aligned.

In Panel C of figure \ref{fig:vf_ode_mcv}, we show a quiver plot on a $25\times25$ subgrid of points of $\mathcal{D}$: arrows with a blue hue show the estimated velocities (computed using the symmetric correction, $\check{v}_\text{sym}(u)$), and arrows with a red hue show the estimated mean curvature vector computed as in equation \eqref{eq:estimated_mcv}.  The arrows have been rescaled by a factor of $0.025$ and $0.2$, respectively, to make them visually comparable (this choice of scaling reflects the fact that $\mathbb{E}[\lVert H(u)\rVert]\approx\tfrac{1}{8}\mathbb{E}[\lVert\mu_\text{ext}(u)\rVert]$).  The strength of the hue of the arrows is determined by their magnitude, with darker colours corresponding to vectors of greater magnitude.  We can see from this plot that, as expected, the MCV is normal and the velocity is tangent to the manifold.  In panels A and B of \ref{fig:vf_ode_velests}, the discrepancy between $L_\varepsilon\iota(u)$ (the blue line) and the analytic velocity (the dashed black line) can be explained by the addition of the MCV vector to the extrinsic velocity.  We also observe in Panel C of figure \ref{fig:vf_ode_mcv}, that the further along the spiral one goes (increasing $u_1$) the smaller the magnitude of the MCV --- which is highly intuitive since the spiral becomes increasingly flat the further away it is from the origin --- thus, the additional normal component in the estimate $\check{v}_\text{md}(u)$ becomes progressively smaller as one traverse the spiral.

\subsection{Numerical Example: Parametrising SDEs}\label{subsec:numerical_results_sde}
To demonstrate how ODKs might be used on data sampled from an SDE, we will take a more practical approach and define two maximum-likelihood objectives that use the rigid structure of ODKs with respect to their dynamic parameters to set up a gradient-based optimisation for exactly these dynamic parameters.  The choice of which objective to use depends on whether or not you are provided with sampling times.

\subsubsection{SDE with known Time Sampling Distribution}
For this example, we generate our data, $\mathcal{D}$, by sampling from a $d$-dimensional Ornstein-Uhlenbeck (OU) process with independent components and constant isotropic noise; the SDE for this process will take the form:
\begin{equation}\label{eq:sde_ou_process}
    dx_t = -\diag(\boldsymbol{\Theta})x_t\mathrm{dt}+\boldsymbol{\sigma}\mathrm{d}w_t,
\end{equation}
where $\Theta\in\mathbb{R}_{>0}^d$, and $w_t$ is a Wiener process.  We define the time sampling distribution for our data to be a mixture of Gaussians truncated over the interval $[0,T]$ with $T<\infty$, i.e.:
\begin{equation}
    t\sim r(t)\coloneqq\text{Trunc}\left(\frac{1}{s_r}\sum_{i=1}^{s_r}\mathcal{N}(t;\bar t_{i},\sigma_{i});[0,T]\right)
\end{equation}
where $\mathcal{N}(\mu,\sigma)$ is a normal distribution with mean $\mu$ and variance $\sigma^2$, and $\{\bar t_{i}\}_{i=1}^{s_r}\subseteq[0,T]$.  We independently sample from the OU process defined in equation \eqref{eq:sde_ou_process} as discussed in section \ref{subsec:setup_sampling_data}, by first sampling a time $t\sim r(t)$, and then sampling $x\sim \int_\mathcal{M}p_0(x')p_t(x|x_0=x')\mathrm{d}x'$; for both this and the next example we assume $p_0(x)=\delta(x-x_0)$ for a fixed and known $x_0$, although we note that we need not make this restriction. In this example we assume that we are not provided with sampling times $\mathcal{T}\coloneqq\{t_i\}_{i=1}^m$, and instead only have knowledge of the sampled states $\mathcal{D}=\{x_i\}_{i=1}^m$, the time sampling distribution $r(t)$, and the fixed initial condition $x_0$.

Given this setup, there is no pointwise dynamical information from which we could estimate the parameters of the underlying SDE; indeed, we do not even have coarser time marginals from which to compare aggregated statistics such as moments.  However, if we model the SDE as a spatially and temporally discrete Markov process, whose state space is the data $\mathcal{D}$, and whose jumps are over a fixed time-interval $\varepsilon>0$, then we can represent this process through the row-stochastic transition probability matrix $P_\varepsilon(\Phi)$, where $\Phi=(\phi_1,\dots,\phi_{|\phi|})\in\mathbb{R}^{|\phi|}$ is our vector of dynamic parameters,\footnote{When we say dynamic parameters, we mean the scalar functions $\tau(x)$, $\rho(x)$ and $l(x)$ that are inputs to an ODK and which we represent as vectors of their evaluations across $\mathcal{D}$.} this will be the normalised backward operator of Theorem \ref{thm:odk_backward_expansion}, and from Theorem \ref{thm:odk_forward_expansion}, we know that the adjoint $P_\varepsilon^*(\Phi)=P_\varepsilon(\Phi)^T$, so we can represent applying the backwards operator to a function by multiplying a vector of this function's evaluations across $\mathcal{D}$ on the right of $P_\varepsilon(\Phi)$, and likewise multiplication of a transposed vector on the left will represent the application of the forward (Fokker-Planck) operator.  Given this transition operator, we can compute the expected transient distribution given the parameterisation $\Phi$, on data sampled through the time-sampling distribution $r(t)$, through:
\begin{equation}\label{eq:sde_transient_distribution}
    \hat q(x_i|\Phi)\coloneqq\sum_{s=1}^S\bar r_\varepsilon(s)p_0^TP_\varepsilon(\Phi)^se_i,
\end{equation}
where: $\bar r_\varepsilon(s)=r(s\varepsilon)/\sum_{s'=1}^Sr(s'\cdot\varepsilon)$ and $S\coloneqq\left\lceil\tfrac{T}{\varepsilon}\right\rceil$; $p_0\in\mathbb{R}_{\geq 0}^m$ is our discretised initial condition vector which could reasonably be set to a $\delta$-distribution on the closest sample in $\mathcal{D}$ to $x_0$, or maybe a small variance symmetric Gaussian around $x_0$; and $[e_i]_j=\delta_{ij}$ is a vector of zeros at every coordinate except the $i$-th in which it is one.  The sum in equation \eqref{eq:sde_transient_distribution}, should be interpreted as a weighted Monte-Carlo sum and is thus approximating the true sampling distribution integral as given in equation \eqref{eq:sampling_distribution_integral}.  By comparing how similar ODK's estimated $\hat q(x_i|\Phi)$ is to the true sampling distribution, $q(x)$, we can formulate an objective that we can attempt to minimise; using the KL-divergence, $D_\text{KL}(\cdot||\cdot)$, as our measure of \textit{closeness}, we define:
\begin{equation}\label{eq:sde_q_objective}
    L_q(\Phi|\mathcal{D},r,x_0,\lambda_\text{reg})\coloneqq D_\text{KL}(q(x)||\hat q(x|\Phi))+\lambda_\text{reg}R(\Phi|\mathcal{D})\approx-\frac{1}{m}\sum_{i=1}^m\log\hat q(x_i|\Phi)+\lambda_\text{reg}R(\Phi|\mathcal{D})+\text{const}
\end{equation}
where $R(\Phi|\mathcal{D})$ denotes some regularisation function on the parameters $\Phi$, and $\lambda_\text{reg}\geq0$ is a hyperparameter controlling the strength of this regularisation. The approximate equality arises when the integral in the KL-divergence is approximated using Monte-Carlo integration.  The approximate solution of this optimisation problem is then:
\begin{equation}
    \Phi^*\coloneqq\argmin_\Phi-\frac{1}{m}\sum_{i=1}^m\log\hat q(x_i|\Phi)+\lambda_\text{reg}R(\Phi)\approx\argmin_{\Phi}L_q(\Phi|\mathcal{D},r,x_0,\lambda_\text{reg})
\end{equation}
The most general and non-parametric approach to parameterise our ODKs is to learn the functions that determine our kernels as vectors of those functions evaluated across $\mathcal{D}$, i.e., we wish to learn $[f]\in\mathbb{R}^m$ where $[f]_i=f(x_i)$ and $f=\tau,\rho$, or $l$. Since $\tau(x)$ controls the speed of the process in the direction $\frac{1}{\kappa(x)}\nabla l(x)$ it does not seem necessary to try and learn $\tau(x)$ and $l(x)$ simultaneously; therefore, depending on how good an ordering function we can initially construct we decide to learn either $\Phi=([\tau],[\rho])$ with $l(x)$ and $\kappa(x)$ fixed at the beginning of our optimisation, or we learn $\Phi=([l],[\rho])$ with $\tau(x)=\kappa(x)\equiv 1$ fixed - this treats $l(x)$ as a potential function that controls the deterministic dynamics of the process. Given this naive parameterisation, it seems sensible to impose some smoothness constraints on our dynamic parameters; for this reason, we define the following two smoothness regularisation functions that we will refer back to throughout this section:
\begin{align}
    R_\text{grad}(f)&=\int_\mathcal{M}\lVert\nabla f(x)\rVert^2\mathrm{d}x\approx\frac{1}{m}\sum_{i=1}^m\lVert G_i^\text{reg}(f_{N_i}-f_i\mathbb{1}_{N_i})\rVert^2\label{eq:sde_gradreg},\\
    R_\text{lap}(f)&=\int_\mathcal{M}|\Delta f(x)|^2\mathrm{d}x\approx\frac{1}{m}\lVert L^\text{reg}f\rVert^2\label{eq:sde_lapreg},
\end{align}
which we refer to as gradient regularisation and Laplacian regularisation, respectively. Both are chosen to promote smoothness of the learned functions.  To solve this optimisation problem we will use gradient based optimisation (we choose to use the L-BFGS algorithm \cite{nocedal1980updating,liu1989limited} as implemented in \cite{NLopt}) where we will exploit the rigid structure of an ODK in terms of the functions $\tau(x)$, $\rho(x)$ and $l(x)$ to compute the gradient analytically.  We give a full description of how we do this, and how we construct the gradient and Laplacian regularisation functions of equations \eqref{eq:sde_gradreg} and \eqref{eq:sde_lapreg}, in Appendix \ref{appBsec:odk_optimisation_scheme}.  

In our first example, we will apply the $L_q$ objective given in equation \eqref{eq:sde_q_objective}, to data sampled from a $2$-dimensional OU process --- with known $x_0=(5,5)^T$ and $r(t)=\text{Trunc}\left(\tfrac{1}{2}\mathcal{N}(0,0.6)+\tfrac{1}{2}\mathcal{N}(2,0.6);[0,2]\right)$ (i.e. we evolve the process over the interval $[0,2]$), and with unknown $\boldsymbol{\Theta}=(1.8,1.2)^T$, and $\boldsymbol{\sigma}=\sqrt{2}I_2$ --- to attempt to learn $\Phi=([l],\rho)\in\mathbb{R}^{m+1}$, where $[l]\in\mathbb{R}^m$ is an ordering function (here effectively viewed as a negative potential) represented as its evaluation across $\mathcal{D}=\{x_i\}_{i=1}^m$, and $\rho\in\mathbb{R}_{>c}$, for some small $c>0$, is a constant diffusion parameter; in total we have $m+1$ parameters to learn.

\begin{figure}
    \centering
    \includegraphics[width=0.98\textwidth]{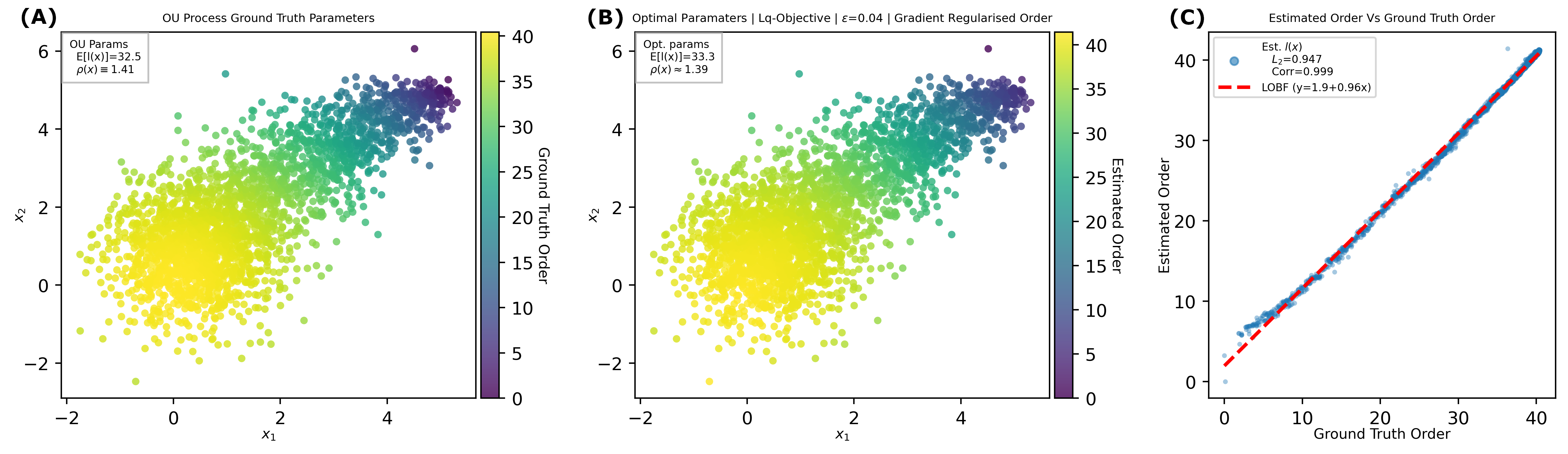}
    \caption{In this figure, we show the results of applying the $L_q$ objective (equation \eqref{eq:sde_q_objective}) to learn the ordering function (equivalently, the negative potential) and constant diffusion of a 2-dimensional OU process with $\Theta=(1.8,1.2)^T$ and $\rho_\text{gt}=\sqrt{2}$; we evolved from process from $x_0=(5,5)^T$ over the time horizon $[0,2]$ with time sampling distribution, $r(t)$, formed as the truncated mixture of two Gaussians centred at $t_1=0$ and $t_2=2$ and with standard deviations both equal to $0.6$.  In \textbf{Panel A}, we visualise the ground truth ordering function $l_\text{gt}(x)=-\frac{1}{2}\sum_{i=1}^d\theta_ix_i^2+C$, with $C$ such that $\min_{x\in\mathcal{D}}l_\text{gt}(x)=0$; as expected, the lowest-order samples lie near $x_0$ and the highest near the steady state at the origin. In \textbf{Panel B}, we show the estimated order and noise function after minimising the $L_q$ objective using L-BFGS \cite{liu1989limited}; we take $\varepsilon=0.04$, and use gradient regularisation on the ordering (equation \eqref{eq:sde_gradreg}, $\lambda_\text{reg}=0.00167$).  We see that the functions represented in Panels A and B are an incredibly close match (Pearson correlation $\approx0.999$), with a slight inflation of the mean order, $\mathbb{E}_{x\sim q}[l_\text{est}]=33.3$ versus $32.5$ ($2.43\%$ relative difference). In \textbf{Panel C}, we plot the estimated order against the ground truth order, with the red dashed line being the Line Of Best Fit (LOBF) between them. The near-perfect collinearity and near-negligible spread perpendicular to the LOBF indicate the estimated order truly recovers the dynamics and thus the foliating structure of $l_\text{gt}(x)$ on $\mathcal{M}$. From Panel C, we can spot two outlying estimates: the minimum-order sample (Panel C, bottom-left) and the maximum-estimated-order sample (Panel C, top-right), which correspond to the samples with the largest and smallest $x_2$, respectively. The estimated diffusion $\rho_\text{est}\approx1.39$ is within $1.38\%$ of the true $\rho_\text{gt}=\sqrt{2}\approx1.41$.}
    \label{fig:ou_sde_example1}
\end{figure}

In figure \ref{fig:ou_sde_example1}, we show the results of applying this optimisation: in panels A and B, we visualise the ground truth order and estimated order, respectively, and observe very little discernable different between the two, the relative error of this estimate averaged across $\mathcal{D}$ being approximately $2.43\%$, and the $L^2$-norm between them being $\lVert l_\text{gt}-l_\text{est}\rVert_{L^2}\approx0.947$. This alignment is more explicitly shown in panel C, where we plot the estimated order against the ground truth order and observe that they are almost perfectly collinear, having a Pearson correlation of 0.999, which illustrates how the optimisation was able to infer even the fine structure of the implied ordering (negative-potential) of this process simply by numerically integrating ODK time-marginals of the process restricted to its samples.  Furthermore, the ground truth and estimated $\rho$'s were equal to $\rho_\text{gt}(x)\equiv\sqrt{2}\approx1.414$ and $\rho_\text{est}(x)\equiv1.395$, respectively, giving a relative error of $1.38\%$, and showing that the $L_q$ objective is capable of isolating a signal left by the underlying noise that generated the data in the geometry (relative spacing and connectivity) of the samples.  What is particularly noteworthy about these results is how little dynamical information we needed: recall that we only provided the model with the time-sampling distribution $r(t)$, the spatial sampling distribution $q(x)$, and an initial state $x_0$.  The initial state implies a distance \textit{from where all samples started}, that surely correlates with the underlying order the process implies, and simply by restricting our representation of the process to only consider sample-to-sample transitions, and constraining the transition probabilities to be like $\varepsilon$-timestep It{\^o} transitions, the $L_q$-objective and ODK representation of the process are capable of performing a non-linear adaptive fitting between the naive implied distance from $x_0$ and the true underlying order.  We can also simultaneously isolate a signal from the intrinsic noise encoded into the geometry of the data.

\subsubsection{SDE with known Sampling Times}
\begin{figure}
    \centering
    \includegraphics[width=0.75\textwidth]{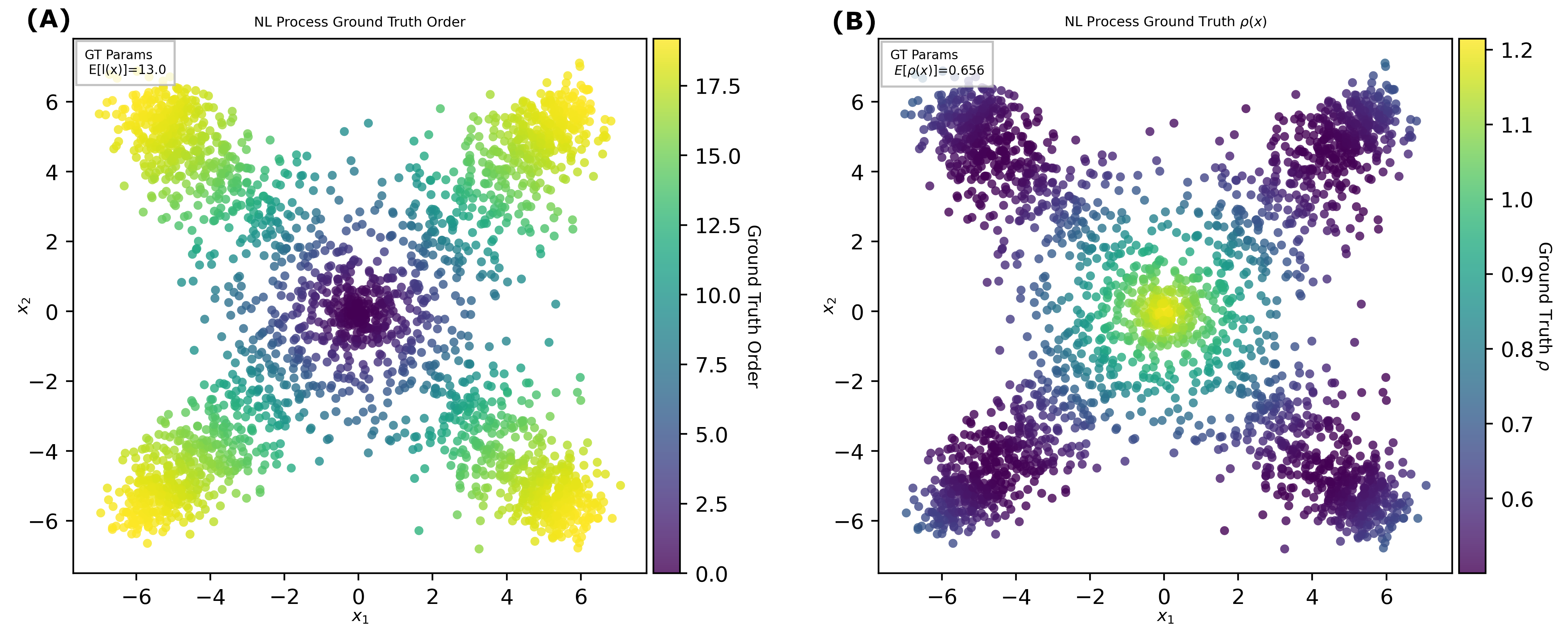}
    \caption{Ground truth ordering function (\textbf{Panel A}) and ground truth noise function (\textbf{Panel B}), as given in equations \eqref{eq:nl_order_function} and \eqref{eq:nl_sde}, respectively.  We evaluate these functions on $2{,}500$ independently sampled points from the NL Process with $x_0=(0,0)^T$ and $r(t)=\text{Unif}([0,4])$.}
    \label{fig:nl_sde_gt}
\end{figure}

For our final example, we consider a more parametrically complex SDE with non-linear drift and diffusion (where the diffusion is now a function of the state variables); however, we now assume that we are provided with the sampling times for each data point. The Non-Linear SDE (that we hereon refer to as the NL process) has the following form:
\begin{equation}\label{eq:nl_sde}
    dx_{t,i}=\left(\frac{3}{2}\sin\left(\tfrac{1}{2}x_{t,i}\right)+\frac{1}{2\pi^2}(4\pi^2x_{t,i}-x_{t,i}^3)\exp\left(-\frac{1}{4}x_{t,i}^2\right)\right)\mathrm{d}t+\frac{1}{2}\sqrt{1+\tfrac{1}{8}\left(\lVert x_t\rVert-2\pi\right)^2}dw_{t,i},
\end{equation}
where $x_{t,i}$ denotes the $i^\text{th}$ component of the state $x$ at time $t$ for $i=1,2$; we use two independent Wiener processes to generate the noise in this SDE, thus corresponding to an isotropic diffusion.  The ordering function (negative potential) corresponding to this process takes the form:
\begin{equation}\label{eq:nl_order_function}
    l(x)=-\sum_{i=1}^23\cos\left(\tfrac{1}{2}x_i\right)+\frac{1}{\pi^2}(4(\pi^2-1)-x_i^2)\exp\left(-\tfrac{1}{4}x_i^2\right),
\end{equation}
To sample our data, we fix $x_0=(0,0)^T$ and define our time-sampling distribution to be uniform on the interval $[0,4]$, i.e. $r(t)=\text{Unif}([0,4])$; as before, for $i=1,\dots,m$, after sampling $t_i\sim r$, we then sample $x_i\in p_t(x|x_0=(0,0)^T)$, and store the sampling times and state samples in the sets $\mathcal{T}=\{t_i\}_{i=1}^m$ and $\mathcal{D}=\{x_i\}_{i=1}^m$, respectively.  In panels A and B of Figure \ref{fig:nl_sde_gt}, we plot the ordering function, as given in equation \eqref{eq:nl_order_function}, and the noise function, as given in the second term of equation \eqref{eq:nl_sde}, respectively, over a dataset of $2{,}500$ independent samples from the NL process.  We can see from panel A that the NL process has 4 steady states located at $(\pm2\pi,\pm2\pi)$, and from panel B we observe that the diffusion is strongest nearer the initial condition and decays once a state appears to have committed to travelling along one of the four possible branches; finally, we note that there is a small increase in the magnitude of the noise when $\lVert x\rVert>2\pi$, this creates a potential well around the circle of radius $2\pi$.

In what follows, we only assume that we are given the sampling times $\mathcal{T}=\{t_i\}_{i=1}^m$, the sampled states $\mathcal{D}=\{x_i\}_{i=1}^m$ (where $x_i$ was sampled at $t_i$), and the initial state $x_0$.  Since we have more dynamical information in this example than in the previous example, we define a new, stronger objective that intrinsically incorporates the time sampling distribution (in a Monte-Carlo sense) and acts as a full maximum likelihood over the joint distribution $p(x,t)=r(t)p_t(x|x_0)$. We define the $L_p$ objective to be:
\begin{equation}\label{eq:Lp_objective}
    L_p(\Phi|\mathcal{D},\mathcal{T},x_0,\lambda_\text{reg})\coloneqq D_\text{KL}(p(x,t)||\hat p(x,t|\Phi))+\lambda_\text{reg}R(\Phi|\mathcal{D})\approx -\frac{1}{m}\sum_{i=1}^m\log\hat p(x_i|t_i,\Phi)+\lambda_\text{reg}R(\Phi|\mathcal{D})+\text{const}.
\end{equation}
This is a standard maximum likelihood objective, where we have used the fact that $r(t)$ does not depend on the parameters of the SDE (by assumption). In the above equation, $\hat p(x_i|t_i,\Phi)$ is the estimated probability of being at the state $x_i$ at time $t_i$ given our $\varepsilon$-timestep ODK transition probabilities, $P_\varepsilon(\Phi)$. The simplest estimate of this function defines $s_i=\left\lceil\tfrac{t_i}{\varepsilon}\right\rceil$ and letting $\hat p(x_i|t_i,\Phi)=\hat p_{s_i}(x_i|\Phi)\coloneqq p_0^TP_\varepsilon(\Phi)^{s_i}e_i$, but we can generalise this further by linearly interpolating between discrete time-steps through:
\begin{equation}
    \hat p(x_i|t_i,\Phi)=a_ip_{s_i-1}(x_i|\Phi)+b_ip_{s_i}(x_i|\Phi)=a_ip_0^TP_\varepsilon(\Phi)^{s_i-1}e_i+b_ip_0^TP_\varepsilon(\Phi)^{s_i}e_i,
\end{equation}
where $a_i=\left(1-\left(s_i-\tfrac{t_i}{\varepsilon}\right)\right)$ and $b_i=\left(s_i-\tfrac{t_i}{\varepsilon}\right)$.  As with the $L_q$ objective used in the previous example we can consider the parameters of this objective to be any of the dynamic functions $l(x)$, $\tau(x)$ or $\rho(x)$ evaluated at each of the data points $x_i$, and it is again possible to compute the gradient of this objective analytically in a single forward pass due to the rigid structure of our ODKs with respect to these parameters, a rigorous description of this computation for the $L_p$ objective with linearly interpolated probability estimates is given in Appendix \ref{appBsec:odk_optimisation_scheme}.  In what follows, we will again optimise for the ordering function, $l(x)$, and the now state-dependent noise, $\rho(x)$; we continue to employ the L-BFGS algorithm to solve our optimisation problem.

It is now worth discussing the issue of parameter identifiability for our proposed objective. To emphasise the issue, suppose the optimiser wished to push probability from one region into another: it could do this by increasing the ordering in this direction and thus creating a deterministic flow towards it; but it could also create a diffusion well around this new region; this would create an effective advection proportional to $-\rho\nabla\rho$ stochastically driving probability into the bottom of this well.  Thus, when one's diffusion depends on the state variable, the optimiser must disentangle where advection is driven deterministically and where it is driven stochastically.  As an example of this issue, it has been shown experimentally that, for a small subset of genes, stem cell differentiation (which is a process most often modelled as being primarily deterministic) can actually be modelled accurately through an entirely diffusive process with state-dependent diffusion creating a well near the differentiated state \cite{hoffmann2008noise}.  The question of identifiability boils down to whether it is always possible for the optimiser to learn this balance between determinism and stochasticity given the $L_p$ (or $L_q$) objective.  Since the $L_p$ objective requires sampling times, the signal of the noise in the data likely comes from two sources: the distribution of samples (as with the $L_q$ objective), where sparser regions can be associated with larger diffusions; and the local variances of sampling times, which we refer to as time uncertainty, this fundamentally will be a function of the total diffusion a state has experienced in getting to any given point.  Thus, if the optimiser overuses a diffusion gradient to push the probability forward, the subsequent evolution will be mostly deterministic, resulting in lower time uncertainties than observed --- recognising this discrepancy should allow the model to identify the true noise values.  The issue is that unless your sampling measure was uniform (highly unlikely for complex dynamics), or you had a tremendous amount of data, you will always have (pointwise) noisy estimates of the time uncertainty and thus one's estimates of the diffusion cannot help but be noisy.  As we will see in the remainder of this section, it does seem possible to simultaneously infer the correct drift and diffusion using the $L_p$ objective, although the estimates of the diffusion can be very noisy and the addition of smoothness regularisation functions (such as those given in equations \eqref{eq:sde_gradreg} \& \eqref{eq:sde_lapreg}) is capable of distorting the local minima away from the true diffusion.  We take this as weak empirical evidence that our optimisation is identifiable, but poorly conditioned for state-dependent noise; we do not prove a result of this sort in this paper but this would be highly relevant future work, we provide a more detailed discussion of identifiability in Appendix \ref{appBsec:identifiability}.

\begin{figure}
    \centering
    \includegraphics[width=0.98\textwidth]{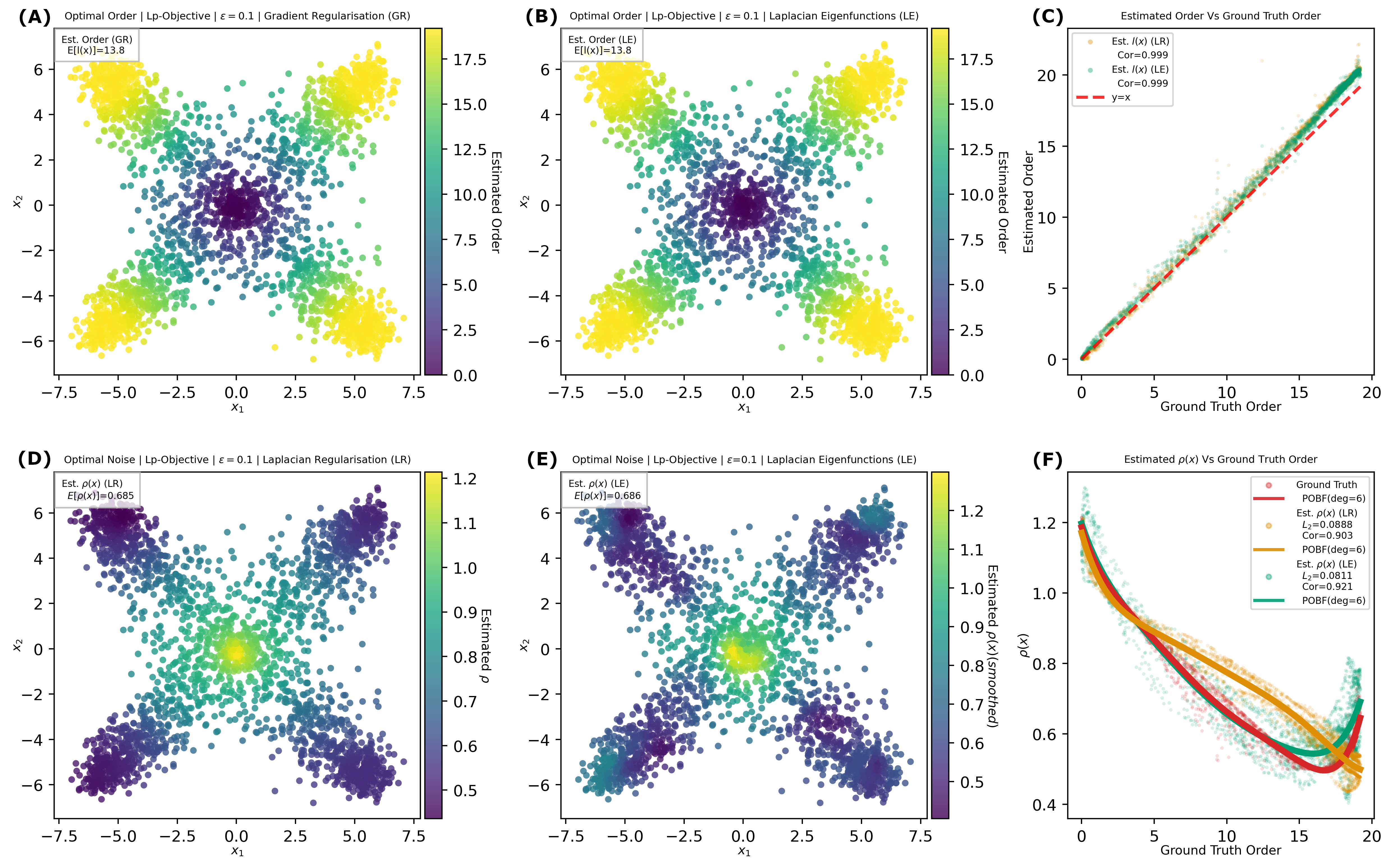}    
    \caption{In this figure we visualise the optimal ordering and noise found by minimising two runs of the $L_p$ objective (equation \eqref{eq:Lp_objective}) with different smoothness regularisations, applied to $2{,}500$ data points sampled from the NL process (equation \eqref{eq:nl_sde}) with $r=\text{Unif}([0,4])$ and from $x_0=(0,0)^T$; the ground truth for both is visualised in figure \ref{fig:nl_sde_gt}.  In \textbf{Panel A} \& \textbf{Panel D}, we show the estimated order and noise functions, respectively, from a run of the $L_p$ objective with gradient smoothness regularisation (GR, equation \eqref{eq:sde_gradreg}, $\lambda_\text{reg}^l=0.03$) applied to $l(x)$, and Laplacian smoothness regularisation (LR, equation \eqref{eq:sde_lapreg}, $\lambda_\text{reg}^\rho=12.5$) applied to $\rho(x)$.  In \textbf{Panel B} and \textbf{Panel E}, we show the estimated order and smoothed (Gaussian smoothing with $k=35$ in equation \eqref{eq:nl_vb_smoothing}) noise functions, respectively, for a run of the $L_p$ objective, with each function being a weighted sum of the top $k=250$ eigenfunctions of the symmetric transition matrix $P_\varepsilon^\text{sym}$ with $\varepsilon=0.1$ (LE; see main text of section \ref{para:laplacian_eigenbasis}, and Appendix \ref{appBsec:odk_optimisation_scheme}) as in equation \eqref{eq:lap_eigenbasis_decomposition}.  \textbf{Panel C} plots the estimates in panel A (orange dots) \& panel B (green dots) against the ground truth order; each estimate is highly correlated with $y=x$ (dashed red line), both having a Pearson correlation of approximately $0.999$ with the ground truth order.  However, we also observe that both runs slightly overestimate the order.  \textbf{Panel F} plots the ground truth noise of panel B in figure \ref{fig:nl_sde_gt} (red dots), the estimated noise of panel D (orange dots) \& the smoothed estimated noise of panel E (green dots), against the ground truth noise; the solid red, orange and green lines are degree-6 Polynomials of Best Fit (POBF) to the scatter plot of each respective colour.  We again observe that both the (GR+LR) run and the LE run capture the rough shape of the true noise. However, the (GR+LR) noise overestimates the noise while travelling along any of the four branches of the NL process and fails to resolve the uptick in the diffusion past $\lVert x\rVert_2=2\pi$; while the LE noise, even after applying additional smoothing, is still more noisy than the regularised estimate, however, the fit of the LE estimate is far more aligned with the groundtruth, seemingly not overestimating the noise along the branches and resolving the uptick past $\lVert x\rVert_2=2\pi$.}
    \label{fig:nl_results}
\end{figure}

A practical issue remains: even with the information contained in our time uncertainty, the signal to the $\rho(x)$ parameters from the $L_p$ objective is much noisier than it is for the ordering function, which makes it more challenging to infer smooth representations of the noise.  Intuitively, this seems to be because from the perspective of the deterministic/stochastic tradeoff, what is most important is the gradient of the order, $\nabla l(x)$, which shows up explicitly in the expansion of ODK, and the gradient of the diffusion, $\nabla \sigma^2(x)$, which does not; thus, the stochastic flux of probability is more obscured than the deterministic flow of the process. In addition, the diffusion will also be more heavily dependent on the observed sampling density of the data, $q_\text{empirical}(x)$, discrepancies between the observed sampling distribution and the process's true transient distribution will then surely manifest in the estimates of the diffusion. We can observe this more noisy signal by running the proposed $L_p$ optimisation without any smoothness regularisation on the ordering or noise function: when we do this, the order estimate is still relatively smooth, although more noisy and less accurate than the results given in figure \ref{fig:nl_results}, while the estimated $\rho$'s appear highly noisy with some pointwise estimates showing very anomalous behaviour.  This shows that for some points (particularly points near the boundary) there is very little signal from the objective, and variations in the sampling density can distort estimates, leading to noisy results.  In addition, if we run the optimisation only for $\rho(x)$ while fixing the order to be equal to the ground truth, we again observe a noisy estimate but we notice that if we - quite heavily ($k\approx 100$) - smooth this estimate with a Gaussian kernel, i.e. we compute $\bar \rho_k(x)$ where for a function $f(x)$ the smoothed estimate is defined as:
\begin{equation}\label{eq:nl_vb_smoothing}
    \bar f_k(x_i)=\frac{1}{Z_i}\sum_{j=1}^m\exp\left(-\frac{\lVert x_i-x_j\rVert^2}{2\lVert x_i-x_i^k\rVert^2}\right)f(x_j),\qquad\text{where:}\qquad Z_i=\sum_{j=1}^m\exp\left(-\frac{\lVert x_i-x_j\rVert^2}{2\lVert x_i-x_i^k\rVert^2}\right),
\end{equation}
and $x_i^k$ is the $k^\text{th}$ nearest neighbour to $x_i$ in $\mathcal{D}$ (not including $x_i$ itself). Then we see that $\bar \rho_k(x)$ is very close to the true noise function (indeed, using $k=90$ the results look very similar to Panel E of figure \ref{fig:nl_results}), which shows that, conditional on the order the true signal of the noise does exist in the data and can be extracted through the $L_p$ objective. Thus, we conclude that when optimising for the order and noise function simultaneously, applying proper smoothness regularisation is important to ensure that variations in the sampling density do not obscure the true signal of the noise in the data and to regularise the loss landscape of the objective so that locally undecidable tradeoffs between deterministic and stochastic dynamics do not skew the estimated order or diffusion.  We propose two methods for introducing a smoothness regularisation and discuss each separately in the next two subsections; the results of applying both methods are shown in figure \ref{fig:nl_results}.
    
\paragraph{Regularisation Functions}
The first method of smoothness regularisation follows from the regularisation we performed in the OU example; we consider adding a gradient regularisation term for the ordering function (i.e. $R_\text{grad}(\hat l)$ as given in equation \eqref{eq:sde_gradreg}) and a Laplacian regularisation term to the noise function (i.e. $R_\text{lap}(\hat\rho)$ as given in equation \eqref{eq:sde_lapreg}), to the $L_p$-objective. We weight the hyperparameters of each term to account for the relative differences in scale between the ground truth order and noise; more precisely, we let $\lambda_\text{reg}^l/\lambda_\text{reg}^\rho\approx \mathbb{E}_{x\sim q}[\rho(x)]^2/\mathbb{E}_{x\sim q}[l(x)]^2$ where we take the square of each function because both the gradient and Laplacian regularisation functions scale multiplicitively as $R_\text{grad}(c\cdot f)=c^2R_\text{grad}(f)$ and $R_\text{lap}(c\cdot f)=c^2R_\text{lap}(f)$. We choose the overall scale by doing a small grid search while constraining the relative strength of both regularisation terms to be small relative to the total objective: ultimately we end up with $\lambda_\text{reg}^l=0.03$ and $\lambda_\text{reg}^\rho=12.5$, which amounts to approximately $3.000\%$ of the total objective coming from these regularisation terms (when evaluated at the ground truth values).  It is desirable that the regularisation is only a small proportion of the total objective, as we don't want to overly influence the functional form of either the ordering or the noise.  

We show the results of running this optimisation in panels A and D of figure \ref{fig:nl_results}. We see very little discernable difference between the ground truth order and estimated order (panel A of figure \ref{fig:nl_sde_gt} and panel A of figure \ref{fig:nl_results}) with average value of the ground truth order being $\mathbb{E}_{x\sim q}[l_\text{gt}(x)]\approx13.05$ and the estimated average value being $\mathbb{E}_{x\sim q}[l_\text{est}^\text{reg}(x)]\approx13.84$, the relative difference of the estimate averaged over the data is $6.1\%$ and shows the scale has been appropriately estimated, and has done so in a uniformly accurate manner across the data.  The $L^2$-norm between the ground truth and the estimated order is $\lVert l_\text{gt}(x)-l_\text{est}^\text{reg}\rVert_{L^2}\approx0.961$, which gives a relative difference ($\lVert l_\text{gt}(x)-l_\text{est}^\text{reg}\rVert_{L^2}/\lVert l_\text{gt}\rVert_{L^2}$) of approximately $6.67\%$, which again shows the estimates are reasonably accurate across the whole dataset.  

In panel C, we plot the estimated orders against the ground truth order, and we observe a high degree of collinearity between the two (where the red dashed line shows the line $y=x$); this is quantified by a Pearson correlation coefficient of $0.999$ between the estimated and ground truth orders.  From panel C, we can also notice that the order is slightly overestimated, indicating a slight preference for moving probability deterministically (a steeper gradient).  Overall, this shows that the first optimisation accurately reconstructs the first-order dynamics of the NL process.  In panel D, we visualise the estimated $\rho_\text{est}^\text{reg}(x)$ given this regularised optimisation and observe that the rough sense of the diffusion has been captured (that being a maxima near the initial condition $x_0=(0,0)^T$ that then decays along the four different branches of the data); however, we note that the ground truth order flicks up once $\lVert x\rVert>2\pi$, which can be seen in the red line of panel F.  The $\rho_\text{est}^\text{reg}(x)$ estimate coming from the regularised optimisation does not seem to capture this, which can be seen in the solid orange line of Panel F. 

We also note that between $l_\text{gt}=5$ and $l_\text{gt}=17.5$ the estimated noise is overestimated, which aligns with our expectations given that the order was slightly overestimated in this region (overestimating the diffusion means the effective drift due to stochasticity is less; overestimating the order then compensates for this).  In spite of this, the relative difference of the scale of the estimated and ground truth $\rho(x)$ stands at $4.48\%$, the $L^2$-norm between the estimate and true values equals $\lVert\rho_\text{gt}-\rho_\text{est}^\text{reg}\rVert_{L^2}\approx 0.0888$, with an average $L^2$-difference of $13.0\%$, and Pearson's correlation coefficient of $0.903$ to the ground truth.

\paragraph{Laplacian Eigenbasis}\label{para:laplacian_eigenbasis}
The second method of smoothness regularisation takes inspiration from Diffusion Maps \cite{coifman2006diffusion}, and works by replacing the pointwise estimate of $l(x)$ and $\rho(x)$ with:
\begin{equation}\label{eq:lap_eigenbasis_decomposition}
    \hat l(x)=\sum_{k=2}^K\alpha_k^l\varphi_k(x),\qquad\text{and:}\qquad\hat\rho(x)=\exp\left(\sum_{k=1}^K\alpha_{k}^\rho\varphi_k(x)\right),
\end{equation}
where $\varphi_k$ satisfies $L_\varepsilon^\text{sym}\varphi_k=\lambda_k\varphi_k$, and $K>1$ is the number of basis functions we use and we choose the eigenfunctions to be those with the $K$ largest eigenvalues.\footnote{In practice, since $\langle\varphi_i,\varphi_j\rangle_{L^2(\mathcal{M},q)}=\delta_{ij}$ are only orthonormal with respect to the sampling measure $q$, they will not be orthonormal in the standard Euclidean $L^2$ sense, so we orthonormalise our eigenfunctions using the Gram-Schmidt process \cite{leon2013gram}.  See appendix \ref{appBsec:odk_optimisation_scheme} for more details.}  We know that $\varphi_1(x)\propto\mathbb{1}(x)$ is the trivial constant eigenfunction, and since shifting the ordering function by a constant does not change the dynamics of our ODK ($\nabla (\hat l(x)+C)=\nabla\hat l(x)$), thus we omit $\varphi_1$ from the eigenspace of the ordering function.  We exponentiate the weighted sum for $\hat \rho$ to enforce strict positivity.

In the example given in figure \ref{fig:nl_results}, we use $K=250$, which, in theory, would be enough to represent almost any smooth function that can reasonably be defined on $\mathcal{D}$ that we might be interested in.  In practice, the discrete eigenvectors we approximate only approximate the eigenfunctions of the Laplacian, and they do so well, but when we start considering higher frequency eigenfunctions (smaller eigenvalues/larger $K$), the discrete eigenvectors depend increasingly more on the specific sampling of the data, i.e. $\mathcal{D}$, that one has --- eventually the eigenfunctions will become almost $\delta$-like.  Thus, there is a trade-off when choosing $K$: too small, and the true function will not be contained in your eigenspace; too large, and you are effectively not enforcing any smoothness regularisation at all.  We see by projecting the true noise onto these eigenfunctions that $\approx99\%$ can be explained by the top 20 eigenfunctions (and even $\approx 90\%$ by the top 10), we also see that the magnitude of the basis weightings are almost negligible after $K>100$.  Using more eigenfunctions than necessary could result in higher-frequency oscillations that make our estimates noisier, but also expand the function space, making it easier for the optimiser to navigate.  The value $K=250$ was chosen after a grid search and is likely too large to enforce smoothness as effectively as we would like.

We show the estimated ordering function and smoothed estimated noise function after running this optimisation in panels B \& E of figure \ref{fig:nl_results}, respectively.  In Panel B, we observe very little discernable difference between the ground truth order (shown in panel A of figure \ref{fig:nl_sde_gt}), with $\mathbb{E}_{x\sim q}[l_\text{est}(x)]\approx13.8$ (an overestimation of $5.5\%$), and with the $L^2$ difference being $\lVert l_\text{gt}(x)-l_\text{est}(x)\rVert_{L^2}=0.897$ which works out to a $5.88\%$ relative error when averaged across the data.  

In Panel C, we again plot the order estimate (green dots) against the ground truth and observe much the same as with the previous estimate (Pearson correlation coefficient $\approx0.999$); both the regularised and Laplacian Eigenbasis optimisations yield very similar optimal orders, with maybe a slight edge to the LE estimate as it seems not to overestimate the true order as much as the GR+LR estimate.  In Panel E, we plot a smoothed estimate of the noise, where we have applied a Gaussian smoothing with $k=35$ as in equation \eqref{eq:nl_vb_smoothing};\footnote{We choose $k=35$ somewhat arbitrarily, using some heuristics for $k<30$ there still appeared to be some noise, and for $k>40$ there was very little marked difference in the estimates before oversmoothing took over.  Due to an inability to decide between $k=30$ and $k=40$, we used $k=35$.} although the smoothed estimate is still highly aligned with the ground truth, we can clearly see that even after applying this additional smoothing procedure, the estimates are still more noisy than those in Panel D.  This indicates to us that the LE smoothness regularisation is not strong enough a constraint to enforce smoothness in our estimates; this is likely due to $K$ being too large and thus includes too many high-frequency sharp eigenfunctions eigenfunctions.  The plot in Panel F provides a clearer picture of the differences between the two estimates: it is immediately clear from looking at the spread of the orange dots (GR+LR noise estimate) compared to the spread of the green dots (LE estimate) that, even after additional smoothing, the LE estimate is still much more noisy than the GR+LR estimate; however, the GR+LR estimate overestimates the diffusion in the central portion of the ordering values, whereas the LE estimates seem to cluster around the true values (red dots).  The solid red, orange and green lines show Polynomials of Best Fit (POBF) of degree-6 that fit a smooth curve to the scatters of the corresponding colours; from this we see that the LE noise estimate really does cluster around the true noise values, even resolving the upwards flick of the diffusion after $\lVert x\rVert_2>2\pi$ --- which we recall the GR+LR estimate does not resolve.  Quantitatively, the smoothed estimates differ in scale from the true values by $4.65\%$, have an absolute $L^2$-difference of approximately $0.0811$, which is equivalent to an $11.9\%$ relative difference; finally, the Pearson correlation coefficient between the estimate and the ground truth stands at $0.921$.

To probe why the estimates of the noise function differ between the (GR+LR) and LE runs of the $L_p$ objective optimisation scheme, we assessed whether the estimated noise from the regularised run, lay within the eigenspace of the LE run: we found that over $99.9\%$ of $\rho_\text{est}^\text{reg}(x)$ could be explained by our chosen eigenbasis, thus it is a feature of the choice of regularisation that causes the specific shape observed in the graph of $\rho_\text{est}^\text{reg}(x)$ as seen in Panel D and F of figure \ref{fig:nl_results}.  But as we noted, the total regularisation for both the order and the noise accounts for only about $3\%$ of the total objective's magnitude.  This indicates that the $L_p$ objective around the true $\rho_\text{gt}(x)$, is very flat and cannot easily distinguish between different functional forms; adding the LR regularisation clearly creates a local minimum around the observed optimal $\rho_\text{est}^\text{reg}(x)$ in the (GR+LR) run.  However, the fact that when we replace the Laplacian regularisation with an eigenbasis decomposition, which doesn't change the objective and merely restricts the parameter space, we get an estimate which clusters strongly around the true values of $\rho_\text{gt}(x)$ indicates to us that the signal of the true noise is present in the dataset, but that this signal is so weak that it cannot easily be distinguished from other nearby functions in parameter space.  We could potentially ameliorate these effects by the addition of a more principled regularisation; for example, we could consider adding a regularisation term of the form $R(\Phi)=\tfrac{1}{m}\sum_{i=1}^m\tfrac{\tau(x_i)^2}{\rho(x_i)^2}$ which is the drift part of the Onsager-Machlup action \cite{onsager1953fluctuations} averaged over the data; in theory, minimising this action will minimise the relative entropy between the inferred path law from a driftless reference law.

The potential unidentifiability of our optimisation problem reduces to how the model balances drift and diffusion, so an objective that depends directly on such a quantity might help amplify the noise signal relative to the drift; furthermore, minimising the Girsanov action is equivalent to maximising the path entropy, thus, in theory, the optimal parameters will assume the least structure beyond what is required by the data.  In this paper we will not explore other possibilities for setting up this optimisation problem, but we note that this is highly relevant and interesting future work.  Assuming the KL-divergence is a good choice of measure for this problem, the $L_p$ objective being flat around $\rho_\text{gt}(x)$ might be empirical evidence that this optimisation problem is theoretically unidentifiable and the fact that either the (GR+LR) or LE run give results that align with the ground truth is the result of some other artifact in the data, that need not be present for all processes or in all samples.

What we have shown with this final example is that both formulations of the optimisation problem given by the $L_p$ objective are capable of accurately reconstructing the dynamics of the NL process given in \eqref{eq:nl_sde}.  Both formulations reconstruct a smooth order that is highly correlated with the true ordering, indicating that we can truly resolve the first-order dynamics in the presence of spatially varying noise arising from a complex non-linear SDE.  The noise estimate from the regularised (GR+LR) optimisation is clearly a smooth function, as desired, and captures the form of the true noise function accurately (maximal initially, and then decays the further from the centre one goes); although, it overestimates the diffusion (underestimates the diffusion gradient) in the central portion of the dynamics and fails to resolve the uptick in diffusion beyond $\lVert x\rVert_2>2\pi$.  In contrast, the estimated noise from the LE optimisation is highly clustered around the true noise function through the entire range of the process and, furthermore, is capable of resolving the uptick in diffusion beyond $\lVert x\rVert_2>2\pi$.  However, the LE noise estimate is highly noisy, likely due both to choosing $K=250$, which is much larger than needed, and to the signal from the noise in the data being very weak, making it difficult to identify a single smooth function within the eigenspace.  Both the $L_q$ and $L_p$ objectives are first attempts to build an optimisation framework that leverages the rigid structure and interpretable design of ODKs to learn unknown first- and second-order dynamics from data generated by a Markov stochastic process and which comes with very little a priori dynamical information.

\section{Conclusion}

In this paper, we have introduced Ordered Diffusion Kernels (ODKs) as a novel class of local kernels for computing transition probabilities and the infinitesimal generator for an arbitrary It{\^o} SDE.  We assume that we have data sampled from a dynamical system not in steady state and with known initial conditions; we do not require any additional dynamical information, such as known rates of evolution, a mechanistic model for the underlying process, exact sampling times, or a known time-sampling distribution.  As a prototypical example of such data, we use single-cell RNA sequencing data (scRNA-seq), which captures a gene expression landscape of a cell population and which we assume is not in steady state, i.e. a stem cell undergoing differentiation.  Representing or learning dynamics on this data is thus a highly underdetermined problem and is often approached by building hugely expressive inference machines to overwhelm the task; these often rely on deep neural networks or complex probabilistic models.  The issues with such models are that, due to the lack of dynamic constraints, one cannot be sure that massively over-parameterised methods have learnt coherent dynamics that faithfully correlate with the true underlying dynamics.  We designed ODKs as a simple alternative to these methods, the expansion of the normalised integral operator of Theorem \ref{thm:odk_backward_expansion} tells us exactly what underlying process we are approximating, that being the It{\^o} SDE with infinitesimal generator $\mathcal{L}=\tfrac{\tau}{\kappa}\nabla l\cdot\nabla+D\Delta$, and the restriction to sample-to-sample transitions constrains the model to not learn unjustifiable out of distribution dynamics.  We took great care in making sure ODKs were as interpretable and modular as possible, allowing practitioners a great deal of freedom to tune the parameters or tailor the ODK to the specific data they are studying.  We emphasised this modularity through the various extensions of ODK that we discussed: we constructed local-Ordered Diffusion Kernels (lODKs) to represent non-gradient flow dynamics (section \ref{subsec:nongradflow_extension}); multi-Ordered Diffusion Kernels (mODKs) are compatible with multiple ordering functions (section \ref{subsec:multi_odk}; anisotropic-Ordered Diffusion Kernels (anODKs) allow for anisotropic diffusions (section \ref{subsec:nonisotropic_extension}; and self-tuning-Ordered Diffusion Kernels (stODKs) allow for symmetric bandwidths of the form $\rho(x)\rho(y)$ (section \ref{subsec:selftuning_odk}).  For this latter extension, we also strengthened the asymptotic results of Variable Bandwidth Kernels \cite{berry2016variable} from weak to strong convergence.  Fundamentally, all these extensions rely on what we call an ordering function, which we conceptualise as a relaxed form of a potential in which only the implied directionality is important, i.e. ODKs will bias transitions towards higher-order samples.  Since we assume the initial conditions are known --- which for the example of scRNA-seq would be some a priori knowledge of the gene expression of stem cells --- we can always construct an ordering as the geodesic distance --- on an undirected connected neighbourhood graph of the data --- of any sample from these initial conditions; so in the abscence of any dynamical information we can still use ODKs to represent dynamics that will correlate with the underlying dynamics that generated the data.  In section \ref{subsec:numerical_results_sde}, we leveraged the rigid structure of these dynamic representations to build two loss functions, namely the $L_q$ and $L_p$ objectives, which we can fit into a simple gradient descent optimisation scheme to learn the drift and diffusion parameters of data sampled on an SDE.  We demonstrated on data sampled from an Ornstein-Uhlenbeck and Non-Linear SDEs that we were capable of learning the true underlying drift and diffusion of the underlying process, even when the diffusion was state-dependent.

Our results follow on from a considerable body of work that formalises the intuition that (un)weighted graph Laplacians produce representations that are coherent with the geometry of the underlying manifold.  Most directly relevant to our work is TMDmap by Banisch et al. \cite{banisch2020diffusion}, who leveraged ideas from the $\alpha$-density regularised Diffusion Maps kernels \cite{coifman2006diffusion} to aim the transitions towards any target measure; however, this creates an inextricable coupling between the first- and second-order terms of the Laplacian.  The structure of ODKs allows us to avoid this coupling, allowing us the freedom to set the drift and diffusion independently of one another.  Also relevant is the work of Berry \& Sauer in their paper on Local Kernels \cite{berry2016local}; here the most general formulation of the asymptotic convergence of (localised) kernels is provided, and is a result we use repeatedly in the proofs given in Appendix \ref{appendix:proofs}.  However, the framing of local kernels was from the perspective of vector fields, which requires abundant dynamical information to infer with confidence.  To summarise, when working with data generated by a dynamical system and in which little dynamical information is available, we believe the reformulation of the unknown dynamics in terms of an ordering function is both a substantial relaxation of the initial inference problem, while still being sufficient to represent the underlying dynamics.  ODKs are then a practical reformulation of Local Kernels within this framework, giving practitioners a tool to convert a minimal representation of the dynamics into a practical tool for further analysis, or a dynamic representation from which to formulate an inference problem.  The well-understood limiting behaviour of ODKs then gives users confidence in the fidelity of this representation and clarity about how each component input to the ODK is resolved in the data.

\paragraph{Limitations:} The choice of an optimal parameterisation for kernel methods is historically a very challenging task since one needs to balance a choice of parameterisation where the kernel is most useful for further analysis, and one must ensure that the kernel is fitting the specific sample of the data that one has.  Our conceptualisation of ODKs as forming a representation of the underlying dynamics that generated the data constrains the notion of what \textit{good} looks like, but for data that is very sparse or inhomogeneously sampled, it is not clear that the \textit{best} representation of the dynamics would be possible given a specific sampling of the data.  Thus, for data of this sort, the intrinsic sensitivity of kernel methods to the sampling of the data cannot be ignored, and must be dealt with carefully.  One of the strongest traits of kernel methods is their non-dependence on the dimensionality of the data, since the only step in the evaluation of a kernel that depends on the dimension is the initial computation of a pairwise distance matrix, this is the essence behind calling the evaluation of a kernel a \textit{representation} of the data; however, what is often left unsaid is how kernels are affected if the distance metric itself is unreliable.  Particularly, for data that is embedded in $\mathbb{R}^n$ with $n\gg 100$, the Euclidean metric is known to struggle in distinguishing distances; particularly, the concentration of measure phenomena \cite{ledoux2001concentration} shows that as $n\to\infty$, the distance to the nearest and furthest samples in a bounded domain will tend to the same value, i.e. the ability for the Euclidean metric to discriminate between distances collapses as the dimension increases --- this is particularly relevant to scRNA-seq data in which one is considering the RNA counts associated to many thousands of genes.  In such situations, one should also take care to ensure that the chosen distance metric is informative enough for the kernel to resolve the data as desired.  We also note that anisotropic diffusions do not fit naturally within the framework of ODKs as they do with local kernels; this can be seen in the expansion of Theorem \ref{thm:noniso_odk_expansion} where the drift is coupled to the anisotropic diffusion tensor, much like the drift is coupled to the diffusion in TMDmaps.  We mollified this coupling by proposing to use the Strang Splitting scheme to compose a drifted ODK with isotropic diffusion and a symmetric anisotropic local kernel; this has the effect of removing this coupling but adds a free parameter $\alpha\in(0,1)$ that one must choose when deciding how to split the diffusion between the two operators.  The choice of $\alpha$ is not arbitrary since if it is made too small or too large one of the two kernels will fail to resolve on a the data; furthermore, if your number of samples is large, $m\gg10^4$, computing the product of three dense matrices of this size will take a not insignificant amount of time, and thus may preclude the use of such a scheme in an iterative framework such as gradient descent in the optimisation schemes of section \ref{subsec:numerical_results_sde}.

\paragraph{Future Directions:} The most obvious future direction is to apply ODKs to the motivating example of scRNA-seq data.  The abundance of so-called \textit{PseudoTime} methods in the scRNA-seq literature \cite{trapnell2014dynamics, haghverdi2016diffusion, luecken2019current} makes the idea of an ordering function almost canonical within the field and allows us to jump straight into exploratory analysis of the dynamics.  Due to the limited dynamical information available, high dimensionality, and inhomogeneous data sampling, this will be a challenging test for ODK; its success or failure will, more generally, indicate where the limits of this and similar kernel methods lie.  However, we note that the method of Diffusion PseudoTime \cite{haghverdi2016diffusion} leverages the work of Coifman \& Lafon \cite{coifman2006diffusion} to construct a diffusion distance on scRNA-seq data constructed through a symmetric bandwidth diffusion kernel \cite{berry2016variable}; this method is considered one of the most reliable and state-of-the-art methods from the PseudoTime literature, thus giving us cause to believe that ODKs can be applied to scRNA-seq data effectively.  Of particular interest is the application of the $L_q$ or $L_p$ objective introduced in section \ref{subsec:numerical_results_sde}; this may allow us to simultaneously infer the potential and the state-dependent diffusion that underlie differentiation. Allowing for state-dependent diffusion would already be a significant advancement in the field, since diffusion is almost always assumed to be constant in the literature.  In this paper, we do not prove that either the $L_q$ or $L_p$ objective is identifiable from our assumed observations; providing such a result, or understanding why it is not possible, would give much-needed context for the sort of scientific statements one can make after running these optimisations (see related work of Lavenant et al. \cite{lavenant2024toward}). Furthermore, we formulated these objectives as simple maximum-likelihoods; it is possible, and indeed probable, that different and better formulations exist that can leverage ODKs' representations of the dynamics to greater effect and potentially improve our chances of identifiability.

\section*{Acknowledgements}
J.H. Soulsby was supported by the EPSRC CASE scholarship in partnership with Pasqal SAS.

\bibliographystyle{plainurl}
\bibliography{refs}

@article{trapnell2014dynamics,
  author    = {Trapnell, Cole and Cacchiarelli, Davide and Grimsby, Jonna 
               and Pokharel, Prapti and Li, Shuqiang and Morse, Michael 
               and Lennon, Niall J and Livak, Kenneth J and Mikkelsen, 
               Tarjei S and Rinn, John L},
  title     = {The dynamics and regulators of cell fate decisions are 
               revealed by pseudotemporal ordering of single cells},
  journal   = {Nature Biotechnology},
  volume    = {32},
  pages     = {381--386},
  year      = {2014},
  doi       = {10.1038/nbt.2859}
}

@article{saelens2019comparison,
  author    = {Saelens, Wouter and Cannoodt, Robrecht and Todorov, Helena 
               and Saeys, Yvan},
  title     = {A comparison of single-cell trajectory inference methods},
  journal   = {Nature Biotechnology},
  volume    = {37},
  pages     = {547--554},
  year      = {2019},
  doi       = {10.1038/s41587-019-0071-9}
}

@article{street2018slingshot,
  author    = {Street, Kelly and Risso, Davide and Fletcher, Russell B 
               and Das, Diya and Ngai, John and Yosef, Nir and Purdom, 
               Elizabeth and Dudoit, Sandrine},
  title     = {Slingshot: cell lineage and pseudotime inference for 
               single-cell transcriptomics},
  journal   = {BMC Genomics},
  volume    = {19},
  pages     = {477},
  year      = {2018},
  doi       = {10.1186/s12864-018-4772-0}
}

@article{haghverdi2016diffusion,
  author    = {Haghverdi, Laleh and B{\"u}ttner, Maren and Wolf, 
               F Alexander and Buettner, Florian and Theis, Fabian J},
  title     = {Diffusion pseudotime robustly reconstructs lineage 
               branching},
  journal   = {Nature Methods},
  volume    = {13},
  pages     = {845--848},
  year      = {2016},
  doi       = {10.1038/nmeth.3971}
}

@article{qiu2022mapping,
  author    = {Qiu, Xiaojie and Zhang, Yan and Martin-Rufino, Jorge D 
               and Weng, Chen and Hosseinzadeh, Shayan and Yang, Dian 
               and Pogson, Angela N and Hein, Marco Y and Min, Kyung Hoi 
               Joseph and Wang, Li and Grody, Emanuelle I and Shurtleff, 
               Matthew J and Yuan, Ruoshi and Xu, Song and Ma, Yian 
               and Replogle, Joseph M and Lander, Eric S and Darmanis, 
               Spyros and Bahar, Ivet and Sankaran, Vijay G and Xing, 
               Jianhua and Weissman, Jonathan S},
  title     = {Mapping transcriptomic vector fields of single cells},
  journal   = {Cell},
  volume    = {185},
  number    = {4},
  pages     = {690--711},
  year      = {2022},
  doi       = {10.1016/j.cell.2021.12.045}
}

@article{luecken2019current,
  author    = {Luecken, Malte D and Theis, Fabian J},
  title     = {Current best practices in single-cell {RNA-seq} analysis: 
               a tutorial},
  journal   = {Molecular Systems Biology},
  volume    = {15},
  number    = {6},
  pages     = {e8746},
  year      = {2019},
  doi       = {10.15252/msb.20188746},
  publisher = {EMBO Press}
}

@article{wang2026regvelo,
  author    = {Wang, Weixu and Hu, Zhiyuan and Weiler, Philipp and 
               Mayes, Sarah and Lange, Marius and Fountain, David M 
               and Haug, Jan Ole and Wang, Jingye and Xue, Zhengyuan 
               and Sauka-Spengler, Tatjana and Theis, Fabian J},
  title     = {{RegVelo}: gene-regulatory-informed dynamics of 
               single cells},
  journal   = {Cell},
  volume    = {189},
  number    = {12},
  pages     = {3773--3800.e44},
  year      = {2026},
  doi       = {10.1016/j.cell.2026.04.022},
  publisher = {Elsevier}
}

@article{ward1954partially, 
  title={Partially ordered topological spaces}, 
  author={Ward, L. E.},
  journal={Proceedings of the American Mathematical Society}, 
  publisher={American Mathematical Society (AMS)},   
  volume={5}, 
  number={1}, 
  ISSN={0002-9939},  
  pages={144–-161},
  year={1954},
  month=Feb,
  DOI={10.1090/s0002-9939-1954-0063016-5},
  url={http://dx.doi.org/10.1090/s0002-9939-1954-0063016-5}
}

@article{berry2016local,
  title={Local kernels and the geometric structure of data},
  author={Berry, Tyrus and Sauer, Timothy},
  journal={Applied and Computational Harmonic Analysis},
  volume={40},
  number={3},
  pages={439--469},
  year={2016},
  month=May,
  publisher={Elsevier BV},
  ISSN={1063-5203}, 
  url={http://dx.doi.org/10.1016/j.acha.2015.03.002}, 
  DOI={10.1016/j.acha.2015.03.002}
}

@article{singer2012vector,
  title={Vector Diffusion Maps and the Connection Laplacian},
  author={Singer, Amit and Wu, H-T},
  journal={Communications on Pure and Applied Mathematics},
  volume={65},
  number={8},
  pages={1067--1144},
  year={2012},
  month=Mar,
  publisher={Wiley},
  ISSN={1097-0312},
  url={http://dx.doi.org/10.1002/cpa.21395}, 
  DOI={10.1002/cpa.21395},
}

@article{lopez2018moore,
  title={Moore-Penrose’s Inverse and Solutions of Linear Systems},
  author={Lopez-Bonilla, J and Lopez-Vanquez, R and Vidal-Beltran, S},
  journal={World Scientific News},
  volume={101},
  pages={246--252},
  year={2018}
}

@article{coifman2006diffusion,
  title={Diffusion Maps},
  author={Coifman, Ronald R. and Lafon, Stéphane},
  journal={Applied and Computational Harmonic Analysis},
  volume={21},
  number={1},
  pages={5--30},
  year={2006},
  month=Jul,
  publisher={Elsevier BV},
  ISSN={1063-5203}, 
  url={http://dx.doi.org/10.1016/j.acha.2006.04.006}, 
  DOI={10.1016/j.acha.2006.04.006}
}

@article{belkin2003laplacian,
  title={Laplacian Eigenmaps for Dimensionality Reduction and Data Representation}, 
  volume={15}, 
  ISSN={1530-888X}, 
  url={http://dx.doi.org/10.1162/089976603321780317}, 
  DOI={10.1162/089976603321780317}, 
  number={6}, 
  journal={Neural Computation}, 
  publisher={MIT Press - Journals}, 
  author={Belkin, Mikhail and Niyogi, Partha}, 
  year={2003}, 
  month=Jun, 
  pages={1373–-1396} 
}

@article{tenenbaum2000global,
  title={A Global Geometric Framework for Nonlinear Dimensionality Reduction}, 
  volume={290}, 
  ISSN={1095-9203}, 
  url={http://dx.doi.org/10.1126/science.290.5500.2319}, 
  DOI={10.1126/science.290.5500.2319}, 
  number={5500}, 
  journal={Science}, 
  publisher={American Association for the Advancement of Science (AAAS)}, 
  author={Tenenbaum, Joshua B. and Silva, Vin de and Langford, John C.}, 
  year={2000}, 
  month=Dec, 
  pages={2319–-2323}
}

@inproceedings{dhillon2004kernel,
  series={KDD04}, 
  title={Kernel k-means: spectral clustering and normalized cuts}, 
  url={http://dx.doi.org/10.1145/1014052.1014118}, 
  DOI={10.1145/1014052.1014118}, 
  booktitle={Proceedings of the tenth ACM SIGKDD international conference on Knowledge discovery and data mining}, 
  publisher={ACM}, 
  author={Dhillon, Inderjit S. and Guan, Yuqiang and Kulis, Brian}, 
  year={2004}, 
  month=Aug, 
  pages={551-–556}, 
  collection={KDD04}
}

@article{sheather1991reliable,
  title={A Reliable Data-Based Bandwidth Selection Method for Kernel Density Estimation}, 
  volume={53}, 
  ISSN={1467-9868}, 
  url={http://dx.doi.org/10.1111/j.2517-6161.1991.tb01857.x}, 
  DOI={10.1111/j.2517-6161.1991.tb01857.x}, 
  number={3}, 
  journal={Journal of the Royal Statistical Society Series B: Statistical Methodology}, 
  publisher={Oxford University Press (OUP)}, 
  author={Sheather, S. J. and Jones, M. C.}, 
  year={1991}, 
  month=Jul, 
  pages={683–-690}
}

@article{hino2017ider,
  title={ider: Intrinsic Dimension Estimation with R}, 
  volume={9}, 
  ISSN={2073-4859}, 
  url={http://dx.doi.org/10.32614/rj-2017-054}, 
  DOI={10.32614/rj-2017-054}, 
  number={2}, 
  journal={The R Journal}, 
  publisher={The R Foundation}, 
  author={Hino, Hideitsu}, 
  year={2017}, 
  pages={329}
}

@book{evans2022partial,
  title={Partial Differential Equations}, 
  ISBN={9781470411442}, 
  ISSN={1065-7339}, 
  url={http://dx.doi.org/10.1090/gsm/019}, 
  DOI={10.1090/gsm/019}, 
  journal={Graduate Studies in Mathematics}, 
  publisher={American Mathematical Society}, 
  author={Evans, Lawrence}, 
  year={2010}, 
  month=Mar
}

@article{berry2016variable,
  title={Variable bandwidth diffusion kernels}, 
  volume={40}, 
  ISSN={1063-5203}, 
  url={http://dx.doi.org/10.1016/j.acha.2015.01.001}, 
  DOI={10.1016/j.acha.2015.01.001}, 
  number={1}, 
  journal={Applied and Computational Harmonic Analysis}, 
  publisher={Elsevier BV}, 
  author={Berry, Tyrus and Harlim, John}, 
  year={2016}, 
  month=Jan, 
  pages={68-–96}
}

@article{schiebinger2015geometry,
    title={The geometry of kernelized spectral clustering}, 
    volume={43}, 
    ISSN={0090-5364}, 
    url={http://dx.doi.org/10.1214/14-aos1283}, 
    DOI={10.1214/14-aos1283}, 
    number={2}, 
    journal={The Annals of Statistics}, 
    publisher={Institute of Mathematical Statistics}, 
    author={Schiebinger, Geoffrey and Wainwright, Martin J. and Yu, Bin}, year={2015}, 
    month=Apr
}

@article{berry2016consistent,
  title={Consistent manifold representation for topological data analysis}, 
  volume={1}, 
  ISSN={2639-8001}, 
  url={http://dx.doi.org/10.3934/fods.2019001}, 
  DOI={10.3934/fods.2019001}, 
  number={1}, 
  journal={Foundations of Data Science}, 
  publisher={American Institute of Mathematical Sciences (AIMS)}, 
  author={Berry, Tyrus and Sauer, Timothy}, 
  year={2019}, 
  pages={1–-38}
}

@article{berry2017density,
  title={Density estimation on manifolds with boundary}, 
  volume={107}, 
  ISSN={0167-9473}, 
  url={http://dx.doi.org/10.1016/j.csda.2016.09.011}, 
  DOI={10.1016/j.csda.2016.09.011}, 
  journal={Computational Statistics \& Data Analysis}, 
  publisher={Elsevier BV}, 
  author={Berry, Tyrus and Sauer, Timothy}, 
  year={2017}, 
  month=Mar, 
  pages={1–-17}
}

@article{banisch2020diffusion,
  title={Diffusion maps tailored to arbitrary non-degenerate Itô processes}, 
  volume={48}, 
  ISSN={1063-5203}, 
  url={http://dx.doi.org/10.1016/j.acha.2018.05.001}, 
  DOI={10.1016/j.acha.2018.05.001}, 
  number={1}, 
  journal={Applied and Computational Harmonic Analysis}, 
  publisher={Elsevier BV}, 
  author={Banisch, Ralf and Trstanova, Zofia and Bittracher, Andreas and Klus, Stefan and Koltai, Péter}, 
  year={2020}, 
  month=Jan, 
  pages={242–-265}
}

@book{chen2011pseudo,
  title={Pseudo-Riemannian Geometry, $\delta$-Invariants and Applications}, 
  ISBN={9789814329644}, 
  url={http://dx.doi.org/10.1142/8003}, 
  DOI={10.1142/8003}, 
  publisher={WORLD SCIENTIFIC}, 
  author={Chen, Bang-Yen}, 
  year={2011}, 
  month=Mar
}

@article{nocedal1980updating,
  title={Updating quasi-Newton matrices with limited storage}, 
  volume={35}, 
  ISSN={0025-5718}, 
  url={http://dx.doi.org/10.1090/s0025-5718-1980-0572855-7}, 
  DOI={10.1090/s0025-5718-1980-0572855-7}, 
  number={151}, 
  journal={Mathematics of Computation}, 
  publisher={American Mathematical Society (AMS)}, 
  author={Nocedal, Jorge}, 
  year={1980}, 
  pages={773-–782}
}

@article{liu1989limited,
  title={On the limited memory BFGS method for large scale optimization}, 
  volume={45}, 
  ISSN={1436-4646}, 
  url={http://dx.doi.org/10.1007/bf01589116}, 
  DOI={10.1007/bf01589116}, 
  number={1-3}, 
  journal={Mathematical Programming}, 
  publisher={Springer Science and Business Media LLC}, 
  author={Liu, Dong C. and Nocedal, Jorge}, 
  year={1989}, 
  month=Aug, 
  pages={503-–528}
}

@misc{NLopt,
  title = {The {NLopt} nonlinear-optimization package},
  author = {Steven G. Johnson},
  year = {2007},
  howpublished = {\url{https://github.com/stevengj/nlopt}}
}

@inbook{shan2022diffusion,
  title={Diffusion Maps: Using the Semigroup Property for Parameter Tuning}, 
  ISBN={9783030458478}, 
  ISSN={2296-5017}, 
  url={http://dx.doi.org/10.1007/978-3-030-45847-8_18}, 
  DOI={10.1007/978-3-030-45847-8_18}, 
  booktitle={Theoretical Physics, Wavelets, Analysis, Genomics}, 
  publisher={Springer International Publishing}, 
  author={Shan, Shan and Daubechies, Ingrid}, 
  year={2022}, 
  month=Dec, 
  pages={409–-424}
}

@article{shi2000normalized,
  title={Normalized cuts and image segmentation}, 
  author={Jianbo Shi and Malik, J.},
  journal={IEEE Transactions on Pattern Analysis and Machine Intelligence}, 
  year={2000},
  volume={22},
  number={8},
  pages={888--905},
  doi={10.1109/34.868688}
}

@inproceedings{ng2001spectral,
  title = {On Spectral Clustering: Analysis and an algorithm},
  author = {Ng, Andrew and Jordan, Michael and Weiss, Yair},
  booktitle = {Advances in Neural Information Processing Systems},
  editor = {T. Dietterich and S. Becker and Z. Ghahramani},
  publisher = {MIT Press},
  url = {https://proceedings.neurips.cc/paper_files/paper/2001/file/801272ee79cfde7fa5960571fee36b9b-Paper.pdf},
  volume = {14},
  year = {2001}
}

@inproceedings{hein2005intrinsic,
  series={ICML ’05}, 
  title={Intrinsic dimensionality estimation of submanifolds in $\mathbb{R}^d$}, 
  url={http://dx.doi.org/10.1145/1102351.1102388}, 
  DOI={10.1145/1102351.1102388}, 
  booktitle={Proceedings of the 22nd International Conference on Machine Learning - ICML ’05}, 
  publisher={ACM Press}, 
  author={Hein, Matthias and Audibert, Jean-Yves},
  year={2005}, 
  pages={289–-296}, 
  collection={ICML ’05}
}

@article{granata2016accurate,
  title={Accurate Estimation of the Intrinsic Dimension Using Graph Distances: Unraveling the Geometric Complexity of Datasets}, 
  volume={6}, 
  ISSN={2045-2322}, 
  url={http://dx.doi.org/10.1038/srep31377}, 
  DOI={10.1038/srep31377}, 
  number={1}, 
  journal={Scientific Reports}, 
  publisher={Springer Science and Business Media LLC}, 
  author={Granata, Daniele and Carnevale, Vincenzo}, 
  year={2016}, 
  month=Aug
}

@article{facco2017estimating,
  title={Estimating the intrinsic dimension of datasets by a minimal neighborhood information}, 
  volume={7}, 
  ISSN={2045-2322}, 
  url={http://dx.doi.org/10.1038/s41598-017-11873-y}, 
  DOI={10.1038/s41598-017-11873-y}, 
  number={1}, 
  journal={Scientific Reports}, 
  publisher={Springer Science and Business Media LLC}, 
  author={Facco, Elena and d’Errico, Maria and Rodriguez, Alex and Laio, Alessandro}, 
  year={2017}, 
  month=Sep
}

@article{di2026scale,
  title={Scale adaptive and robust intrinsic dimension estimation via optimal neighbourhood identification}, 
  volume={16}, 
  ISSN={2045-2322}, 
  url={http://dx.doi.org/10.1038/s41598-026-48005-4}, 
  DOI={10.1038/s41598-026-48005-4}, 
  number={1}, 
  journal={Scientific Reports}, 
  publisher={Springer Science and Business Media LLC}, 
  author={Di Noia, Antonio and Macocco, Iuri and Glielmo, Aldo and Laio, Alessandro and Mira, Antonietta}, 
  year={2026}, 
  month=Apr
}

@article{allegra2020data,
  title={Data segmentation based on the local intrinsic dimension}, 
  volume={10}, 
  ISSN={2045-2322}, 
  url={http://dx.doi.org/10.1038/s41598-020-72222-0}, 
  DOI={10.1038/s41598-020-72222-0}, 
  number={1}, 
  journal={Scientific Reports}, 
  publisher={Springer Science and Business Media LLC}, 
  author={Allegra, Michele and Facco, Elena and Denti, Francesco and Laio, Alessandro and Mira, Antonietta}, 
  year={2020}, 
  month=Oct
}

@article{scholkopf1998nonlinear,
  title={Nonlinear Component Analysis as a Kernel Eigenvalue Problem}, 
  volume={10}, 
  ISSN={1530-888X}, 
  url={http://dx.doi.org/10.1162/089976698300017467}, 
  DOI={10.1162/089976698300017467}, 
  number={5}, 
  journal={Neural Computation}, 
  publisher={MIT Press - Journals}, 
  author={Schölkopf, Bernhard and Smola, Alexander and Müller, Klaus-Robert}, 
  year={1998}, 
  month=Jul, 
  pages={1299–-1319}
}

@article{couillet2016kernel,
  title={Kernel spectral clustering of large dimensional data}, 
  volume={10}, 
  ISSN={1935-7524}, 
  url={http://dx.doi.org/10.1214/16-ejs1144}, 
  DOI={10.1214/16-ejs1144}, 
  number={1}, 
  journal={Electronic Journal of Statistics}, 
  publisher={Institute of Mathematical Statistics}, 
  author={Couillet, Romain and Benaych-Georges, Florent}, 
  year={2016}, 
  month=Jan
}

@article{article,
  title = {Kernel Spectral Clustering},
  author = {Giulini, Ilaria},
  journal={arXiv preprint arXiv:1606.06519},
  year = {2016},
  month = Jun,
  doi = {10.48550/arXiv.1606.06519},
  url={https://arxiv.org/abs/1606.06519}
}

@inproceedings{vankadara2020optimality,
  title={On the optimality of kernels for high-dimensional clustering},
  author={Vankadara, Leena C and Ghoshdastidar, Debarghya},
  booktitle={International Conference on Artificial Intelligence and Statistics},
  pages={2185--2195},
  year={2020},
  organization={PMLR}
}

@inproceedings{pillaud2023kernelized,
  title={Kernelized diffusion maps},
  author={Pillaud-Vivien, Loucas and Bach, Francis},
  booktitle={The Thirty-Sixth Annual Conference on Learning Theory},
  pages={5236--5259},
  year={2023},
  organization={PMLR}
}

@article{belkin2008towards,
  title={Towards a theoretical foundation for Laplacian-based manifold methods}, 
  volume={74}, 
  ISSN={0022-0000}, 
  url={http://dx.doi.org/10.1016/j.jcss.2007.08.006}, 
  DOI={10.1016/j.jcss.2007.08.006}, 
  number={8}, 
  journal={Journal of Computer and System Sciences}, 
  publisher={Elsevier BV}, 
  author={Belkin, Mikhail and Niyogi, Partha}, 
  year={2008}, 
  month=Dec, 
  pages={1289-–1308} 
}

@article{schiebinger2019optimal,
  title={Optimal-Transport Analysis of Single-Cell Gene Expression Identifies Developmental Trajectories in Reprogramming}, 
  volume={176}, 
  ISSN={0092-8674}, 
  url={http://dx.doi.org/10.1016/j.cell.2019.01.006}, 
  DOI={10.1016/j.cell.2019.01.006}, 
  number={4}, 
  journal={Cell}, 
  publisher={Elsevier BV}, 
  author={Schiebinger, Geoffrey and Shu, Jian and Tabaka, Marcin and Cleary, Brian and Subramanian, Vidya and Solomon, Aryeh and Gould, Joshua and Liu, Siyan and Lin, Stacie and Berube, Peter and Lee, Lia and Chen, Jenny and Brumbaugh, Justin and Rigollet, Philippe and Hochedlinger, Konrad and Jaenisch, Rudolf and Regev, Aviv and Lander, Eric S.}, 
  year={2019}, 
  month=Feb, 
  pages={928-–943}
}

@inproceedings{guan2026gradientflow,
    title={Gradient-Flow {SDE}s Have Unique Transient Population Dynamics},
    author={Vincent Guan and Joseph Janssen and Nicolas Lanzetti and Antonio Terpin and Geoffrey Schiebinger and Elina Robeva},
    booktitle={The 29th International Conference on Artificial Intelligence and Statistics},
    year={2026},
    url={https://arxiv.org/abs/2505.21770},
    journal={arXiv preprint arXiv:2505.21770},
    doi={10.48550/arXiv.2505.21770}
}

@book{bogachev2022fokker,
  title={Fokker–Planck–Kolmogorov Equations}, 
  ISBN={9781470427931}, 
  ISSN={2331-7159}, 
  url={http://dx.doi.org/10.1090/surv/207}, 
  DOI={10.1090/surv/207}, 
  journal={Mathematical Surveys and Monographs}, 
  publisher={American Mathematical Society}, 
  author={Bogachev, Vladimir and Krylov, Nicolai and Röckner, Michael and Shaposhnikov, Stanislav}, 
  year={2015}, 
  month=Dec
}

@article{cristofol2017simultaneous,
  title={Simultaneous determination of the drift and diffusion coefficients in stochastic differential equations}, 
  volume={33}, 
  ISSN={1361-6420}, 
  url={http://dx.doi.org/10.1088/1361-6420/aa7a1c}, 
  DOI={10.1088/1361-6420/aa7a1c}, 
  number={9}, 
  journal={Inverse Problems}, 
  publisher={IOP Publishing}, 
  author={Cristofol, Michel and Roques, Lionel}, 
  year={2017}, 
  month=Aug, 
  pages={095006}
}

@article{bergen2020generalizing,
  title={Generalizing RNA velocity to transient cell states through dynamical modeling}, 
  volume={38}, 
  ISSN={1546-1696}, 
  url={http://dx.doi.org/10.1038/s41587-020-0591-3}, 
  DOI={10.1038/s41587-020-0591-3}, 
  number={12}, 
  journal={Nature Biotechnology}, 
  publisher={Springer Science and Business Media LLC}, 
  author={Bergen, Volker and Lange, Marius and Peidli, Stefan and Wolf, F. Alexander and Theis, Fabian J.}, 
  year={2020}, 
  month=Aug, 
  pages={1408–-1414}
}

@article{jin2018scepath,
  title={scEpath: energy landscape-based inference of transition probabilities and cellular trajectories from single-cell transcriptomic data}, 
  volume={34}, 
  ISSN={1367-4811}, 
  url={http://dx.doi.org/10.1093/bioinformatics/bty058}, 
  DOI={10.1093/bioinformatics/bty058}, 
  number={12}, 
  journal={Bioinformatics}, 
  publisher={Oxford University Press (OUP)}, 
  author={Jin, Suoqin and MacLean, Adam L and Peng, Tao and Nie, Qing}, 
  editor={Berger, Bonnie}, 
  year={2018}, 
  month=Feb, 
  pages={2077–-2086}
}

@article{bengio2013representation,
  title={Representation Learning: A Review and New Perspectives}, 
  volume={35}, 
  ISSN={2160-9292}, 
  url={http://dx.doi.org/10.1109/tpami.2013.50}, 
  DOI={10.1109/tpami.2013.50}, 
  number={8}, 
  journal={IEEE Transactions on Pattern Analysis and Machine Intelligence}, 
  publisher={Institute of Electrical and Electronics Engineers (IEEE)}, 
  author={Bengio, Y. and Courville, A. and Vincent, P.}, 
  year={2013}, 
  month=Aug, 
  pages={1798–-1828}
}

@article{fefferman2016testing,
  title={Testing the manifold hypothesis}, 
  volume={29}, 
  ISSN={0894-0347}, 
  url={http://dx.doi.org/10.1090/jams/852}, 
  DOI={10.1090/jams/852}, 
  number={4}, 
  journal={Journal of the American Mathematical Society}, 
  publisher={American Mathematical Society (AMS)}, 
  author={Fefferman, Charles and Mitter, Sanjoy and Narayanan, Hariharan}, 
  year={2016}, 
  month=Feb, 
  pages={983-–1049}
}

@article{strang1968construction,
  title={On the Construction and Comparison of Difference Schemes}, 
  volume={5}, 
  ISSN={1095-7170}, 
  url={http://dx.doi.org/10.1137/0705041}, 
  DOI={10.1137/0705041}, 
  number={3}, 
  journal={SIAM Journal on Numerical Analysis}, 
  publisher={Society for Industrial & Applied Mathematics (SIAM)}, 
  author={Strang, Gilbert}, 
  year={1968}, 
  month=Sep, 
  pages={506-–517}
}

@book{hairer2006numerical,
  series={Springer Series in Computational Mathematics},
  title={Geometric Numerical Integration: Structure-Preserving Algorithms for Ordinary Differential Equations},
  author={Hairer, Ernst and Wanner, Gerhard and Lubich, Christian},
  ISBN={3540306633}, 
  url={http://dx.doi.org/10.1007/3-540-30666-8}, 
  DOI={10.1007/3-540-30666-8}, 
  publisher={Springer Berlin, Heidelberg}, 
  year={2006},
  edition={2}
}

@article{trotter1959product,
  title={On the product of semi-groups of operators}, 
  volume={10}, 
  ISSN={0002-9939}, 
  url={http://dx.doi.org/10.1090/s0002-9939-1959-0108732-6}, 
  DOI={10.1090/s0002-9939-1959-0108732-6}, 
  number={4}, 
  journal={Proceedings of the American Mathematical Society}, 
  publisher={American Mathematical Society (AMS)}, 
  author={Trotter, H. F.}, 
  year={1959}, 
  month=Aug, 
  pages={545–-551}
}

@inproceedings{belkin2012toward,
  title={Toward Understanding Complex Spaces: Graph Laplacians on Manifolds with Singularities and Boundaries},
  author={Belkin, Mikhail and Que, Qichao and Wang, Yusu and Zhou, Xueyuan},
  booktitle={Proceedings of the 25th Annual Conference on Learning Theory},
  pages={36.1--36.26},
  year={2012},
  editor={Mannor, Shie and Srebro, Nathan and Williamson, Robert C.},
  volume={23},
  series={Proceedings of Machine Learning Research},
  month=Jun,
  publisher={PMLR},
  url={https://proceedings.mlr.press/v23/belkin12.html},
}

@article{hoffmann2008noise,
  title={Noise-Driven Stem Cell and Progenitor Population Dynamics}, 
  volume={3}, 
  ISSN={1932-6203}, 
  url={http://dx.doi.org/10.1371/journal.pone.0002922}, 
  DOI={10.1371/journal.pone.0002922}, 
  number={8}, 
  journal={PLoS ONE}, 
  publisher={Public Library of Science (PLoS)}, 
  author={Hoffmann, Martin and Chang, Hannah H. and Huang, Sui and Ingber, Donald E. and Loeffler, Markus and Galle, Joerg}, 
  editor={Secomb, Timothy}, 
  year={2008}, 
  month=Aug, 
  pages={e2922}
}

@book{chung1997spectral,
  title={Spectral graph theory},
  author={Chung, F. R. K},
  volume={92},
  year={1997},
  publisher={American Mathematical Society},
  DOI={10.1090/cbms/092}
}

@article{kolmogoroff1931analytischen,
  title={{\"U}ber die analytischen Methoden in der Wahrscheinlichkeitsrechnung},
  volume={104}, 
  ISSN={1432-1807}, 
  url={http://dx.doi.org/10.1007/bf01457949}, 
  DOI={10.1007/bf01457949}, 
  number={1}, 
  journal={Mathematische Annalen}, 
  publisher={Springer Science and Business Media LLC}, 
  author={Kolmogoroff, A.}, 
  year={1931}, 
  month=Dec, 
  pages={415-–458}
}

@book{ledoux2001concentration,
  title={The concentration of measure phenomenon},
  author={Ledoux, Michel},
  number={89},
  year={2001},
  publisher={American Mathematical Society},
  DOI={10.1090/surv/089}
}

@book{grandis2009directed,
  title={Directed Algebraic Topology: Models of Non-Reversible Worlds}, 
  ISBN={9780511657474}, 
  url={http://dx.doi.org/10.1017/cbo9780511657474}, 
  DOI={10.1017/cbo9780511657474}, 
  publisher={Cambridge University Press}, 
  author={Grandis, Marco}, 
  year={2009}, 
  month=Sep
}

@article{krishnan2009convenient,
  title={A Convenient Category of Locally Preordered Spaces}, 
  volume={17}, 
  ISSN={1572-9095}, 
  url={http://dx.doi.org/10.1007/s10485-008-9140-9}, 
  DOI={10.1007/s10485-008-9140-9}, 
  number={5}, 
  journal={Applied Categorical Structures}, 
  publisher={Springer Science and Business Media LLC}, 
  author={Krishnan, Sanjeevi}, 
  year={2008}, 
  month=Jun, 
  pages={445–-466}
}

@article{haucourt2012streams,
  title={Streams, d-Spaces and Their Fundamental Categories}, 
  volume={283}, 
  ISSN={1571-0661}, 
  url={http://dx.doi.org/10.1016/j.entcs.2012.05.008}, 
  DOI={10.1016/j.entcs.2012.05.008}, 
  journal={Electronic Notes in Theoretical Computer Science}, 
  publisher={Elsevier BV}, 
  author={Haucourt, Emmanuel}, 
  year={2012}, 
  month=Jun, 
  pages={111–-151}
}

@article{zhou2011behavior,
  title={Behavior of Graph Laplacians on Manifolds with Boundary}, 
  author={Xueyuan Zhou and Mikhail Belkin},
  year={2011},
  DOI={10.48550/arXiv.1105.3931},
  journal={arXiv preprint arXiv:1105.3931},
  url={https://arxiv.org/abs/1105.3931},
}

@book{robert2004monte,
  title={Monte Carlo Statistical Methods},
  author={Christian P. Robert and George Casella},
  ISSN={1431-875X},
  DOI={10.1007/978-1-4757-4145-2},
  edition={2},
  year={2004},
  series={Springer Texts in Statistics},
  publisher={Springer New York, NY}
}

@inproceedings{belkin2006convergence,
  title={Convergence of Laplacian Eigenmaps}, 
  ISBN={9780262256919}, 
  url={http://dx.doi.org/10.7551/mitpress/7503.003.0021}, 
  DOI={10.7551/mitpress/7503.003.0021}, 
  booktitle={Advances in Neural Information Processing Systems}, 
  editor = {B. Sch\"{o}lkopf and J. Platt and T. Hoffman},
  volume={19},
  publisher={The MIT Press}, 
  author={Belkin, Mikhail and Niyogi, Partha}, 
  year={2007}, 
  month=Sep, 
  pages={129-–136}
}

@article{gine2006empirical,
  title={Empirical Graph Laplacian Approximation of Laplace–Beltrami Operators: Large Sample Results}, 
  journal = {Lecture Notes-Monograph Series},
  ISBN={9780940600676}, 
  ISSN = {07492170},
  url={http://dx.doi.org/10.1214/074921706000000888}, 
  DOI={10.1214/074921706000000888}, 
  booktitle={High Dimensional Probability}, 
  publisher={Institute of Mathematical Statistics}, 
  author={Giné, Evarist and Koltchinskii, Vladimir}, 
  year={2006}, 
  pages={238-–259},
  volume={51}
}

@article{singer2006graph,
  title={From graph to manifold Laplacian: The convergence rate}, 
  volume={21}, 
  ISSN={1063-5203}, 
  url={http://dx.doi.org/10.1016/j.acha.2006.03.004}, 
  DOI={10.1016/j.acha.2006.03.004}, 
  number={1}, 
  journal={Applied and Computational Harmonic Analysis}, 
  publisher={Elsevier BV}, 
  author={Singer, A.}, 
  year={2006}, 
  month=Jul, 
  pages={128-–134}
}

@article{hein2007graph,
  author = {Hein, Matthias and Audibert, Jean-Yves and Luxburg, Ulrike von},
  title = {Graph Laplacians and their Convergence on Random Neighborhood Graphs},
  year = {2007},
  issue_date = {12/1/2007},
  publisher = {JMLR},
  volume = {8},
  ISSN = {1532-4435},
  journal = {Journal of Machine Learning Research},
  month = Dec,
  pages = {1325–-1370}
}

@article{wang2015spectral,
  title={Spectral Convergence Rate of Graph Laplacian},
  author={Wang, Xu},
  journal={arXiv preprint arXiv:1510.08110},
  year={2015},
  DOI={10.48550/arXiv.1510.08110},
  url={https://arxiv.org/abs/1510.08110}
}

@inproceedings{ting2011analysis,
  author = {Ting, Daniel and Huang, Ling and Jordan, Michael I.},
  title = {An analysis of the convergence of graph Laplacians},
  year = {2010},
  isbn = {9781605589077},
  publisher = {Omnipress},
  address = {Madison, WI, USA},
  booktitle = {Proceedings of the 27th International Conference on International Conference on Machine Learning},
  pages = {1079-–1086},
  numpages = {8},
  location = {Haifa, Israel},
  series = {ICML'10}
}

@article{tang2018limit,
  title={Limit theorems for eigenvectors of the normalized Laplacian for random graphs}, 
  volume={46}, 
  ISSN={0090-5364}, 
  url={http://dx.doi.org/10.1214/17-aos1623}, 
  DOI={10.1214/17-aos1623}, 
  number={5}, 
  journal={The Annals of Statistics}, 
  publisher={Institute of Mathematical Statistics}, 
  author={Tang, Minh and Priebe, Carey E.}, 
  year={2018}, 
  month=Oct
}

@article{inagaki2025spectral,
  title={Spectral convergence of graph Laplacians with Ricci curvature bounds and in non-collapsed Ricci limit spaces}, 
  author={Masato Inagaki},
  year={2025},
  url={https://arxiv.org/abs/2506.07427}, 
  journal={arXiv preprint arXiv:2506.07427},
  DOI={10.48550/arXiv.2506.07427}
}

@inproceedings{zelnik2004self,
  author = {Zelnik-manor, Lihi and Perona, Pietro},
  booktitle = {Advances in Neural Information Processing Systems},
  editor = {L. Saul and Y. Weiss and L. Bottou},
  publisher = {MIT Press},
  title = {Self-Tuning Spectral Clustering},
  url = {https://proceedings.neurips.cc/paper_files/paper/2004/file/40173ea48d9567f1f393b20c855bb40b-Paper.pdf},
  volume = {17},
  year = {2004}
}

@inbook{elworthy1998stochastic, 
  author={Elworthy, K. D.},
  title={Stochastic Differential Equations on Manifolds},
  booktitle={Probability Towards 2000},
  series={Lecture Notes in Statistics},
  year={1998},
  publisher={Springer New York},
  address={New York, NY},
  pages={165--178},
  isbn={978-1-4612-2224-8},
  doi={10.1007/978-1-4612-2224-8_10},
  url={https://doi.org/10.1007/978-1-4612-2224-8_10},
  volume={128}
}

@article{hinton2006reducing,
  title={Reducing the Dimensionality of Data with Neural Networks}, 
  volume={313}, 
  ISSN={1095-9203}, 
  url={http://dx.doi.org/10.1126/science.1127647}, 
  DOI={10.1126/science.1127647}, 
  number={5786}, 
  journal={Science}, 
  publisher={American Association for the Advancement of Science (AAAS)}, 
  author={Hinton, G. E. and Salakhutdinov, R. R.}, 
  year={2006}, 
  month=Jul, 
  pages={504–-507}
}

@article{vincent2010stacked,
  author  = {Pascal Vincent and Hugo Larochelle and Isabelle Lajoie and Yoshua Bengio and Pierre-Antoine Manzagol},
  title   = {Stacked Denoising Autoencoders: Learning Useful Representations in a Deep Network with a Local Denoising Criterion},
  journal = {Journal of Machine Learning Research},
  year    = {2010},
  volume  = {11},
  number  = {110},
  pages   = {3371--3408},
  url     = {http://jmlr.org/papers/v11/vincent10a.html}
}

@article{kingma2013auto,
  title={Auto-encoding variational bayes},
  author={Kingma, Diederik P and Welling, Max},
  journal={arXiv preprint arXiv:1312.6114},
  year={2022},
  url={https://arxiv.org/abs/1312.6114},
  DOI={10.48550/arXiv.1312.6114}
}

@article{edelsbrunner2002topological,
  title={Topological Persistence and Simplification}, 
  volume={28}, 
  ISSN={1432-0444}, 
  url={http://dx.doi.org/10.1007/s00454-002-2885-2}, 
  DOI={10.1007/s00454-002-2885-2}, 
  number={4}, 
  journal={Discrete \& Computational Geometry}, 
  publisher={Springer Science and Business Media LLC}, 
  author={Edelsbrunner and Letscher and Zomorodian}, 
  year={2002}, 
  month=Nov, 
  pages={511-–533}
}

@inproceedings{zomorodian2004computing,
  series={SoCG04}, 
  title={Computing persistent homology}, 
  url={http://dx.doi.org/10.1145/997817.997870}, 
  DOI={10.1145/997817.997870},
  booktitle={Proceedings of the twentieth annual symposium on Computational geometry}, 
  publisher={ACM}, 
  author={Zomorodian, Afra and Carlsson, Gunnar}, 
  year={2004}, 
  month=Jun, 
  pages={347-–356}, 
  collection={SoCG04}
}

@article{carlsson2009topology,
  title={Topology and data}, 
  volume={46}, 
  ISSN={0273-0979}, 
  url={http://dx.doi.org/10.1090/s0273-0979-09-01249-x}, 
  DOI={10.1090/s0273-0979-09-01249-x}, 
  number={2}, 
  journal={Bulletin of the American Mathematical Society}, 
  publisher={American Mathematical Society (AMS)}, 
  author={Carlsson, Gunnar}, 
  year={2009}, 
  month=Jan, 
  pages={255-–308}
}

@article{loftsgaarden1965nonparametric,
  title={A Nonparametric Estimate of a Multivariate Density Function}, 
  volume={36}, 
  ISSN={0003-4851}, 
  url={http://dx.doi.org/10.1214/aoms/1177700079}, 
  DOI={10.1214/aoms/1177700079}, 
  number={3}, 
  journal={The Annals of Mathematical Statistics}, 
  publisher={Institute of Mathematical Statistics}, 
  author={Loftsgaarden, D. O. and Quesenberry, C. P.}, 
  year={1965}, 
  month=Jun,
  pages={1049-–1051}
}

@article{leon2013gram,
  title={Gram‐Schmidt orthogonalization: 100 years and more}, 
  volume={20}, 
  ISSN={1099-1506}, 
  url={http://dx.doi.org/10.1002/nla.1839}, 
  DOI={10.1002/nla.1839}, 
  number={3}, 
  journal={Numerical Linear Algebra with Applications}, 
  publisher={Wiley}, 
  author={Leon, Steven J. and Björck, {\r A}ke and Gander, Walter}, 
  year={2012}, 
  month=Jun, 
  pages={492–-532}
}

@article{peleg1970utility,
  title={Utility Functions for Partially Ordered Topological Spaces}, 
  volume={38}, 
  ISSN={0012-9682}, 
  url={http://dx.doi.org/10.2307/1909243}, 
  DOI={10.2307/1909243}, 
  number={1}, 
  journal={Econometrica}, 
  publisher={JSTOR}, 
  author={Peleg, Bezalel}, 
  year={1970}, 
  month=Jan, 
  pages={93}
}

@article{evren2011multi,
  title={On the multi-utility representation of preference relations}, 
  volume={47}, 
  ISSN={0304-4068}, 
  url={http://dx.doi.org/10.1016/j.jmateco.2011.07.003}, 
  DOI={10.1016/j.jmateco.2011.07.003}, 
  number={4-5}, 
  journal={Journal of Mathematical Economics}, 
  publisher={Elsevier BV}, 
  author={Evren, {\"O}zgür and Ok, Efe A.}, 
  year={2011}, 
  month=Aug, 
  pages={554-–563}
}

@article{troue1969topology,
  title={Topology and Order by Leopoldo Nachbin. Van Nostrand Mathematical Studies, Van Nostrand 1965 (Princeton)},
  author={Trou{\'e}, J},
  journal={Canadian Mathematical Bulletin},
  volume={12},
  number={2},
  pages={242--244},
  year={1969},
  publisher={Cambridge University Press}
}

@article{bosi2012continuous,
  title={Continuous multi-utility representations of preorders}, 
  volume={48}, 
  ISSN={0304-4068}, 
  url={http://dx.doi.org/10.1016/j.jmateco.2012.05.001}, 
  DOI={10.1016/j.jmateco.2012.05.001}, 
  number={4}, 
  journal={Journal of Mathematical Economics}, 
  publisher={Elsevier BV}, 
  author={Bosi, Gianni and Herden, Gerhard}, 
  year={2012}, 
  month=Aug, 
  pages={212-–218}
}

@article{fajstrup2006algebraic,
  title = {Algebraic topology and concurrency},
  journal = {Theoretical Computer Science},
  volume = {357},
  number = {1},
  pages = {241--278},
  year = {2006},
  note = {Clifford Lectures and the Mathematical Foundations of Programming Semantics},
  issn = {0304-3975},
  doi = {https://doi.org/10.1016/j.tcs.2006.03.022},
  url = {https://www.sciencedirect.com/science/article/pii/S030439750600274X},
  author = {Lisbeth Fajstrup and Martin Raußen and Eric Goubault},
}

@phdthesis{belkin2004problems,
  author = {Belkin, Mikhail},
  advisor = {Niyogi, Partha},
  title = {Problems of learning on manifolds},
  year = {2003},
  school = {The University of Chicago},
  note = {AAI3097083}
}

@book{lehmann1998theory,
  title={Theory of Point Estimation}, 
  ISBN={9781461580713}, 
  url={http://dx.doi.org/10.1007/978-1-4615-8071-3}, 
  DOI={10.1007/978-1-4615-8071-3}, 
  publisher={Springer US}, 
  author={Lehmann, E. L.}, 
  year={1991}
}

@inproceedings{hashimoto2016learning,
  title = 	 {Learning Population-Level Diffusions with Generative RNNs},
  author = 	 {Hashimoto, Tatsunori and Gifford, David and Jaakkola, Tommi},
  booktitle = 	 {Proceedings of The 33rd International Conference on Machine Learning},
  pages = 	 {2417--2426},
  year = 	 {2016},
  editor = 	 {Balcan, Maria Florina and Weinberger, Kilian Q.},
  volume = 	 {48},
  series = 	 {Proceedings of Machine Learning Research},
  address = 	 {New York, New York, USA},
  month = 	 {20--22 Jun},
  publisher =    {PMLR},
  url = 	 {https://proceedings.mlr.press/v48/hashimoto16.html},
}

@article{lavenant2024toward,
  title={Toward a mathematical theory of trajectory inference}, 
  volume={34}, 
  ISSN={1050-5164}, 
  url={http://dx.doi.org/10.1214/23-aap1969}, 
  DOI={10.1214/23-aap1969}, 
  number={1A}, 
  journal={The Annals of Applied Probability}, 
  publisher={Institute of Mathematical Statistics}, 
  author={Lavenant, Hugo and Zhang, Stephen and Kim, Young-Heon and Schiebinger, Geoffrey},
  year={2024}, 
  month=Feb
}

@inproceedings{chizat2022trajectory, 
  series={NeurIPS 2022}, 
  title={Trajectory Inference Via Mean-Field Langevin in Path Space}, 
  url={http://dx.doi.org/10.52202/068431-1217}, 
  DOI={10.52202/068431-1217}, 
  booktitle={Advances in Neural Information Processing Systems 35}, 
  publisher={Neural Information Processing Systems Foundation, Inc. (NeurIPS)}, 
  author={Chizat, Lénaïc and Zhang, Stephen and Heitz, Matthieu and Schiebinger, Geoffrey}, 
  year={2022}, 
  pages={16731–-16742}, 
  collection={NeurIPS 2022} 
}

@inproceedings{neklyudov2023action,
  title={Action Matching: A Variational Method for Learning Stochastic Dynamics from Samples},
  author={Kirill Neklyudov and Daniel Severo and Alireza Makhzani},
  booktitle={NeurIPS 2022 Workshop on Score-Based Methods},
  year={2022},
  url={https://openreview.net/forum?id=rXiZMBJBdB},
}

@article{onsager1953fluctuations,
  title={Fluctuations and Irreversible Processes}, 
  volume={91}, 
  ISSN={0031-899X}, 
  url={http://dx.doi.org/10.1103/physrev.91.1505}, 
  DOI={10.1103/physrev.91.1505}, 
  number={6}, 
  journal={Physical Review}, 
  publisher={American Physical Society (APS)}, 
  author={Onsager, L. and Machlup, S.}, 
  year={1953}, 
  month=Sep, 
  pages={1505-–1512}
}

@article{dijkstra2022note,
  title={A note on two problems in connexion with graphs}, 
  volume={1}, 
  ISSN={0945-3245}, 
  url={http://dx.doi.org/10.1007/bf01386390}, 
  DOI={10.1007/bf01386390}, 
  number={1}, 
  journal={Numerische Mathematik}, 
  publisher={Springer Science and Business Media LLC}, 
  author={Dijkstra, E. W.}, 
  year={1959}, 
  month=Dec, 
  pages={269–-271}
}

@article{warshall1962theorem,
  title={A Theorem on Boolean Matrices}, 
  volume={9}, 
  ISSN={1557-735X}, 
  url={http://dx.doi.org/10.1145/321105.321107}, 
  DOI={10.1145/321105.321107}, 
  number={1}, 
  journal={Journal of the ACM}, 
  publisher={Association for Computing Machinery (ACM)}, 
  author={Warshall, Stephen}, 
  year={1962}, 
  month=Jan, 
  pages={11-–12}
}

@article{floyd1962algorithm,
  title={Algorithm 97: Shortest path}, 
  volume={5}, 
  ISSN={1557-7317}, 
  url={http://dx.doi.org/10.1145/367766.368168}, 
  DOI={10.1145/367766.368168}, 
  number={6}, 
  journal={Communications of the ACM}, 
  publisher={Association for Computing Machinery (ACM)}, 
  author={Floyd, Robert W.}, 
  year={1962}, 
  month=Jun, 
  pages={345-–345}
}

@book{bird2014transport,
  title = {Introductory {Transport} {Phenomena}},
  isbn = {978-1-118-95371-6},
  url = {https://books.google.co.uk/books?id=wwXYBgAAQBAJ},
  publisher = {Wiley},
  author = {Bird, R.B. and Stewart, W.E. and Lightfoot, E.N. and Klingenberg, D.J.},
  year = {2015}
}

@article{wilson1931distribution,
  author = {Edwin B. Wilson  and Margaret M. Hilferty },
  title = {The Distribution of Chi-Square},
  journal = {Proceedings of the National Academy of Sciences},
  volume = {17},
  number = {12},
  pages = {684--688},
  year = {1931},
  doi = {10.1073/pnas.17.12.684},
  URL = {https://www.pnas.org/doi/abs/10.1073/pnas.17.12.684}
}

@article{rosenblatt1956,
  title={Remarks on Some Nonparametric Estimates of a Density Function}, 
  volume={27}, 
  ISSN={0003-4851}, 
  url={http://dx.doi.org/10.1214/aoms/1177728190}, 
  DOI={10.1214/aoms/1177728190}, 
  number={3}, 
  journal={The Annals of Mathematical Statistics}, 
  publisher={Institute of Mathematical Statistics}, 
  author={Rosenblatt, Murray}, 
  year={1956}, 
  month=Sep, 
  pages={832-–837}
}

@article{parzen1962estimation,
  title={On Estimation of a Probability Density Function and Mode}, 
  volume={33}, 
  ISSN={0003-4851}, 
  url={http://dx.doi.org/10.1214/aoms/1177704472}, 
  DOI={10.1214/aoms/1177704472}, 
  number={3}, 
  journal={The Annals of Mathematical Statistics}, 
  publisher={Institute of Mathematical Statistics}, 
  author={Parzen, Emanuel}, 
  year={1962}, 
  month=Sep, 
  pages={1065-–1076}
}

@book{silverman2018density,
  title={Density Estimation for Statistics and Data Analysis}, 
  ISBN={9781315140919}, 
  url={http://dx.doi.org/10.1201/9781315140919}, 
  DOI={10.1201/9781315140919},
  publisher={Routledge}, 
  author={Silverman, B.W.}, 
  year={2018}, 
  month=Feb
}

@book{scott2015multivariate,
  title={Multivariate Density Estimation: Theory, Practice, and Visualization}, 
  ISBN={9781118575574}, 
  ISSN={1940-6347}, 
  url={http://dx.doi.org/10.1002/9781118575574}, 
  DOI={10.1002/9781118575574}, 
  journal={Wiley Series in Probability and Statistics}, 
  publisher={Wiley}, 
  author={Scott, David W.}, 
  year={2015}, 
  month=Mar
}

@inbook{debreu1983representation, 
  title={Representation of a preference ordering by a numerical function}, 
  ISBN={9781139052092}, 
  url={http://dx.doi.org/10.1017/ccol052123736x.007}, 
  DOI={10.1017/ccol052123736x.007}, 
  booktitle={Mathematical Economics}, 
  publisher={Cambridge University Press}, 
  author={Debreu, Gerard}, 
  pages={105-–110},
  chapter={6},
  year={1983}
}

\newpage
\appendix
\renewcommand{\theequation}{\thesection.\arabic{equation}}
\numberwithin{equation}{section}
\section{Proofs of Main Results}\label{appendix:proofs}
To complete the proofs of section \ref{sec:main_results}, we will use two Lemmas that we take from \cite{coifman2006diffusion}, which were initially formulated in the PhD Thesis of Belkin \cite{belkin2004problems}; however, we state them in the more convenient form that Coifman \& Lafon provide.  We define and use two systems of coordinates on our manifold $d$-dimensional Riemannian manifold $\mathcal{M}$: the normal coordinates --- constructed implicitly using the exponential map --- which we denote by $(s_1,\dots,s_d)$; and the projection coordinates --- constructed extrinsically using the implicit function theorem to define a neighbourhood where we can project the embedded manifold onto a $d$-dimensional affine subspace of $\mathbb{R}^n$ --- which we denote by $(u_1,\dots,u_d)$.  For a much more rigorous discussion of the construction of both of these coordinate charts, please consult appendix \ref{appBsec:prelims}.

The first Lemma gives an expression for the components of the normal coordinates, in terms of the components of the projection coordinates:
\begin{lemma}[Conversion of Coordinates \cite{coifman2006diffusion}]\label{lemma:coordinate_transformation}
    For any $x\in \mathcal{M}$, for sufficiently small $\varepsilon$ and $y\in B(x,\varepsilon^{1/2})$ (the closed ball centred as $x$ with radius $\varepsilon^{1/2}$), we can write the normal coordinate $s_i$ as:
    \begin{equation}\label{appAeq:coordinate_transformation}
        s_i=u_i+p_{x,3}(u)+\bigO(\varepsilon^2)
    \end{equation}
    for $i=1,\dots,d$; and where $p_{x,k}(u)$ is a homogeneous polynomial of degree $k$ whose coefficients depend on $x$.\footnote{To be precise $p_{x,k}(u)$ will be some polynomial of the form $\sum_{\alpha:|\alpha|=k} c_\alpha\prod_{i=1}^d u_i^{\alpha_i}$, where $\alpha=(\alpha_1,\dots,\alpha_d)\in\mathbb{N}_{>0}^d$ is a multi-index and $|\alpha|=\alpha_1+\dots+\alpha_d$.}
\end{lemma}
The second Lemma has two parts: the first gives an expression for the Euclidean distance $\lVert y-x\rVert^2$ in terms of projection coordinates, where we interpret this not as a geodesic distance on $\mathcal{M}$ but a literal Euclidean distance in $\mathbb{R}^n$, i.e. if $\iota:\mathcal{M}\to\mathbb{R}^n$ is an embedding of our manifold into ambient space then for two nearby points $x,y\in\mathcal{M}$, we mean $\lVert y-x\rVert=\lVert \iota(y)-\iota(x)\rVert$; the second, gives us a change of coordinates formula for the volume form $\mathrm{d}y$ (which we simply write as $\mathrm{d}y$) in terms of projection coordinates.
\begin{lemma}[Metric \& Volume Comparison\cite{coifman2006diffusion}]
    For any $x\in \mathcal{M}$, for sufficiently small $\varepsilon$ and $y\in B(x,\varepsilon^{1/2})$, we have:
    \begin{equation}\label{appAeq:metric_comparison}
        i).\quad \lVert y-x\rVert^2=\lVert u\rVert^2+p_{x,4}(u)+p_{x,5}(u)+\bigO(\varepsilon^3)
    \end{equation}
    \begin{equation}\label{appAeq:volume_comparison}
        ii).\quad \det\left(\frac{dy}{du}\right)=1+p_{x,2}(u)+p_{x,3}(u)+\bigO(\varepsilon^2)
    \end{equation}
\end{lemma}
If $f:\mathcal{M}\to\mathbb{R}$ is some scalar function defined on the manifold, then the Taylor Expansion about and $x\in\mathcal{M}$ for $y$ in a sufficiently small ball $B(x,\varepsilon^{1/2})$ around $x$, would read (using Einstein summation notation) as:
\begin{equation}
    f(y)=\tilde f(0)+s_i\frac{\partial \tilde{f}}{\partial s_i}(0)+\frac{1}{2}s_is_j\frac{\partial^2\tilde{f}}{\partial s_i\partial s_j}(0)+\frac{1}{6}s_is_js_k\frac{\partial^3\tilde{f}}{\partial s_i\partial s_j\partial s_k}(0)+R_{x,4}(s)
\end{equation}
where $\tilde{f}$ views $f$ as a function of the normal coordinates $(s_1,\dots,s_d)$ around $x$ which means $\tilde{f}(0)=f(x)$; and $R_{x,4}(s)$ is the remainder term for $f\in C^3(\mathcal{M})$ which we note is of order $\bigO(\lVert s\rVert^4)$.  Using \ref{lemma:coordinate_transformation}, we can write this in projection coordinates as:
\begin{equation}
    f(y)=f(x)+u_i\frac{\partial \tilde{f}}{\partial s_i}(0)+\frac{1}{2}u_iu_j\frac{\partial^2\tilde{f}}{\partial s_i\partial s_j}(0)+p_{x,3}(u)+\bigO(\lVert u\rVert^4)
\end{equation}
Using this, we prove the following Lemma, which puts the expression inside of $h(\cdot)$ in the definition of a standard ODK, given in equation \eqref{eq:general_odk}, in terms of projection coordinates:
\begin{lemma}\label{lemma:odk_metric_expansion}
    For any $x\in\mathcal{M}$, for sufficiently small $\varepsilon$ and $y\in B(x,\varepsilon^{1/2})$, we have:
    \begin{align}
        (\star)\coloneqq \lVert y-x\rVert^2&+\frac{(l(y)-l(x)-\varepsilon\tau(x)\kappa(x))^2}{\kappa(x)^2}-\frac{(l(y)-l(x))^2}{\kappa(x)^2}\label{appAeq:odk_prehexpression}\\
        &=\lVert u-\varepsilon\mu(x)\rVert^2+\Theta_4(u)+\Theta_5(u)+\bigO(\varepsilon^3,\varepsilon\lVert u\rVert^4)
    \end{align}
    where $\mu(x)=\frac{\tau(x)}{\kappa(x)}\nabla l(x)$, and $\Theta_4(u)$ and $\Theta_5(u)$ are expressions of order $\bigO(\varepsilon^2)$ and $\bigO(\varepsilon^{5/2})$, respectively.
\end{lemma}
\begin{proof}
    We note that the expression given in equation \eqref{appAeq:odk_prehexpression}, which we now denote with $(\star)$, can be written as:
    \begin{equation}
        (\star)=\lVert y-x\rVert^2-2\varepsilon\tau(x)\left(\frac{l(y)-l(x)}{\kappa(x)}\right)+\varepsilon^2\tau(x)^2
    \end{equation}
    We can take a Taylor expansion of $l(y)$ around $x$ in projection coordinates to arrive at:
    \begin{equation}
        l(y)=l(x)+u_i\frac{\partial \tilde{l}}{\partial u_i}(0)+\frac{1}{2}u_iu_j\frac{\partial^2\tilde{l}}{\partial u_i\partial u_j}+p_{x,3}(u)+\bigO(\lVert u\rVert^4)
    \end{equation}
    where we have used the coordinate transformation \eqref{appAeq:coordinate_transformation} to transform from $s\mapsto u$, and have noted that $R_{x,3}(s)=p_{x,3}(u)+\bigO(\lVert u\rVert^4)$.  Now using the metric comparison equation of \eqref{appAeq:metric_comparison} we can expand:
    \begin{align}
        (\star)&=\lVert u\rVert^2+p_{x,4}(u)+p_{x,5}(u)-2\varepsilon\frac{\tau(x)}{\kappa(x)}\left(u_i\frac{\partial \tilde{l}}{\partial u_i}(0)+\frac{1}{2}u_iu_j\frac{\partial^2\tilde{l}}{\partial u_i\partial u_j}+p_{x,3}(u)\right) +\varepsilon^2\tau(x)^2+\bigO(\varepsilon^3,\varepsilon\lVert u\rVert^4)\\
        &=\lVert u\rVert^2-2\varepsilon\frac{\tau(x)}{\kappa(x)}u_i\frac{\partial\tilde{l}}{\partial s_i}(0)+\varepsilon^2\frac{\tau(x)^2}{\kappa(x)^2}\lVert\nabla l(x)\rVert^2+\Theta_4(u)+\Theta_5(u)+\bigO(\varepsilon^3,\varepsilon\lVert u\rVert^4)\label{appAeq:l_lemma_intermediate}
    \end{align}
    where we have added an subtracted $\varepsilon^2\tfrac{\tau(x)^2}{\kappa(x)^2}$, and defined:
    \begin{align}
        \Theta_4(u)&=p_{x,4}(u)+\varepsilon^2\tau(x)^2\left(1-\frac{\lVert \nabla l(x)\rVert^2}{\kappa(x)^2}\right) - 2\varepsilon\frac{\tau(x)}{\kappa(x)}u_iu_j\frac{\partial^2\tilde{l}}{\partial s_i\partial s_j}\\
        \Theta_5(u)&=p_{x,5}(u)-2\varepsilon\frac{\tau(x)}{\kappa(x)}p_{x,3}(u)
    \end{align}
    Continuing from \eqref{appAeq:l_lemma_intermediate} we see that:
    \begin{equation}
        \lVert u\rVert^2-2\varepsilon\frac{\tau(x)}{\kappa(x)}u_i\frac{\partial\tilde{l}}{\partial s_i}(0)+\varepsilon^2\frac{\tau(x)^2}{\kappa(x)^2}\lVert\nabla l(x)\rVert^2=\left\langle u-\varepsilon\frac{\tau(x)}{\kappa(x)}\nabla l(x), u-\varepsilon\frac{\tau(x)}{\kappa(x)}\nabla l(x)\right\rangle
    \end{equation}
    and thus setting $\mu(x)=\tfrac{\tau(x)}{\kappa(x)}\nabla l(x)$ completes the proof.
\end{proof}
With Lemma \ref{lemma:odk_metric_expansion} in hand, we will now complete the proof of the main result:
\begin{theorem}[Expansion of Ordered Diffusion Kernels]\label{appAthm:odk_backward_expansion}
    Suppose $\mathcal{M}$ is a $d$-dimensional Riemannian manifold embedded in $\mathbb{R}^n$ and with metric inherited from the standard Euclidean metric in ambient space. Let $l\in C^2(\mathcal{M})$ be an ordering function, and let $K_\varepsilon(x,y)$ be the corresponding ODK as defined in \eqref{eq:general_odk}. Then for any $f\in C^3(\mathcal{M})$ the associated family of backward integral operators $\{G_\varepsilon\}_{\varepsilon> 0}$ admits the following asymptotic expansion as $\varepsilon\to0$ when applied to $f$:
    \begin{equation}
        G_\varepsilon f(x) = m(x)f(x) + \varepsilon\left( \omega_l(x)f(x) + m(x)\mathcal{L}f(x) \right)  + \varepsilon^{3/2}\Omega_l(x)+ \bigO(\varepsilon^2),
    \end{equation}
    at a point $x\in\mathcal{M}$ a distance larger than $\varepsilon^\gamma$ from $\partial\mathcal{M}$, for a fixed $\gamma\in(0,1/2)$; and where $\omega_l(x)$ depends on the manifold and the kernel $K_\varepsilon(x,\cdot)$ at $x$, and $\Omega_l(x)$ depends on the manifold, the kernel $K_\varepsilon(x,\cdot)$ at $x$, and also the test function $f(x)$. The differential operator $\mathcal{L}$ is a second-order elliptic operator of the form:
    \begin{equation}\label{appAeq:backward_differential_operator}
        \mathcal{L}f(x) = \mu(x)\cdot\nabla f(x) + D(x)\Delta f(x)
    \end{equation}
    where $m(x)$, $\mu(x)$ and $D(x)$, satisfy:
    \begin{align}
        m(x)&=\int_{T_x\mathcal{M}}h\left(\tfrac{\lVert z\rVert^2}{\rho(x)^2}\right)\mathrm{d}z\\
        \mu(x)&=\frac{\tau(x)}{\kappa(x)}\nabla l(x)\label{appAeq:odk_backward_expansion_drift}\\
        D(x)&=\ \frac{1}{2}m(x)^{-1}\int_{T_x\mathcal{M}}z_1^2\ h\left(\tfrac{\lVert z\rVert^2}{\rho(x)^2}\right)\mathrm{d}z\backwardscoloneqq \tfrac{1}{2}\sigma_h(x)^2\rho(x)^2\label{appAeq:odk_backward_expansion_diffusion}\\
        &\ \ \quad\text{where}\quad\sigma_h(x)^2=\frac{1}{\int_{T_x\mathcal{M}}h(\lVert z\rVert^2)\dif z }\int_{T_x\mathcal{M}}z_1^2h(\lVert z\rVert^2)\dif z.\label{appAeq:sigma_h_def}
    \end{align}
    The normalised backward integral operator can thus be expanded in terms of $\varepsilon$ as:
    \begin{equation}
        P_\varepsilon f(x) = f(x)+\varepsilon\mathcal{L}f(x)+\bigO(\varepsilon^{3/2}).
    \end{equation}
\end{theorem}
\begin{proof}[Proof of Theorem \ref{appAthm:odk_backward_expansion}]
    Take any $\gamma\in (0,1/2)$, and let $x\in\mathcal{M}$ be a distance of atleast $\varepsilon^\gamma$ from the boundary of $\mathcal{M}$.  Since, by assumption, an ODK is bounded by an exponential, there will exist parameters $C_h,\sigma>0$ such that $K_\varepsilon(x,y)=h(-\star)\leq C_h\exp(-\star/\sigma)$.  As a result, the contribution to integral in $G_\varepsilon f(x)$ that is outside the ball of radius $\varepsilon^\gamma$ can be bounded as:
    \begin{align}
        &\left|\varepsilon^{-d/2}\int_{y\in\mathcal{M}:\lVert y-x\rVert>\varepsilon^\gamma} K_\varepsilon(x,y)f(y)\ \mathrm{d} y\right|\leq\varepsilon^{-d/2}\lVert f\rVert_\infty\int_{y\in\mathcal{M}:\lVert y-x\rVert>\varepsilon^\gamma}K_\varepsilon(x,y)\ \mathrm{d} y\\
        &\leq \varepsilon^{-d/2}C_h \lVert f\rVert_\infty\int_{y\in\mathcal{M}:\lVert y-x\rVert>\varepsilon^\gamma}\exp\left(-\frac{\lVert y-x\rVert^2}{\varepsilon\sigma\rho(x)^2}+\frac{2(l(y)-l(x))\tau(x)}{\sigma\rho(x)^2\kappa(x)}-\frac{\varepsilon\tau(x)^2}{\sigma\rho(x)^2}\right)\ \mathrm{d} y\\
        &\leq MC_h\lVert f\rVert_\infty \int_{\lVert z\rVert>\varepsilon^{\gamma-1/2}}\exp\left(-\frac{\lVert z\rVert^2}{\varepsilon\sigma\rho(x)^2}\right)\ \mathrm{d}y=MC_h\lVert f\rVert_\infty Q\left(\tfrac{d}{2},c\varepsilon^{2\gamma-1}\right)\label{appA:proofint_sub}\\
        &\leq MC_h\lVert f\rVert_\infty e^{-c\varepsilon^{2\gamma -1}}\mathcal{P}(\varepsilon^{2\gamma-1})=\bigO(\varepsilon^2)\label{appA:proofint_exppoly}
    \end{align}
    where $M=\exp(4\lVert l\rVert_\infty\tau(x)/(\sigma\rho(x)^2\kappa(x)))$, $c=1/(\rho(x)^2\sigma)$, and in \eqref{appA:proofint_sub} we made the substitution $y-x=\sqrt{\varepsilon} z$; the function $Q(s,r)$ is the regularised upper incomplete gamma function.  In the final line, we use the fact that when $r>s-1$, we have $Q(s,r)\leq e^{-r}r^s[\Gamma(s)(r-(s-1))]^{-1}$; supposing $\varepsilon<\left(\tfrac{2c}{d}\right)^{1/(1-2\gamma)}$ then we can set $\mathcal{P}(\varepsilon^{2\gamma-1})=\Gamma\left(\tfrac{d}{2}\right)^{-1}c^{d/2}\left(\varepsilon^{2\gamma-1}\right)^{d/2}$ thus equation \eqref{appA:proofint_exppoly} can be bounded by any power of $\varepsilon$ so we take this to be $\varepsilon^2$.  
    
    Using the previous Lemma, we have:
    \begin{equation}\label{appAeq:odk_backward_intermed_kerexp}
        h\left(\frac{(\star)}{\varepsilon\rho(x)^2}\right)=h\left(\frac{\lVert u-\varepsilon\tau\kappa^{-1}\nabla l\rVert^2}{\varepsilon\rho(x)^2}+\frac{\Theta_4(u)+\Theta_5(u)}{\varepsilon\rho(x)^2}+\bigO(\varepsilon^2,\lVert u\rVert^4)\right)
    \end{equation}
    Now expanding $h$ about $(\varepsilon\rho^2)^{-1}\lVert u-\varepsilon\tau\kappa^{-1}\nabla l\rVert^2$, we have:
    \begin{equation}
        h\left(\frac{\lVert u-\varepsilon\tau\kappa^{-1}\nabla l\rVert^2}{\varepsilon\rho^2}\right)+\frac{\Theta_4(u)+\Theta_5(u)}{\varepsilon\rho^2}h'\left(\frac{\lVert u-\varepsilon\tau\kappa^{-1}\nabla l\rVert^2}{\varepsilon\rho^2}\right)+\bigO(\varepsilon^2,\lVert u\rVert^4)
    \end{equation}
    For convenience, let $\xi$ denote $(\varepsilon\rho^2)^{-1}\lVert u-\varepsilon\tau\kappa^{-1}\nabla l\rVert^2$, we can then write:
    \begin{equation}\label{proof:odk_three_term_expansion}
        \varepsilon^{d/2}G_\varepsilon f(x)=\int h\left(\xi\right)f(y)\mathrm{d}y+\int\frac{\Theta_4(u)}{\varepsilon\rho^2}h'(\xi)f(y)\mathrm{d}y+\int\frac{\Theta_5(u)}{\varepsilon\rho^2}h'(\xi)f(y)\mathrm{d}y
    \end{equation}
    where integrals are over the domain $B(x,\varepsilon^\gamma)$.  The first term is a local kernel, and taking a Taylor expansion of $f$ about $x$, the second term becomes:
    \begin{equation}
        f(x)\int \frac{\Theta_4(u)}{\varepsilon\rho(x)^2}h'(\xi)\mathrm{du}+\int\frac{\Theta_4(u)}{\varepsilon\rho(x)^2}h'(\xi)u_i\frac{\partial \tilde{f}}{\partial s_i}\mathrm{du}+\bigO(\varepsilon^2,\lVert u\rVert^4)
    \end{equation}
    and the term becomes:
    \begin{equation}
        f(x)\int\frac{\Theta_5(u)}{\varepsilon\rho(x)^2}h'(\xi)\mathrm{du}+\bigO(\varepsilon^2,\lVert u\rVert^4)
    \end{equation}
    respectively.  Since the tails of the integral $\varepsilon^{-d/2}\int K_\varepsilon(x,y)f(y)\mathrm{d}y$ when extended to $\mathbb{R}^d$ are of order $\bigO(\varepsilon^2)$, we extend the three integrals above to the full tangent space.  Since the first term is a local kernel, we will have:
    \begin{equation}\label{appA:odk_proof_local_kernel_expansion}
        \varepsilon^{-d/2}\int_{\mathbb{R}^d} h\left(\frac{\lVert u-\varepsilon\tau\kappa^{-1}\nabla l\rVert^2}{\varepsilon\rho^2}\right)f(y)\mathrm{d}y=m(x)f(x)+\varepsilon(\omega(x)f(x)+\tilde{\mathcal{L}}f(x))+\varepsilon^{3/2}\Omega(x)+\bigO(\varepsilon^2)
    \end{equation}
    where $\tilde{\mathcal{L}}=\tilde{\mu}(x)\cdot\nabla+\tilde{D}(x)\Delta$ and the zeroth-, first- and second-order moments are as given in equation (1) of \cite{berry2016local} and take the form:
    \begin{align}
        m(x)&=\int_{T_x\mathcal{M}}h\left(\tfrac{\lVert z\rVert^2}{\rho(x)^2}\right)\mathrm{dz}\\
        \tilde{\mu}(x)&=m(x)\tfrac{\tau(x)}{\kappa(x)}\nabla l(x)\backwardscoloneqq m(x)\mu(x)\\
        \tilde{D}(x)&=m(x)\cdot m(x)^{-1}\frac{1}{2}\int_{T_x\mathcal{M}}z_1^2\cdot h\left(\tfrac{\lVert z\rVert^2}{\rho(x)^2}\right)\mathrm{dz}\backwardscoloneqq m(x)D(x)
    \end{align}
    From this we see that $\tilde{\mathcal{L}}=\tilde{\mu}(x)\cdot\nabla+\frac{1}{2}\tilde{D}(x)\Delta=m(x)\left[\mu(x)\cdot\nabla+D(x)\Delta\right]\backwardscoloneqq m(x)\mathcal{L}$.  Now making the substitution $u\mapsto\varepsilon^{1/2}z$ we have:
    \begin{align}
        \Theta_4(u)&\mapsto \varepsilon^2\left(p_{x,4}(z)+\tau(x)^2\left(1-\frac{\lVert\nabla l(x)\rVert^2}{\kappa(x)^2}\right)-2\frac{\tau(x)}{\kappa(x)}z_iz_j\frac{\partial^2 \tilde{l}}{\partial s_i\partial s_j}\right)\\
        \Theta_5(u)&\mapsto \varepsilon^{5/2}\left(p_{x,5}(z)-2\frac{\tau(x)}{\kappa(x)}p_{x,3}(z)\right)
    \end{align}
    Defining $\Theta'_4(z)$ and $\Theta'_5(z)$ to be the two terms in parentheses, respectively, and noting that $h'(\xi)$ contributes no additional powers of $\varepsilon$ after this substitution, we can define:
    \begin{align}
        \omega_l(x)&\coloneqq\omega(x)+\int\frac{\Theta'_4(z)}{\rho(x)^2}h'(\xi)\mathrm{dz}\\
        \Omega_l(x)&\coloneqq \Omega(x)+\int\frac{\Theta'_4(z)}{\rho(x)^2}h'(\xi)z_i\frac{\partial\tilde{f}}{\partial s_i}\mathrm{dz}+f(x)\int\frac{\Theta'_5(z)}{\rho(x)^2}h'(\xi)\mathrm{dz}
    \end{align}
    Putting this together with the local kernels expansion of equation \eqref{appA:odk_proof_local_kernel_expansion}, we arrive at the expansion:
    \begin{equation}
        G_\varepsilon f(x)=m(x)f(x)+\varepsilon(\omega_l(x)f(x)+m(x)\mathcal{L}f(x))+\varepsilon^{3/2}\Omega_l(x)+\bigO(\varepsilon^2)
    \end{equation}
    
    This completes the first part of the proof.  To verify the expansion of the normalised integral operator $P_\varepsilon f(x)=(G_\varepsilon\mathbb{1}(x))^{-1}G_\varepsilon f(x)$, we note that:
    \begin{equation}
        G_\varepsilon \mathbb{1}(x) = m(x)\left(1+\varepsilon\frac{\omega_l(x)}{m(x)}+\varepsilon^{3/2}\frac{\Omega_l(x)}{m(x)}+\bigO(\varepsilon^2)\right)
    \end{equation}
    where we have used the fact that $\mathcal{L}\mathbb{1}\equiv0$.  For sufficiently small $\varepsilon$ we can use the expansion $\frac{1}{1+x}=1-x+x^2-\dots$, and find that, up to order $\varepsilon^{3/2}$:
    \begin{align}
        P_\varepsilon f(x)&=\frac{G_\varepsilon f(x)}{G_\varepsilon\mathbb{1}(x)}=m(x)^{-1}\cdot(m(x)f(x)+\varepsilon\omega_l(x)f(x)+\varepsilon m(x)\mathcal{L}f(x))\left(1-\varepsilon\frac{\omega_l(x)}{m(x)}\right)+\bigO(\varepsilon^{3/2})\\
        &=f(x)+\varepsilon\frac{\omega_l(x)}{m(x)}f(x)+\varepsilon\mathcal{L}f(x)-\varepsilon\frac{\omega_l(x)}{m(x)}f(x)+\bigO(\varepsilon^{3/2})\\
        &=f(x)+\varepsilon\mathcal{L}f(x)+\bigO(\varepsilon^{3/2})
    \end{align}
    which completes the proof.
\end{proof}

The associated expansion for the adjoint operator is done in the weak sense, which other shows that up to the resolution of the $L^2(\mathcal{M})$ inner product, the operator $A$ is equal to the operator $B$ when applied to some $f\in L^2(\mathcal{M})\cap\text{Dom}(A)\cap\text{Dom}(B)$ is for any $g\in L^2(\mathcal{M})\cap\text{Dom}(A)\cap\text{Dom}(B)$, we have $\langle g,Af\rangle_{L^2(\mathcal{M})}=\langle g,Bf\rangle_{L^2(\mathcal{M})}$.  To do this, we must restrict the domain of $f$ to be in $C^3(\mathcal{M})\cap L^2(\mathcal{M})$, and thus the Theorem and proof are as follows:
\begin{theorem}\label{appAthm:odk_forward_expansion}
    Let $K_\varepsilon(x,y)$ be an Ordered Diffusion Kernel under the same aassumptions as Theorem \ref{appAthm:odk_backward_expansion}, then, in the weak sense, the family of adjoint integral operator $\{G_\varepsilon^*\}_{\varepsilon>0}$ admit the following asymptotic expansion as $\varepsilon\to 0$ when applied to a test function $f\in C^3(\mathcal{M})\cap L^2(\mathcal{M})$:
    \begin{equation}
        G_\varepsilon^*f(x)\coloneqq\varepsilon^{-d/2}\int_\mathcal{M}K_\varepsilon(y,x)f(y)\mathrm{d}y=m(x)f(x)+\varepsilon(\omega_l(x)f(x)+\mathcal{L}^*(fm)(x))+\bigO(\varepsilon^{3/2}),
    \end{equation}
    where $\mathcal{L}^*$ is the formal adjoint of $\mathcal{L}$ given in Theorem \ref{appAthm:odk_backward_expansion}, this operator takes the form (\cite{evans2022partial}, p320):
    \begin{equation}
        \mathcal{L}^*f(x)=-\text{div}\left(f(x)\mu(x)\right)+\Delta(D(x)f(x)),
    \end{equation}
    where $\mu(x)$ and $D(x)$ are the normalised first and second moments of $K_\varepsilon$, as given in Theorem \ref{appAthm:odk_backward_expansion}. The normalised adjoint integral operator can thus be expanded in terms of $\varepsilon$ as:
    \begin{equation}
        P_\varepsilon^*f(x)=f(x)+\varepsilon\mathcal{L}^*f(x)+\bigO(\varepsilon^{3/2}).
    \end{equation}
\end{theorem}
\begin{proof}[Proof of Theorem \ref{thm:odk_forward_expansion}]
    Let $g\in C^3(\mathcal{M})\cap L^2(\mathcal{M})$ be an arbitrary test function, then:
    \begin{align}
        \langle g,G_\varepsilon^*f\rangle &=\int_\mathcal{M}g(x)G_\varepsilon^*f(x)\mathrm{d}x=\int_\mathcal{M}\int_\mathcal{M}g(x)K_\varepsilon(y,x)f(y)\mathrm{d}y\mathrm{d}x\\
        &=\int_\mathcal{M}f(y)\int_\mathcal{M}K_\varepsilon(y,x)g(x)\mathrm{d}x\mathrm{d}y=\int_\mathcal{M}f(y)G_\varepsilon g(y)\mathrm{d}y\label{appAeq:adjoint_intermediary1}\\
        &=\int_\mathcal{M}f(y)\left[m(y)g(y)+\varepsilon(\omega_l(y)g(y)+\tilde{\mathcal{L}}g(y))+\bigO(\varepsilon^{3/2})\right]\mathrm{d}y\\
        &=\langle f,mg+\varepsilon\omega_lg\rangle+\varepsilon\langle\tilde{\mathcal{L}}g,f\rangle+\bigO(\varepsilon^{3/2})=\left\langle g,mf+\varepsilon(\omega_lf+\mathcal{L}^*(mf))+\bigO(\varepsilon^{3/2})\right\rangle\label{appAeq:adjoint_intermediary2}
    \end{align}
    where in equation \eqref{appAeq:adjoint_intermediary1} we have switched the order of integration and in equation \eqref{appAeq:adjoint_intermediary2} we have used the fact $\langle f,gh\rangle=\langle fh,g\rangle$ and, recalling that $\tilde{\mathcal{L}}g(x)=m(x)\mathcal{L}g(x)$, we have used the fact that:
    \begin{equation}
        \langle \tilde{\mathcal{L}} g,f\rangle=\langle m\mathcal{L}g,f\rangle=\langle\mathcal{L}g,mf\rangle=\langle g,\mathcal{L}^*(mf)\rangle
    \end{equation}
    Therefore $\langle g,G_\varepsilon^*f\rangle=\langle g,mf+\varepsilon(\omega_l f+\mathcal{L}^*(mf))+\bigO(\varepsilon^{3/2})\rangle$, so in the weak sense:
    \begin{equation}
        G_\varepsilon^*f(x)=mf+\varepsilon(\omega_lf+\mathcal{L}^*(mf)\rangle+\bigO(\varepsilon^{3/2})
    \end{equation}
    which completes the first part of the proof.  

    To prove the expansion for the associated normalised adjoint operator, we again employ the expansion $(G_\varepsilon\mathbb{1}(x))^{-1}=m(x)^{-1}\left(1-\varepsilon\frac{\omega_l(x)}{m(x)}\right)+\bigO(\varepsilon^{3/2})$, and:
    \begin{align}\label{appAeq:proof_adjoint_laplacian_intermediary}
        G_\varepsilon^*((G_\varepsilon \mathbb{1})^{-1}f)(x)&=(G_\varepsilon\mathbb{1}(x))^{-1}\left[m(x)f(x)+\varepsilon\omega_l(x)f(x)\right]+\varepsilon\mathcal{L}^*\left(\frac{mf}{G_\varepsilon\mathbb{1}}\right)(x)+\bigO(\varepsilon^{3/2})\\
        &=f(x)+\varepsilon\mathcal{L}^*f(x)+\bigO(\varepsilon^{3/2})
    \end{align}
    which completes the proof.
\end{proof}

It is worth making a short practical comment about $\kappa(x)$, which we introduced as a function of utility, acting more like a hyperparameter of the kernel than a parameter needing to be inferred; we did this because it allows one to decouple the scale of the ordering function (which has no effect on the LQO, $\leq_l$, that it implies) from the magnitude of the drift-term, thus separating the task of inferring direction and magnitude cleanly.  If one wishes to make the kernel entirely agnostic to the scale of $l(x)$, then one should set $\kappa(x)=\lVert\nabla l(x)\rVert$.  In the proofs of Theorem \ref{thm:odk_backward_expansion}, we assumed that $\kappa(x)$ was bounded away from zero, but if $l(x)$ is chosen arbitrarily and $\kappa(x)=\lVert\nabla l(x)\rVert$, then this might not always be the case, thus creating a singularity within the kernel.  If this occurs let $\mathcal{M}_c$ denote a connected critical submanifold of $l(x)$ - where we note $\text{dim}(\mathcal{M}_c)$ can equal any value from $0$ (an isolated critical point) to $d$ (an entire critical region) - then there are three possible courses of action: firstly, one could attempt to mollify $l(x)$ into a new function $\tilde l(x)$ that implies an identical ordering but avoids the criticality, this is only guaranteed to be possible when the critical submanifold is of dimension $(d-1)$ and $l(x)$ is strictly monotonic on the leaves of the implied foliation of some open region containing the critical manifold.  More precisely, there must exist a $d$-dimensional neighbourhood $\mathcal{U}$ of the $(d-1)$-dimensional critical submanifold, $\mathcal{M}_c$, such that there exists a submersion $\varphi:\mathcal{U}\to\mathbb{R}$ with $\varphi^{-1}(0)=\mathcal{M}_c$, and a smooth strictly monotonic function $f:\varphi(\mathcal{U})\to l(\mathcal{U})$ such that $l|_U=f\circ\varphi$; the mollified $\tilde l(x)$ is then constructed as a smooth gluing of $l|_{\mathcal{M}\setminus\mathcal{U}}$ and the submersion $\varphi$.  In all other cases, one can either enforce a uniform lower bound on $\kappa(x)$ (i.e. define $\tilde \kappa(x)=\max(\kappa(x),c)$ for $c>0$ as small as the data allows), this will result in a new drift-term $\tau(x)\tfrac{\nabla l(x)}{\tilde\kappa(x)}$ which strictly underestimates the magnitude of the velocity in the set $\{x\in\mathcal{M}:0<\kappa(x)<c\}$, but on the critical manifold - since we know the drift term will be zero\footnote{The only critical manifold in which the ordering function $l(x)$ implies a direction is the aforementioned monotonic $(d-1)$-dimensional critical submanifold case.} - the dynamics will be as desired; if $\text{dim}(\mathcal{M}_c)=d$ this is the only available correction.  The final course of action that can be applied if $\text{dim}(\mathcal{M}_c)<d$ is to simply cut out the critical manifold(s) and treat them like boundaries of the new manifold $\mathcal{M}\setminus\mathcal{M}_c$, thus referring us back to the previous discussion.  Unless the critical manifold has very specific dynamical significance, and thus underestimating the dynamics would be damaging, it is likely more practical to simply apply a uniform lower bound on $\kappa(x)$ or to mollify away the critical submanifold if possible or to change the form of $\kappa(x)$ altogether.  When working with real data, in all cases except where $l(x)$ has been constructed to have a $d$-dimensional critical submanifold, it is highly unlikely that we would ever numerically compute the derivative of $l(x)$ to be so close to zero as to cause a singularity in our kernel.  This is simply due to the fact that whenever $\text{dim}(\mathcal{M}_c)<d$ we will have $\text{vol}(\mathcal{M}_c)=0$, is a set of measure zero and thus is unlikely to bother us in practice; indeed depending on how one numerically computes $\lVert\nabla l(x)\rVert$ is may already be reasonable to clamp $\kappa(x)$ between a lower and upper value to avoid anomalous estimates skewing ones results.

\subsection{Local-Ordered Diffusion Kernels}
To prove the expansion result of the locally Ordered Diffusion Kernel (lODK), we only assume that the local ordering function, $l(x,y)$, is defined on a restricted domain in the second argument; thus, there is an additional technical consideration that we will discuss first before completing the proof without ambiguity.  In the proof of Theorem \ref{thm:odk_backward_expansion} we take $\gamma\in(0,\tfrac{1}{2})$ and consider the integral of the kernel about $x$ over $\mathcal{M}$ in a ball of radius $\varepsilon^\gamma$.  But the local ordering function at $x\in\mathcal{M}$, written $l_x:U_x\to\mathbb{R}:y\mapsto l(x,y)$, is only assumed to be defined on a restricted domain; it is always possible to extend $l_x:U_x\to\mathbb{R}$ smoothly to all of $\mathcal{M}$, so the integral in $G_\varepsilon f$ outside of $B(x,\varepsilon^\gamma)$ are still well-defined and of order $\bigO(\varepsilon^2)$, as discussed in the proof of Theorem \ref{appAthm:odk_backward_expansion}, but we need to make sure that the integral inside the ball is independent of this extension.

In light of the discussion of the local ordering function given in Appendix \ref{appBsubsec:ordering_functions}, we know that the domain of $l_x$ is an open neighbourhood of $x$ over which the exponential map must be a diffeomorphism; we can frame this in terms of the injectivity radius of $\exp_x$, which is defined as:
\begin{equation}
    \text{inj}(x)=\sup\{r>0:\exp_x|_{B(0,r)}\text{ is a diffeomorphism onto its image.}\}   
\end{equation}
Thus, we can require that $U_x\subseteq \exp_x^{-1}(\text{inj}(x))$, although we note that the exponential map will be diffeomorphic on a \textit{star-shaped} open set and thus we are taking a restriction by requiring $U_x$ to be a ball.  We assume in the statement of this Theorem that $\mathcal{M}$ has bounded geometry, which means there exists some $c>0$ such that $\inf_{x\in\mathcal{M}}\text{inj}(x)=c$, and thus we can take the domain of each $U_x=\exp_x^{-1}(B(0,c))$ and by nature of the exponential map being an isometry we need to choose $\varepsilon$ small enough such that $\varepsilon\leq c^{1/\gamma}$; such a $\varepsilon$ guarantees that the asymptotic expansion of $G_\varepsilon f(x)$ is valid uniformly across all $x$.  When we compute the integral in $G_\varepsilon f(x)$ we are implicitly taking any bounded $C^2(\mathcal{M}\times\mathcal{M})$ extension of $l(x,y)$, and taking the integral over all $y\in\mathcal{M}$.

After making this consideration, we can state and prove the expansion of lODKs in much the same manner as we have done previously:
\begin{theorem}[Expansion of locally Ordered Diffusion Kernels]\label{appAthm:nongradflow_odk_expansion}
    Let $\mathcal{M}$ be a $d$-dimensional Riemannian manifold with bounded geometry embedded in $\mathbb{R}^n$, and with the metric inherited from the standard Euclidean metric.  Let $l(x,y)$ be a $C^2(\text{Dom}(l))$ local ordering function of two inputs and let $b(x)=\nabla l_x(x)$.  For the lODK $K_\varepsilon(x,y)$, as defined in \ref{def:nongradflow_odk}, and for any $f\in C^3(\mathcal{M})$ we have the following asymptotic expansion of the associated family of backward integral operators $\{G_\varepsilon\}_{\varepsilon>0}$ as $\varepsilon\to 0$, applied to $f$:
    \begin{equation}
        G_\varepsilon f(x) = m(x)f(x) + \varepsilon\left( \omega_l(x)f(x) + m(x)\mathcal{L}_bf(x) \right)  + \varepsilon^{3/2}\Omega_l(x)+ \bigO(\varepsilon^2)
    \end{equation}
    at a point $x\in\mathcal{M}$ a distance larger than $\varepsilon^\gamma$ from $\partial\mathcal{M}$, for a fixed $\gamma\in(0,1/2)$.  The second-order differential operator $\mathcal{L}_b$ takes the form:
    \begin{equation}
        \mathcal{L}_b f(x)=\tfrac{\tau(x)}{\kappa(x)}b(x)\cdot\nabla f(x)+\frac{1}{2}D(x)\Delta f(x)
    \end{equation}
    In the above, $m(x)$ and $D(x)$ are as in Theorem \ref{appAthm:odk_backward_expansion}.  The normalised integral operators and their Laplacians, defined in the usual way, can be expanded as:
    \begin{equation}
        P_\varepsilon f(x) =f(x)+\varepsilon\mathcal{L}_b f(x)+\bigO(\varepsilon^{3/2})\implies L_\varepsilon f(x)=\mathcal{L}_b f(x)+\bigO(\varepsilon^{1/2})
    \end{equation}
    Analogous to the previous theorems, the family of adjoint operators, $\{G_\varepsilon^*\}_{\varepsilon>0}$, satisfies the following asymptotic expansion, in the weak sense, when applied to $f\in C^3(\mathcal{M})\cap L^2(\mathcal{M})$:
    \begin{equation}
        G_\varepsilon^*f(x)=m(x)f(x)+\varepsilon(\omega_l(x)f(x)+\mathcal{L}_b^*(fm)(x))+\bigO(\varepsilon^{3/2})
    \end{equation}
    where $\mathcal{L}_b^*$ is the formal adjoint of $\mathcal{L}_b$ and takes the form:
    \begin{equation}
        \mathcal{L}_b^* f(x)=-\text{div}\left(f(x)\tfrac{\tau(x)}{\kappa(x)}b(x)\right)+\Delta(D(x)f(x))
    \end{equation}
    The normalised adjoint integral operator $P_\varepsilon^*f=G_\varepsilon^*((G_\varepsilon \mathbb{1})^{-1} f)$ and its Laplacian can be expanded as:
    \begin{equation}
        P_\varepsilon^* f(x)=f(x)+\varepsilon\mathcal{L}_b^* f(x)+\bigO(\varepsilon^{3/2})\implies L_\varepsilon^* f(x)=\mathcal{L}_b^* f(x)+\bigO(\varepsilon^{1/2})
    \end{equation}
\end{theorem}
\begin{proof}[Proof of Theorem \ref{thm:nongradflow_odk_expansion}]
    Fix some $\gamma\in (0,1/2)$ and take any $x\in\mathcal{M}$ that is at least a distance $\varepsilon^\gamma$ from the boundary of $\mathcal{M}$; noting that a local ordering function is defined such that $l(x,x)=0$ for any $x\in\mathcal{M}$, we can write:
    \begin{equation}
        l(x,y)=l_x(y)=u_i\frac{\partial\tilde{l_x}}{\partial u_i}+\frac{1}{2}u_iu_j\frac{\partial^2\tilde{l_x}}{\partial u_i\partial u_j}+p_{x,3}(u)+\bigO(\lVert u\rVert^4)
    \end{equation}
    by definition of $b(x)$ we have $\frac{\partial \tilde{l_x}}{\partial u_i}(0)=b_i(x)$; therefore our kernel can be written as:
    \begin{equation}
        h\left(\frac{\lVert y-x\rVert^2}{\varepsilon\rho(x)^2}-\frac{\tau(x)l(x,y)}{\rho(x)^2\kappa(x)}+\frac{\varepsilon\tau(x)^2}{\rho(x)^2}\right)=h\left(\frac{\lVert u-\varepsilon b(x)\rVert^2}{\varepsilon\rho(x)^2}+\frac{\Theta_4(u)+\Theta_5(u)}{\varepsilon\rho(x)^2}+\bigO(\varepsilon^2,\lVert u\rVert^4)\right)
    \end{equation}
    Which we see is identical to equation \eqref{appAeq:odk_backward_intermed_kerexp} in the proof of Theorem \ref{appAthm:odk_backward_expansion}; so after taking any arbitrary $C^2$ extension of $l_x(y)$ to all of $\mathcal{M}$ we can conclude immediately that:
    \begin{equation}
        G_\varepsilon f(x) = m(x)f(x) + \varepsilon\left( \omega_l(x)f(x) + m(x)\mathcal{L}_bf(x) \right)  + \varepsilon^{3/2}\Omega_l(x)+ \bigO(\varepsilon^2)
    \end{equation}
    and:
    \begin{equation}
        P_\varepsilon f(x) =f(x)+\varepsilon\mathcal{L}_b f(x)+\bigO(\varepsilon^{3/2})\implies L_\varepsilon f(x)=\mathcal{L}_b f(x)+\bigO(\varepsilon^{1/2})
    \end{equation}
    where $\mathcal{L}_b=\frac{\tau}{\kappa}b\cdot \nabla+D(x)\Delta$, and $D(x)=\tfrac{1}{2}\sigma_h(x)^2\rho(x)^2$ as before.  Setting $\tau(x)=\kappa(x)\equiv 1$ gives the general form of the infinitesimal generator of an It{\^o} SDE with drift $b(x)$ and isotropic diffusion $\frac{1}{2}D(x)I_d$.  The proof of the adjoint expansion is identical to the proof of Theorem \ref{appAthm:odk_forward_expansion} and is done in the weak sense, to give:
    \begin{equation}
        G_\varepsilon^*f(x)=m(x)f(x)+\varepsilon(\omega_l(x)f(x)+\mathcal{L}_b^*(fm)(x))+\bigO(\varepsilon^{3/2})
    \end{equation}
    and:
    \begin{equation}
        P_\varepsilon^* f(x)=f(x)+\varepsilon\mathcal{L}_b^* f(x)+\bigO(\varepsilon^{3/2})\implies L_\varepsilon^* f(x)=\mathcal{L}_b^* f(x)+\bigO(\varepsilon^{1/2})
    \end{equation}
    which completes the proof.
\end{proof}

\subsection{Anisotropic-Ordered Diffusion Kernels}

In the definition of an anisotropic ODK, we have the term $(y-x)^TA(x)^{-1}(y-x)$, where $A$ is a symmetric, uniformly elliptic $\binom{2}{0}$-tensor field on $\mathcal{M}$.  However, the meaning of this weighted inner product is ambiguous: suppose $\iota:\mathcal{M}\to\mathbb{R}^n$ is our isotropic embedding of $\mathcal{M}$ --- for now, let us denote for any $x\in\mathcal{M}$, $\tilde x=\iota(x)$; before we stated Lemma \ref{lemma:odk_metric_expansion} we noted that $\lVert y-x\rVert\coloneqq \lVert \tilde y-\tilde x\rVert$, which is sloppy notation but doesn't cause any issues.  However, in an anODK, for a fixed $x\in\mathcal{M}$, and denoting $B=A(x)^{-1}\in\mathbb{R}^{d\times d}$, then $(\tilde y-\tilde x)^TB(\tilde y-\tilde x)$ does not make any sense as the dimensions are incompatible.  So we either have to view $B$ as being in $\mathbb{R}^{n\times n}$, or $y-x\in\mathbb{R}^d$, if we suppose that $A$ really is a $(2,0)$-tensor field on $\mathcal{M}$, which would make $B$ a $\binom{0}{2}$-tensor field similar to the Riemannian metric tensor, then can resolve both interpretations through:
\begin{equation}\label{appAeq:aniso_metric_eq1}
    (y-x)^TB(x)(y-x)\coloneqq (\tilde y-\tilde x)^TD\iota(x) B(x)D\iota(x)^T(\tilde y-\tilde x)
\end{equation}
Since we construct the normal coordinates and projection coordinates with respect to the same orthonormal basis, we will have $D\iota^TD\iota =I_d$ by definition of the isometry, which makes the columns of $D\iota(x)\in\mathbb{R}^{n\times d}$ the orthonormal basis of the affine tangent space $T_{\tilde x}\iota(\mathcal{M})\subset\mathbb{R}^n$, and thus the coordinate chart of the projection coordinates around $x$ will equal $\phi_x:U\to\mathbb{R}^d:y\mapsto D\iota(x)^T(\iota(y)-\iota(x))$.  Therefore, assuming $B(x)$ has been written in projection coordinates, equation \eqref{appAeq:aniso_metric_eq1} will equal $u^TA^{-1}u$ exactly.  We can also view this equation as before with $\tilde y-\tilde x$ unchanged and with the symmetric quadratic form being with respect to $\tilde B(x)\coloneqq D\iota(x)B(x)D\iota(x)^T\in\mathbb{R}^{n\times n}$, where naturally $\tilde B$ will be a rank-$d$ $n\times n$ matrix.  

The more complex case arises if we don't consider $\tilde B(x)$ to be a projection onto the tangent space, i.e. we view $y-x$ as $\iota(y)-\iota(x)$ and $\tilde B(x)\in\mathbb{R}^{n\times n}$ is a symmetric uniformly elliptic matrix with rank greater than $d$.  Theoretically, we can avoid such cases by restricting to the tangent space, but in practice, one might encounter a situation such as this if one constructs $\tilde{B}$ from the local covariance structure of the embedded data $\mathcal{D}\subset\mathbb{R}^n$.  To handle this we will have to establish the following notation: for $\iota:\mathcal{M}\to\iota(M)\backwardscoloneqq\tilde{\mathcal{M}}\subset\mathbb{R}^n$, we denote by $T_{\tilde x}\tilde{\mathcal{M}}\subset\mathbb{R}^n$ and $N_{\tilde x}\tilde{\mathcal{M}}\subset\mathbb{R}^n$ the ambient tangent and normal bundles of $\mathcal{M}$ induced by $\iota$.  We note that for every $\tilde x\in\tilde{\mathcal{M}}$, we will have $T_{\tilde x}\mathbb{R}^n=T_{\tilde x}\tilde{\mathcal{M}}\oplus N_{\tilde x}\tilde{\mathcal{M}}$; we will assume that $A$ is symmetric positive semi-definite $\binom{2}{0}$-tensor field on $\mathbb{R}^n$; furthermore, we assume $A$ can be decomposed through the direct sum $A_\top\oplus A_\perp$, where $A_\top$ is a symmetric rank-$d$ uniformly elliptic tensor-field on the subbundle $T_0^2\tilde{\mathcal{M}}\leq T_0^2\mathbb{R}^n$, and $A_\perp$ is a symmetric, potentially degenerate tensor field on the subbundle $N_0^2\tilde{\mathcal{M}}\leq T_0^2\mathbb{R}^n$.  We note that in a much more canonical setting of equation \eqref{appAeq:aniso_metric_eq1} --- which we note is also the setting of local kernels --- this is equivalent to this latter scenario but with $B_\perp|_{\tilde{\mathcal{M}}}\equiv\mathbb{0}_{(n-d)\times(n-d)}$.

Since $A_\top$ is everywhere rank-$d$ on a $d$-dimensional vector space, its inverse will always exist, so let us denote this with $B_\top$ (i.e. $B_\top A_\top=A_\top B_\top=I_d$).  Given this setup, we provide the following generalisation of the metric comparison Lemma given in \ref{lemma:odk_metric_expansion}:
\begin{lemma}[Generalised Metric Comparison]\label{appAlemma:genmetcomp}
    Let $\mathcal{M}$ be a $d$-dimensional Riemannian manifold that can be embedded in $\mathbb{R}^n$, and let $\tilde B$ be a field of matrices over $\mathbb{R}^n$ such that $\tilde B|_{\tilde{\mathcal{M}}}=B_\top\oplus B_\perp$ where $B_\top\in \mathcal{T}_2^0\tilde{\mathcal{M}}$ is a symmetric and uniformly elliptic rank-$d$ 2-tensor field acting on tangent vectors of $\tilde{\mathcal{M}}$, and $B_\perp\in \mathcal{N}_2^0\tilde{\mathcal{M}}$ is a symmetric, and potentially degenerate, 2-tensor field acting on normal vectors of $\tilde{\mathcal{M}}$.  For any $x\in\mathcal{M}$, and $\varepsilon$ sufficiently small, then $\forall y\in B(0,\varepsilon^{1/2})$, we will have:
    \begin{equation}\label{appAlemma:generalised_metric_comparison}
        (\tilde y-\tilde x)^T\tilde B(\tilde y-\tilde x)=u^TB_\top u+p_{x,4}(u)+p_{x,5}(u) +\bigO(\varepsilon^3)
    \end{equation}
\end{lemma}
\begin{proof}
    The projection coordinates, $u$, about any $x\in\mathcal{M}$ are defined on a neighbourhood guaranteed by the implicit function theorem; furthermore, there will exist a diffeomorphism $g:\mathbb{R}^d\to\mathbb{R}^{n-d}$ such that in this neighbourhood, $\tilde y-\tilde x=(u,g(u))^T$.  This function will satisfy $g(0)=\mathbb{0}_{n-d}$ and for any component of $g$ $\tfrac{\partial g_j(0)}{\partial u_i}=0$, since $T_{\tilde x}\tilde{\mathcal{M}}$ is tangent to $g$ at $0$.

    One should think of this $g$ as encoding all the curvature of the manifold coming from the embedding $\iota:\mathcal{M}\to\tilde{\mathcal{M}}$.  By the assumed decomposition $\tilde B$ on $T_{\tilde x}\tilde{\mathcal{M}}$ we may write $\tilde B$ in these coordinates as the block matrix:
    \begin{equation}
        \tilde B(x)=\begin{pmatrix}
            B_\top(x) & \mathbb{0}_{d\times (n-d)}\\
            \mathbb{0}_{(n-d)\times d} & B_\perp(x)
        \end{pmatrix}
    \end{equation}
    We can now expand the matrix-vector product in $\mathbb{R}^n$ as:
    \begin{equation}\label{appAeq:genlem_int}
        (\tilde y-\tilde x)^T\tilde B(\tilde y-\tilde x)=\begin{pmatrix}u\\g(u)\end{pmatrix}^T\tilde B\begin{pmatrix}u\\g(u)\end{pmatrix}=u^TB_\top u+g(u)^TB_\perp g(u)
    \end{equation}
    For each component function, $g_k:\mathbb{R}^{d}\to\mathbb{R}$, we can take a Taylor expansion about $0$:
    \begin{equation}
        g_l(u)=\frac{1}{2}\sum_{i,j=1}^du_iu_j\frac{\partial^2g_l(0)}{\partial s_i\partial s_j}+\frac{1}{6}\sum_{i,j,k=1}^du_iu_ju_k\frac{\partial^3g_l(0)}{\partial s_i\partial s_j\partial s_k}+p_{x,4}(u)+\bigO(\lVert u\rVert^5)
    \end{equation}
    where we have used the localised and tangential properties of $g(u)$; we can represent this in vector form through:
    \begin{equation}
        g(u)=\boldsymbol{p}_{x,2}(u)+\boldsymbol{p}_{x,3}(u)+\boldsymbol{p}_{x,4}(u) +\bigO(\lVert u\rVert^5)
    \end{equation}
    where $\boldsymbol{p}_{x,k}$ is a vector of dimension $n-d$ in which each component is a homogeneous polynomial of degree $k$.  Plugging this into equation \eqref{appAeq:genlem_int} and collecting higher order terms we deduce:
    \begin{align}
        &u^TB_\top u+\boldsymbol{p}_{x,2}(u)^TB_\perp\boldsymbol{p}_{x,2}(u)+2\boldsymbol{p}_{x,2}(u)B_\perp \boldsymbol{p}_{x,3}(u)+\bigO(\lVert u\rVert^6)\\
        &\qquad=u^TB_\top u+\sum_{i,j=1}^{n-d}(B_\perp)_{ij}(\boldsymbol{p}_{x,2}(u))_i(\boldsymbol{p}_{x,2}(u))_j+2\sum_{i,j=1}^{n-d}(B_\perp)_{ij}(\boldsymbol{p}_{x,2}(u))_i(\boldsymbol{p}_{x,3}(u))_j\\
        &\qquad=u^TB_\top u+p_{x,4}(u)+p_{x,5}(u)+\bigO(\lVert u\rVert^6)
    \end{align}
    completing the proof.
\end{proof}
We can deduce from this Lemma that if we chose $\tilde B$ not to have this block diagonal structure, then the cross terms between the tangent and normal bundles would manifest itself in a $p_{x,3}(u)$ error term.  Given the schematic of the proof for Theorem \ref{appAthm:odk_backward_expansion}, we will end up with an integral term equal to:
\begin{equation}
    \frac{\sqrt{\varepsilon}}{\rho^2}\int_{\mathbb{R}^d}p_{x,3}(z-\sqrt{\varepsilon}\mu_A)(f+\sqrt{\varepsilon}\langle (z-\sqrt{\varepsilon}\mu_A),\nabla f\rangle)h'\left(\frac{z^TB_\top z}{\rho^2}\right)\mathrm{d}z+\bigO(\varepsilon^{3/2})
\end{equation}
we note that $p_{x,3}(z-\sqrt{\varepsilon}\mu_A)=p_{x,3}(z)+\sqrt{\varepsilon}p_{x,2}(z)+\varepsilon p_{x,1}(z)+\bigO(\varepsilon^{3/2})$; since we are integrating over a symmetric domain and $h'(z^TB_\top z/\rho^2)$ is an even function, after expanding the integrand we can set some of these integrals to zero and others we can collect in the $\omega_Af$ and $\bigO(\varepsilon^{3/2})$ terms, however the following term does not collapse into any other:
\begin{align}
    \frac{\varepsilon}{\rho^2}\nabla f\cdot&\int_{\mathbb{R}^d}z z^TB_{\top\perp}\begin{pmatrix}
        z^T\nabla^2g_1(0)z\\ \vdots\\z^T\nabla^2g_{n-d}(0)z
    \end{pmatrix}h'\left(\frac{z^TB_\top z}{\rho^2}\right)\mathrm{d}z\\
    &\frac{\varepsilon}{\rho^2}\sum_{i,j,k,l=1}^d\sum_{m=1}^{n-d}(B_{\top\perp})_{jm}\frac{\partial^2g_m(0)}{\partial s_k\partial s_l}\partial_i f\int_{\mathbb{R}^d}z_iz_jz_kz_lh'\left(\frac{z^TB_\top z}{\rho^2}\right)\mathrm{d}z
\end{align}
Since this depends on the integral of a positive even function over a symmetric domain, it will not vanish; and thus we would resolve an additional advection term in the Laplacian, which would depend on the curvature of the manifold at $\tilde x$ and something like the fourth moment of our kernel (which is a measure of peakedness/heavy-tails).  If $h(u)=\exp(-u/2)$, and thus $h'(u)=-\tfrac{1}{2}\exp(-u/2)$, then dividing this integral by $m_A(x)$ would yield:
\begin{equation}
    -\frac{\varepsilon}{2\rho^6}\sum_{i,j,k,l=1}^d\sum_{m=1}^{n-d}(B_{\top\perp})_{jm}\frac{\partial^2g_m(0)}{\partial s_k\partial s_l}\partial_i f\left[(B_\top)_{ij}(B_\top)_{kl}+(B_\top)_{ik}(B_\top)_{jl}+(B_\top)_{il}(B_\top)_{jk}\right]
\end{equation}
where we have applied Isserlis's Theorem to compute the fourth moments.  Now that we understand why we have constrained our matrix $\tilde B$, and recalling that when $B_\perp|_{\tilde{\mathcal{M}}}\equiv \mathbb{0}_{(n-d)\times(n-d)}$ we return to the case of equation \eqref{appAeq:aniso_metric_eq1} --- which is the theoretically proper setting --- we can continue by generalising the ODK metric equation, given in Lemma \ref{lemma:odk_metric_expansion}, to anisotropic diffusions; we denote an anODK as $K_\varepsilon(x,y)=h\left(\tfrac{\star\star}{\varepsilon\rho(x)^2}\right)$:
\begin{lemma}
    Under the same conditions as Lemma \ref{appAlemma:genmetcomp}, and denoting an anODK as $K_\varepsilon(x,y)=h\left(\tfrac{\star\star}{\varepsilon\rho(x)^2}\right)$, we have:
    \begin{align}
        (\star\star)\coloneqq (y-x)^T\tilde B(y-x)&+\frac{(l(x,y)-\varepsilon\tau(x)\kappa(x))^2}{\kappa(x)^2}-\frac{l(x,y)^2}{\kappa(x)^2}\\
        &=(u-\varepsilon\mu_A(x))^TB_\top (u-\varepsilon\mu_A(x))+\Theta_4(u)+\Theta_5(u) + \bigO(\varepsilon^3,\varepsilon\lVert u\rVert^4)
    \end{align}
    where $\mu_A=\tfrac{\tau(x)}{\kappa(x)}A_\top\nabla l_x(x)$, and $\Theta_4(u)$ and $\Theta_5(u)$ are expressions of order $\bigO(\varepsilon^2)$ and $\bigO(\varepsilon^{5/2})$, respectively.
\end{lemma}
\begin{proof}
    Using Lemma \ref{appAlemma:genmetcomp} and expanding $l(x,y)$ about $x$, and letting $b(x)=\nabla l_x(x)$ we arrive at:
    \begin{equation}
        2\frac{\tau(x)}{\kappa(x)}l(x,y)=2\frac{\tau(x)}{\kappa(x)}(u^Tb(x)+p_{x,2}(u)+p_{x,3}(u))=2\frac{\tau(x)}{\kappa(x)}u^TB_\top A_\top b(x)+2\frac{\tau(x)}{\kappa(x)}(p_{x,2}+p_{x,3}(u))
    \end{equation}
    where we have used the fact that $B_\top^{-1}$ exists (by assumption of uniformly elliptic on $T_{\tilde x}\tilde{\mathcal{M}}$) and this inverse is $A_\top$.  Combining this with the expansion of other terms gives us:
    \begin{align}
        &u^TB_\top u-2\varepsilon\frac{\tau(x)}{\kappa(x)}u^TB_\top A_\top b(x)+\varepsilon^2\frac{\tau(x)^2}{\kappa(x)^2}\lVert B_\top^{1/2}A_\top b(x)\rVert^2+\Theta_{4}(u)+\Theta_5(u)+\bigO(\varepsilon^3,\varepsilon\lVert u\rVert^4)\\
        &\qquad=(u-\varepsilon\mu_A(x))^TB_\top (u-\varepsilon\mu_A(x))+\Theta_4(u)+\Theta_5(u) + \bigO(\varepsilon^3,\varepsilon\lVert u\rVert^4)
    \end{align}
    where $\mu_A(x)\coloneqq \tfrac{\tau(x)}{\kappa(x)}A_\top(x) b(x)$, and:
    \begin{align}
        \Theta_4(u)&=p_{x,4}(u)+\varepsilon^2\tau(x)^2\left(1-\frac{\lVert B_\top^{1/2}A_\top b(x)\rVert^2}{\kappa(x)^2}\right)-2\varepsilon\frac{\tau(x)}{\kappa(x)}p_{x,2}(u)\\
        \Theta_5(u)&=p_{x,5}(u)-2\varepsilon\frac{\tau(x)}{\kappa(x)}p_{x,3}(u)
    \end{align}
    completing the proof.
\end{proof}

\begin{theorem}[Expansion of Anisotropic Ordered Diffusion Kernels, Revised]\label{appAthm:noniso_odk_backward_expansion}
    Under the same assumptions as theorem \ref{thm:odk_backward_expansion}, the anisotropic ODK has the expansion:
    \begin{equation}
        G_\varepsilon f(x)=m_A(x)f(x)+\varepsilon(\omega_l(x)f(x)+m_A(x)\mathcal{L}_Af(x))+\varepsilon^{3/2}\Omega_l(x)+\bigO(\varepsilon^2),
    \end{equation}
    where $m_A(x)=m(x)\sqrt{|A(x)|}$ where $m(x)$ is the zeroth-moment from theorem \ref{thm:odk_backward_expansion}.  The second-order differential operator, $\mathcal{L}_A$, takes the form:
    \begin{equation}\label{appAeq:non_isotropic_laplacian}
        \mathcal{L}_Af(x)=\mu_A(x)\cdot\nabla f(x)+\tr_{D_A}\nabla^2f(x)=\mu_A^i(x)\partial_i f(x)+D_A^{ij}(x)\partial_i\partial_j f(x)
    \end{equation}
    where in the second equality we have written this operator in Euclidean normal coordinates (and we note that the trace of a $\binom{0}{2}$-tensor $T$ with respect to a $\binom{2}{0}$-tensor $B$ is defined to equal $\tr_B T=B^{ij}T_{ji}$).  The first and second order moments $\mu_A(x)$ and $D_A(x)$ satisfy:
    \begin{align}
        \mu_A(x)&=\frac{\tau(x)}{\kappa(x)}A(x)\nabla l(x)=A(x)\mu(x)\\
        D_A(x)&=\frac{1}{2}m_A(x)^{-1}\int_{T_x\mathcal{M}}zz^Th\left(\frac{z^TA(x)^{-1}z}{\rho(x)^2}\right)\dif z=\tfrac{1}{2}\sigma_h(x)^2\rho(x)^2A(x)=A(x)D(x)
    \end{align}
    where $\mu(x)$ and $D(x)$ are the same as in theorem \ref{thm:odk_backward_expansion}.  The normalised backward integral operator and its Laplacian can be expanded in terms of $\varepsilon$ as:
    \begin{equation}
        P_\varepsilon f(x) = f(x)+\varepsilon\mathcal{L}_A f(x)+\bigO(\varepsilon^{3/2})\implies L_\varepsilon f(x)=\mathcal{L}_A f(x)+\bigO(\varepsilon^{1/2})
    \end{equation}
    analogously to theorem \ref{thm:odk_forward_expansion}, in the weak sense the adjoint operator satisfies:
    \begin{equation}
        G_\varepsilon^*f(x)=m_A(x)f(x)+\varepsilon(\omega_l(x)f(x)+\mathcal{L}_A^*(fm_A)(x))+\bigO(\varepsilon^{3/2}),
    \end{equation}
    where $\mathcal{L}_A^*$ is the adjoint of the operator $\mathcal{L}_A$ in theorem \ref{thm:noniso_odk_expansion}, satisfying $\langle f,\mathcal{L}_A^*g\rangle_{L^2}=\langle \mathcal{L}_Af,g\rangle_{L^2}$ for any $f,g\in L^2(\mathcal{M})$, this operator takes the form (\cite{evans2022partial}, p320):
    \begin{equation}
        \mathcal{L}_A^*f(x)=-\diver\left(f(x)\mu_A(x)\right)+\Delta(D_A(x)f(x))=-\partial_i(f(x)\mu_A^i(x))+\partial_i\partial_j(D_A^{ij}(x)f(x))
    \end{equation}
    and the normalised adjoint integral operator and its Laplacian can be expanded as:
    \begin{equation}
        P_\varepsilon^*f(x)=f(x)+\varepsilon\mathcal{L}_A^*f(x)+\bigO(\varepsilon^{3/2})\implies L_\varepsilon^* f(x)=\mathcal{L}_A^* f(x)+\bigO(\varepsilon^{1/2}).
    \end{equation}
\end{theorem}
\begin{proof}[Proof of Theorem \ref{thm:noniso_odk_expansion}]
    Fix some $\gamma\in (0,1/2)$ and take any $x\in\mathcal{M}$ that is at least a distance $\varepsilon^\gamma$ from the boundary of $\mathcal{M}$; from Lemma \ref{appAlemma:genmetcomp} we see that the only difference between our expansion in the proof of Theorem \ref{appAthm:odk_backward_expansion} and what we have currently is the terms inside $h(\cdot)$ are of the form $(u-\varepsilon\mu_A)^TB_\top(u-\varepsilon\mu_A)$ which is still a local kernel and thus we can still apply Lemma 3.9 from \cite{berry2016local} to arrive at:
    \begin{equation}
        G_\varepsilon f(x)=m_A(x)f(x)+\varepsilon(\omega_A(x)f(x)+m_A(x)\mathcal{L}_A f(x))+\varepsilon^{3/2}\Omega_A(x)+\bigO(\varepsilon^2)
    \end{equation}
    Where $m_A\mathcal{L}_A=\tilde{\mathcal{L}}_A$ and $\tilde{\mathcal{L}}_A=\tilde{\mu}_A\cdot\nabla+\text{tr}_{\tilde{D}_A}\nabla^2$.  The zeroth, first- and second-order moments of this local kernel are:
    \begin{align}
        m_A(x)&=\int_{T_x\mathcal{M}}h\left(\frac{z^TB_\top(x)z}{\rho(x)^2}\right)\mathrm{d} z=\sqrt{|A_\top(x)|}m(x)\\
        \tilde{\mu}_A(x)&=m_A(x)\frac{\tau(x)}{\kappa(x)}A(x)\nabla l(x)=m_A(x)\mu_A(x)\\
        (\tilde{D}_A)_{ij}(x)&=\frac{1}{2}\int_{T_x\mathcal{M}}z_iz_jh\left(\frac{z^TB_\top(x)^{-1}z}{\rho(x)^2}\right)\mathrm{d} z\backwardscoloneqq m_A(x)D_A(x)
    \end{align}
    We note that, just as before, our $\Theta_4(u)$ and $\Theta_5(u)$ terms are well-behaved and thus have no bearing on the expansion we derive.  In much the same way as before, we can now immediately derive the expression for the normalised backwards integral operator and the Laplacian:
    \begin{equation}
        P_\varepsilon f(x)=f(x)+\varepsilon\mathcal{L}_Af(x)+\bigO(\varepsilon^{3/2})\implies L_\varepsilon f(x)=\mathcal{L}_Af(x)+\bigO(\varepsilon^{1/2})
    \end{equation}
    The expansion of the adjoint is computed in the weak sense as:
    \begin{equation}
        G_\varepsilon^* f(x)=m_A(x)f(x)+\varepsilon(\omega_A(x)f(x)+\mathcal{L}_A^*(m_Af)(x))+\bigO(\varepsilon^{3/2})
    \end{equation}
    and to finish the proof:
    \begin{equation}
        P_\varepsilon^* f(x)=f(x)+\varepsilon\mathcal{L}_A^*f(x)+\bigO(\varepsilon^{3/2})\implies L_\varepsilon^* f(x)=\mathcal{L}_A^*f(x)+\bigO(\varepsilon^{1/2})
    \end{equation}
\end{proof}

\subsubsection{Viewing $A(x)$ as a non-trivial Riemannian Metric}\label{appAsec:noniso_metric}
At the beginning of section \ref{subsec:nonisotropic_extension} we explained that in definition \ref{def:nonisotropic_odk} we split the diffusion, $D(x)=\tfrac{1}{2}\rho(x)^2A(x)$, into terms involving $\rho(x)^2$ and $A(x)$ so that we can discriminate between viewing $A(x)$ as a diffusion and a non-trivial Riemannian metric --- we will now consider $A(x)$ to be as in this latter case.  The Riemannian metric is a $\binom{0}{2}$-symmetric tensor that defines an inner product on $\mathcal{M}$, if $u,v\in T_x\mathcal{M}$ for some $x\in\mathcal{M}$, then $g(x)(u,v)=g_{ij}(x)u^iv^j$ will be their inner product with respect to some coordinate system.  We see from equation \eqref{eq:non_isotropic_odk} that $A^{-1}(x)$ acts like a metric on $\mathcal{M}$, thus we set $(g(x))_{ij}=(A^{-1}(x))_{ij}$ --- we will use the convention that $g^{ij}$ equals the $(i,j)^\text{th}$ element of the inverse metric and since $g=A^{-1}$ we will denote $g^{ij}=A^{ij}$.  The infinitesimal generator of an It{\^o} SDE with gradient-flow drift given by $\mu(x)=\tfrac{\tau(x)}{\kappa(x)}\nabla l(x)$, and isotropic diffusion given by $\rho(x)^2I_d$ defined on the manifold $(\mathcal{M},g)$ would equal (dropping function inputs):
\begin{align}
    \mathcal{L}_gf&=\tfrac{\tau}{\kappa}g(\nabla_g l,\nabla_g f)+\tfrac{1}{2}\rho^2\Delta_g f=\tfrac{\tau}{\kappa}g_{ij}(g^{ik}\partial_k l)(g^{jm}\partial_m f)+\tfrac{1}{2}\rho^2\text{div}_g(\nabla_g f)\\&=\tfrac{\tau}{\kappa}g^{ij}\partial_i l\partial_j f+\tfrac{1}{2}\tfrac{\rho^2}{\sqrt{|g|}}\partial_i\left(\sqrt{|g|}g^{ij}\partial_j f\right)\label{eq:nontrivialmetric_laplacian}
\end{align}
where $\sqrt{|g|}=\sqrt{|\text{det}(g)|}$ is the volume element with respect to the metric $g$; $\Delta_gf\coloneqq \text{div}_g(\nabla f)$ is the Laplace-Beltrami operator; $\text{div}_g(X)=|g|^{-1/2}\partial_i(|g|^{1/2} X^i)$ the divergence; and $\nabla_g f\coloneqq \text{grad}_g(f)=g^{ij}\partial_j f$ is the gradient defined to be such that $g(x)(\text{grad}(f)(x),v)=df_x(v)$ for all $v\in T_x\mathcal{M}$ where $df_x\in T_x^*\mathcal{M}$ is the differential of $f$ at $x\in\mathcal{M}$.  This follows directly from $\mathcal{L}f=\tfrac{\tau}{\kappa}\nabla l\cdot\nabla f+\tfrac{1}{2}\rho^2\Delta f$ when the metric is the standard Euclidean metric; all we have done is reparameterise all derivatives and inner products in terms of the new metric.  We first note that the observed coupling between $\mu=\tfrac{\tau}{\kappa}\nabla l$ and $\mu_A=A\mu=g^{ij}\mu_j$ is now working in our favour, since by the ordering function not being aware of the rescaling of coordinates by $A(x)^{-1}$, it has implicity adapted to these coordinates.  We see that the advection term in equation \eqref{appAeq:non_isotropic_laplacian} satisfies $\mu_A\cdot\nabla f=A^{ij}\mu_j\partial_if=A^{-1}_{ij}(A^{ik}\mu_k)(A^{jl}\partial_l f)$.  However, it is not so simple for the second-order term in equation \eqref{appAeq:non_isotropic_laplacian}, since $\tr_{\rho^2A}\nabla^2 f=\rho^2A^{ij}\partial_i\partial_j f$, but we can expand the Laplace-Beltrami operator as given in equation \eqref{eq:nontrivialmetric_laplacian} to be:
\begin{equation}
    \Delta_g f=\tfrac{1}{\sqrt{|g|}}\partial_i\left(\sqrt{|g|}g^{ij}\partial_j f\right)=g^{ij}\partial_i\partial_j f+\tfrac{1}{\sqrt{|g|}}\partial_i\left(\sqrt{|g|}g^{ij}\right)\partial_j f=\tr_g\nabla^2f+\nabla f\cdot\text{div}_g(g^{-1})
\end{equation}
where $\text{div}_g T$ is usually defined to act on the rows of a $\binom{2}{0}$-tensor $T$, i.e. $(\text{div}_g T)_i=\partial_jT^{ij}$, since $g^{-1}$ is symmetric we need not worry about this convention.  Hence, we see that the second-order term of the Laplacian in equation \eqref{appAeq:non_isotropic_laplacian}, does not capture the effective drift due to variations in the volume element.  Reconstructing the Laplace-Beltrami operator with respect to a non-trivial Riemannian metric was consider by Berry \& Sauer in Local Kernels \cite{berry2016local}, where they navigated this same issue by defining a symmetrised local kernel $\bar K_\varepsilon(x,y)=K_\varepsilon(x,y)+K_\varepsilon(y,x)$ whose Laplacian approximates the Laplace-Beltrami operator with respect to some metric $\tilde g$ under the assumptions that $\mathcal{D}=\{x_i\}_{i=1}^m$ can been sampled according to the volume measure on $\mathcal{M}$ defined by $\tilde g$.  This works because the adjoint kernel $K_\varepsilon(y,x)$ depends on $A(y)$ as $y$ varies, thus allowing the symmetrised kernel to be aware of volumetric variations; the downside of their construction is that they lose all other advection terms.  We proceed by rewriting the missing advection term from the Laplacian as:
\begin{equation}
    \tfrac{1}{\sqrt{|g|}}\partial_i\left(\sqrt{|g|}g^{ij}\right)\partial_j f=\left(\partial_i g^{ij}+g^{ij}\tfrac{1}{\sqrt{|g|}}\partial_i\sqrt{|g|}\right)\partial_j f=\left(\partial_i g^{ij}+g^{ij}\partial_i\log\sqrt{|g|}\right)\partial_j f
\end{equation}
which we can see is defining a vector field $\mu_\text{eff}=\left(\partial_i g^{ij}+\tfrac{1}{2}g^{ij}\partial_i\log|g|\right)\partial_j$ which we advect $f$ along.  In order for us to be able to recapture this term we need it be true that there exists functions $a,\varphi:\mathcal{M}\to\mathbb{R}$ such that $\mu_\text{eff}(x)=a(x)\nabla\varphi(x)$, let $\omega=\mu_\text{eff}^\flat\in \mathcal{T}^*\mathcal{M}$ be this vector-fields associated 1-form, then we know that $d\omega=da\wedge d\omega$, and therefore we know that $\omega\wedge d\omega =0$ (since $ad\varphi\wedge da\wedge d\varphi\equiv 0$), but the converse is not necessarily true.  The intuition behind this constraint (that one might notice is the Frobenius integrability condition) is that for any $x\in\mathcal{M}$ where $\omega\neq 0$, then $B_x\{v\in T_x\mathcal{M}:\omega_x(v)=0\}$ defines a $(d-1)$-dimensional subspace of $T_x\mathcal{M}$, we can stich this subspaces together to get $(d-1)$-dimensional submanifolds of $\mathcal{M}$, very much like the level sets of an ordering function.  The differential $d\omega$, effectively measures how close $\omega$ is to being (locally) generated by the gradient of a scalar function --- we can see this from the previous argument that $\omega=a d\varphi\implies \omega\wedge d\omega=0$ --- strictly, $d\omega$ is measuring the circulations of $\omega$ around $x$ (i.e. the rotational dynamics).  The condition $\omega\wedge d\omega=0$ is weaker than $\omega=a d\varphi$, because it allows for singularities when $\mu_\text{eff}(x)=0$, just as in our discussion of handling poorly conditioned ordering functions earlier.  Thus, quite poetically, for us to represent the infinitesimal generator of the It{\^o} SDE $dx_t=\tfrac{\tau}{\kappa}\nabla l\ \mathrm{d}t+\rho \mathrm{d}w_t$ on a Riemannian manifold with non-trivial metric $(\mathcal{M},A^{-1})$ as given in equation \eqref{eq:nontrivialmetric_laplacian}, and where $l(x)$ is an ordering function (almost) foliating the manifold, then the effective drift implied by the metric $A^{-1}$ in the form of the Laplace-Beltrami operator must also (almost) foliate the manifold into sets perpendicular to its advection.

We will now assume that $\text{div}(g^{-1})=0$, which immediately implies $\mu_\text{eff}=A^{ij}|A|^{1/2}\partial_i|A|^{-1/2}=\tfrac{1}{2}A^{ij}\partial_i\log|A|^{-1}$ and thus is effectively generated by a potential; we will now show how you can reconstruct the desired limiting operator using: mODKs from section \ref{subsec:multi_odk}; $\alpha$-regularised ODKs from section \ref{subsec:discrete_continuous}; and stODKs from section \ref{subsec:selftuning_odk}.  The condition that $\text{div}(A)(x)\equiv0$ makes $A$ into a stress-tensor field and roughly means that there are no sources of sinks of weight in the coordinate functions of $A$; effectively, $A(x)$ is the diffusion equivalent of a gradient-flow drift, but in a more structure way than if we only have $\omega\wedge d\omega=0$.  In all of what follows, we assume that $\sigma_h^2\equiv 1$.

Firstly, if using mODKs we can simply add another ordering function, $\tilde l(x)=\tfrac{1}{2}\log|A(x)|^{-1}$, with $\tilde\tau$ and $\tilde\kappa$ satisfying $\tilde{\kappa}^{-1}\tilde\tau=\tfrac{1}{2}\rho^2$, merging this with the anisotropic ODK will give a limiting operator:
\begin{align}
    \mathcal{L}_{A}^mf&=A\left(\tfrac{\tau}{\kappa}\nabla l+\tfrac{1}{2}\rho^2\nabla \log\tfrac{1}{\sqrt{|A|}}\right)\cdot\nabla f+\tfrac{1}{2}\rho^2\tr_{A}\nabla^2 f\\
    &=\tfrac{\tau}{\kappa}A^{ij}\partial_j l\partial_if+\tfrac{1}{2}\rho^2\left(\sqrt{|A|}A^{ij}\partial_j\left(\tfrac{1}{\sqrt{|A|}}\right)\partial_i f+A^{ij}\partial_i\partial_j f\right)
\end{align}
where, given that $\partial_i A^{ij}\equiv 0$, the second term in brackets equals $\Delta_{A^{-1}} f$.  This produces the correct limiting operator, but the associated transition probabilities treat the effective drift due to variations in the volume element as an actual drift.  We can achieve a similar expansion by assuming that the sampling measure, $q(x)$, is uniform with respect to the volume form, or in other words $q(x)\propto |A(x)|^{-1/2}$, and setting $\alpha=\tfrac{1}{2}$ in the $\alpha$-regularised ODK (where the $\alpha$-regularisation was introduced in \cite{coifman2006diffusion}), this results in:
\begin{align}
    \mathcal{L}_A^qf&=\tfrac{\tau}{\kappa}A\nabla l\cdot\nabla f+\tfrac{1}{2}\rho^2A\frac{\nabla q}{q}\cdot\nabla f+\tfrac{1}{2}\rho^2\tr_{A}\nabla^2 f\\
    &=\tfrac{\tau}{\kappa}A^{ij}\partial_j l\partial_i f+\tfrac{1}{2}\rho^2\left(\sqrt{|A|}A^{ij}\partial_j\left(\tfrac{1}{\sqrt{|A|}}\right)\partial_j f+A^{ij}\partial_i\partial_j f\right)
\end{align}
where again the expression inside the brackets equals $\Delta_{A^{-1}} f$ when $\partial_iA^{ij}\equiv 0$.  We note that the assumption on the sampling density can be turned around and thought of as choosing $A(x)$ such that this covariance matrix is the best fit to the data; furthermore, given we assume our data is generated by some dynamical process, it is likely true that the covariance that best fits the data and the underlying covariance $\sigma(x)\sigma(x)^T$ will be highly correlated at any $x\in\mathcal{M}$.  We also add that if $A(x)=\rho^2I_d$, then $|A(x)|^{-1/2}=\rho^{-d}$, so if $\rho(x)\propto\lVert x-x^k\rVert$ where $x^k$ is the $k^\text{th}$ closest sample in $\mathcal{D}$ to any $x\in\mathcal{M}$, then by the simplest possible approximation to the sampling density we would already have $q(x)\propto \lVert x-x^k\rVert^{-d}$ \cite{loftsgaarden1965nonparametric}; having $q(x)\propto|A(x)|^{-1/2}$ is thus the simplest anisotropic generalisation of this result.  

In a similar vein to the previous reasoning we could define an stODK, as described in section \ref{subsec:selftuning_odk}, with $\rho_1(x)\propto |A(x)|^{\eta/2}\rho^2$ and $\rho_2(y)\propto|A(y)|^{-\eta/2}$ with proportionality such that $\rho(x)^{-2}\left[\rho_1(x)\rho_2(x)\right]=\left[\eta(d+2)\right]^{-1}$, then after combining the anisotropic ODK with the self-tuning ODK (replace all occurences of $\rho(x)^2$ in equation \eqref{eq:non_isotropic_odk} with $\rho_1(x)\rho_2(y)$), we will find the limiting operator to be:
\begin{align}
    \mathcal{L}_A^{st}f&=\tfrac{\tau}{\kappa}A\nabla l\cdot\nabla f+\tfrac{d+2}{2}\rho_1A\nabla\rho_2\cdot\nabla f+\tfrac{1}{2}\tr_{A}\nabla^2 f\\
    &=\tfrac{\tau}{\kappa}A^{ij}\partial_j l\partial_i f+\eta\tfrac{d+2}{2}\rho_1\rho_2\sqrt{|A|}A^{ij}\partial_j\left(|A|^{-1/2}\right)\partial_i f+\frac{1}{2}A^{ij}\partial_i\partial_j f\\
    &=\tfrac{\tau}{\kappa}A^{ij}\partial_j l\partial_i f+\tfrac{1}{2}\rho^2\left(\sqrt{|A|}A^{ij}\partial_j\left(|A|^{-1/2}\right)\partial_i f+A^{ij}\partial_i\partial_j f\right)
\end{align}
where the expression in brackets equals $\Delta_{A^{-1}} f$ when $\partial_i A^{ij}\equiv 0$.  To gain some intuition, let us define $\rho_A(x)=\sqrt{|A(x)|}^{1/d}$ to be the $d^\text{th}$ root of the inverse volume form, and which can be thought of as a length scale (geometric mean of length scales along different dimensions) with respect to the metric $A^{-1}$, then in the above formulation $\rho_2\propto\rho_A^{-\eta d}$ and, assuming $\rho\propto\rho_A$, we will also have $\rho_1\propto \rho_A^{\eta d+2}$, with the product $\rho_1\rho_2\propto\rho_A^2$ combining nicely into a factor that is of the desired dimensionality.  If $\eta=1$, then $\rho_2$ is proportional to the volume form and, and thus we are biasing rescaling the terms of the ODK inside $h(\cdot)$ by a factor proportional to $\text{vol}_{A^{-1}}(x)/\text{vol}_{A^{-1}}(y)$, thus biasing towards transitioning to higher volume (or lower density if $q(x)\propto|A(x)|^{-1/2}$) states; in addition $\rho_1\propto\rho_A^{d+2}$ which fits nicely with the factor of $(d+2)$ in the self-tuning bandwidth formulation of ODK, recalling that this factor initially comes from making a substitution within an integral and thus it makes sense why taking $\rho_2$ to be a volume form would recover this.  If we take $\eta=-1/d$, then $\rho_1\propto\rho_A$ and $\rho_2\propto\rho_A$, which is a much more familiar symmetric bandwidth formulation (i.e. $\beta=1/2$ regularisation in section \ref{subsec:selftuning_odk}), thus with the correct constants of proportionality we can also use stODKs to recapture the Laplace-Beltrami operator when $\partial_iA^{ij}\equiv 0$.

To summarise, we have extended ODKs to anisotropic diffusions at the cost of coupling the drift and diffusion in the infinitesimal generator; analogous to how TMDmap kernels couple drift and diffusion unavoidably.  We then proposed to approximate the uncoupled anisotropic transition operator by using the Strang Splitting scheme, and as a result we are capable (especially upon combining anODKs with lODKs) of approximating the infinitesimal generator of an arbitrary It{\^o} process.  We then discussed how one can interpret the anisotropic covariance matrix $A(x)$ as a non-trivial Riemannian metric; in this scenario, the coupling we observed between the first- and second-order terms now works in our favour, but the second-order term only captures the trace of the Hessian and thus does not approximate the Laplace-Beltrami operator with respect to $A^{-1}$.  We showed that under some mild constraints on $A^{-1}$, namely that the effective drift it implies is capable of (almost) foliating $\mathcal{M}$ with its tranverse submanifolds; this is an interesting condition as it is effectively identical to the discussion we had in section \ref{subsec:ordering_functions} about what makes $l(x)$ a good ordering function.  Under this constraint, we show that one can reconstruct the full infinitesimal generator (advection term and Laplace-Beltrami term) by applying multi-Ordered Diffusion Kernels discussed in section \ref{subsec:multi_odk}; $\alpha$-regularised ODKs as discussed in section \ref{subsec:discrete_continuous}; and self-tuning ODKs as discussed in the next section.  Using an mODK makes the effective drift very real from the perspective of the transition probabilities, which is a technical inaccuracy but is easy to apply; $\alpha$-normalised ODKs are the most satisfying method but they require that the sampling distribution $q(x)$ is uniform with respect to the volume form $|A|^{-1/2}$; stODKs reconstruct the effective drift in a truly effective manner, and do so with no extra constraints, however one does need to be careful that proportionality constants satisfy $\eta(d+2)\rho_1\rho_2=\rho^2$.

\subsubsection{Self-Tuning Ordered Diffusion Kernels}

In this appendix, we will provide the proof for the expansion of self-tuning Ordered Diffusion Kernels (stODKs) as defined in section \ref{subsec:selftuning_odk}.  All the difficulty in this proof is making a detailed accounting of all the new integral terms that occur once you expand $\rho_2(y)$ about $\rho_2(x)$; this task can be greatly simplified through the following Lemma:
\begin{lemma}\label{appAlemma:simplifying_lemma}
    Let $(\star)$ be as in Lemma \ref{lemma:odk_metric_expansion} written with either a local or global ordering function.  Let $\xi=\tfrac{(\star)}{\varepsilon}$, the following kernels will have the same limiting operator as $h(\xi)$:
    \begin{enumerate}[nolistsep,noitemsep]
        \item[(i).] $h\left(\xi+\varepsilon^{k/2} c(x)\right)$ for any integer $k\geq 2$ and any bounded $c(x)$,
        \item[(ii).] $h\left(\xi+\tfrac{p_{x,a}(u)}{\varepsilon}\right)$ for any integer $a\geq 4$,
        \item[(iii).] $h\left(\xi +p_{x,b}(u)\xi\right)$ for any integer $b\geq 2$.
    \end{enumerate}
    
\end{lemma}
\begin{proof}
    All of these proofs will rely on expanding $h(\cdots)$ about $\xi$, so let $\zeta$ be any of the additional terms, and we will have:
    \begin{equation}
        h(\xi+\zeta)=h(\xi)+\zeta h'(\xi)+\tfrac{1}{2}\zeta^2h''(\xi)+\bigO(\zeta^3)
    \end{equation}
    The validity of (ii) and (iii) will ultimately be determined by the validity of (i); thus, let us expand $h(\xi+\varepsilon^rc(x))$ about $h(\xi)$ in the backward integral operator:
    \begin{equation}
        \varepsilon^{d/2}G_\varepsilon f(x)=\int h(\xi)f(y)\mathrm{d}y+\varepsilon^{k/2} c(x)\int h'(\xi)f(y)\mathrm{d}y+\varepsilon^{k}c(x)^2\int h''(\xi)f(y)\mathrm{d}y+\bigO(\varepsilon^{3k/2})
    \end{equation}
    Since $k\geq 2$, the third term in the above expression will be $\bigO(\varepsilon^2)$; we also see that once we change the integral to projection coordinates and substitute $u\mapsto \varepsilon^{1/2}z$, the only integrals that are of lower order than $\bigO(\varepsilon^2)$ will be:
    \begin{equation}
        \varepsilon^{k/2}f(x)c(x)\int h'(\xi)\mathrm{d}z+\varepsilon^{(k+1)/2}c(x)\int h'(\xi)\langle u,\nabla f\rangle\mathrm{d}u+\bigO(\varepsilon^2)
    \end{equation}
    If $k\geq 3$, we see that the second term is now $\bigO(\varepsilon^2)$ and the first term can be absorbed into the $\varepsilon^{3/2}\Omega_l(x)$, and thus has no effect on the limiting operator; if $k=2$, then similarly the first term can be absorbed into $\varepsilon\omega_l(x)f(x)$ and the second term can be absorbed into $\varepsilon^{3/2}\Omega_l(x)$, and again the limiting operator is unchanged.  This finishes the proof of (i).

    To prove (ii) and (iii), we note that under the substitution $u\mapsto \varepsilon^{1/2}z$, we will have:
    \begin{align}
        \frac{p_{x,a}(u)}{\varepsilon}&\longmapsto \varepsilon^{(a-2)/2}p_{x,a}(z)\\
        p_{x,b}(u)\xi&\longmapsto \varepsilon^{b/2}p_{x,b}(z)\lVert z-\varepsilon^{1/2}\mu\rVert^2
    \end{align}
    Thus, so long as $a\geq 4$ in (ii), and $b\geq 2$ in (iii), the same reasoning as was applied in the proof of (i) will apply, and the limiting operator will be unchanged.
\end{proof}
This Lemma effectively says that if a term is $\bigO(\varepsilon,\lVert u\rVert^2)$ inside of $h(\cdots)$, then it will not have any effect on the limiting operator for the Laplacian of an ODK.  Taking the Taylor expansion of $\rho_2(y)$ about $\rho_2(x)$ introduces many such terms, and thus the proof of the asymptotic expansion for stODKs will utilise this Lemma multiple times.  The statement and proof of the theorem is as follows:
\begin{theorem}[Expansion of Self-Tuning Ordered Diffusion Kernels]\label{appAthm:self_tuning_expansion}
    Suppose $\mathcal{M}$ is a $d$-dimensional Riemannian manifold embedded in $\mathbb{R}^n$ and with metric inherited from the standard Euclidean metric in ambient space. Let $l\in C^2(\mathcal{M})$ be an ordering function (or a local ordering function), and let $K_\varepsilon^{st}(x,y)$ be the corresponding stODK, as defined in \eqref{def:selftuning_odk}, with self-tuning bandwidth of the form $\rho_1(x)\rho_2(y)$ for smooth functions $\rho_1,\rho_2\in C^2(\mathcal{M})$ bounded uniformly away from zero. Then for any $f\in C^3(\mathcal{M})$ the associated family of backward integral operators $\{G_\varepsilon\}_{\varepsilon> 0}$ admits the following asymptotic expansion as $\varepsilon\to0$ when applied to $f$:
    \begin{equation}
        G_\varepsilon f(x) = m(x)f(x) + \varepsilon\left( \omega_l(x)f(x) + m(x)\mathcal{L}_{st}f(x) \right)  + \varepsilon^{3/2}\Omega_l(x)+ \bigO(\varepsilon^2),
    \end{equation}
    at a point $x\in\mathcal{M}$ a distance larger than $\varepsilon^\gamma$ from $\partial\mathcal{M}$, for a fixed $\gamma\in(0,1/2)$; and where $\omega_l(x)$ depends on the manifold and the kernel $K_\varepsilon(x,\cdot)$ at $x$, and $\Omega_l(x)$ depends on the manifold, the kernel $K_\varepsilon(x,\cdot)$ at $x$, and also the test function $f(x)$. The differential operator $\mathcal{L}$ is a second-order elliptic operator of the form:
    \begin{equation}\label{appAeq:st_backward_differential_operator}
        \mathcal{L}_{st}f(x) = \left[\mu(x)+(d+2)D(x)\frac{\nabla \rho_2(x)}{\rho_2(x)}\right]\cdot\nabla f(x) + D(x)\Delta f(x)
    \end{equation}
    where $m(x)$, $\mu(x)$ and $D(x)$, satisfy:
    \begin{equation}
        m(x)=\int_{T_x\mathcal{M}}h\left(\tfrac{\lVert z\rVert^2}{\rho_1(x)\rho_2(x)}\right)\mathrm{d}z,\qquad\mu(x)=\frac{\tau(x)}{\kappa(x)}b(x),\qquad D(x)=\ \tfrac{1}{2}\sigma_h(x)^2\rho_1(x)\rho_2(x)\\
    \end{equation}
    The normalised backward integral operators and their Laplacians can thus be expanded in terms of $\varepsilon$ as:
    \begin{equation}
        P_\varepsilon f(x) = f(x)+\varepsilon\mathcal{L}_{st}f(x)+\bigO(\varepsilon^{3/2})\implies L_\varepsilon^{st}f(x)=\mathcal{L}_{st}f(x)+\bigO(\varepsilon^{1/2})
    \end{equation}
    Analogous to previous theorems, the family of adjoint operators $\{G_\varepsilon^*\}_{\varepsilon>0}$ satisfies the following asymptotic expansion, in the weak sense, when applied to $f\in C^3(\mathcal{M})\cap L^2(\mathcal{M})$:
    \begin{equation}
        G_\varepsilon^*f(x)=m(x)f(x)+\varepsilon(\omega_l(x)f(x)+\mathcal{L}_{st}^*(fm)(x))+\bigO(\varepsilon^{3/2})
    \end{equation}
    where $\mathcal{L}_{st}^*$ is the formal adjoint of $\mathcal{L}_{st}$ and takes the form:
    \begin{equation}
        \mathcal{L}_{st}^* f(x)=-\text{div}\left(f(x)\left[\mu(x)+(d+2)D(x)\frac{\nabla \rho_2(x)}{\rho_2(x)}\right]\right)+\Delta(D(x)f(x))
    \end{equation}
    The normalised adjoint integral operator $P_\varepsilon^*f=G_\varepsilon^*((G_\varepsilon \mathbb{1})^{-1} f)$ and its Laplacian can be expanded as:
    \begin{equation}
        P_\varepsilon^* f(x)=f(x)+\varepsilon\mathcal{L}_{st}^* f(x)+\bigO(\varepsilon^{3/2})\implies L_\varepsilon^* f(x)=\mathcal{L}_{st}^* f(x)+\bigO(\varepsilon^{1/2})
    \end{equation}
\end{theorem}
\begin{proof}[Proof of Theorem \ref{appAthm:self_tuning_expansion}]
    Fix some $\gamma\in (0,1/2)$ and take any $x\in\mathcal{M}$ that is at least a distance $\varepsilon^\gamma$ from the boundary of $\mathcal{M}$.  The Taylor Expansion of $\rho_2(y)$ at $x$ yields:
    \begin{equation}
        \rho_2(y)=\rho_2(x)+\langle u,\nabla \rho_2(x)\rangle+p_{x,2}(u)+p_{x,3}(u)+\bigO(\lVert u\rVert^4)
    \end{equation}
    Utiliting the expansion of $\tfrac{1}{1+x}=1-x+x^2-\dots$, we deduce that:
    \begin{equation}
        \frac{1}{\rho_2(y)}=\tfrac{1}{\rho_2(x)}\left(1-\langle u,\rho_2^{-1}\nabla\rho_2\rangle+p_{x,2}(u)+p_{x,3}(u)\right) +\bigO(\lVert u \rVert^4)
    \end{equation}
    Thus, letting $\xi=\lVert u-\varepsilon \mu\rVert^2/(\varepsilon\rho_1(x)\rho_2(x))$ and using the results of lemma \ref{lemma:odk_metric_expansion}, the terms inside of $h(\cdots)$ will equal:
    \begin{equation}
        \xi-\xi\langle u,\rho_2^{-1}\nabla \rho_2\rangle +\xi \left[p_{x,2}(u)+p_{x,3}(u)\right]+\frac{\Theta_4(u)}{\varepsilon\rho_1\rho_2}+\left(\frac{\Theta_5(u)}{\varepsilon \rho_1\rho_2}-\frac{\Theta_4(u)}{\varepsilon\rho_1\rho_2}\langle u,\rho_2^{-1}\nabla\rho_2\rangle\right) +\bigO(\varepsilon^2,\lVert u\rVert^4)
    \end{equation}
    Applying Lemma \ref{appAlemma:simplifying_lemma}, we see that the third term is of type (iii), and given the forms of $\Theta_4(u)$ and $\Theta_5(u)$ given in the proof of Theorem \ref{appAthm:odk_backward_expansion}, these will be a mixture of type (i) and (ii) terms, and thus none of these terms have any effect on the limiting operator of the self-tuning ODK.  As a result, we will ignore these terms entirely as if they weren't there.  Expanding $h(\cdots)$ about $h(\xi)$ are backward operator becomes:
    \begin{equation}
        \varepsilon^{d/2}G_\varepsilon f(x)=\int h(\xi)f(y)\mathrm{d}y-\int \langle u,\rho_2^{-1}\nabla\rho_2\rangle \xi h'(\xi)f(y)\mathrm{d}y+\frac{1}{2}\int \langle u,\rho_2^{-1}\nabla\rho_2\rangle^2 \xi^2 h''(\xi)f(y)\mathrm{d}y + \bigO(\lVert u\rVert^5)
    \end{equation}
    The first integral is a standard local kernel, and after setting $\rho(x)=\sqrt{\rho_1(x)\rho_2(x)}$ we know the limiting operator will equal $\mathcal{L}=\mu\cdot\nabla+\tfrac{1}{2}\sigma_h^2\rho_1\rho_2\Delta$, where $\sigma_h(x)^2$ is as in equation \eqref{eq:sigma_h_def}.  The third integral, by the same reasoning as the proof given in Lemma \ref{appAlemma:simplifying_lemma}, will be absorbed into $f(x)\omega_l(x)$ and $\Omega_l(x)$ and thus does not affect the limiting operator.  Let $(\dagger)=-\varepsilon^{-d/2}\int \langle u,\rho_2^{-1}\nabla\rho_2\rangle \xi h'(\xi)f(y)\mathrm{d}y$, changing this integral to projection coordinates and substituting $u\mapsto \varepsilon^{1/2}w$, we find:
    \begin{align}
        -(\dagger)&=\varepsilon^{-d/2}\int\langle u,\frac{\nabla\rho_2}{\rho_2}\rangle \xi h'(\xi)f(y)\mathrm{d}y=\varepsilon^{-d/2}f\int\langle u,\frac{\nabla\rho_2}{\rho_2}\rangle\xi h'(\xi)\mathrm{d}u+\varepsilon^{-d/2}\int \langle u,\nabla f\rangle\langle u,\frac{\nabla\rho_2}{\rho_2}\rangle\xi h'(\xi)\mathrm{d}u\\
        &=\sqrt{\varepsilon}f\int_{\mathbb{R}^d}\langle w,\frac{\nabla\rho_2}{\rho_2}\rangle \xi h'(\xi)\mathrm{d}w+\varepsilon\int_{\mathbb{R}^d}\langle w,\nabla f\rangle\langle w,\frac{\nabla\rho_2}{\rho_2}\rangle\xi h'(\xi)\mathrm{d}w+\bigO(\varepsilon^{3/2})\label{appAeq:stproof_intermediary1}
    \end{align}
    where we have extended the integral to the whole tangent space by using the fact that the tails are $\bigO(\varepsilon^2)$; we also note that $\xi=\tfrac{\lVert u-\varepsilon\mu\rVert^2}{\varepsilon\rho_1\rho_2}$ now takes the form $\xi=\frac{\lVert w-\sqrt{\varepsilon}\mu\rVert^2}{\rho_1\rho_2}$.  If we now make the substitution $z=w-\sqrt{\varepsilon}\mu$, then the first term of equation \eqref{appAeq:stproof_intermediary1} will break into two integrals of order $\sqrt{\varepsilon}$ and $\varepsilon$ in the following manner:
    \begin{equation}
        \sqrt{\varepsilon}f\int_{\mathbb{R}^d}\langle w,\rho_2^{-1}\nabla\rho_2\rangle \xi h'(\xi)\mathrm{d}w=\sqrt{\varepsilon}f\int_{\mathbb{R}^d}\langle z,\rho_2^{-1}\nabla\rho_2\rangle\xi h'(\xi)\mathrm{d}z+\varepsilon f\langle\mu,\rho_2^{-1}\nabla\rho_2\rangle\int_{\mathbb{R}^d}\xi h'(\xi)\mathrm{d}z
    \end{equation}
    We notice that since $\xi h'(\xi)$ is symmetric about the origin in $z$-coordinates, the first integral will equal zero since it is the integral of an odd function (since an odd function, $z$, multiplied by an even function, $\xi h'(\xi)$, is itself odd) over a symmetric domain.  The second term is of order $\bigO(\varepsilon)$ and equals some $x$-dependent value (not depending on $f$) multiplied by $f$, and thus can be absorbed into the $\omega_l(x)f(x)$ term of the expansion given in Theorem \ref{appAthm:odk_backward_expansion}, and as we know these $\omega_l(x)$ terms cancel upon normalisation by $G_\varepsilon\mathbb{1}(x)$.  Thus we can entirely ignore the first integral in equation \eqref{appAeq:stproof_intermediary1} and focus solely on the second integral, which, upon making the same $z$ substitution and collecting all higher-order terms into $\bigO(\varepsilon^{3/2})$, can be written as:
    \begin{equation}
        \varepsilon\int_{\mathbb{R}^d}\langle z,\nabla f\rangle\langle z,\rho_2^{-1}\nabla\rho_2\rangle\xi h'\left(\xi\right)\mathrm{d}z=\varepsilon\rho_2^{-1}\partial_i\rho_2\partial_j f\int_{\mathbb{R}^d}z_iz_j\xi h'\left(\xi\right)\mathrm{d}z=\varepsilon\rho_2^{-1}\partial_i\rho_2\partial_i f\int_{\mathbb{R}^d}z_i^2\xi h'\left(\xi\right)\mathrm{d}z
    \end{equation}
    where we have expanded the inner product terms and again used the fact that when $i\neq j$ $z_iz_j$ is an odd function, and integrating an odd function over a symmetric domain cancels the integral.  Since $\xi h'(\xi)$ is symmetric in any component of $z$ we need only compute the value for $z_1$, thus we define the integral $I(t)$ parametrised by $t>0$ to be:
    \begin{equation}
        I(t)=\int_{\mathbb{R}^d}z_1^2h\left(\xi t\right)\mathrm{d}z=t^{-\tfrac{d+2}{2}}\int_{\mathbb{R}^d}\tilde{z}_1^2h(\xi)\mathrm{d}\tilde z=t^{-\tfrac{d+2}{2}}m\sigma_h^2\rho_1\rho_2
    \end{equation}
    where in the second equality we have substituted $\tilde z=\sqrt{t}z$, and in the final equality we have applied the fact that $m^{-1}\int z_1^2h\left(\tfrac{\lVert z\rVert^2}{\rho_1\rho_2}\right)\mathrm{d}z=\sigma_h^2\rho_1\rho_2$ as before.  Taking the derivative of $I(t)$ with respect to $t$ gives:
    \begin{equation}
        I'(t)=\int_{\mathbb{R}^d}z_1^2\xi h'(\xi t)\mathrm{d}z=-\frac{d+2}{2}t^{-\tfrac{d+4}{2}}m\sigma_h^2\rho_1\rho_2
    \end{equation}
    and we notice that $\int_{\mathbb{R}^d}z_1^2\xi h'(\xi)\mathrm{d}z=I'(1)=-\tfrac{1}{2}(d+2)m\sigma_h^2\rho_1\rho_2$, and thus:
    \begin{align}
        (\dagger)&=-\varepsilon\rho_2^{-1}\partial_i\rho_2\partial_i f\int_{\mathbb{R}^d}z_i^2\xi h'\left(\xi\right)\mathrm{d}z+\bigO(\varepsilon^{3/2})=-\varepsilon\rho_2^{-1}\langle \nabla\rho_2,\nabla f\rangle I'(1)+\bigO(\varepsilon^{3/2})\\
        &=\varepsilon m(d+2)D_{st}\frac{\nabla\rho_2}{\rho_2}\cdot\nabla f+\bigO(\varepsilon^{3/2})
    \end{align}
    where $D_{st}\coloneqq\tfrac{1}{2}\sigma_h^2\rho_1\rho_2$.  Putting this together with the local kernels expansion, we have:
    \begin{equation}
        G_\varepsilon f(x)=m(x)f(x)+\varepsilon\left(\omega_l^{st}(x)f(x)+m(x)\mathcal{L}_{st}f(x)\right)+\varepsilon^{3/2}\Omega_l^{st}(x)+\bigO(\varepsilon^2)
    \end{equation}
    where $\mathcal{L}_{st}f=\tfrac{\tau}{\kappa}\mu\cdot\nabla f+(d+2)D_{st}\frac{\nabla \rho_2}{\rho_2}\cdot\nabla f+D_{st}\Delta f$, and $\mu(x)=\nabla l(x)$ or $\mu(x)=\nabla l_x(x)$ depending on whether a local or global ordering function has been used.  By normalising by $G_\varepsilon\mathbb{1}(x)$ as before we find:
    \begin{equation}
        P_\varepsilon f(x)=f(x)+\varepsilon\mathcal{L}_{st}f(x)+\bigO(\varepsilon^{3/2})\implies L_\varepsilon f(x)=\mathcal{L}_{st}f(x)+\bigO(\varepsilon^{1/2})
    \end{equation}
    The proof of the adjoint expansion is identical to the proof of Theorem \ref{appAthm:odk_forward_expansion} and is done in the weak sense, to give:
    \begin{equation}
        G_\varepsilon^*f(x)=m(x)f(x)+\varepsilon(\omega_l(x)f(x)+\mathcal{L}_{st}^*(fm)(x))+\bigO(\varepsilon^{3/2})
    \end{equation}
    and:
    \begin{equation}
        P_\varepsilon^* f(x)=f(x)+\varepsilon\mathcal{L}_{st}^* f(x)+\bigO(\varepsilon^{3/2})\implies L_\varepsilon^* f(x)=\mathcal{L}_{st}^* f(x)+\bigO(\varepsilon^{1/2})
    \end{equation}
    which completes the proof.
\end{proof}

An immediate corollary of this result is that it strengthens the proof given by Berry \& Harlim in ``Variable Bandwidth Kernels" (VBKs, \cite{berry2016variable}) since in their Appendix A.3 they prove the expansion of the \textit{right-formulation} in the weak sense, whereas we do not take such a step, and by setting $\tau(x)\equiv 0$, $l(x)\equiv 0$ and $\kappa(x)\equiv 1$ we recover the form of a VBK.  The reason Berry \& Harlim cannot make such a step is because of the order of operations they take in performing their normalisation.  Particularly, they would need to make a substitution of the form $\hat y=\tfrac{y-x}{\sqrt{\rho_2(y)}}+x$, but such a substitution need not be invertible and thus cannot be directly applied.  If instead they took the expansion as we do, and then make the substitution $z\mapsto \tilde z/\sqrt{\rho_1(x)\rho_x(x)}$, the integrals defining $\mu(x)$ and $D_{st}(x)$ would provide the necessary factors of $\rho_1^{d/2}\rho_2^{d/2}$ that VBKs use in their normalisation.

\subsection{Density Regularised ODKs}\label{appAsec:densreg}
In the literature, particularly in Diffusion Maps \cite{coifman2006diffusion}, they introduce an $\alpha\in[0,1]$ parameter to tune the strength of density regularisation, with $\alpha=0$ being no density regularisation, $\alpha=1/2$ being their \textit{Fokker-Planck} normalisation, and $\alpha=1$ being fully regularised.  This amounts to defining $K_{\varepsilon,\alpha}(x,y)\coloneqq K_\varepsilon(x,y)q(y)^{-\alpha}$ and only makes sense if we change the integrals for the backward and forward integral operators to be with respect to the underlying sampling measure $q(x)$.  In other words, we now consider the operator:
\begin{equation}
    G_{\varepsilon,\alpha}f(x)\coloneqq \varepsilon^{-d/2}\int_\mathcal{M}K_{\varepsilon,\alpha}(x,y)f(y)q(y)\ \mathrm{d}y=\varepsilon^{-d/2}\int_\mathcal{M} K_\varepsilon(x,y) f(y)q(y)^{1-\alpha}\ \mathrm{d}y=G_\varepsilon(fq^{1-\alpha})(x)
\end{equation}
from the results of Theorem \ref{thm:odk_backward_expansion} we can deduce the expansion of the non-uniform normalised transition operator $P_{\varepsilon,\alpha}f\coloneqq (G_{\varepsilon,\alpha}\mathbb{1})^{-1}G_{\varepsilon,\alpha}f$, we can write this as:
\begin{align}
    P_{\varepsilon,\alpha}f&= \frac{G_{\varepsilon}(fq^{1-\alpha})}{G_{\varepsilon}(q^{1-\alpha})}=\frac{P_\varepsilon(fq^{1-\alpha})}{P_\varepsilon(q^{1-\alpha})}=\frac{fq^{1-\alpha}+\varepsilon\mathcal{L}(fq^{1-\alpha})+\bigO(\varepsilon^{3/2})}{q^{1-\alpha}+\varepsilon\mathcal{L}(q^{1-\alpha})+\bigO(\varepsilon^{3/2})}\\
    &=q^{\alpha-1}(fq^{1-\alpha}+\varepsilon\mathcal{L}(fq^{1-\alpha}))(1-\varepsilon q^{\alpha-1}\mathcal{L}(q^{1-\alpha}))+\bigO(\varepsilon^{3/2})\\
    &= f+\varepsilon q^{\alpha-1}\left(\mathcal{L}(fq^{1-\alpha})-f\mathcal{L}(q^{1-\alpha})\right)+\bigO(\varepsilon^{3/2})
\end{align}
since $\mathcal{L}f=\mu\cdot\nabla f+D\Delta f$, we can compute:
\begin{align}
    \mathcal{L}(q^{1-\alpha})&=(1-\alpha)q^{-\alpha}\mu\cdot\nabla q+D\Delta q^{1-\alpha}\\
    \mathcal{L}(fq^{1-\alpha})&=q^{1-\alpha}\mu\cdot\nabla f+(1-\alpha)fq^{-\alpha}\mu\cdot\nabla q+D\left(f\Delta q^{1-\alpha}+q^{1-\alpha}\Delta f+2(1-\alpha)q^{-\alpha}\nabla q\cdot\nabla f\right)
\end{align}
putting this together, we deduce:
\begin{align}\label{eq:alphanormalised_TP}
    P_{\varepsilon,\alpha}f&=f+\varepsilon q^{\alpha-1}\left(q^{1-\alpha}\mu\cdot\nabla f+D(q^{1-\alpha}\Delta f+2(1-\alpha)q^{-\alpha}\nabla q\cdot\nabla f)\right)+\bigO(\varepsilon^{3/2})\\
    &=f+\varepsilon\left(\mu\cdot\nabla f+2D(1-\alpha)q^{-1}\nabla q\cdot\nabla f+D\Delta f\right)+\bigO(\varepsilon^{3/2})\\
    &=f+\varepsilon(\mathcal{L}f+2D(1-\alpha)q^{-1}\nabla q\cdot\nabla f)+\bigO(\varepsilon^{3/2})
\end{align}
This gives us the following form for the $\alpha$-regularised Laplacian:
\begin{equation}
    \mathcal{L}_\alpha f=\mu\cdot\nabla f+2(1-\alpha)D\frac{\nabla q}{q}\cdot\nabla f+D\Delta f
\end{equation}
We see that if one doesn't fully regularise for the sampling density, then there will be another advection term proportional to $\nabla q$ that is weighted by a factor of $(1-\alpha)$.  If one progressively increases $\alpha$ from $0$ to $1$, then the contribution of this additional term decreases.  If one wants an exact reconstruction of the underlying dynamics then one should always consider the $\alpha=1$ case; but if one's data is highly noisy (many outliers) then pushing the effective drift towards the core concentration of the data may be wise depending on your ultimate goal.  The $\alpha$-normalised adjoint has an expansion analogous to the previous results given in Theorem \ref{thm:odk_forward_expansion}.

\subsection{Discrete Operator Convergence}\label{appAsec:disc2cont}

The final proof we provide completes the theory that underlies most kernel methods, and gives convergence in probability of the discrete operators to the continuous operators.  In what follows, we will use the Bernstein inequality as given below:
\begin{lemma}[I.I.D. Bernstein Inequality]\label{appA:bernstein}
    Let $X_1,X_2,\dots,X_m$ be independent and identically distributed zero-mean real-valued random variables, and suppose $\exists a,b>0$ such that $|X_i|\leq a$ and $\mathbb{E}[X_i^2]\leq b$.  Then for any $\delta>0$:
    \begin{equation}\label{appAeq:bernstein}
        P\left(\left|\frac{1}{m}\sum_{i=1}^mX_i\right|>\delta\right)\leq2\exp\left(-\frac{\tfrac{1}{2}m\delta^2}{b+\tfrac{1}{3}a\delta}\right)
    \end{equation}
\end{lemma}

\begin{theorem}[Discrete Operator Convergence]\label{appAthm:discrete_conv_proof}
    Let $\mathcal{M}$ be a $d$-dimensional Riemannian manifold embedded in $\mathbb{R}^n$ and with metric inherited from the standard Euclidean metric in ambient space.  Suppose $\mathcal{D}=\{x_i\}_{i=1}^m$ is a dataset of $m$ independent samples from the sampling distribution $q(x)$ which we assume satisfies $\inf_{x\in\mathcal{M}}q(x)\backwardscoloneqq c_q>0$, i.e. $\text{supp}(q)=\mathcal{M}$ and $q$ doesn't decay to zero anywhere on $\mathcal{M}$ - this condition necessitates that $\mathcal{M}$ be bounded.  Suppose $K_\varepsilon(x,y)$ is an Ordered Diffusion Kernel and $L_\varepsilon$ its associated Laplacian, with the sampling density fully regularised away ($\alpha=1$), and suppose, after fixing an indexing of $\mathcal{D}$, $L_{m,\varepsilon}\in\mathbb{R}^{m\times m}$ is the discretisation of $L_\varepsilon$ onto $\mathcal{D}$ (i.e. $[L_\varepsilon f]\approx L_{m,\varepsilon}[f]=\frac{1}{m}[L_\varepsilon]\cdot[f]$).  Then for any bounded $f\in L^2(\mathcal{M},q)\cap C^3(\mathcal{M})$, i.e. $\exists C_f>0$ such that $\sup_{x\in\mathcal{M}}|f(x)|=C_f$, and any $\delta>0$:
    \begin{equation}
        P\left(\left\lVert L_{m,\varepsilon}[f]-[L_\varepsilon f]\right\rVert_\infty>\delta\right)\leq2m\exp\left(-\frac{\tfrac{1}{2}m\varepsilon^{d/2+2}\delta^2}{C_f(C_K+\varepsilon^{d/2})(C_f+\tfrac{1}{3}\varepsilon\delta)}\right)
    \end{equation}
    where $C_K>0$ is a constant that depends on $K_\varepsilon$ and $q$ and is independent of $x_i$.
\end{theorem}
\begin{proof}
    We will proceed rather directly, for any $x\in\mathcal{D}$ we define $V_x(y)=P_\varepsilon(x,y)f(y)$, and $W_x(y)=V_x(y)-\mathbb{E}_{y\sim q}[V_x(y)]$ - it is clear that $\mathbb{E}_{y\sim q}[W_x(y)]=0$. We recall that $P_\varepsilon(x,y)\coloneqq \frac{K_\varepsilon(x,y)}{Z(x)q(y)}$, where $Z(x)=\int_\mathcal{M}K_\varepsilon(x,y)\mathrm{d}y$.  Since $K_\varepsilon(x,y)=h(u(x,y))\leq C_h\exp(-u(x,y)/\sigma^2)$ we know that $K_\varepsilon(x,y)$ is bounded from above by $C_h>0$; we also know that, due to the $\varepsilon$ scaling in the defining of an ODK, $Z(x)=\varepsilon^{d/2}m(x)$ where $m(x)$ is as in Theorem \ref{thm:odk_backward_expansion}, let us assume that $\inf_{x\in\mathcal{M}}m(x)=c_m>0$ - this is a sort of non-degeneracy condition.  By combining these with the assumptions that $|f(x)|\leq C_f$ and $q(x)\geq c_q$, then we can deduce:
    \begin{align}
        |W_x(y)|&\leq |V_x(y)|+|\mathbb{E}_{y\sim q}[V_x(y)]|=|P_\varepsilon(x,y)f(y)|+\left|\int_\mathcal{M}P_\varepsilon(x,z)f(z)q(z)\mathrm{d}z\right|\\
        &\leq \frac{C_hC_f}{\varepsilon^{d/2}c_mc_q}+C_f^2\int_\mathcal{M}\frac{K_\varepsilon(x,z)}{Z(x)}\mathrm{d}z\backwardscoloneqq \varepsilon^{-d/2}C_KC_f+C_f^2
    \end{align}
    where we have defined $C_K=C_hc_m^{-1}c_q^{-1}$.  We can also deduce:
    \begin{align}
        |\mathbb{E}_{y\sim q}[W_x(y)^2]|&\leq|\mathbb{E}_{y\sim q}[V_x(y)^2]|+|\mathbb{E}_{y\sim q}[V_x(y)]^2|=\int_\mathcal{M}P_\varepsilon(x,y)^2f(y)^2q(y)\mathrm{d}y+\left(\int_\mathcal{M}P_\varepsilon(x,y)f(y)q(y)\mathrm{d}y\right)\\
        &\leq \frac{C_K C_f^2}{\varepsilon^{d/2}}\int_\mathcal{M}\frac{K_\varepsilon(x,y)}{Z(x)}\mathrm{d}y+C_f^2\left(\int_\mathcal{M}\frac{K_\varepsilon(x,y)}{Z(x)}\mathrm{d}y\right)^2=\varepsilon^{-d/2}C_KC_f^2+C_f^2
    \end{align}
    We can now apply Bernstein's inequality, as given in Lemma \ref{appA:bernstein}, to any $x_i\in\mathcal{D}$:
    \begin{align}
        P\left(\left|\frac{1}{m}\sum_{j=1}^mW_{x_i}(x_j)\right|>\delta\right)&=P\left(\left|\frac{1}{m}\sum_{j=1}^mP_\varepsilon(x_i,x_j)f(x_j)-\mathbb{E}_{y\sim q}[P_\varepsilon(x_i,y)f(y)]\right|>\delta\right)\\
        &\leq 2\exp\left(-\frac{\tfrac{1}{2}m\varepsilon^{d/2}\delta^2}{C_f^2(C_K+\varepsilon^{d/2})+\tfrac{1}{3}C_f(C_K+\varepsilon^{d/2})\delta}\right)
    \end{align}
    We note that if we added and subtracted $f(x_i)$ and replaced $\delta\mapsto\varepsilon\delta$ in the second expression, then all $P_\varepsilon$ terms could be replaced by $L_\varepsilon$ terms; thus:
    \begin{equation}
        P\left(\left|\frac{1}{m}\sum_{j=1}^mL_\varepsilon(x_i,x_j)f(x_j)-L_\varepsilon f(x_i)\right|>\delta\right)\leq2\exp\left(-\frac{\tfrac{1}{2}m\varepsilon^{d/2+2}\delta^2}{C_f(C_K+\varepsilon^{d/2})(C_f+\tfrac{1}{3}\varepsilon\delta)}\right)
    \end{equation}
    The first term within the modulus expression is equal to $(L_{m,\varepsilon}[f])_i$ and the second term equals $[L_\varepsilon f]_i$.  The $L^\infty$-norm across all $i=1,\dots,m$ will equal $\left\lVert L_{m,\varepsilon}[f]-[L_\varepsilon f]\right\rVert_\infty=\max_{i=1,\dots,m}|(L_{m,\varepsilon}[f])_i-[L_\varepsilon f]_i|$, by the union bound inequality:
    \begin{align}
        P\left(\left\lVert L_{m,\varepsilon}[f]-[L_\varepsilon f]\right\rVert_\infty>\delta\right)&\leq \sum_{i=1}^mP\left(\left|(L_{m,\varepsilon}[f])_i-[L_\varepsilon f]_i\right|>\delta\right)\\
        &\leq2m\exp\left(-\frac{\tfrac{1}{2}m\varepsilon^{d/2+2}\delta^2}{C_f(C_K+\varepsilon^{d/2})(C_f+\tfrac{1}{3}\varepsilon\delta)}\right)
    \end{align}
\end{proof}

\newpage
\section{Supplementary Notes}\label{appendix_B}

\subsection{Ordering Functions}\label{appBsubsec:ordering_functions}
In this supplementary note, we seek to clarify what is meant by an \textit{Ordering Function}.  It is important to say that an ordering function, as we conceptualise it, is not really something that can be defined independently based on some distinguished set of properties, and is more like a representation of some mathematical ordering of elements (i.e. a notion of something being larger or smaller, or happening later or earlier, than something else)

Orderings of mathematical spaces are studied in order theory, and many deep results come from this field; however, there are many different axiomatic constructions to choose from when defining an ordering (\cite{ward1954partially,krishnan2009convenient,haucourt2012streams}).  The idea of an ordering function is not novel in principle, since ultimately what an ordering function does is try to map the natural ordering of the real numbers onto some abstract mathematical space; in the literature, concepts like our ordering function are referred to as the utility functions or utility representations \cite{troue1969topology,peleg1970utility,evren2011multi,bosi2012continuous}.  We will give a discussion of how an ordering function fits within the more well-studied field of representations reflexive and transitive quasi-orders, sometimes called pre-orders, although we note that it may be more appropriate to consider a more specialised construction such as d-spaces \cite{grandis2009directed}, which consider orderings defined by a set of directed paths, or streams \cite{krishnan2009convenient}, which considers a \textit{circulation} that associates a quasi-order to the open sets of a topological space, both of which generalise quasi-orders.
\begin{definition}[Quasi/Pre-Order]
    Let $X$ be some set.  A binary relation $\leq:X\times X\to\{0,1\}$ is called a quasi-order if it satisfies the following two conditions:
    \begin{enumerate}[noitemsep,nolistsep]
        \item[i).] (Reflexivity): $\forall x\in X,\ x\leq x$
        \item[ii).] (Transitivity): $\forall x,y,z\in X,\ (x\leq y)\land(y\leq z)\implies (x\leq z)$
    \end{enumerate}
    If, in addition, the quasi-order satisfies the following linearity (totality) condition, we call $\leq$ a linear (total) quasi-order (LQO):
    \begin{enumerate}[noitemsep,nolistsep]
        \item[iii).] (Linearity): $\forall x,y\in X,\ (x\leq y)\lor (y\leq x)$
    \end{enumerate}
\end{definition}
Let $\leq$ be a quasi-order: then if $x\leq y$ or $y\leq x$, we say that $x$ and $y$ are comparable --- a Linear Quasi-Order (LQO) is thus a quasi-order where all elements of $X$ are comparable; we can define a strict quasi-order, $<$, as a restriction of $\leq$, through $x<y\iff x\leq y\land y\nleq x$; we can define an equivalence relation $x\simeq y\iff x\leq y\land y\leq x$, which we call (quasi-)equality; for two subsets $U,V\subset X$, if for every $x\in U$ and $y\in V$ we have $x\leq y$ ($x< y\text{ or }x\simeq y)$, then we write $U\leq V$ ($U<V$ or $U\simeq V$); for any $A\subset X$ we can define the lower, equal and upper sets as:
\begin{equation}
    L(A)=\{x\in X|\exists y\in A:x\leq y\},\quad U(A)=\{x\in X|\exists y\in A:y\leq x\},\quad E(A)=\{x\in X|\exists y\in A:x\simeq y\}
\end{equation}
A quasi-order is called spacious if for every $x,y\in X$ such that $x\leq y$, then $L(x)\subseteq L(y)$ (if you are given a strict quasi-order this condition would be $\text{cl}(L(x))\subseteq L(y)$ and it is true that the strict restriction of a quasi-order satisfying the former property, will also satisfy the latter property.  A strict quasi-order is said to be separable if for any $x,y\in X$ such that $x<y$, there will exist a $p\in X$ such that $x<p<y$, likewise $\leq$ is separable if its strict restriction is separable; if $X$ is a topological space that is itself separable, and $\leq$ is continuous (see definition \ref{appBdef:cts_qo} below) then $<$ will be separable.  Assume from now on that $X$ is a topological space, currently a quasi-order on $X$ will be entirely independent of the topology on $X$, but we may connect them with notions of continuity that mirror $T_1$ and $T_2$ (Hausdorff) spaces: effectively, if $\leq$ is lower-(upper-)semicontinuous then for all $x,y\in X$ such that $x\nleq y$ ($y\nleq x$) there exists an open set of $U\subset X$ containing $a$ such that $U\nleq b$ ($b\nleq U$); a result of this is that if $\leq$ is lower (upper) semicontinuous then its lower (upper) sets will be topologically closed.  A quasi-order is called semicontinuous if it is both upper- and lower-semicontinuous.  A quasi-order is said to be continuous if:
\begin{definition}\label{appBdef:cts_qo}
    A quasi-order, $\leq$, of a topological space $X$, is said to be continuous if $\forall x,y\in X$ such that $x\nleq y$, there exists open neighbourhoods $U$ and $V$ of $x$ and $y$, respectively, such that $U\cap V=\varnothing$ and $\forall u\in U$ and $\forall v\in V$ then $u\nleq v$.
\end{definition}
It is true that continuity implies semicontinuity, and indeed, if $\leq$ is linear, semicontinuity and continuity are equivalent \cite{ward1954partially}.  For a continuous quasi-order, the lower and upper sets will be topologically closed, and furthermore, the lower and upper sets of the strict restriction will be open; a continuous quasi-order will have a closed graph in the product topology of $X\times X$, whereas semicontinuity does not guarantee this.  For what follows, we will always assume that $\leq$ is continuous.

The notion of an order we have built so far is very abstract and rather difficult to study and visualise, but we realise that for any reasonable ordering we might consider we are effectively mapping the only canonical ordering we have, that is the ordering of the real numbers, onto local regions of $X$.  In other words, under mild assumptions on $X$ and $\leq$, we have for every $x\in X$ a neighbourhood $U$ of $x$ and an order-preserving (i.e. strictly increasing/isotone) map $f:U\to\mathbb{R}$; the question becomes when will these local functions be \textit{nice} (i.e. continuous) and can we glue these local $f$'s together to form a global continuous function $f:(X,\leq)\to(\mathbb{R},\leq)$ that satisfies:
\begin{equation}\label{appBeq:richter_peleg_rep}
    \forall x,y\in X: x< y\implies f(x)< f(y)
\end{equation}
where we must use the strict relation, as otherwise any constant $f$ would be satisfactory.  In the literature, these functions are referred to as utility functions, and more broadly they form a utility representation of quasi-ordered topological spaces.\footnote{More precisely, this is a Richter-Peleg representation.}  We note that in the representation given in equation \eqref{appBeq:richter_peleg_rep}, $f$ does not entirely determine/characterise the ordering since the implication only goes one way; however, generalising this to an if and only if statement would immediately imply that $\leq$ is linear, which is a strong condition.  This motivated the development of so-called multi-utility representations of quasi-orders, in which we take a (potentially uncountable) family of functions $\boldsymbol{F}=\{f_i\}_{i\in I}$ such that $\forall x,y\in X: x\leq y\implies \left(f(x)\leq f(y)\text{ for all }f\in\boldsymbol{F}\right)$.  Thus, the family $\boldsymbol{F}$ entirely characterises the ordering, which feels more useful; there is still the issue of existence that we tackle in the following proposition:
\begin{proposition}\label{appBprop:utility_representations}
    Let $(X,\leq)$ be a continuous quasi-ordered topological space.  The following will hold:
    \begin{enumerate}[noitemsep,nolistsep]
        \item[(i).] If the associated strict quasi-order, $<$, is separable and spacious, then there exists a continuous utility function, $f:X\to\mathbb{R}$, such that $x\leq y\implies f(x)\leq f(y)$, for $(X,\leq)$.
        \item[(ii).] If the associated strict quasi-order, $<$, is separable and spacious, then there exists a continuous utility function, $f:X\to\mathbb{R}$, such that $x\leq y\iff f(x)\leq f(y)$ if and only if $\leq$ is linear.
        \item[(iii).] If $X$ is a locally-compact Hausdorff space that is also $\sigma$-compact,\footnote{A topological space $X$ is Hausdorff if any two points are jointly distinguished by disjoint open sets; $X$ is locally compact if around every point there exists an open set whose closure is compact; a locally compact space is $\sigma$-compact if it can be written as a countable union of these compact sets.} then there always exists a continuous multi-utility representation of $(X,\leq)$.  If in addition $X$ is second-countable ($\iff C(X)$ is separable), then this representation can be taken to be countable.
    \end{enumerate}
\end{proposition}
\begin{proof}
    (i) This is a known result of Peleg (Theorem 3.1 of \cite{peleg1970utility}) that shows for a separable and spacious, continuous strict quasi-order, $<$, there will exist a continuous utility function, $f:X\to\mathbb{R}$, in the form $x<y\implies f(x)<f(y)$.  If we first take the quotient of $X$ by the equality relation, $\simeq$, then we can define the relation $\preceq$, on $X/\simeq$ through: $[x]\preceq [y]\iff \exists \tilde x\in[x],\tilde{y}\in[y]:\tilde{x}\leq \tilde y$.  This relation will be well-defined, reflexive and transitive and thus a quasi-order; furthermore, if $<$ is separable and spacious, so too will $\prec$.  Therefore by the Theorem of Peleg, there will exist $\tilde f:(X/\simeq)\to\mathbb{R}$ such that $[x]\prec [y]\implies \tilde f([x])<\tilde f([y])$ that is constant on the equal sets of $(X,\leq)$ (i.e. if $x\simeq y$ then $f([x])=f([y])$).  Thus, by properties of quotient spaces, there will be a continuous function $f:X\to\mathbb{R}$ such that $x<y\implies f(x)< f(y)$ and $x\simeq y\implies f(x)=f(y)$.
    
    (ii) The forward implication is trivial, and (i) shows that there exists a continuous $f$ such that $x< y\implies f(x)< f(y)$; if $x\simeq y$, then $f(x)= f(y)$, if $x<y$ then $f(x)< f(y)$, and if $y<x$ $f(y)<f(x)$, thus $f$ entirely characterises $\leq$ and by linearity we are done.
    
    (iii) This is a known result and derives from a result of Evren \& Ok (see Theorem 1 of \cite{evren2011multi}), which states the existence of a continuous multi-utility representation.  The countability of this representation comes from locally compact Hausdorff spaces that are also $\sigma$-compact and second-countable, having a countable basis that makes $C(X)$ separable.
\end{proof}
There are several takeaways here: first, utility representations are abundant for \textit{nice} quasi-orders (continuous, separable and spacious) defined on \textit{nice} topological spaces (effectively, we want $X$ to be metrisable, and we note that second-countable, locally compact Hausdorff spaces that are $\sigma$-compact are always metrisable).  Thus, if we are focused on continuous linear quasi-orders, we may as well consider the function $l:X\to\mathbb{R}$, that generates the LQO, $\leq_l$, through:
\begin{equation}
    x\leq_l y\iff l(x)\leq l(y)
\end{equation}
as being the fundamental object of study, and thus we call this an ordering function (although, as we mentioned at the beginning, it doesn't seem correct to define the ordering function outright); we could also view the multiple-Ordered Diffusion Kernels (mODKs) of section \ref{subsec:multi_odk} as being analogous to multi-utility representations.  If we now restrict ourselves to defining ordering functions on a $d$-dimensional topological manifold, i.e. $l:\mathcal{M}\to\mathbb{R}$, then we see that if $dl_x\neq 0$ for any $x$ in an open set $U\subset \mathcal{M}$, then the level sets of $l$ form a codimension-one foliation of $U$, i.e. we have sliced $U$ along the level sets of $l$, which are equivalently the equal sets of the order $\leq_l$, and this slicing is such that each level set is a $(d-1)$-dimensional submanifold.  We can then view the second term of the ODK:
\begin{equation}
    K_\varepsilon(x,y)=h\left(\frac{\lVert y-x\rVert^2}{\varepsilon\rho(x)^2}+\frac{(l(y)-l(x)-\varepsilon\tau(x)\kappa(x))^2}{\varepsilon\rho(x)^2\kappa(x)^2}-\frac{(l(y)-l(x))^2}{\varepsilon\rho(x)^2\kappa(x)^2}\right)
\end{equation}
as an advection across the leaves of the foliation, and the first term as a diffusion applied after this initial advection.  By Theorem \ref{thm:odk_backward_expansion}, we see that the first-order term of the expansion of $K_\varepsilon$ is proportional to $\nabla l$, thus, (spacious and separable) continuous linear quasi-orders (whose equal sets are sufficiently smooth) defined on differentiable manifolds is entirely characterised by a gradient flow; we could say that the order is generated by this flow.  Thus these $C^2(\mathcal{M})$ ordering functions characterise all sufficiently smooth linear quasi-orders; this motivates us to think of an ordering as some representation of a vector field, and thus the existence of a suitable ordering function depends on some notion of integrability of this vector field.  Furthermore, if we allow for our vector fields to be more general than gradient-flows, then statement (ii) in proposition \ref{appBprop:utility_representations}, gives us another characterisation relation (an if and only if statement) between orderings and multi-utility representations, thus taking an mODK with a finite ($k=\min(2d,n)$) number of ordering functions would allow us to represent arbitrary vector-field proportional to $\sum_{i=1}^k\tau_i(x)\nabla l_i(x)$, thus allowing us to generalise beyond (sufficiently smooth) linear quasi-orders to (sufficient smooth) non-linear (partial) quasi-orders.

However, if we are now trying to view an ordering that is both nice and as generally applicable as possible, then we notice that (global) transitivity limits us from representing cyclic orders.  As an example, suppose we wished to order the hour hand of a clock: clearly $1<2$, but we also have that $2<3$, and $3<4$, and so on, until we have $12<1$; by transitivity, this implies $2<1$ and thus our ordering is contradictory.  To order the clock, or more generally to order cyclic domains, we need to relax the assumption of transitivity; there have been several formulations of such an idea of \textit{locally ordered topological spaces}, these include: a very intuitive formulation in terms of open sets through local partial orders, whereby partial we mean that $\leq$ is also antisymmetric\footnote{Where antisymmetry is taken to mean $\forall x,y\in X$ if $x\leq y\land y\leq x\implies x=y$ where equality is of elements not of orders, i.e. not as with $\simeq$ as we have had with quasi-orders.}, such a space is called a \textit{local pospace} (see \cite{fajstrup2006algebraic} for more details); a simple removal of transitivity, where $\leq$ is simply reflexive and its graph, $\{(x,y)\in X\times X:x\leq y\}$, is closed with respect to the product topology, we call such spaces weakly ordered topological spaces (or \textit{wospaces}); a characterisation of an ordering based on a set of distinguished paths, $dX=\{\gamma:I\to X\}$, that satisfies all constant path are in $dX$, $dX$ is closed under reparametisation of the interval $I$ by weakly increasing functions, and $dX$ is closed under concaternation of consecutive paths in $dX$, thus we can define a relation $x\leq y$ if there is a path in $dX$ originating at $x$ and ending at $y$, we call such spaces \textit{d-spaces} (see \cite{grandis2009directed} for more details); a characterisation of an ordering based on a circulation that is effectively a pre-sheaf, i.e. we can restrict the circulation to a preorder on any open set, we call such an ordering a \textit{stream} (see \cite{krishnan2009convenient} for more details).

All of the above methods (other than weakly ordered spaces) will generalise easily to cyclic orders; the study of local pospaces naturally includes the study of global pospaces (which is similar to the discussion we give above), while d-spaces and streams naturally can represent any continuous quasi-ordered topological space --- therefore, these last two constructions are most aligned with the generalisation we seek to make.  Just as the discussion of continuous quasi-orders defined on a $d$-dimensional manifold led to a discussion of vector fields and integrability, so too will the notion of a local continuous quasi-order lead to a discussion of locally defined vector fields and their local integrability.  Furthermore, since if our construction allows us to say that a subset $U\subset X$ has a local ordering $\leq_U$, then any restriction of this ordering to $V\subset U$ will be a local ordering on $V$, which fits nicely into the description of sheafs and thus functors between categories, local orderings are often studied in this light; intrestingly, it was proven by Haucourt that streams and d-spaces are almost isomorphic in a category theoretic sense \cite{haucourt2012streams}.  

In this paper we will not generalise an ordering function (or utility representation) to these more general notions of an ordering, but will rather define a simple construction --- that mirrors the definition of local pospaces given in \cite{fajstrup2006algebraic} but differs in some important ways --- that fits naturally within the framework of a local ordering function, $l(x,y)$, as introduced in section \ref{subsec:nongradflow_extension}.  

We define a Locally Linear Quasi-Ordered Topological Space (LLQOTS) as follows:
\begin{definition}[LLQOTS]
    A Locally Linear Quasi-Ordered Topological Space (LLQOTS) is a pair $(X,O_X)$, where $X$ is a topological space and $O_X=\{(U,\leq_U)\}$ is called the order set of $X$, and is composed of open linearly ordered subsets of $X$ satisfying the following properties:
    \begin{enumerate}[noitemsep,nolistsep]
        \item[1).] Each $U\in O_X$ is open, connected and has compact closure.
        \item[2).] Every $\leq_U$ is a reflexive, transitive and linear binary relation on $U$ that is continuous and separable with respect to the subspace topology.
        \item[3).] There is no element $U\in O_X$ that is strictly contained in another.
        \item[4).] For every $(U,\leq_U),(V,\leq_V)\in O_X$, the orderings will agree on their intersection.
        \item[5).] There are no extremal (maximal or minimal) elements in any $(U,\leq_U)\in O_X$.
    \end{enumerate}
    Let us denote by $X_r=\bigcup_{U\in O_X}U$ the regular set of the order (or the regularly ordered set), and let $X_c=X\setminus X_r$ be the critical set of the order.  If $\text{cl}(X)=\text{cl}(X_r)$, then we say the local ordering is complete; if, in addition, $\text{cl}(X)=\bigcup_{U\in O_X}\text{cl}(U)$, we say that $O_X$ is locally complete.
\end{definition}
We note that if $O_X$ is complete and $(X_r,\leq_{X_r})\in O_X$, then by assumption (4) we will have returned to the case of continuous LQOs on $X$.  This construction is very similar to how local pospaces are defined, except we do not require partiality and we do not assume that $\bigcup_{U\in O_X}U=X$ --- this relaxation is important when introducing local ordering functions (which we could reasonably call local utility-representations).  To help our intuition we could imagine that we are somewhere on Earth looking out into the distance, our eyes will only allow us to see so far before the curvature of earth breaks our line of sight, or a mountain blocks the view of something behind it; the locally ordered neighbourhoods encode that we can only see so far, and the peaks of mountains would be contained in the closures of these neighbourhoods.

We do not develop the theory of LLQOTs fully in this paper as it does not seem suitable, but we will make the following general claims: the restriction to local orderings allows us to represent cyclic orders (this can be seen in the numerical example given in section \ref{subsec:numerical_results_operators}); formalising the idea of a maximal extension\footnote{An extension, $(W,\leq_W)$, over $O_X$ is any open connected ordered subset of $X$, s.t. $W$ is covered by sets in $O_X$, and $\leq_W$ agrees with $\leq_{U}$ on $W\cap U$ for any $U\in O_X$.  A maximal extension is an extension that is not strictly contained in any other extension.} of a connected ordered set will allow us to study the topology of $X$ through the \textit{criticalities} (the parts we have to cut out in order to define an LQO) of the ordering, and furthermore, the set of all maximal extensions of a set can itself be ordered;\footnote{Furthermore, if we formalise a notion of \textit{path-extensions} it is likely possible to make a direct connection to d-spaces, as described in \cite{grandis2009directed}.} the elements of $X$ that are not contained in $\bigcup_{U\in O_X}U$ will conveniently be the local maximal and minimal elements of the local order, which as we recall from the discussion of criticalities of the ordering function (i.e. $x\in \mathcal{M}$ such that $\nabla l(x)=0$), will allow us to assume the local ordering functions form proper foliations of each $U\in O_X$, thus making it easier to discuss and characterise a local ordering function.

By the assumption (2) above and statement (ii) in proposition \ref{appBprop:utility_representations}, we can associate a continuous function $f_U$ to every $(U,\leq_U)\in O_X$, and by assumption (1) and (5) we know that $f_U$ does not have any extremal values in $U$, but does take both a maximal and minimal value on $\partial U$ (the boundary of $U$, defined as $\partial U=\text{cl}(U)\setminus\text{int}(U)$), since these utility representations are invariant under composition with strictly increasing functions (i.e. if $\varphi:f_U(U)\to\mathbb{R}$ is strictly increasing then $\varphi\circ f_U$ is also a utility function), we can assume that $f_U:U\to I=(0,1)$.  Assumption (4) implies for every $(U,\leq_U),(V,\leq_V)\in O_X$ there exists a strictly increasing continuous bijection $\varphi_U^V:f_U(U\cap V)\to f_V(U\cap V)$ such that $f_V=\varphi_U^V\circ f_U$.

Assume now that $X$ is Hausdorff and $O_X$ is locally complete, i.e. every point $x$ is in the closure of at least one of the ordered sets.  We seem to construct a local ordering function (local utility-representation), $f(x,y)$, that satisfies $x\leq_{U_x} y\iff f(x,y)\geq 0$, where $U_x$ is some open neighbourhood of $x$.  Following from the previous discussion, we see that the only difficulty is defining the open sets $U_x$ for elements of $\mathcal{M}_c$; suppose $x\in\partial U$ for some $U\in O_X$ is not contained in any other element of $O_X$, then this can break down in two distinct situations: the closure (with sequences) of the order, $\leq_U$ on $\text{cl}(U)$, forms a cyclic ordering; there is another $(V,\leq_V)\in O_X$ such that $x\in\partial V$ and the closure of the ordering $\leq_U$ does not agree with the closure of the ordering $\leq_V$.  In the latter case we can either define the set $U_x=(U\cup V)\cup(\partial U\cap \partial V)$ and allow for $f_x(y)\coloneqq f(x,y)$ to be discontinuous or we must assume the stronger version of assumption (4) that being the transverse closure by sequences of any $(U,\leq_U),(V,\leq_V)\in O_X$ must agree on intersection --- we could call this a transversly closed local ordering.

Resolving the cyclic order is more technical: for this to occur $x$ must be both a maximal and minimal element of $U$ (since $\text{cl}(U)$ is compact), thus we must take an equal set, $E$, in $U$ that is strictly less (equivalently strictly greater) than $x$, letting $E_x=[x]\cap U$ denote the equal set of $x$ in the closure of $U$, we can and define $U_x=(U\cup E_x)\setminus E$, this will be an open subset of $X$ that is linearly orderable given $O_X$.  The difficulty comes from choosing which $E$ to cut out and is effectively equivalent to choosing which point on a circle to cut out to make $S^1\setminus\{p\}\cong (0,1)$, to do this unambiguously you would need, at least, that $X$ be metrisable so that you can take $E$ such that it is in the \textit{midpoint} of $U$.  For all other points of closure we can construct $U_x$ as follows: if $x\in\partial X$, then simply taking $U_x=\text{cl}(U)$ is satisfactory (since this will be open in the topology of $X$); if $x$ is a transverse or extremal point of closure, we could take either $U_x=(U\cup V)\cup(\partial U\cap \partial V)$, or $U_x=\text{int}(U^\uparrow\cup V^\downarrow)$ where $U^\uparrow$ is the upward closure (closure in terms of increasing sequences) and $V^\downarrow$ is the downward closure (closure in terms of decreasing sequences), respectively.  

Therefore, if $X$ is Hausdorff, and $O_X$ is locally complete and transversely closed, then for every $x\in X$ we can construct an open set $U_x$ containing $x$ such that $O_X$ implies a continuous LQO on $U_x$.  By statement (ii) of proposition \ref{appBprop:utility_representations} there will exist a continuous $f_x$ for every such neighbourhood; thus on the domain $D=\bigsqcup_{x\in X} U_x$, we can construct a local ordering function as $f:D\to\mathbb{R}:(x,y)\mapsto f_x(y)$ and which satisfies $x\leq_{U_x} y\iff f(x,y)\geq 0$.  We note that the $D$ that we construct is not necessarily the \textit{best} (i.e. largest) domain possible; to construct a maximal domain for $f$ we would need each $U_x$ to be the maximal extension over $O_X$ containing $x$; since we do not formally introduce this idea in this paper, we will end that discussion here.  In the above discussion, we have proven the following proposition:
\begin{proposition}\label{appBprop:local_ordering_existence}
    Let $X$ be a Hausdorff topological space, and let $O_X$ be a locally complete and transversely closed continuous local linear quasi-order.  Then there always exists a local utility-representation, $f:\text{Dom}(f)\to\mathbb{R}$ of $O_X$, where $\text{Dom}(f)=\bigsqcup_{x\in X}U_x$ for open neighbourhoods $U_x$ of $x$ satisfying for any $x\in X$ and $y\in U_x$, $x\leq_{U_x} y\iff f(x,y)\geq 0$.  If, in addition, $X$ is metrisable (e.g. locally compact, second-countable, and $\sigma$-compact) and the maximal extensions of $O_X$ are not arbitrarily small, then each $U_x$ can be chosen to be not arbitrarily small.
\end{proposition}
\begin{proof}
    The first statement has been shown by construction in the previous discussion.  If $X$ is metrisable then it is homeomorphic to a metric space $(M,d)$, let us denote this homeomorphism as $\phi:(X,\tau)\to(M,d)$.  Therefore, \text{not arbitrarily small} means there is a uniform lower bound on the diameter of the sets $\phi(U_x)\subset M$ with respect to $d$.  As mentioned, we can construct $U_x$ to be a maximal extension over $O_X$ containing $x$, and by assumption of the proposition, the maximal extensions will not be arbitrarily small, so we are done.
\end{proof}

A local ordering function allows us to unify Richter-Peleg utility functions for continuous LQOs, and multi-utility representations for non-linear (incomplete) quasi-orders:
\begin{align}
    \text{Utility Function:}&\qquad f(x,y)=f(y)-f(x)\\
    \text{Multi-Utility Representations:}&\qquad f(x,y)=\inf_{i\in I}\left[f_i(y)-f_i(x)\right]
\end{align}
where $f:X\to\mathbb{R}$ is a Richter-Peleg utility representation, and $\{f_i\}_{i\in I}$ is a multi-utility representation.  Trivially in both cases $f(x,x)=0$ and by the definition of the utility representations, we will have $x\leq y\iff f(x,y)\geq 0$.  Furthermore, can represent the more complex notion of continuous LLQOs and even local pospaces naturally within the local ordering function.

If we restrict ourselves to smooth d-dimensional Riemannian manifolds, which we denote by $\mathcal{M}$, then we can frame a local ordering function in terms of vector fields again.  We note that it is true that any arbitrary smooth vector field $b:\mathcal{M}\to T\mathcal{M}$ can be written as $b(x)=\nabla l_x(x)$ for some smooth scalar function $l_x:\mathcal{M}\to\mathbb{R}$.  We can see this by simply writing $l(x,y)=g_x(\exp_x^{-1}(y),b(x))$ where $\exp_x:\tilde U\to U$ where $U$ is an open neighbourhood of $x$, and $\tilde U$ is an open ball of radius $\text{inj}(x)$ centred at the origin of $\mathbb{R}^d$; the exponential map defines the normal coordinates about $x$ and in these normal coordinates $g_x=I_d$ and therefore $l(x,y)=y_ib_i(x)$, it is then trivial true that $\nabla l_x(x)=b(x)$.  The assumption of bounded geometry simply allows us to bound the injectivity radius away from zero, i.e. there exists some $c>0$ such that $\inf_{x\in\mathcal{M}}\text{inj}(x)=c$; thus, $\tilde l:\mathcal{M}\times B(0,c)\to\mathbb{R}:(x,u)\to l(x,\exp_x(u))$ is a well-defined smooth function.  The condition $x\leq y\iff l(x,y)\geq 0$ is equivalent to $y-x$ being in the same half-space as $b(x)$, which is simply a measure of directionality - exactly as we desire for our ordering function to imply.  Furthermore, if $\omega:\mathcal{M}\to T^*\mathcal{M}$ was the dual covector field of $b:\mathcal{M}\to T\mathcal{M}$, then by Frobenius' Theorem, we know that the kernel of $\omega$ will generate a foliation if and only if $\omega\wedge d\omega=0$.  This equivalence allows us to construct an ordering function from any codimension-one transversally orientable\footnote{A codimension-one foliation is said to be transversally orientable if the normal vector to the leaves can be chosen to have a globally consistent unit normal vector field.} foliation as opposed to a foliation defined from an ordering.

Assume that $O_\mathcal{M}$ is locally complete, and let $l_U$ be an ordering function for some $U\in O_\mathcal{M}$, then we will have an atlas of codimension-one foliated charts that cover $\mathcal{M}_r$ and whose foliations agree on their intersections, thus $\mathcal{M}_r$ admits a global codimension-one foliation, and this is generated by a nowhere vanishing 1-form $\omega$ satisfying $\omega_{U}=h_Udf_U$ --- this follows from the fact that transversal orientability being immediate from the definition of an LLQOTS.  The scalar functions $h_U$ are constructed such that $h_U/h_V=D\varphi_U^V(f_U)$ where $\varphi_U^V:l_U(U\cap V)\to l_V(U\cap V)$ are the continuous strictly increasing functions that we us to represent the agreement of the order on a LLQOTS,\footnote{This condition relates back to our discussion about the continuity of a local ordering function in the first argument given in section \ref{subsec:nongradflow_extension}.} with this condition arising from:
\begin{equation}
    \omega_{U\cap V}=h_Udf_U=h_Vdf_V=h_VD\varphi_U^V(f_U)df_U\implies\frac{h_U}{h_V}=D\varphi_U^V(f_U)>0
\end{equation}
We can construct each $h_U$ as follows: take a partition of unity $\{\rho_U\}_{U\in O_{\mathcal{M}}}$ on a subcover $\mathcal{V}$ of $O_\mathcal{M}$, for every $U\in O_\mathcal{M}$ define $g_U=\sum_{V\in \mathcal{V}}\rho_V\lambda_U^V$, where $\lambda_U^V=\log D\varphi_U^V(l_U)$.  Since $\varphi_U^W=\varphi_V^W\circ\varphi_U^V$ we will have $\lambda_U^W=\lambda_U^V + \lambda_V^W$ and thus $g_U-g_V=\sum\rho_W(\lambda_U^W-\lambda_V^W)=\lambda_x^y$ and thus $h_U\coloneqq\exp(g_U)$ we will have $h_U/h_V=\varphi_U^V$.  Therefore, given a locally complete LLQOTS, the associated local ordering function will be equivalent to $l(x,y)=\omega(x)(\exp_x^{-1}(y))$, thus connecting the structures of vector fields, foliations and locally linear quasi-orders, and representing them all with a local ordering function.

To summarise, if $\mathcal{M}$ is a smooth $d$-dimensional Riemannian manifold, and $O_\mathcal{M}$ is an order set that is locally complete and transversely closed, then the regular set will admit a codimension-one foliation with leaves being the equal contours of $O_\mathcal{M}$.  Furthermore, $O_\mathcal{M}$ being locally complete implies $\mathcal{M}_c\subseteq \partial\mathcal{M}_r$ which can be at most a $(d-1)$-dimensional submanifold which means $\text{vol}_g(\mathcal{M}_c)=0$, i.e. almost everything is foliated.  By proposition \ref{appBprop:local_ordering_existence} there will always exist a local ordering function, $l(x,y)$, and by the previous discussion there will exist a smooth 1-form, $\omega:\mathcal{M}\to T^*\mathcal{M}$ such that $l(x,y)=\omega(x)(\exp_x^{-1}(y))$ that satisfies $\omega\neq 0$ and $\omega\wedge d\omega =0$ on $\mathcal{M}_r$.  Although we resisted providing a definition of an ordering function, if we were forced to give one, it would be:
\begin{definition}[Ordering Functions]
    Let $\mathcal{M}$ be a $d$-dimensional Riemannian manifold with bounded geometry, and let $\omega:\mathcal{M}\to T^*\mathcal{M}$ be a $C^2$ 1-form satisfying $\omega\wedge d\omega=0$ and $\text{supp}(\omega)\coloneqq\text{cl}(\{x\in\mathcal{M}:\omega(x)\neq 0\})=\mathcal{M}$.  Then we call the function $l:\text{Dom}(l)\to\mathbb{R}:(x,y)\mapsto \omega(x)(\exp_x^{-1}(y))$ a local ordering function, where $\text{Dom}(l)=\bigsqcup_{x\in \mathcal{M}}U_x$ is such that for any $x\in\mathcal{M}$ the exponential map, $\exp_x^{-1}$, is a diffeomorphism on $U_x$.  If, in addition, there exists a scalar function $\tilde l\in C^2(\mathcal{M})$ such that $l(x,y)=\tilde l(y)-\tilde l(x)$ (equivalently $\omega=d\tilde l$) then we call $\tilde l(x)$ a global ordering function.
\end{definition}
From this definition, we know that a local ordering function defines transversally orientable codimension-one foliations almost everywhere, i.e. the leaves of the foliation have a well-defined normal vector.\footnote{We know this because the foliation is represented by a nowhere-vanishing 1-form.}  Our restriction of the domain ensures that $l(x,y)$ is well-defined, and the assumption of bounded geometry ensures this domain is not arbitrarily small, and trivially we have that $l(x,x)=\omega(x)(0)=0$ for every $x\in\mathcal{M}$.  We also make the important practical note that the restriction of $\mathcal{M}$ to be such that it admits a transversely orientable codimension-one foliation is not at all limiting since we are assuming the manifold is generated by some stochastic process, i.e. there really will be a well-defined normal vector for every leaf.

\subsection{Estimating $\boldsymbol{\lVert\nabla l(x)\rVert}$ on Discrete Samples}\label{appBsubsec:gradest}
In our discussion of the form of the order regularisation function $\kappa(x)$, we proposed setting $\kappa(x)=\lVert\nabla l(x)\rVert$ to make ODK agnostic to the scale of $l(x)$.  In practice, it is unlikely that you will have an analytic expression for $l(x)$, so one must compute $\lVert\nabla l(x)\rVert$ numerically from evaluations of $l(x)$ on the data $\mathcal{D}$.  There are many methods that one could use towards this end, but we find that one in particular feels very canonical given our assumption that the data clusters on or near an $d$-dimensional Riemannian submanifold $\mathcal{M}$ of $\mathbb{R}^n$.  Let us first set up the problem: the gradient of a function $f:\mathbb{R}^n\to\mathbb{R}$ at a point $x\in\mathbb{R}^n$, denoted $\nabla f(x)$, is the unique vector satisfying:
\begin{equation}\label{appBeq:localpca1}
    \lim_{y\to x}\frac{f(y)-f(x)-\nabla f(x)^T(y-x)}{\lVert y-x\rVert}=0
\end{equation}
In other words, $\nabla f(x)$ defines the best possible linear approximation of $f$ around $x$.  For some $y\in B(x,\eta(x))$, where $\eta(x)$ is a sufficiently small radius such that this linear approximation is reasonable, we can ignore the limit in equation \eqref{appBeq:localpca1} and rearrange to find:
\begin{equation}\label{appBeq:localpca2}
    f(y)=f(x)+\nabla f(x)^T(y-x)
\end{equation}
Let $\mathcal{D}=\{x_i\}_{i=1}^m\subset\mathbb{R}^n$ be a finite sample of points on $\mathcal{M}$, for every $i=1,\dots,m$, let $N_i=(B(x_i,\eta_i)\cap\mathcal{D})\setminus\{x_i\}$, where again $\eta_i$ is some sufficiently small radius.  Then we can see that equation \eqref{appBeq:localpca2} defines a linear regression problem; more precisely, it defines the problem $\delta_i f=X_i^T\beta_i$, where $\beta_i$ denotes our approximation to $\nabla f(x_i)$, and $(\delta_i f)_j=f(x_j)-f(x_i)$ and $(X_i)_{:,j}=x_j-x_i$, for $x_j\in N_i$ --- we are assuming for convenience that these samples are indexed from $j=1,\dots,|N_i|$ in $N_i$.  If our goal is to minimise the squared residuals, i.e. we wish to minimise $R^2(\beta_i)=\lVert\delta_i f-X_i^T\beta_i\rVert^2$, then this is a least squares problem and if $|N_i|\gg 1$ will have a unique solution given by:
\begin{equation}
    \beta_i=(X_iX_i^T)^{-1}X_i\delta_if
\end{equation}
However, if $|N_i|< n$ --- which would often be the case if we are choosing a small neighbourhood $N_i$ of $x_i$ and $n\gg 1$ --- then the system will be underdetermined and there will be infinitely many solutions that minimise the squared residuals.  To break this tie, the usual method is to take the solution that minimises $R^2(\beta_i)+\lVert \beta_i\rVert^2$, which will select the solution of smallest magnitude that minimises $R^2(\beta_i)$.  This can also be expressed exactly through \cite{lopez2018moore}:
\begin{equation}
    \beta_i=X_i(X_i^TX_i)^{-1}\delta_i f
\end{equation}
It turns out that both of these expressions are equivalent if written in terms of components of the Singular Value Decomposition (SVD) of $X_i^T=U_iS_iV_i^T$, where $U_i\in\mathbb{R}^{|N_i|\times |N_i|}$ and $V_i\in\mathbb{R}^{n\times n}$ are unitary matrices (i.e. $V^{-1}=V^T$), and $S_i\in\mathbb{R}^{|N_i|\times n}$ with positive values along the diagonal.  In this case, our solution can be expressed as:
\begin{equation}\label{appBeq:localpca3}
    \beta_i=V_iS_i^{-1}U_i^T\delta_i f
\end{equation}
This discussion justifies solving for the gradient by setting up a regression problem if the function $f$ is defined in $\mathbb{R}^n$, but in our setup, the ordering function, $l:\mathcal{M}\to\mathbb{R}$ is imagined to be be defined on a $d$-dimensional Riemannian manifold $\mathcal{M}$ embedded in $\mathbb{R}^n$, with $n\gg d$, and thus we are including far too many dimensions in equation \eqref{appBeq:localpca3} which can only be adding errors since the $l(x)$ does not vary in these directions.  Furthermore, if $n\gg 1$, the Euclidean metric itself begins to become less informative as the concentration of measure phenomenon starts to take effect \cite{ledoux2001concentration}.  Given this consideration, the most principled next step would be to try and restrict our regression problem to the tangent space $T_x\mathcal{M}$ in which $\nabla l(x)$ actually lives; furthermore, we know that $T\mathcal{M}$ is a rank-$d$ subbundle of $T\mathbb{R}^n$ and thus at at $x_i\in\mathcal{D}$ we need to restrict to a $d$-dimensional affine space centred about $x_i$.  This is exactly the sort of problem that Principal Component Analysis (PCA) is designed to solve, and thus, in what follows, we will apply PCA locally on every difference matrix $X_i$ for $i=1,\dots,m$, to construct an approximation of the tangent bundle on $\mathcal{M}$.  This method is referred to as local-PCA (for obvious reasons) and is familiar to many people working in the field of manifold learning, as it is one of the simplest constructive methods to access the tangent spaces of a manifold \cite{singer2012vector}.

Before proceeding, we remark that our restriction to the neighbourhood $N_i$ needs to include enough points to resolve the function locally, but not so many points that the linear approximation of $\beta_i$ to $\delta_il$ breaks down, and also must not be so large that the curvature of the manifold invalidates using a linear approximation to the tangent space.  The determination of these neighbourhoods is either done by setting $\eta_i\equiv\eta$ to be constant, as is done by Singer \& Wu in \cite{singer2012vector}, and which allows for some cleaner theoretical results; or by setting $\eta_i=\lVert x_i-x_i^k\rVert$ where $x_i^k$ is the $k^\text{th}$ nearest neighbour to $x_i$ not including $x_i$, which uniformises our neighbourhoods (i.e. $|N_i|\equiv k$) and thus gives each estimate the same chance to fail or succeed.  If the data is very inhomogeneously sampled, then choosing the latter parameterisation is essentially necessary, and the burden will then fall on how large you allow your neighbourhood to be; we have found that for manifolds of relatively low dimensions $k\in[10,30]$ often produces good results.  We also note that it is standard to centre $X_i$ before applying PCA, but since we are only interested in the projection, we skip this step.

It is important to recall that we do not a priori assume we know the value of the intrinsic dimensionality $d$.  We could now employ one of the many intrinsic dimensionality estimation methods in the literature (see, for example, \cite{hein2005intrinsic, granata2016accurate, hino2017ider, facco2017estimating, allegra2020data, di2026scale}), but instead we will keep within the domains of PCA and define a local dimensionality $d_i$ at every $x_i\in\mathcal{D}$ which will equal the minimum number of dimensions required to explain a preset proportion of the variance in our neighbourhood $N_i$.  This is equivalent to choosing some $\gamma\in (0,1)$ and defining:
\begin{equation}
    d_i=\min \left\{k\in\{1,\dots,|N_i|:\frac{\sum_{j=1}^k\sigma_j^2}{\sum_{j=1}^{|N_i|}\sigma_j^2}>\gamma\right\}
\end{equation}
where $\{\sigma_j\}_{j=1}^{|N_i|}$ are the singular values of $X_i$, i.e. the diagonal elements of $S_i$.\footnote{We usually set $\gamma=0.9$, as was done in \cite{singer2012vector}, and found it suitable in almost all cases.}  It would be most rigorous to now define:
\begin{equation}
    \hat d=\text{median}\{d_i\}_{i=1}^m
\end{equation}
so that we have a uniform estimate of the intrinsic dimensionality across the manifold; however, in practice we do not see much benefit to our regression problem from making such a restriction, so we often approximate $T_{x_i}\mathcal{M}$ with respect to each $d_i$.  Now we are ready to apply PCA on each difference matrix $X_i$ to form the rank-$d_i$ projection $\tilde X_i\in\mathbb{R}^{d_i\times |N_i|}$,  whose span determines our approximation to the tangent space.  More precisely, we compute the SVD of $X_i^T=U_iS_iV_i^T$, and then construct $\tilde{X}_i=\tilde V_i^TX_i=\tilde S_iU_i^T$ as the projection onto the first $d_i$ principal components, where $\tilde{V}_i=(V_i)_{:,1:d_i}$, $\tilde{S}_i=(S_i)_{1:d_i,:}$.\footnote{Recall that the principal components are matrices are the columns of $Q$ in the eigendecomposition $XX^T=Q\Lambda Q^{-1}$, where symmetry implies $Q^{-1}=Q^T$; and if $X^T=USV^T$ then $XX^T=VS^2V^T$, thus $Q=V$.}  Finally, from equation \eqref{appBeq:localpca3} we define our estimate to be:
\begin{equation}
    \hat{\nabla} l(x_i)\coloneqq \tilde{S}_i^{-1}U_i^T\delta_i l
\end{equation}
Define $G_i\coloneqq \tilde{S}_i^{-1}\tilde U_i^T$, then $\lVert\nabla l(x_i)\rVert^2\approx \delta_i l^TG_i^TG_i\delta_i l$.  It is exactly these $G_i$ that we denoted as $G_i^\text{reg}$ in the gradient regularisation of section \ref{subsec:numerical_results_sde}, equation \eqref{eq:sde_gradreg}.

There are other gradient estimation methods that we could apply --- most notably, a finite-difference scheme tailored to non-uniform grids --- but we believe that local-PCA is the most principled approach given our assumptions about the data.  Furthermore, while testing different schemes, we found that local-PCA was the most reliable on inhomogeneously sampled data.

\subsection{Heuristic Parameterisations of ODK}\label{appBsubsec:heuristic_parametrisations}

In this appendix, we will briefly discuss heuristic parametrisations that one might use for ODKs; we should view these as data-adaptive parameterisations that are not related to the underlying dynamics that generated the data but only seek to fit ODKs effectively to the data.  It is worth noting that all of these heuristics will have tangible drawbacks, and some may not work at all on certain datasets; it is important to tune --- or design new --- heuristics based on the realities of the data you are studying.  In what follows, $\mathcal{D}=\{x_i\}_{i=1}^m\subset\mathbb{R}^n$ will be our data, which we assume clusters around a $d$-dimensional submanifold, $\mathcal{M}$, of $\mathbb{R}^n$.

\subsubsection{Geodesic Construction of an Ordering Function}\label{appBsubsec:heuristic_orderingfunction}
We will begin by discussing the geodesic construction of the ordering function that we have mentioned several times in the main text.  First, we construct a symmetric (undirected) neighbourhood graph $\mathcal{G}=(\mathcal{V},\mathcal{E})$ where a fixed indexing of our data forms the vertex set, $\mathcal{V}$; and $\mathcal{E}$ is our edge set, i.e. $i\to j\iff (i,j)\in\mathcal{E}$, and where $i\to j$ reads \textit{i connects to j}; and we note that the property of being symmetric is equivalent to $\mathcal{G}$ being undirected.  It is not particularly important which graph construction method you use to construct your graph so long as it satisfies the following three properties:
\begin{enumerate}[nolistsep,noitemsep]
    \item[1).] $\mathcal{G}$ is symmetric: $\forall i,j\in\mathcal{V}$ such that $(i,j)\in \mathcal{E}$, then we must also have $(j,i)\in\mathcal{E}$;
    \item[2).] $\mathcal{G}$ is connected: $\forall i,j\in\mathcal{V}$ there exist vertices $i_0,i_1,\dots,i_k\in\mathcal{V}$ such that $i_0=i$ and $i_k=j$ and $(i_l,i_{l+1})\in\mathcal{E}$ for all $l=0,1,\dots,k-1$;
    \item[3).] $\mathcal{G}$ is a neighbourhood graph: there exists some $r>0$, that is sufficiently small with respect to the maximum length scale of $\mathcal{D}$ (i.e. $r\ll\max_{i,j\in \mathcal{V}}\lVert x_i-x_j\rVert$), such that $\forall i\in\mathcal{V}$, we have $\max_{j\in N_i}d(x_i,x_j)<r$, where $N_i=\{j\in\mathcal{V}:(i,j)\in\mathcal{E}\}$ are the descendants of $i$, and $d:\mathcal{D}\times\mathcal{D}\to\mathbb{R}$ is some distance metric.
\end{enumerate}
As a disclaimer, we have invented the definition of a neighbourhood graph ourselves purely for convenience; but we hope that it is clear that we are simply requiring that $\mathcal{G}$ encodes neighbourhood structure that is consistent with the neighbourhood structure of the underlying manifold; and thus doesn't allow long-distance connections.  A convenient method to represent such graphs is to use a distance threshold function $f(i,j)$ that defines the edge set through $(i,j)\in\mathcal{E}\iff d(x_i,x_j)<f(i,j)$.  The most common graph is a k-Nearest Neighbours (kNN) graph where $f(i,j)=\lVert x_i-x_i^k\rVert$, where $x_i^k$ denotes the $k^\text{th}$ nearest neighbour to $x_i$ (not including $x_i$ itself), but we notice that $f(i,j)\neq f(j,i)$ so this graph is not necessarily connected; we could then consider the following two variations: the kNN-AND graph $f(i,j)=\min(\lVert x_i-x_i^k\rVert,\lVert x_j-x_j^k\rVert)$, or the kNN-OR graph $f(i,j)=\max(\lVert x_i-x_i^k\rVert,\lVert x_j-x_j^k\rVert)$.  Our personal preference is to use the Continuous-kNN (CkNN) graph, introduced by Berry \& Sauer in \cite{berry2016consistent}, which sets $f(i,j)=\delta\sqrt{\lVert x_i-x_i^k\rVert\cdot\lVert x_j-x_j^k\rVert}$ for some \textit{sensitivity} $\delta>0$; the CkNN graph has much more stable geometric properties as the dataset size increases, and on inhomogeneously sampled data it is easier to construct a connected CkNN graph for a modest value of $k$.

After one has constructed a graph, $\mathcal{G}$, and given that we know a distribution, $p_0(x)$, on the initial states, then we can construct the ordering function as:
\begin{equation}\label{appBeq:geodesic_order}
    l(x_j)=\sum_{i\in\text{supp}(p_0)}\bar p_0(x_i)\min_{P\in \mathcal{P}_{ij}}\sum_{k=0}^{|P|-1}\lVert x_{i_k}-x_{i_{k+1}}\rVert
\end{equation}
where $\bar p_0(x_i)=p_0(x_i)/\sum_{j\in\mathcal{V}}p_0(x_j)$ and $\mathcal{P}_{ij}=\bigcup_{k=0}^m\{(i_0,i_1,\dots,i_k)\in\mathcal{V}^{k+1}:i_0=i\land i_k=j\}$ is the set of all possible paths originating at $i$ and terminating at $j$.  If we $\bar p_0(x)=e_{i_0}$ for some $i_0\in\mathcal{V}$ (i.e. $p_0(x)=\delta(x-x_{i_0})$), then the outer sum would disappear and equation \eqref{appBeq:geodesic_order} would be the standard shortest path length formula on a graph, denoted as $d_\mathcal{G}(i_0,j)$; otherwise, we are computing the average of $d_{\mathcal{G}}(i,j)$ over $i$.  For symmetric and connected graphs, the shortest path will always exist and be well defined, and if our discrete graph is a neighbourhood graph, the shortest path will asymptotically approximate the geodesic distance on $\mathcal{M}$.  In practice, if $|\text{supp}(p_0)|\ll m$, you should use Dijkstra's Algorithm \cite{dijkstra2022note} to compute the shortest path; otherwise, or if at some future point you may need a matrix of pairwise geodesic distances, then you should use the Floyd-Warshall algorithm \cite{floyd1962algorithm, warshall1962theorem}.  

The geodesic construction of the ordering function should yield $\lVert\nabla l(x)\rVert\equiv 1$, so in theory there is no need to estimate the scale of the ordering.  However, in practice, inhomogeneities in in the sampling can case the effective local scales to vary and thus one still might wish to compute $\kappa(x)=\lVert\nabla l(x)\rVert$; in this case, and if one has used a different construction of the ordering, we refer the reader to Appendix \ref{appBsubsec:gradest} for a detailed discussion of how we may approximate this quantity.

\subsubsection{Fundamental Constraints on ODKs Parameters} 
We will begin by discussing $\varepsilon$, which in some sense is the most difficult parameter to choose.  We will not review all relevant literature in this field, only what we consider important to mention, and we refer interested readers to the following works \cite{rosenblatt1956, parzen1962estimation, silverman2018density, scott2015multivariate}.  In the classical literature, a bandwidth is denoted as $h$, but to avoid confusion with our kernel function, let us denote it as $\bar h$.  For a symmetric Diffusion Maps kernel of the form $h\left(\tfrac{\lVert y-x\rVert^2}{\varepsilon}\right)$, $\bar h=\varepsilon$ must act simultaneously as a fixed bandwidth (fitting the kernel to the data) and as the asymptotic variable which we wish to take to zero.  Let $\bar h_\text{min}(i)$ be the smallest value of the bandwidth such that the kernel $K_\varepsilon(x_i,\cdot)$ would not be degenerate (we will clarify this later), then for a Diffusion Maps kernel, we would choose $\varepsilon^*=\max_{i=1,\dots,m}\bar h_\text{min}(i)$ which would be the smallest epsilon which guarantees that the kernel is non-degenerate everywhere.  If we generalise to the left formulation of a variable bandwidth, i.e. kernels of the form $h\left(\tfrac{\lVert y-x\rVert^2}{\varepsilon\rho(x)^2}\right)$, then the same discussion applies and we will have $\varepsilon^*\rho^*(x_i)^2=\bar h_\text{min}(i)$, so we need to compute some heuristic length scale $L_i$, for every $x_i\in\mathcal{D}$, such that $\exists c>0$ such that $L_i\approx c^{-1/2}\bar h_\text{min}(i)^{1/2}$, then setting $\rho^*(x_i)=L_i$ implies $\varepsilon^*=c$ is our optimal parametrisation, but since our heuristic estimate of $L_i$ will not satisfy the proportionality relation exactly, in reality we would still need $\varepsilon^*=\max_{i=1,\dots,m}\bar h_\text{min}(x_i)/L_i^2$.  When we make the connection to the underlying stochastic process, as we do often in the main text, we add an additional constraint for $\varepsilon$, in that it must also be small enough such that, relative to the underlying dynamics, the short time-step approximation of a Gaussian transition operator is valid; this constraint is then equivalent to choosing a small increment $dt$ such that $\int_0^Tp_t(x)\mathrm{d}t\approx\sum_{s=1}^Sp_{s\cdot dt}(x)dt$, where each $p_t(x)$ is a time-marginal of the underlying process, and $S\approx T/dt$.  The error in this approximation is $\bigO(dt T\sup_{t\in[0,T]}\partial_t p_t)$; let $E(x)$ denote the error of this approximation, substitution $\varepsilon$ for $dt$, we know that:
\begin{equation}\label{appBeq:integral_error}
    E(x)=\left|\int_0^Tp_t(x)\mathrm{d}t-\sum_{s=1}^Sp_{s\cdot \varepsilon}(x)\varepsilon\right|\leq \frac{\varepsilon T}{2}\sup_{t\in[0,T]}|\partial_tp_t|\implies \lVert E\rVert_\infty\leq\frac{\varepsilon T}{2}\sup_{t\in[0,T]}\lVert\mathcal{L}^*p_t\rVert_\infty
\end{equation}
where in the first inequality you use the mean value theorem and the fact that the integral of a linear function is the area of a triangle, and the inequality for the supremum norm of the error is immediate after commuting the suprema in time and space.  The supremum on the right-hand side of the finally inequality in equation \eqref{appBeq:integral_error} depends on the underlying dynamics, and thus if the dynamics are highly irregular, our approximation would only be valid for very small $\varepsilon$; but if the true dynamics are highly irregular and you have sampled so few points that $\varepsilon$ cannot be made sufficiently small, then solving the inverse-problem of inferring the dynamics is ill-posed in general, thus it isn't surprising that a proper parameterisation exists.

The integer $S$ is the number of segments in our uniformly partitioned Riemann sum approximation of the integral; the minimum number of these segments required to approximate the integral up to a fixed resolution will be a function of the underlying dynamics, denote this minimum value as $s(\mu,D;\delta_T)$ where $\delta_T$ is the fixed resolution we require for the time resolved integrals.  Furthermore, we could also consider the number of samples that we need to draw from $q$ such that we can approximate the integral:
\begin{equation}
    \int_\mathcal{M}K_\varepsilon(x_i,y)\mathrm{d}y\approx \frac{1}{|N_i|}\sum_{j\in N_i}K_\varepsilon(x_i,x_j)
\end{equation}
for every $x_i\in\mathcal{M}$ given a value of $\varepsilon$, denote this as $m(\varepsilon,q;\delta_\mathcal{M})$ where $\delta_\mathcal{M}$ is the fixed resolution we require for the spatially resolved integrals.  Since we assume $\varepsilon$ is linear in the denominator of $h(\cdot)$, we can show that $m(\varepsilon,q;\delta_\mathcal{M})\approx\varepsilon^{-d/2}m(q;\delta_\mathcal{M})$ for some function $m(q)$ whose value depends on the inhomogeneity of $q(x)$.\footnote{The dependence on the dimension $d$ in this formula for $m(\varepsilon,q)$ is a realisation of the so-called \textit{curse-of-dimensionality}, whereby to resolve our statistics we must sample exponentially more points as the dimensionality increases.}  Putting this together:
\begin{equation}
    S\geq s(\mu,D;\delta_T)\implies \varepsilon\leq \frac{T}{s(\mu,D;\delta_T)},\qquad\text{and}\qquad m\geq m(\varepsilon,q;\delta_\mathcal{M})\implies\varepsilon\geq \left(\frac{m(q;\delta_\mathcal{M})}{m}\right)^{2/d}
\end{equation}
The right-hand side of the lower bound on $\varepsilon$ is effectively our $\bar h_\text{min}$, so for ODKs, or more generally any local kernel, that seeks to approximate the underlying dynamics of a dynamical system, a given sample of $m$ samples from a distribution $q(x)$ is insufficient for this task if:
\begin{equation}
    \sqrt{\frac{T}{s(\mu,D;\delta_T)}}<\left(\frac{m(q;\delta_\mathcal{M})}{m}\right)^{1/d}
\end{equation}
Which is an additional constraint that doesn't need to be satisfied for kernels that only seek to represent the geometry of $\mathcal{M}$.

\subsubsection{Heuristics for $\boldsymbol{\varepsilon}$}
Now that we have made this consideration, we will begin providing heuristic methods that one can use to determine a good parametrisation of $(\varepsilon,\tau(x),\rho(x)^2)$.  We will proceed by first constructing a heuristic that quantifies the degree of degeneration that occurs as the bandwidth decreases; for convenience, we will assume that $\rho(x)$ has already been set.  We define the maximum transition probability from a sample $x_i$ as:
\begin{equation}\label{appBeq:pmax}
    p_\text{max}^i(\varepsilon)=\max_{x_j\in\mathcal{D}}P_\varepsilon(x_i,x_j)
\end{equation}
where $P_\varepsilon$ is to be viewed as our row-stochastic transition matrix, not the continuous operator; we also define $p_\text{max}(\varepsilon)\coloneqq\max_{i=1,\dots,m}p_\text{max}^i(\varepsilon)$.  For diffusion maps kernels with monotonically decreasing choices of $h(\cdot)$, we will have $p_\text{max}^i(\varepsilon)=P_\varepsilon(x_i,x_i)$, which is the self-transition probability, or, alternatively, the probability of not doing anything.  Since $\sum_{j=1}^mP_\varepsilon(x_i,x_j)=1$, for any fixed set of parameters we will have $p_\text{max}(\varepsilon)\to 0$ as $m\to\infty$, so it may seem as if this metric doesn't tell us much, but for a fixed dataset $p_\text{max}(\varepsilon)\to 1$ as $\varepsilon\to 0$.  So if we set some threshold value $\eta_p\in(0,1)$, for example $\eta_p=0.9$, then for a given value of $\varepsilon$ $p_\text{max}(\varepsilon)>\eta_p$ implies there is some $x_i\in\mathcal{D}$ where a diffusion kernel only diffuses $10\%$ of the time --- which seems like a poor representation of diffusion on $\mathcal{D}$.  Thus, we can use $p_\text{max}$ to determine whether our bandwidth is failing to resolve somewhere in $\mathcal{D}$; setting a more reasonable threshold, such as $\eta_p=0.5$, can allow us to define a lower bound on our bandwidth based on a heuristic notion of what we feel determines a good fit for us; i.e. we can define $\varepsilon_\text{min}=\inf\{\varepsilon>0:p_\text{max}(\varepsilon)<\eta_p\}$.\footnote{In practice, there is some flexibility in how harshly we apply this threshold, it is likely true that not all samples are equally important to us, so we could relax this to $Q_{0.9}(\{p_\text{max}(i):i=1,\dots,m\})>\eta_p$ where $Q_{\theta}$ is the $\theta^\text{th}$ quantile of a set (assuming a uniform distribution over the values).}

Supposing we have restricted the transitions of our kernel (ODK or otherwise) to a symmetric, connected, neighbourhood graph\footnote{Under certain circumstance it might be reasonable to relax the symmetry constraint.} of the data, denoted $\mathcal{G}_\mathcal{D}=(\mathcal{V},\mathcal{E}_\mathcal{D})$, and noting that this should not be the same graph you use to compute the geodesic order of equation \eqref{appBeq:geodesic_order} and should be significantly larger.  We could think of $\mathcal{G}_\mathcal{D}$ as only being computed for the practical purposes of reducing the computation load of the kernels and avoiding floating-point issues when we have transition probabilities being approximately zero (long-distance transitions).  Then we can construct a heuristic metric for a maximum bandwidth (and indeed a maximum value of $\tau(x)$) by computing the resolution as:
\begin{equation}
    \text{res}_i(\varepsilon)=\sum_{j:(i,j)\in\mathcal{E}_\mathcal{D}}P_\varepsilon^\mathcal{D}(x_i,x_j)
\end{equation}
where $P_\varepsilon^\mathcal{D}$ has not been restricted to the graph $\mathcal{G}_\mathcal{D}$ and is applied to the whole dataset; as before we can also define $\text{res}(\varepsilon)=\max_{i=1,\dots,m}\text{res}_i(\varepsilon)$.  This computes how much detail we are losing by restricting transitions to $\mathcal{G}_\mathcal{D}$.  Thus, if we are happy with the validity of our constructed graph, then we can set some threshold $\eta_r\in(0,1)$ and if $\text{res}(\varepsilon)>\eta_r$ we know that the bandwidth is too large to be represented on our graph in some region of the data.  This allows us to define $\varepsilon_\text{max}=\sup\{\varepsilon>0:\text{res}(\varepsilon)>\eta_r\}$.  Putting this together with $\varepsilon_\text{min}=\inf\{\varepsilon>0:p_\text{max}(\varepsilon)<\eta_p\}$, we can define an interval $(\varepsilon_\text{min},\varepsilon_\text{max})$ such that we are confident our kernel is well defined for any value of $\varepsilon$ in this range.

We would also like to bring attention to a nice method from Shan \& Daubechies \cite{shan2022diffusion}, that provides a principled and easily applicable metric to select on \textit{optimal} bandwidth.  It is true that the family of transition operators $\{P_\varepsilon\}_{\varepsilon>0}$ will form a semi-group under composition,\footnote{A semi-group is a set with an associative binary operation (a group without identity of inverses).} and will satisfy the relation $P_sP_t=P_{s+t}$.  Our discrete realisations of $P_\varepsilon$ will not satisfy this property exactly, but choosing the value of $\varepsilon$ where this discrepancy --- referred to as the Semi-Group Error (SGE) --- is minimised gives us a transition operator that best respects the underlying algebra of the set.  Since it seems reasonable to assume that the SGE will be minimised for values of $\varepsilon$ which best approximate $\exp(\varepsilon\mathcal{L})$, since these operators will trivially satisfy the composition rule, we can assume that the minimiser of the SGE will also give us a measure of when our approximation to the continuous $P_\varepsilon(x,y)$ is the best given our restricting to the data.  For clarity, we will state the SGE as defined in \cite{shan2022diffusion}:
\begin{equation}
    \text{SGE}(\varepsilon)=\lVert P_\varepsilon P_\varepsilon-P_{2\varepsilon}\rVert
\end{equation}
where in their paper they used the operator norm, but we found simpler norms, such the Frobenius norm, worked just as well.  When working with real data (as opposed to synthetic data), we have often used $p_\text{max}(\varepsilon)$ and $\text{res}(\varepsilon)$ to determine a maximal range of valid $\varepsilon$'s, and then used $\text{SGE}(\varepsilon)$ to localise us within this range, and within this local region we can then perform a grid search to determine our optimal choice of $\varepsilon$.

\subsubsection{Redundancy of ODK Parametrisations}
The above discussion required $\rho(x)$ to be fixed, as otherwise the above would not relate to $\varepsilon$ but rather $\varepsilon\rho(x)^2$, and thus we would still have to choose how to balance this weight between the two parameters.  To proceed without ambiguity, we must note that there is some redundancy in the parameter space of an ODK.  Precisely, denote by $\Phi=(\varepsilon,\tau(x),\rho(x))\in\mathbb{R}_{>0}\times C^2(\mathcal{M})\times C^2(\mathcal{M})$, and element of this parameter space, and let $K_\Phi$ be the ODK with this parameterisation; assuming $l(x)$ and $\kappa(x)$ are fixed, then we can define the following equivalence relationship $\Phi_1\sim\Phi_2\iff K_{\Phi_1}(x,y)=K_{\Phi_2}(x,y),\ \forall x,y\in\mathcal{M}$, by inspecting the definition of an isotropic ODK as given in equation \eqref{eq:general_odk}, we see that this is equivalent to the following:
\begin{align}
    \Phi_1\sim\Phi_2\iff&\forall x\in\mathcal{M},\  \varepsilon_1\tau_1(x)=\varepsilon_2\tau_2(x)\ \land\ \varepsilon_1\rho_1(x)^2=\varepsilon_2\rho_2(x)^2\label{appBeq:param_equivalence1}\\
    \iff&\forall x\in\mathcal{M},\ \frac{\tau_1(x)}{\rho_1(x)^2}=\frac{\tau_2(x)}{\rho_2(x)^2}\ \land\ \frac{\varepsilon_1\tau_1(x)^2}{\rho_1(x)^2}=\frac{\varepsilon_2\tau_2(x)^2}{\rho_2(x)^2}\label{appBeq:param_equivalence2}
\end{align}
Equation \eqref{appBeq:param_equivalence1} shows us that discussing the heuristic parametrisation of $\tau(x)$ and $\rho(x)^2$ without fixing $\varepsilon$ creates redundancy in our parameter space since the equivalence class (set of all parameters which produce the same kernel) of any parameterisation would satisfy:
\begin{equation}\label{appBeq:param_equivalence3}
    \left[(\varepsilon,\tau(x),\rho(x)^2)\right]=\{(\varepsilon',\tau'(x),\rho'(x)^2):\exists c>0,\ \varepsilon'=\varepsilon/c,\ \tau'(x)=c\tau(x),\ \rho'(x)^2=c\rho(x)^2\}
\end{equation}
If we view $\varepsilon$ as a time-scale again, then this equation is saying that all ODKs are equivalent up to a linear time-change.  It is worth dwelling on this: in the main text we always consider $\varepsilon$ as being a time-scale of the underlying process that generated the data, relative to the time horizon $[0,T]$, with $T<\infty$.  Supposing we have parametrised our ODK with ground truth drift and diffusion, and some $\varepsilon>0$, then we must apply the adjoint transition operator to $p_0$ approximately $S=\tfrac{T}{\varepsilon}$ times before our time marginals have spanned the data entirely.  Equation \eqref{appBeq:param_equivalence3} tells us that if we applied a linear time change $\varepsilon'= c\varepsilon$, for some $c>0$, we will get an identical kernel, and thus we must still apply our reparametrised ODK $S$ times to span our data, which implies the new effective time horizon of the underlying process is $T'=S\varepsilon'=cT$; and thus this redundancy in ODKs parameter space is structurally equivalent to equating all SDEs that are equal up to a linear time-change.

This observation provides us much more freedom in our choice of $\varepsilon$ when we are parametrising our kernel heuristically, since we can always transform back to a regime where $\varepsilon$ is a real time-scale of the underlying process.  As an example, suppose $T=5$, then we might decide that given the rough scale of the data (maximum length scale $\max_{x_i,x_j\in\mathcal{D}}\lVert x_i-x_j\rVert$, and the spread of the data $\mathbb{V}_{x\sim q}(x)$) the most likely order of magnitude of $\tau(x)$ and $\rho(x)^2$ would allow us to set $\varepsilon^*=0.1$ as our reference time-scale.  We could then arbitrarily set $\varepsilon=1$ and then heuristically estimate the values for $\tau(x_i)$ and $\rho(x_i)$ --- using techniques already mentioned above and/or those that we will shortly discuss --- without any consideration for the value of $\varepsilon$.  After we are done, we can simply transform back to the physically relevant reference time-scale regime, and thus our kernel would be unchanged, and our new heuristic estimates would be $\tau^*(x_i)=\tau(x_i)/\varepsilon^*$ and $\rho^*(x_i)=\rho(x_i)/\sqrt{\varepsilon^*}$.

It is important to note that we are not getting anything for free here; we cannot simply rescale to $\varepsilon'\ll 1$ and claim our convergence results are better, since only $\varepsilon$'s that resolve on the data and can reasonably be interpreted as a time-scale will give the asymptotics any definite meaning with respect to the underlying process.  Furthermore, if one simply rescaled $\varepsilon$ and $T$ and then applied the $L_q$ or $L_p$ optimisation schemes, we would get identical results as we would if we didn't rescale, since the time horizons change identically.  The utility of this rescaling would only be realised if we intended to estimate the drift and diffusion in a \textit{dimensionless} manner based on some geometric or statistical heuristics, and where, to this end, getting rid of the extra degree of freedom is convenient.  Another scenario where this rescaling invariance is of great use, is if we didn't know the time-horizon of the underlying process, in this case we cannot ascribe any meaning to $\varepsilon$ as a timescale, making it a truly dimensionless quantity that may as well be set equal to $1$.\footnote{We remark that this is the case because the concept of a short or long time only makes sense relative to the global time-horizone which you are considering; i.e. a decade is a long time relative to the life of a human, but it is an exceptionally short time relative to the age of the universe.}.  If one inferred the values of $\tau(x)$ and $\rho(x)$ in a such a manner as to match the true underlying dynamics that generated the data (for example, you applied the $L_q$ or $L_p$ objective of section \ref{subsec:numerical_results_sde}); then the fact that it is possible to rescale ODKs to a more physical time-scale allows you to still claim you have learnt the true underlying dynamics, just up to a multiplicative constant.  Furthermore, quantities of the form $\tau(x)^k\rho(y)^{-2k}$ for any non-zero $k$ will be invariant to linear time-changes; thus, one can talk very precisely about such relative quantities.

We can group the equivalence conditions given in equations \eqref{appBeq:param_equivalence1} \& \eqref{appBeq:param_equivalence2} into the following list:
\begin{equation}\label{appBeq:param_equivalence4}
    \Psi(x)=\frac{\tau(x)}{\rho(x)^2},\qquad \psi(x)=\frac{\varepsilon\tau(x)^2}{\rho(x)^2},\qquad b(x)=\varepsilon\tau(x),\qquad a(x)=\varepsilon\rho(x)^2
\end{equation}
We call $\Psi(x)$ the drift-to-diffusion ratio, $b(x)$ the drift-scale, and $a(x)$ the diffusion scale; the dimensions of each can be determined from: $\varepsilon\tau\sim [L]$ - a length scale; $\tau\sim[L][t]^{-1}=[u]$ - rate of advection; and $\rho^2\sim [D]=[L]^2[t]^{-1}$ - rate of diffusion; thus $\psi(x)\sim\tfrac{[L][u]}{[D]}=[1]$, is a dimensionless number in the form of a P{\'e}clet number \cite{bird2014transport}.  A P{\'e}clet number can be seen as a measure of the relative difference in the rate of transport due to advection, against the rate of transport due to diffusion.  We can also see that $\Psi(x)\sim[L]^{-1}$ is an inverse length scale, and particularly, it is invariant under linear time changes; so fixing a functional form for one regime will fix it for all others.  This property makes $\Psi(x)$ difficult to set heuristically, unless one had some additional information that resulted in a restricted functional form of $\Psi(x)$; if one were running multiple heuristic parametrisations of ODK, then setting the functional form of $\Psi(x)$ to be different across each would ensure each ODK is fundamentally distinct.

It is the quantities given in equation \eqref{appBeq:param_equivalence4} that we will construct heuristics for estimating.  We note that there are only two degrees of freedom in this choice, since fixing any two of these quantities determines the others; thus, of the below-listed heuristics, we only need to apply at most two.

\subsubsection{Constructing a Length Scale}
If we wished to fix $\varepsilon\tau(x)\sim[L]$ or $\varepsilon\rho(x)^2\sim[L]^2$, we need a local estimate of the length scale in each neighbourhood of $x_i\in\mathcal{D}$, this length scale will be adaptive to the sampling density thus almost guarantees a good fit for the kernel by ensuring it doesn't degenerate.  The most common method to extract such a length scale is to choose some $k\in\mathbb{N}_{>0}$ and define $r_{i,k}=\lVert x_i-x_i^k\rVert\sim[L]$.  It has been shown that such bandwidths significantly improve convergence of the discrete Laplacian to the continuous Laplacian when the data is sampled non-uniformly \cite{berry2016variable}.  If one assumes that $k$ is large enough as to mitigate the variance in computing local estimates, but small enough so that over $B_{i,k}\coloneqq B(x_i,\lVert x_i-x_i^k\rVert)$ we can assume that $q(x)$ is approximately uniform, then up to the leading order we will have $q(x)\propto \lVert x_i-x_{i,k}\rVert^{-d}$.  This can be seen through the following argument:
\begin{equation}\label{appBeq:lengthscale_density}
    \frac{k}{m}\approx q(B_{i,k})=\int_{B_{i,k}}q(y)\mathrm{d}y\approx q(x_i)\text{Vol}(B_{i,k})\implies q(x_i)\propto \frac{1}{\text{Vol}(B_{i,k})}=\frac{\Gamma\left(\tfrac{d}{2}+1\right)}{\pi^{d/2}r_{i,k}^d}\propto \frac{1}{\lVert x_i-x_{i,k}\rVert^d}
\end{equation}
The sketch we provide does not produce a consistent estimator of the density, for this one needs $\hat q(x)=\left((k(m)-1)/m\right)\text{Vol}(B(x,r_{k(m)}))^{-1}$, where $k(m)$ satisfies $\lim_{m\to\infty}k(m)=\infty$ and $\lim_{m\to\infty}k(m)/m=0$, and $r_{k(m)}=\lVert x_i-x_i^{k(m)}\rVert$ is a random variable depending on the sample $\mathcal{D}$, and where we use the fact that $\mathbb{P}(B(x,r_{k(m)}))\sim\text{Beta}(k(m),m-k(m)+1)$; to see how this all comes together we refer the reader to \cite{loftsgaarden1965nonparametric} for more details.  To us the important detail from the proof is that the neighbourhood parameter $k(m)=o(m)$ and $k(m)\to\infty$ as $m\to\infty$, so we know that arbitrarily choosing $k=30$ for every sample is not prudent.  Equation \eqref{appBeq:lengthscale_density} is effectively saying that the neighbourhood one defines in a kNN-graph is specifically adapted to the underlying sampling density; therefore, choosing $L_i= \lVert x_i-x_i^k\rVert$, and setting $\varepsilon\rho(x_i)^2\propto L_i^2$ will also make you kernel adaptive to the underlying distribution, and thus it is very unlikely one will have issues with kernel degeneration in this setting.
Other estimators of a local length scale are possible, such as:
\begin{equation}
    \log L_i=\frac{1}{k}\sum_{j=1}^k\log \lVert x_i-x_i^j\rVert \iff L_i=\left(\prod_{j=1}^k\lVert x_i-x_i^j\rVert\right)^{1/k}
\end{equation}
By averaging in log-space, we reduce the variance of our length scale estimate and stabilise the absolute differences $\log L_i-\log L_j$, which can be equivalently seen as stabilising the relative differences $L_i/L_j$, making our length scales more comparable over the data.  Furthermore, since the length scale is estimated from samples in a $d$-dimensional volume around $x_i$, it is may be more justified to consider the geometric mean as opposed to the arithmetic mean, since the relationship between length and volume is multiplicative; this is indicated by the fact that the geometric mean is far less dampened estimator of the radius of data sampled uniformly in a $d$-dimensional hypersphere.

Since $\varepsilon\tau(x_i)\sim[L]$, fixing a length scale can allow us to set one or both of $\rho(x_i)=c_\rho L_i$ and $\tau(x_i)=c_\tau L_i$, for appropriate constants of proportionality $c_\tau,c_\rho>0$.  It is important to note that the length scales of the drift and diffusion need not be the same, and setting them as such enforces a hard constraint on the drift-to-diffusion ratio $\Psi(x)\propto L_i^{-1}$, i.e. the drift-to-diffusion ratio will be adapted to the underlying density.

\subsubsection{Determing $\boldsymbol{\psi(x)}$ through a notion of Advection Rate}
Consider an idealised $d$-dimensional Gaussian density $p\sim\mathcal{N}(\varepsilon\tau e_1,\varepsilon\rho^2I_d)$, where $e_1$ is a basis vector in the first coordinate dimension of $\mathbb{R}^d$.  We define the Forward Bias as $\text{FB}(\varepsilon,\tau,\rho)\coloneqq\mathbb{P}(\{x\in\mathbb{R}^d:x_1>0\}$, if we suppose $e_1=\frac{\nabla l(x)}{\kappa(x)}$, $\text{FB}$ then the forward bias tells us the probability, given our parameterisation, of transitioning to a higher-ordered state.  We can derive an expression for this through:
\begin{align}
    \text{FB}(\varepsilon,\tau,\rho)&=\frac{1}{(2\pi\varepsilon\rho^2)^{d/2}}\int_0^\infty\int_{\mathbb{R}^{d-1}}\exp\left(-\frac{\lVert x-\varepsilon\tau e_d\rVert^2}{2\varepsilon\rho^2}\right)\mathrm{d}x_{2:d}\mathrm{d}x_1\\
    &=\frac{1}{(2\pi\varepsilon\rho^2)^{1/2}}\int_0^\infty\exp\left(-\frac{(x_1-\varepsilon\tau)^2}{2\varepsilon\rho^2}\right)\mathrm{d}x_1\prod_{i=2}^d\frac{1}{(2\pi\varepsilon\rho^2)^{1/2}}\int_{-\infty}^\infty\exp\left(-\frac{x_i^2}{2\varepsilon\rho^2}\right)\mathrm{d}x_i\\
    &=\frac{1}{\sqrt{\pi}}\int_{-\tfrac{\sqrt{\varepsilon}\tau}{\sqrt{2}\rho}}^\infty e^{-u^2}\mathrm{d}u=\frac{1}{2}\text{erfc}\left(\sqrt{\frac{\psi}{2}}\right)
\end{align}
Hence the forward bias is entirely a function of the P{\'e}clet-like ratio $\varphi(x)$.  So if we preset some value for the Forward Bias, say $\text{FB}=\eta_f$, then we will have $\psi=2\text{erfc}^{-1}(2\eta_f)^2$, and thus fixes one of the four quantities of equation \eqref{appBeq:param_equivalence4}.

\subsubsection{Additional Comment on the Effects of Increasing Dimensionality}
It is true that in high dimensions the covariance of a Gaussian does not act as our intuitions expect; consider a zero mean Gaussian distribution on $\mathbb{R}^d$, with covariance $\sigma^2 I_d$, then it is known that $\mathbb{P}(B(0,r))=P\left(\tfrac{d}{2},\tfrac{r^2}{2\sigma^2}\right)$, where $P(s,x)$ is the regularised lower-incomplete gamma function, and can be written as:
\begin{equation}
    P(s,x)=\frac{1}{\Gamma(s)}\int_0^xt^{s-1}e^{-t}\mathrm{d}t=\frac{x^se^{-x}}{\Gamma(s+1)}\sum_{k=0}^\infty \frac{x^k}{s^{(k)}}
\end{equation}
where $s^{(k)}=\prod_{i=0}^{k-1}(s+i)$ is the Pochhammer symbol.  By applying Stirlings formula for the gamma function, we can show that $\log P(s,x)\approx -s(\log s-\log x -1)-x$ which implies up to the leading order $\log P(s,x)\sim -s\log s$ and thus $P(s,x)\to 0$ as $s\to\infty$; we show evaluations of the $P\left(\tfrac{d}{2},\tfrac{r^2}{2\sigma^2}\right)$ across $r^2$ being integer multiples of the covariance and $d$ varying from $1-20$ in table \ref{appBtab:reg-lower-inc-gamma}.

\begin{table}[htbp]
    \centering
    \footnotesize
    \setlength{\tabcolsep}{4pt}
    \renewcommand{\arraystretch}{1.35}
    \begin{tabular}{|lccccccccccc|}
        \hline
        $r\ \backslash\ d$ & $1$ & $2$ & $4$ & $6$ & $8$ & $10$ & $12$ & $14$ & $16$ & $18$ & $20$ \\
        \hline
        $\sigma$  & $0.6827$ & $0.3935$ & $0.0902$ & $0.0144$ & $0.0018$ & $0.0002$ & $\sim 10^{-5}$ & $\sim10^{-6}$ & $\sim10^{-8}$ & $\sim10^{-9}$ & $\sim10^{-10}$ \\
        $2\sigma$ & $0.9545$ & $0.8647$ & $0.5940$ & $0.3233$ & $0.1429$ & $0.0527$ & $0.0166$ & $0.0045$ & $0.0011$ & $0.0002$ & $\sim10^{-5}$ \\
        $3\sigma$ & $0.9973$ & $0.9889$ & $0.9389$ & $0.8264$ & $0.6577$ & $0.4679$ & $0.2971$ & $0.1689$ & $0.0866$ & $0.0403$ & $0.0171$ \\
        $4\sigma$ & $0.9999$ & $0.9997$ & $0.9970$ & $0.9862$ & $0.9576$ & $0.9004$ & $0.8088$ & $0.6866$ & $0.5470$ & $0.4075$ & $0.2834$ \\
        $5\sigma$ & $1.0000$ & $1.0000$ & $0.9999$ & $0.9997$ & $0.9984$ & $0.9947$ & $0.9852$ & $0.9654$ & $0.9302$ & $0.8751$ & $0.7986$ \\
        \hline
    \end{tabular}
    \caption{Numerical evaluation of the regularised lower incomplete gamma function $P\left(\tfrac{d}{2},\tfrac{r^2}{2\sigma^2}\right)=\mathbb{P}\left(\{x\in\mathbb{R}^d:\lVert x\rVert<r\}\right)$, where the latter expression gives the probability of sample $x\sim \mathcal{N}(0,\sigma^2I_d)$ being inside the ball or radius $r$.  Each row represents gives $r=k\sigma$ for $k=1,2,3,4,5$, and we vary the dimensionality along the columns, $d=1,2,4,\dots,20$.}
    \label{appBtab:reg-lower-inc-gamma}
\end{table}

Our intuition from lower dimensions is that being a distance of $3\sigma$ away from the mean of a Gaussian distribution is \ textit {far away} from the perspective of one's likelihood to be sampled there.  Table \ref{appBtab:reg-lower-inc-gamma} shows us that the probability of being sampled at a distance of less than $3\sigma$ from the mean decays rapidly as $d$ increases; indeed, when $d=10$ it is now more probable that we sample a state to be further than $3\sigma$ from the mean than it is to sample one closer.  This is a clear manifestation of the concentration of measure \cite{ledoux2001concentration}, where it is shown that high-dimensional Euclidean metrics begin to lose their discriminative power as $d\to\infty$.  We can be more rigorous in this discussion by letting $P_{d/2}(x)=P(d/2,x)$ and noting that this is a Cumulative Distribution Function (CDF) and thus is monotonically increasing in $x$, and its inverse exists; define the radius to be:
\begin{equation}\label{appBeq:radius_definition}
    \text{rad}(d,\eta_d;\sigma)=\sigma\sqrt{2P_{d/2}^{-1}(\eta_d)}\approx\sigma\sqrt{d}\left(1+\Phi_{0,1}^{-1}(\eta_d)\sqrt{\tfrac{2}{9d}}-\tfrac{2}{9d}\right)^{3/2}
\end{equation}
where $\Phi_{0,1}$ is the CDF of a standard univariate normal distribution with mean $0$ and variance $1$.  To deduce the approximate form in equation \ref{appBeq:radius_definition}, we note that if $x\sim\mathcal{N}(0,I_d)$ then $\lVert x\rVert^2\sim\chi_d^2$, by definition, and therefore $\eta_d=P(d/2,r^2/2)=\mathbb{P}(\chi_d^2\leq r^2)$ which allows us to rewrite the formula for the radius as $\text{rad}(d,\eta_d)=\sqrt{Q_{\chi_d^2}(\eta_d)}$, where $Q_{\chi_d^2}$ is the quantile function (inverse cdf) of a $\chi_d^2$ random variable.  In equation 5 of \cite{wilson1931distribution}, Wilson gives the following approximation: $W\coloneqq(\chi_d^2/d)^{1/3}\approx\mathcal{N}(\mu=1-\tfrac{2}{9d},\sigma^2=\tfrac{2}{9d})$, putting this together:
\begin{equation}
    \eta_d=\mathbb{P}(\chi_d^2\leq r^2)=\mathbb{P}\left(\sigma^{-1}(W-\mu)\leq \sigma^{-1}\left(\left(\tfrac{r^2}{d}\right)^{1/3}-\mu\right)\right)=\Phi_{0,1}(s)
\end{equation}
Therefore, setting $z_{\eta_d}\coloneqq\Phi_{0,1}^{-1}(\eta_d)$ we can easily show that $r=\sqrt{d}(\mu+\sigma z_{\eta_d})^{3/2}$, which, after rescaling by $\sigma$, is the exact form of equation \eqref{appBeq:radius_definition} after substituting the values for $\mu$ and $\sigma$.  Thus, the radius defining the boundary beyond which we don't expect to see many samples scales with $\sqrt{d}$ as $d$ increases.  We can get the same result by applying the central limit theorem to $\lVert x\rVert^2$, in which we see:
\begin{equation}
    \lVert x\rVert^2= d\left(\frac{\sum_{i=1}^dx_i^2}{d}-\sigma^2\right)+d\sigma^2\sim\mathcal{N}(d\sigma^2,2d\sigma^4)
\end{equation}
where we have used the fact that $\mathbb{E}[x_i^2]=\sigma^2$ and $\mathbb{E}[x_i^4]-\mathbb{E}[x_i^2]^2=(\kappa-1)\sigma^4=2\sigma^4$ where $\kappa$ is the kurtosis of the Gaussian random variable $x_i$, and which we know is equal to 3.  We can then write $\lVert x\rVert^2=d\sigma^2+\sqrt{2d}\sigma^2z$ where $z\sim\mathcal{N}(0,1)$, the leading order term of which is $d\sigma^2$ which is exactly the leading order term of $\text{rad}(d,\eta;\sigma)^2$.  

Thus, whether we introduce a dimension-dependent rescaling would depend partly on whether $L_i\sim \sigma$ or $L_i\sim \sigma\sqrt{d}$, i.e. whether $L_i$ is a statistical or geometric scale.  It is true that the best isotropic approximation to a $d$-dimensional Gaussian with arbitrary covariance $\Sigma$ is given by $\left(\frac{1}{d}\text{tr}\ \Sigma\right)I_d$, and we can see that:
\begin{equation}
    \tfrac{1}{d}\text{tr}\ \Sigma=\tfrac{1}{d}\sum_{i=1}^d\mathbb{E}[(x_i-\mathbb{E}[x_i])^2]=\tfrac{1}{d}\mathbb{E}[\left\lVert x-\mathbb{E}[x]\right\rVert^2]\propto \tfrac{1}{d}\lVert x-x^k\rVert^2\sim \sigma^2
\end{equation}
Hence if $L_i\sim \sigma_i$, then setting $\varepsilon\rho^2(x_i)\sim d\cdot L_i^2$ will remove some of the dimensionality dependence in the norm and thus such a scaling will make the resulting kernel more intuitive from the perspective of this length scale; however, if $L_i\sim \sqrt{d}\sigma$, i.e. you have used $L_i=\lVert x_i-x_i^k\rVert$, then it would not be wise to rescale by the dimensionality and thus leaving $\varepsilon\rho^2\sim L_i^2$ is more prudent.

Ultimately, if you use the Euclidean metric in high dimensions, you will run into issues due to the concentration of measure: distances will become very large and roughly symmetric about their mean.  No amount of rescaling can fix the loss of discrimination, but if you have constructed a length scale that depends on intrinsic distributional scales, then introducing a dimension-dependent rescaling may help make one's kernels more interpretable.

\subsection{Strang Splitting Scheme}\label{appBsubsec:strang_splitting}
In this supplementary note, we will give some theoretical context for the Strang splitting scheme that we introduce in section \ref{subsec:nonisotropic_extension} to generalise ODKs to anisotropic diffusions; we will also discuss the choice of the parameter $\alpha\in(0,1)$ that splits the diffusion between our two operators.

To recap, we have applied two kernels that approximate the following two operators:
\begin{align}
    \mathcal{L}^\mu &= \mu^i\partial_i +\frac{1}{2}D^{ij}_\mu\partial_i\partial_j & \text{where:}\qquad D_\mu&=\alpha\sigma_h^2\rho^2\lambda_\text{min}(A)I_d\\
    \mathcal{L}^A &= \frac{1}{2}D^{ij}_A\partial_i\partial_j & \text{where:}\qquad D_A&=\sigma_h^2\rho^2(A-\alpha\lambda_\text{min}(A)I_d)=\sigma_h^2\rho^2 A-D_\mu
\end{align}
where $A(x)$ is a field of symmetric positive definite $\binom{2}{0}$-tensors over $\mathcal{M}$; $\mu(x)$ is an arbitrary vector field since we could construct $\mathcal{L}^\mu$ using a local order; $\sigma_h(x)$ is a feature of the choice of ODK; $\tfrac{1}{2}\rho^2 A(x)$ is the diffusion tensor we wish to reconstruct; $\lambda_\text{min}(A(x))>0$ is the smallest eigenvalue of $A(x)$ which must be greater than zero due to the assumption of $A$ being positive definite; finally, $\alpha\in(0,1)$ is a parameter we choose that decides how to split the diffusion contributions between the two kernels.  Our goal is to approximate the operator $\mathcal{L}=\mu^i\partial_i+\frac{1}{2}D^{ij}\partial_i\partial_j$ where $D=\sigma_h^2\rho^2A$ and we note that indeed $\mathcal{L}=\mathcal{L}^\mu+\mathcal{L}^A$ so we could stop here if $\mathcal{L}$ was all we wanted.  However, by doing this we are breaking the usual chain of constructions $G_\varepsilon\to P_\varepsilon\to L_\varepsilon$, since we do not know the $\varepsilon$ transition matrix of the Markov process we are building, and indeed if the errors in estimating $\mathcal{L}^\mu$ and $\mathcal{L}^A$ with $L_\varepsilon^\mu$ and $L_\varepsilon^A$ add up unfavourably, there may exist no valid matrix $P_\varepsilon$ such that $\varepsilon^{-1}(P_\varepsilon-I)=L_\varepsilon^\mu+L_\varepsilon^A$.  More fundamentally, the family of transition operators $\{P_\varepsilon\}_{\varepsilon>0}\cup\{I\}$ forms a semi-group (no inverses) under composition (which is matrix multiplication in the discrete case) and $\mathcal{L}$ is the (operator) right derivative of this group at $\varepsilon=0$, and acts (when it exists) as the unique generator of the transition group through exponentiation, i.e. $P_\varepsilon=\exp(\varepsilon\mathcal{L})$, where the operator exponential is defined in terms of its power-series $\exp(\varepsilon\mathcal{L})=\sum_{k=0}^\infty\tfrac{\varepsilon^k}{k}\mathcal{L}^k$.  Since composition of transition operators is a multiplicative operation, it seems natural to construct our desired $P_\varepsilon$ in such a manner; we see that:
\begin{equation}
    P_\varepsilon=\exp(\varepsilon\mathcal{L})=\exp(\varepsilon\mathcal{L}^\mu+\varepsilon\mathcal{L}^A)
\end{equation}
we also have that $P_\varepsilon^\mu=\exp(\varepsilon \mathcal{L}^\mu)$ and $P_\varepsilon^A=\exp(\varepsilon\mathcal{L}^A)$, the crux is that unless $\mathcal{L}^\mu$ and $\mathcal{L}^A$ commute, $\exp(\varepsilon\mathcal{L}^\mu+\varepsilon\mathcal{L}^A)\neq\exp(\varepsilon \mathcal{L}^\mu)\exp(\varepsilon \mathcal{L}^A)$.  For two operators $A,B$, we denote the \textit{commutator} as $[A,B]\coloneqq AB-BA$; the commutator trivially satisfies that $A$ and $B$ commute if and only if $[A,B]=0$.  The Baker-Campbell-Hausdorff (BCH) formula (as given in \cite{hairer2006numerical}, Chapter 3.4) gives as expansion in terms of $\varepsilon$ for an operator $C_\varepsilon$ such that $\exp(C_\varepsilon)=\exp(\varepsilon A)\exp(\varepsilon B)$, the CBH gives this expansion as:
\begin{equation}\label{appB:CBH_formula}
    C_\varepsilon = \varepsilon(A+B)+\frac{1}{4}\varepsilon^2[A,B]+\frac{1}{12}\varepsilon^3\left([A,[A,B]]+[B,[B,A]]\right)+\frac{1}{24}\varepsilon^4[A,[B,[B,A]]]+\dots
\end{equation}
with all terms of order $\varepsilon^2$ and higher being in terms of the commutator of $A$ and $B$ - a quick sanity check then shows that if $[A,B]=0$ that $C_\varepsilon=\varepsilon(A+B)$.  Thus, the error of approximating $P_\varepsilon$ with $P_\varepsilon^AP_\varepsilon^\mu$, is $\bigO(\varepsilon^2)$ in the log-space of $P_\varepsilon$, where $\log P_\varepsilon = \log\exp(\varepsilon\mathcal{L})=\varepsilon\mathcal{L}$ the inverse of the operator exponential (defined in terms of its power series).  This approximation scheme is often called Lie-Trotter Splitting (\cite{trotter1959product}), and is a generalisation of the first-order Euler method for numerical integration.  Naturally, we can reduce the error in our approximation by considering generalisations of higher-order integration methods; a generalisation of the Verlet Scheme gives:
\begin{equation}\label{appBeq:appB_strang_splitting}
    P_\varepsilon\approx P_{\varepsilon/2}^AP_\varepsilon^\mu P_{\varepsilon/2}^A
\end{equation}
which, in the context of operators, is called the Strang Splitting Scheme (\cite{strang1968construction}).  Applying the BCH expansion twice, we find the Laplacian approximation error to be:
\begin{equation}\label{appBeq:appB_strang_error}
    \log\left(P_{\varepsilon/2}^AP_\varepsilon^\mu P_{\varepsilon/2}^A\right)-\varepsilon\mathcal{L}=\frac{1}{12}\varepsilon^3\left(\left[\mathcal{L}^\mu,[\mathcal{L}^\mu,\mathcal{L}^A]\right]-\tfrac{1}{2}\left[\mathcal{L}^A,[\mathcal{L}^A,\mathcal{L}^\mu]\right]\right)+\bigO(\varepsilon^5)
\end{equation}
where every even power of $\varepsilon$ has been cancelled.  We note that this cancellation is a property of the symmetry of the product, not of the halving of $\varepsilon$ in the outer terms; the $\varepsilon/2$ ensures the first-order term equals $\varepsilon\mathcal{L}$, and indeed matches the Verlet Scheme for integrating functions.  Naturally, we are faced with a choice of whether to have $P_\varepsilon^\mu$ or $P_\varepsilon^A$ as the central term, but we found during numerical experimentation that having the drifted term in the centre (as in equation \eqref{appBeq:appB_strang_splitting}) is more accurate and more robust to the choice of $\alpha$ than the other way around.  This intuitively makes sense, as if the drift were anything other than linear, it would push the probability in the direction of $\mu(x)$ and then $\mu(x+\tfrac{\varepsilon}{2}\mu(x))$ which are not necessarily colinear for a practical value of $\varepsilon$; the robustness to $\alpha$ is likely due to the eigenstructure of $A(x)$, if the smallest eigenvalue is too small then $\varepsilon$ cannot be made small without degenerating the kernel.

\subsubsection{Choosing $\boldsymbol{\alpha\in(0,1)}$}\label{appBsubsec:choosing_alpha}
To apply the Strang splitting scheme we introduced a parameter $\alpha\in (0,1)$ to determine how the diffusion is split between the two operators $\mathcal{L}_\alpha^\mu=\mu\cdot\nabla+\tfrac{1}{2}\alpha\rho^2\lambda_\text{min}(A)\Delta$ and $\mathcal{L}_\alpha^A=\tfrac{1}{2}\rho^2(A-\alpha\lambda_\text{min}(A)I_d)_{ij}\nabla_i\nabla_j$, where $\lambda_\text{min}(A-\alpha\lambda_\text{min}(A)I_d)=(1-\alpha)\lambda_\text{min}(A)$.  In theory, $\alpha$ is a free parameter for you to choose, but in practice, we need to choose $\alpha$ such that both $P_{\varepsilon/2}^A$ and $P_{\varepsilon}^\mu$ are best resolved on the data.  We will proceed to make two rather Heuristic arguments, before ultimately contradicting ourselves.

Suppose in a local neighbourhood of any given $x\in\mathcal{D}$, samples are roughly uniformly distributed; then the kernel becomes degenerate when it decays almost entirely (along its weakest direction) before being evaluated on a neighbouring data point.  We could write this as $\varepsilon\lambda_\text{min}(\Sigma)\geq c(\mathcal{D})$ where $\Sigma$ is a positive-definite covariance matrix, and $c(\mathcal{D})$ will depend on the dimension of the underlying manifold and the number of data points, as we are not being rigorous we simply assume $c$ to be constant across all $x\in\mathcal{D}$.  In the case of $P_{\varepsilon/2,\alpha}^A$ and $P_{\varepsilon,\alpha}^\mu$ we have the constraints:
\begin{align}
P_{\varepsilon/2,\alpha}^A:&\qquad \tfrac{1}{2}\varepsilon\rho^2\lambda_\text{min}(A-\alpha\lambda_\text{min}(A)I_d)&\geq c &\qquad&\implies& \tfrac{1}{2}\varepsilon(1-\alpha)&\geq\frac{c}{\rho^2\lambda_\text{min}(A)} \\
P_{\varepsilon,\alpha}^\mu:&\qquad\qquad\qquad\qquad\  \varepsilon\rho^2\alpha\lambda_\text{min}(A)&\geq c &\qquad&\implies& \varepsilon\alpha &\geq\frac{c}{\rho^2\lambda_\text{min}(A)} 
\end{align}
If we say that the optimal $\alpha$ is the one which allows us to make $\varepsilon$ the smallest:
\begin{align}
    \varepsilon^*(\alpha)&=\max\left(\frac{2c}{(1-\alpha)\lambda_\text{min}(A)\rho^2},\frac{c}{\alpha\lambda_\text{min}(A)\rho^2}\right)=\max\left(\frac{2}{1-\alpha},\frac{1}{\alpha}\right)\frac{c}{\rho^2\lambda_\text{min}(A)}\\
    &=\min\left(\frac{1-\alpha}{2},\alpha\right)^{-1}\frac{c}{\rho^2\lambda_\text{min}(A)}
\end{align}
i.e. choosing $\alpha$ such that $\min((1-\alpha)/2,\alpha)$ is minimal allows us to choose the smallest value of $\varepsilon$.  It is not hard to show that $\alpha=\frac{1}{3}$ is this value.

Our second argument looks at the $\varepsilon^3$ error term in equation \eqref{appBeq:appB_strang_error}.  This argument is incredibly hand-wavy, but if we suppose $\mathcal{L}^\mu\propto \alpha$ and $\mathcal{L}^A\propto (1-\alpha)$, then, after completely ignoring the application of the commutator, we would have:
\begin{equation}
    2\left[\mathcal{L}^\mu,[\mathcal{L}^\mu,\mathcal{L}^A]\right]-\left[\mathcal{L}^A,[\mathcal{L}^A,\mathcal{L}^\mu]\right]\propto 2\alpha^2(1-\alpha)-\alpha(1-\alpha)^2
\end{equation}
the absolute value of which is minimised at $\alpha=\frac{1}{3}$.  Putting this unreasonable argument together with the heuristic argument involving $\varepsilon^*$ we feel emboldened in our belief that $\alpha_\text{heuristic}^*=\frac{1}{3}$ is the optimal value of $\alpha$.  In spite of our confidence, when we tested this splitting method on synthetic data (as shown in the operator results of section \ref{subsec:numerical_results_operators}), we found that $\alpha^*_\text{data}=\frac{2}{3}$ produced the best approximation of $L_\varepsilon$ to $\mathcal{L}$ (quantified via the $L^2$-norm).  

Furthermore, if we reversed the order of composition, both of the above heuristic arguments would result in a new optimum of $1-\alpha^*_\text{heuristic}$, but by switching the order of composition in the numerical example of \ref{subsec:numerical_results_operators} we see that the optimum remains fixed at $\tfrac{2}{3}$.  This tells us that the two heuristics we gave above seem to have no bearing on the optimal $\alpha$, and that the composition scheme is not the limiting factor in making the approximation $\bar P_{\varepsilon,\alpha}$.  In hindsight this latter point is not surprising, since by Theorem \ref{thm:odk_backward_expansion}, we know that the leading order error in approximating any $P_{\varepsilon}$ can be as high as $\bigO(\varepsilon^{3/2})$, and thus the $\bigO(\varepsilon^3)$ error coming from the composition would be insignificant in comparison.  This leads us to conclude that the choice of $\alpha$ depends on the eigenstructure of $A(x)$ and one should choose $\alpha$ as large as is possible without $P_{\varepsilon/2,\alpha}$ degenerating on the data.  We also believe that the optimal order of composition for the Strang splitting scheme will depend on which of $A(x)$ and $\mu(x)$ is changing most rapidly on the data --- and since in practice variations of $\mu$ are likely to cause more pronounced errors we would generally recommend using $P_{\varepsilon/2,\alpha}^AP_{\varepsilon,\alpha}^\mu P_{\varepsilon/2,\alpha}^A$.

\subsection{Riemannian Geometry Preliminaries}\label{appBsec:prelims}
In this appendix, we will give some background on the tools from differential geometry that we use throughout this paper; we hope that this provides clarity on which differential geometric concepts are important to us and how and why we use them.

We assume that we are working on a $d$-dimensional Riemannian manifold embedded in $\mathbb{R}^n$, where $n\gg d$; we usually denote this manifold as $\mathcal{M}$, but when we do so, we are talking about the manifold intrinsically.  It would be more proper to say that $(\mathcal{M},g_\mathcal{M})$ is a $d$-dimensional Riemannian manifold with metric $g_\mathcal{M}$, and there is an embedding $\iota:\mathcal{M}\to\mathbb{R}^n$ whose image we denote $\mathcal{N}\coloneqq \iota(\mathcal{M})$.  In differential geometry, an embedding is an immersion that is also a differentiable homeomorphism onto its image (i.e. a diffeomorphism); the latter condition means that the topological structure (the open sets) of $\mathcal{M}$ and $\mathcal{N}$ are identical, while being an immersion (guarantees that the manifold structure is preserved through the mapping; therefor, $\iota$ being an embedding means that $\mathcal{N}$ is a $d$-dimensional submanifold of the $n$-dimensional manifold $\mathbb{R}^n$.  In addition, $(\mathcal{M},g_\mathcal{M})$ is (locally) a metric space, so for us to consider $(\mathcal{M},g_\mathcal{M})$ and $(\mathcal{N},g_\mathcal{M})$ is be isomorphic under $\iota$ we need the additional condition that $\iota:(\mathcal{M},g_\mathcal{M})\to(\mathcal{N},g_\mathcal{N})$ be an isometry --- in other words, distances on $\mathcal{M}$ and $\mathcal{N}$ are also the same --- this is equivalent to $g_\mathcal{M}=\iota^*g_\mathcal{N}$ where $\iota^*$ is the pullback and is defined as $\iota^*g_\mathcal{N}(u,v)=g_\mathcal{N}(D\iota(u),D\iota(v))$.  Any calculations we make we do so on $\mathcal{N}\subset\mathbb{R}^n$, but we assume there is an isometric embedding $\iota$ that allows us to consider all quantities computed on $\mathcal{N}$ as being computed on $\mathcal{M}$ and we do so without ambiguity.  Furthermore, we assume that $g_\mathcal{N}$ is inheritted from the Euclidean metric on $\mathbb{R}^n$, i.e. for any $x\in\mathcal{M}$ let $u,v\in T_x\mathcal{M}$, then $\tilde u=D\iota(x)(u)$ and $\tilde v=D\iota(x)(v)$ will be two vectors in $\mathbb{R}^n$ such that $\iota(x)+\tilde u$ and $\iota(x)+\tilde v$ are tangent to $\mathcal{N}$ in $\mathbb{R}^n$ at $\iota(x)$, we assume that $g_\mathcal{N}(\tilde u,\tilde v)=\langle \tilde u,\tilde v\rangle_2=\tilde u^T\tilde v=u(D\iota^T D\iota)v$, and thus the metric $g_\mathcal{M}=D\iota^T D\iota$ when written in orthonormal coordinates (projection coordinates as we will soon call them) on $\mathcal{N}$.  From now on we hope that you are convinced that we may work unambiguously on the manifold $\mathcal{M}$ with metric $g_\mathcal{M}$.

A geodesic through a point $x\in \mathcal{M}$ is a curve $\gamma:I\to\mathcal{M}$, where $I\subset\mathbb{R}$ is a closed interval of the form $[0,T]$ such that $\gamma(0)=x$ and whose rate of change is constant (constant velocity), and which we formulate as $\gamma$ having no acceleration.  To formalise this we need to have a notion of differentiation vector fields which require us to be able to take differences of the form $\gamma'(t+dt)-\gamma'(t)$, where $\gamma'(t)$ is the standard derivative of $\gamma:I\to\varphi(U)$ where $(U,\varphi)$ is a coordinate chart containing $\gamma(t)$ and which has a coordinate independent meaning up to composition with the transition functions of compatible charts.  The issue is that $\gamma'(t+dt)\in T_{\gamma(t+dt)}\mathcal{M}$ and $\gamma'(t)\in T_{\gamma(t)}\mathcal{M}$ exist in different vector spaces and thus cannot be naively subtracted.  This results in the notion of a connection which allows us to \text{connect} different vector spaces defined over $\mathcal{M}$ together; more formally, suppose that $E$ is a vector bundle over $\mathcal{M}$, and $\mathcal{E}(\mathcal{M})$ is the space of smooth sections of $E$\footnote{A rank-n vector bundle, $E$, over $\mathcal{M}$ is a set of the form $E=\bigsqcup_{x\in\mathcal{M}} D_x$ where each $D_x$ has the structure of a vector space and all $D_x$ have dimension $n$; a smooth section of a bundle is a smooth function $V:\mathcal{M}\to E:x\mapsto (x,v)$ where $v\in D_x$, i.e. a section of a bundle maps points of $\mathcal{M}$ to vectors in $E$.}, then a connection in $E$ is a map $\nabla:\mathcal{T}\mathcal{M}\times\mathcal{E}(\mathcal{M})\to\mathcal{E}(\mathcal{M})$ which takes a section of the tangent bundle (i.e. a vector field) and a section of the $E$ and outputs another section of $E$, and must be linear over $C^\infty(\mathcal{M})$ in the first argument, linear in $\mathbb{R}$ in the second argument, and for $f\in C^\infty(\mathcal{M})$ must satisfy the product rule $\nabla_X(fY)=f\nabla_X(Y)+(Xf)Y$ where $Xf(p)=df(X(p))$ is the directional derivative of $f$ in the direction $X(p)$.  This gives the connection the structure of a directional derivative and thus $\nabla_X Y$ is called the covariant derivative of $Y$ in the direction of $X$.  If $\mathcal{E}(\mathcal{M})=\mathcal{T}(\mathcal{M})$, i.e. $E=T\mathcal{M}$ is the tangent bundle, then we call $\nabla:\mathcal{T}(\mathcal{M})\times\mathcal{T}(\mathcal{M})\to\mathcal{T}(\mathcal{M})$ a linear connection on $\mathcal{M}$; without going into detail, there is a one-to-one correspondence between linear connections and choices of $d^3$ smooth functions $\Gamma_{ij}^k$ on $\mathcal{M}$, and they are related through:
\begin{equation}
    \nabla_X Y=(X^i\partial_iY^k+X^iY^j\Gamma_{ij}^k)\partial_k
\end{equation}
The smooth functions $\Gamma_{ij}^k$ are called the Christoffel symbols, and although they are written with tensor indices, they are not a tensor.  It is also true that: every smooth manifold will admit a linear connection and this can be proved constructively; the existence of a linear connection induces (under the addition of constraints) a unique connection on each tensor bundle $T_l^k\mathcal{M}$.\footnote{This notation defines the bundle of all type $\binom{k}{l}$-tensors across $\mathcal{M}$.  A $\binom{k}{l}$-tensor can be viewed as a multilinear map acting on $k$-covariant vectors and $l$-contravariant vectors; we recall that co- and contravariant vectors are determined by the transformation properties of the vector with respect to a change in coordinates, where covariant vectors transform identically to the coordinate basis and contravariant vectors inversely to the coordinates.}  The connection is the mathematical machinery that allows us to take the derivative of a vector field along a curve: supposing $\nabla$ is a linear connection and $\gamma:I\to\mathcal{M}$ a curve, then $\nabla$ determines a unique operator $D_t:\mathcal{T}(\gamma)\to\mathcal{T}(\gamma)$ which maps smooth vector fields to smooth vector fields over $\gamma(I)$, that satisfies linearity over $\mathbb{R}$, the product rule and (if $V\in\mathcal{T}(\gamma)$ is a restriction of some vector field $\tilde V\in\mathcal{T}(\mathcal{M})$) $D_tV(t)=\nabla_{\gamma'(t)}\tilde V$; we call $D_t$ the covariant derivative along $\gamma$.  A vector field $V$ is said to be parallel along $\gamma$ if $D_t V\equiv 0$, if $V$ is parallel along any curve then we just say $V$ is parallel, it is true that for any curve $\gamma$ and $V_0\in T_{\gamma(0)}\mathcal{M}$ that there is a unique parallel vector field $V$ along $\gamma$ such that $V(\gamma(0))=V_0$; the uniqueness in our statement allows us to unambiguously transport a vector $V_0$ along $\gamma$ to a vector $V(\gamma(T))\in T_{\gamma(T)}\mathcal{M}$. 

We already have enough machinery to describe a geodesic, but we will deal with two further technicalities first.  The first issue is that if $X,Y\in T_p\mathcal{M}$ and we take a short path parallel to $X$, taking us to $x$, and then from $x$ take another short path parallel to $Y$ to get to $q$, supposing the short paths we the same, would $r\in\mathcal{M}$ --- the point reached after taking the short path parallel to $Y$ to get to $y$, and then the short path parallel $X$ at $y$ --- be the same as $q$?  It turns out that this need not be true, and the extent of this failure can be quantified by the torsion tensor $\tau(X,Y)=\nabla_XY-\nabla_YX-[X,Y]$, where $[\cdot,\cdot]$ is the Lie-Bracket of vector fields (or equivalently the commutator when viewing a vector field as an operator on smooth functions).  For any connection on a Euclidean space $\tau(X,Y)\equiv 0$, so we will assume that our chosen connection $\nabla$ satisfies this as well --- in such a case the linear connection is called symmetric, and this is equivalent to the associated Christoffel symbols satisfying $\Gamma_{ij}^k=\Gamma_{ji}^k$.  The final constraint we need is that our chosen connection is compatible with the metric, i.e. $\nabla_Xg( Y,Z)=g(\nabla_XY,Z)+g(Y,\nabla_XZ)$, this is equivalent to: $\nabla g\equiv 0$; if $V,W$ are parallel vector fields along a curve $\gamma$, then $g(V,W)$ is constant; parallel translation is an isometry.  Thus, we can invoke the fundamental lemma of Riemannian geometry, which states: for a Riemannian manifold $(\mathcal{M},g)$, there is a unique linear connection $\nabla$ that is compatible with $g$ and symmetric --- we call this unique connection the Levi-Civita connection.  

As a result, we can define a geodesic through a point $x\in\mathcal{M}$ as a curve $\gamma:I\to\mathcal{M}$ originating at $x$ and satisfying $D_t\dot\gamma\equiv 0$ where $D_t$ is with respect to the unique symmetric connection compatible with $g_\mathcal{M}$.  Given $x\in\mathcal{M}$ and $v\in T_x\mathcal{M}$ there exists a unique maximal geodesic $\gamma_v:I\to\mathcal{M}$ such that $\gamma(0)=x$ and $\dot\gamma(0)=v$.  We can define a subset of the tangent bundle $\mathcal{E}\coloneqq \{V\in T\mathcal{M}:\gamma_V\text{ is defined on the interval }[0,1]\}\subset T\mathcal{M}$ in which the so-called exponential map can be naturally defined $\exp:\mathcal{E}\to\mathcal{M}:V\mapsto\gamma_V(1)$, with $\exp_x$ being the exponential map restricted to $T_x\mathcal{M}$ (i.e. $\mathcal{E}\cap T_x\mathcal{M}$).  Two important properties of the exponential map are: it commutes over isometries, i.e. $\iota\circ\exp_x=\exp_{\iota(x)}\circ D\iota$; there exists a neighbourhood $V\subset T_x\mathcal{M}$ of the origin, and a neighbourhood $U$ of $x$ in $\mathcal{M}$ such that $\exp_x:V\to U$ is a diffeomorphism, any neighbourhood $U$ in which this property holds is called a normal neighbourhood.  If $\{e_i\}_{i=1}^d$ form an orthonormal basis of $T_x\mathcal{M}$ with respect to the metric $g=g_\mathcal{M}$, then we can define the isomorphism $E:\mathbb{R}^d\to T_x\mathcal{M}:(x^1,\dots,x^d)\mapsto x^ie_i$, if $U$ is a normal neighbourhood of $x$ then $\varphi\coloneqq E^{-1}\circ\exp_x^{-1}:U\to\mathbb{R}^d$ defines a coordinate chart, the respective coordinates are called the normal coordinates at $x$.  Expressing things in terms of normal coordinates flattens $U$ in the most canonical manner; thus, many Euclidean intuitions hold true in these coordinates. In particular, if we write the metric in these coordinates, we will find $g_{ij}=\delta_{ij}$, i.e. the metric appears flat.

We can define another set of coordinates on $\mathcal{M}$ that flattens the metric using the fact that $\mathcal{N}=\iota(\mathcal{M})\subset \mathbb{R}^n$, thus we can view $T_{x}\mathcal{M}$, for any $x\in\mathcal{M}$, as being an affine subspace of $\mathbb{R}^n$, i.e. a $d$-dimensional subspace $W_x\leq\mathbb{R}^n$ and $T_{x}\mathcal{M}\cong \iota(x)+W_x$, and let us denote $T_x\mathcal{N}=x+W_x$.  Taking any orthonormal basis $\{\bar{e}_i\}_{i=1}^d$ of $W_x$, noting that $\bar{e}_i\in\mathbb{R}^n$ and with orthonormality being in the sense of the Euclidean inner product, then defining the projection operator $P_x:\mathbb{R}^n\to\mathbb{R}^d:v^1\bar{e}_1+\dots v^d\bar{e}_d\mapsto (v^1,\dots,v^d)$ we can define the coordinate chart $\phi_x:U\to\mathbb{R}^d:y\mapsto P_x(\iota(y)-\iota(x))$, where $U$ is an open neighbourhood of $x$ such that $P_x$ is injective on $\iota(U)$; the existence of such a neighbourhood is guaranteed by the implicit function theorem.  We call these coordinates the projection coordinates, and note that unlike the normal coordinates, they depend on the embedding $\iota$.

If we define the normal coordinates, $(s^1,\dots,s^d)$, and projection coordinates, $(u^1,\dots,u^d)$, with respect to the same orthonormal basis --- i.e. $D\iota(e_1)=\bar e_1$ using the previous notation --- then these two coordinate charts are naturally related.  The projection coordinates take the best possible affine approximation to $\mathcal{N}$ at some $\iota(x)$, whereas the normal coordinates, by nature of their implicit definition on $\mathcal{M}$, will curve along $\mathcal{N}$, we could think of the normal coordinates as being the shadow of the projection coordinates projected down onto $\mathcal{N}$; hence they are highly related but not identical equal to one another.  We now make a very practical note that is relevant to us: when we talk about derivatives of functions or write the norm $\lVert x-y\rVert^2$ we are imagining to be doing this directly on the manifold $\mathcal{M}$ with respect to normal coordinates; but it is difficult to perform some of these computations exactly (such as the expansion of $G_\varepsilon f(x)$ given in section \ref{sec:main_results}) in these coordinates since the exponential map is non-linear, but since the projection coordinates are truly flat on $\mathcal{N}$, it is much more preferable to change from normal coordinates to projection coordinates.

To change to coordinates in the integral of $G_\varepsilon f(x)$ given in the proof of Theorem \ref{thm:odk_backward_expansion} in Appendix \ref{appendix:proofs}, we must expand $\lVert y-x\rVert^2$ and $\mathrm{d}y$ in terms of projection coordinates, and we must convert the Taylor expansion of $f(y)$ and $l(y)$ in normal coordinates, to being in terms of projection coordinates.  Doing the former allows us to compute the integrals we need, and the latter combines intrinsic derivatives (calculated using normal coordinates) with projection coordinates, allowing us to derive our Laplacian $\mathcal{L}$.  In the remainder of this subsection, we will overview some other relevant concepts from differential geometry.

In section \ref{subsec:nonisotropic_extension}, we discuss a non-isotropic extension of ODK and consider the metric on $\mathcal{M}$ to be non-trivial; and in section \ref{subsec:numerical_results_ode}, we consider a non-isometric embedding, which also results in having to consider a non-trivial metric.  There will still exist normal coordinates on $\mathcal{M}$ such that $g_{ij}=\delta_{ij}$, but these will not align, as before, with the projection coordinates we implicitly use; therefore, the metric must remain non-trivial.  This adds complexity to the derivatives, where now the differential of a scalar function $df\in \mathcal{T}^*\mathcal{M}$, and the gradient $\text{grad}(f)\backwardscoloneqq\nabla_g f$ must be compatible with the metric through the following relation 
\begin{equation}
    df(v)=g(\nabla_g f,v)\iff\partial_ifv^i=g_{ij}\nabla_g^i fv^j\iff \nabla_g^if=g^{ij}\partial_jf
\end{equation}
where $\partial_i f$ are the standard partial derivatives of $f$ with respect to some coordinates, and $g^{ij}$ is the inverse metric tensor satisfying $g^{ik}g_{kj}=\delta_j^i$.  This introduces the idea of using the metric to raise or lower and index; this is often denoted using musical notation with ``X-flat" being the lowering of the index of contravariant vector $X$, with $X^\flat=g_{ij}X^jdx^i$, and ``$\omega$-sharp" being the raising of the index of a covariant vector $\omega$, with $\omega^\sharp=g^{ij}\omega_j\partial_i$; thus, we can see that $\nabla_g f=df^\sharp$.  The divergence will also change as is again defined implicitly through $d(i_XdV)=(\text{div}_g X)dV$, where $dV=\sqrt{|\det g|}dx_1\wedge\dots\wedge dx_d$ is the volume form and $i_X\omega(V_1,\dots,V_{k-1})=\omega(X,V_1,\dots,V_{k-1})$ is in the interior multiplication operator taking $k$-forms to $(k-1)$-forms, so if $X$ is your vector field then $i_XdV$ will take the parallelepiped defined by the volume for $dV$ break it into a sum of each side of it not including $dx^i$ scaled by $X^i$, so is the sum of scaled area elements; $d(i_XdV)$ asks how all these area elements are changing and glues it back together into another volume form; the divergence is a measure of how much rescaling by $X$ managed to change the volume form and can be found to be equal to $\text{div}_g(X)=|g|^{-1/2}\partial_i(|g|^{1/2}X^i)$.  The Laplace-Beltrami operator is defined to be $\Delta_g f=\text{div}_g(\nabla_g f)$ and thus we can write this in coordinates as:
\begin{equation}
    \Delta_g f=\tfrac{1}{\sqrt{|g|}}\partial_i\left(\sqrt{|g|}g^{ij}\partial_j f\right)
\end{equation}
When our metric is trivial, we will have $g_{ij}=\delta_{ij}$ and $|g|\equiv 1$ and thus $\Delta_g f=\partial_i\partial_j f$, so we arrive back at the regular Euclidean Laplacian, but if we cannot trivialise our metric, then we must account for how the metric affects our derivatives. In general, the non-trivial backwards Laplacian will equal:
\begin{equation}
    \mathcal{L}_gf=\mu(f)+\tfrac{1}{2}\rho^2\Delta_g f=\mu^i\partial_if+\tfrac{\rho^2}{2\sqrt{|g|}}\partial_i\left(\sqrt{|g|}g^{ij}\partial_j f\right)
\end{equation}
where we note that the first term is equivalent to $g(\nabla_g f,\mu)$.  The adjoint operator will take the form:
\begin{equation}
    \mathcal{L}_g^* f=-\text{div}_g(\mu f)+\tfrac{1}{2}\Delta_g(\rho^2 f)=-\tfrac{1}{\sqrt{|g|}}\partial_i(\sqrt{|g|}\mu^if)+\tfrac{1}{2\sqrt{|g|}}\partial_i\left(\sqrt{|g|}g^{ij}\partial_j(\rho^2f)\right)
\end{equation}
Since we still implicitly use the projection coordinates, the Laplacians that we approximate will no longer respect the more complex Laplacians above, as they will not incorporate the change of the metric of the manifold.  Local Kernels introduced a method to approximate, under mild constraints, $\Delta_g$ by constructing a symmetrised kernel $\bar K_\varepsilon(x,y)=K_\varepsilon(x,y)+K_\varepsilon(y,x)$; and we, in section \ref{subsec:nonisotropic_extension}, construct the full Laplacians given above under harder constraints, namely, we ask that $\mu_\text{eff}\coloneqq (\partial_i g^{ij}+\tfrac{1}{2}g^{ij}\partial_i\log |g|)\partial_j$, which is the effective advection term introduced by $\Delta_g$ be such that $\exists a,\varphi\in C^\infty(\mathcal{M})$ satisfying $\mu_\text{eff}=a\nabla\varphi$ which allows us to use extensions of ODK such as mODK, stODK, or the density $\alpha$-regularised ODK, which introduce terms in this form and which can then combine with the non-isotropic ODK to form the above two Laplacian operators.

The last piece of machinery from differential geometry that we discuss is the second fundamental form, which allows us to quantify the curvature of a manifold within the context of an embedding into another, larger manifold, which for us would be $\iota:\mathcal{M}\to\mathbb{R}^n$.  Suppose $\bar \nabla$ is the standard covariant derivative on $\mathbb{R}^n$, then the second fundamental form defined as $II:\mathcal{T}(\mathcal{M})\times\mathcal{T}(\mathcal{M})\to\mathcal{N}(\mathcal{M}):(X,Y)\mapsto (\bar\nabla_XY)^\perp$, i.e. we take the normal part of the covariant derivative of $Y$ in the direction of $X$ in $\mathbb{R}^n$ and is naturally then a measure of curvature.  The trace of the second fundamental form is called the mean curvature vector and is denoted as $H(x)=\sum_{i=1}^dII(x)(\partial_i,\partial_i)$.  This is relavant to us because the mean displacement $\check v_\text{md}(x_i)=\sum_{j=1}^mP_\varepsilon(x_i,x_j)(x_j-x_i)=L_\varepsilon\iota(x_i)$, where $\mathcal{L}\iota=(\mathcal{L}\iota_1,\dots,\mathcal{L}\iota_n)^T$, and thus $\check{v}_\text{md}=L_\varepsilon\iota\approx v+\tfrac{1}{2}\rho^2\Delta\iota$, where $v$ is the underlying velocity and $\tfrac{1}{2}\rho^2$ a diffusion and $\Delta \iota(x)=H(x)$ the mean curvature vector of the embedding (this is proven on Pg. 41 of \cite{chen2011pseudo}).  By the definition of the second fundamental form giving values in the normal bundle of $\iota(\mathcal{M})$, the trace --- which is then a sum of normal vectors --- will also point normal to $\iota(\mathcal{M})$ in $\mathbb{R}^n$ and thus the mean displacement has a normal component that encodes the curvature of the manifold; in section \ref{subsec:numerical_results_ode} we introduce two estimators that remove this normal component, allowing you to reconstruct the underlying dynamics.

\subsection{ODK Parameter Optimisation Scheme}\label{appBsec:odk_optimisation_scheme}
In this appendix, we will describe how to analytically compute the gradients of the transient objective, $L_q$, and the negative log-likelihood objective, $L_p$, both computed as:
\begin{align}
    L_q(\Phi;\lambda_\text{reg})&=-\frac{1}{m}\sum_{i=1}^m\log\hat q(x_i)+\lambda_\text{reg}R(\Phi)+\text{constant}\label{appB:Lq_objective}\\
    L_p(\Phi;\lambda_\text{reg})&=-\frac{1}{m}\sum_{i=1}^m\log \hat{p}(t_i,x_i)++\lambda_\text{reg}R(\Phi)+\text{constant}\label{appB:Lp_objective}
\end{align}
where $m=|\mathcal{D}|$, and $\lambda_\text{reg}$ is a hyperparameter controlling the strength of the regularisation function $R(\Phi)$.  We will not talk much about how one might compute the gradient of the regularisation function since we leave this free for the reader to choose.  However, for the gradient and Laplacian regularisations we introduce in section \ref{subsec:numerical_results_sde}, they are computed as:
\begin{equation}
    R_\text{grad}(f)=\frac{1}{m}\sum_{i=1}^m(f_{N_i}-f_i\mathbb{1}_{k_\text{reg}})^TG_i(f_{N_i}-f_i\mathbb{1}_{k_\text{reg}}),\qquad R_\text{lap}(f)=\frac{1}{m}\lVert L_\text{reg}f\rVert^2
\end{equation}
where: $N_i$ is a precomputed set of indices for the $k_\text{reg}$ nearest neighbours to $x_i$ (not including $x_i$ itself), and $\mathbb{1}_{k_\text{reg}}$ is a vector of $k_\text{reg}$ 1's; and $L_\text{reg}$ is an appropriately normalised and fixed Laplacian matrix (thus $Lf\approx\Delta f$).  We state the gradients of these regularisations as:
\begin{align}
    (\nabla R_\text{grad}(f))_i&=\sum_{j|\exists k:i=(N_j)_k}2\left(G_j(f_{N_j}-f_j\mathbb{1}_{N_j})\right)_k-\sum_{k\in N_i}2\left(G_i(f_{N_i}-f_i)\right)_k\\
    \nabla R_\text{lap}(f)&=2L_\text{reg}f
\end{align}
Considering the q-objective first, we recall that we compute the estimated transient distribution through:
\begin{equation}
    \hat q(x_i) = \sum_{s=1}^S\bar r(s)p_s(x_i)=\sum_{s=1}^S\bar r(s)p_0^TP_\varepsilon^se_i
\end{equation}
where: $S=\left\lceil\frac{T}{\varepsilon}\right\rceil$, $\bar r(s)=r(s\varepsilon)/\sum_{s=1}^Sr(s\varepsilon)$, and $p_s(x_i)=p_0^TP_\varepsilon^se_i$ where $(e_i)_j=\delta_{ij}$.  Plugging this into the derivative of $L_q$ we see that:
\begin{equation}\label{appB:dLq_dphi1}
    -m\frac{\partial L_q}{\partial \phi_k}=\sum_{i=1}^m\frac{\partial}{\partial \phi_k}\log\hat q(x_i)=\sum_{i=1}^m\frac{1}{\hat q(x_i)}\frac{\partial}{\partial \phi_k}\sum_{s=1}^S\bar r(s)p_0^TP_\varepsilon^se_i=\sum_{s=1}^S\bar r(s)p_0^T\frac{\partial P_\varepsilon^s}{\partial\phi_k}\sum_{i=1}^m\frac{e_i}{\hat q(x_i)}=\sum_{s=1}^S\bar r(s)p_0^T\frac{\partial P_\varepsilon^s}{\partial\phi_k}\hat q^{-1}
\end{equation}
where we write $\hat q^{-1}=\sum_{i=1}^m\frac{1}{\hat q(x_i)}e_i$.  The gradient of the matrix $P_\varepsilon$ with respect to the parameter $\phi_k$ is done componentwise, i.e. $\left(\frac{\partial P_\varepsilon}{\partial \phi_k}\right)_{ij}=\frac{\partial}{\partial \phi_k}p_\varepsilon(x_i,x_j)$, as a result we can expand the derivative of the matrix exponent, $P_\varepsilon^s$, as:
\begin{equation}
    \frac{\partial P_\varepsilon^s}{\partial \phi_k}=P_\varepsilon^{s-1}\frac{\partial P_\varepsilon}{\partial \phi_k}+\frac{\partial P_\varepsilon^{s-1}}{\partial \phi_k}P_\varepsilon=P_\varepsilon^{s-1}\frac{\partial P_\varepsilon}{\partial \phi_k}+P_\varepsilon^{s-2}\frac{\partial P_\varepsilon}{\partial \phi_k}P_\varepsilon+\frac{\partial P_\varepsilon^{s-2}}{\partial \phi_k}P_\varepsilon^2=\sum_{t=0}^{s-1}P_\varepsilon^t\frac{\partial P_\varepsilon}{\partial \phi_k}P_\varepsilon^{s-1-t}
\end{equation}
which is in terms of powers of $P_\varepsilon$ and $\frac{\partial P_\varepsilon}{\partial \phi_k}$, only.  We can plug this back into equation \eqref{appB:dLq_dphi1} to get:
\begin{equation}
    -m\frac{\partial L_q}{\partial \phi_k}=\sum_{s=1}^S\bar{r}(s)\sum_{t=0}^{s-1}p_0^TP_\varepsilon^t\frac{\partial P_\varepsilon}{\partial \phi_k}P_\varepsilon^{s-1-t}\hat{q}^{-1}=\sum_{s=1}^S\bar r(s)\sum_{t=0}^{s-1}p_t^T\frac{\partial P_\varepsilon}{\partial \phi_k}\hat w_{s-1-t}
\end{equation}
where $\hat w_{t}\coloneqq P_\varepsilon^t\hat q^{-1}$.  We see from this expression that to compute the gradient analytically we need only compute $p_t$ and $v_t$ for $t=0,1,\dots,S-1$, which can be done in a single forward pass (and indeed the $p_t$'s will already have been calculated to compute the objective function) and then factor these vectors through $\frac{\partial P_\varepsilon}{\partial \phi_k}$.  

To compute the analytic gradient of the $L_p$ objective is slightly more complicated but follows the same structure.  We proceed initially by defining $\hat p(t_i,x_i)=p_0^TP_{\varepsilon}^{s_i}e_i$ where $s_i=\left\lceil\frac{t_i}{\varepsilon}\right\rceil$, then we write the derivative of the objective as:
\begin{equation}
    -m\frac{\partial L_p}{\partial \phi_k}=\sum_{i=1}^m\frac{\partial}{\partial\phi_k}\log \hat p(t_i,x_i)=\sum_{i=1}^m\frac{1}{\hat p(t_i,x_i)}\frac{\partial}{\partial \phi_k}\hat p(t_i,x_i)=\sum_{i=1}^mp_0^T\frac{\partial P_\varepsilon^{s_i}}{\partial\phi_k}\frac{e_i}{p_{s_i}(x_i)}=\sum_{i=1}^m\sum_{t=0}^{s_i-1}p_t^T\frac{\partial P_\varepsilon}{\partial\phi_k}P_\varepsilon^{s_i-1-t}\frac{e_i}{p_{s_i}(x_i)}
\end{equation}
The difficulty here is that we cannot swap the order of the state sum and the time sum as we did before since the limits of the time sum depend on $s_i$.  The trick to get around this problem is to sum from $0$ to $S-1$ regardless but mask away unwanted terms with a characteristic function, i.e.:
\begin{equation}\label{appB:dLp_dphi1}
    -m\frac{\partial L_p}{\partial \phi_k}=\sum_{i=1}^m\sum_{t=0}^{S-1}p_t^T\frac{\partial P_\varepsilon}{\partial\phi_k}P_\varepsilon^{s_i-1-t}\frac{e_i}{p_{s_i}(x_i)}\mathbb{1}(t\leq s_i-1)=\sum_{t=0}^{S-1}p_t^T\frac{\partial P_\varepsilon}{\partial\phi_k}w_t
\end{equation}
where $w_t$ is defined as:
\begin{align}
    w_t\coloneqq\sum_{i=1}^mP_\varepsilon^{s_i-1-t}\frac{e_i}{p_{s_i}(x_i)}\mathbb{1}(t\leq s_i-1)=\sum_{s=1}^{S}\sum_{i:s_i=s}P_\varepsilon^{s-1-t}\frac{e_i}{p_{s_i}(x_i)}\mathbb{1}(t\leq s-1)=\sum_{s=t+1}^{S}\sum_{i:s_i=s}P_\varepsilon^{s-1-t}\frac{e_i}{p_{s_i}(x_i)}
\end{align}
and satisfies the following recursive relationship:
\begin{equation}
    w_t=P_\varepsilon w_{t+1}+\sum_{i:s_i-1=t}\frac{e_i}{p_{s_i}(x_i)}
\end{equation}
From this recursion relation we see that if we take the sum given in equation \eqref{appB:dLp_dphi1} backwards (from $S-1$ to $0$) then we can compute $\frac{\partial L_p}{\partial \phi_k}$ through a single forward pass (as with the $L_q$ objective).  

We can generalise the $L_p$ objective to work with a linearly interpolated probability, $\hat p(t_i,x_i)=a_ip_{s_i-1}(x_i)+b_ip_{s_i}(x_i)$ where $a_i=1-\left(s_i-\frac{t_i}{\varepsilon}\right)$ and $b_i=s_i-\frac{t_i}{\varepsilon}$.  From equation \eqref{appB:dLp_dphi1} that:
\begin{align}
    -m\frac{\partial L_p}{\partial \phi_k}=\sum_{i=1}^m\frac{a_i}{\hat p(t_i,x_i)}\frac{\partial}{\partial \phi_k}p_{s_i-1}(x_i)+\frac{b_i}{\hat p(t_i,x_i)}\frac{\partial}{\partial \phi_k}p_{s_i}(x_i)=\sum_{i=1}^m&\sum_{t=0}^{s_i-2}p_t^T\frac{\partial P_\varepsilon}{\partial\phi_k}P_\varepsilon^{s_i-2-t}\frac{a_i}{\hat p(t_i,x_i)}e_i\\
    &+\sum_{t=0}^{s_i-1}p_t^T\frac{\partial P_\varepsilon}{\partial\phi_k}P_\varepsilon^{s_i-1-t}\frac{b_i}{\hat p(t_i,x_i)}e_i
\end{align}
using the same trick as before we define:
\begin{align}
    w_t^a&=\sum_{i=1}^mP_\varepsilon^{s_i-2-t}\frac{a_i}{\hat p(t_i,x_i)}e_i\mathbb{1}(t\leq s_i-2)=\sum_{s=t+2}^S\sum_{i:s_i=s}P_\varepsilon^{s-2-t}\frac{a_i}{\hat p(t_i,x_i)}e_i\\
    w_t^b&=\sum_{i=1}^mP_\varepsilon^{s_i-1-t}\frac{b_i}{\hat p(t_i,x_i)}e_i\mathbb{1}(t\leq s_i-1)=\sum_{s=t+1}^S\sum_{i:s_i=s}P_\varepsilon^{s-1-t}\frac{b_i}{\hat p(t_i,x_i)}e_i
\end{align}
which again satisfy the recursion relations:
\begin{align}
    w_t^a&=P_\varepsilon w_{t+1}^a+\sum_{i:s_i-2=t}\frac{a_i}{\hat p(t_i,x_i)}e_i\\
    w_t^b&=P_\varepsilon w_{t+1}^b+\sum_{i:s_i-1=t}\frac{b_i}{\hat p(t_i,x_i)}e_i
\end{align}
Letting $w_t^{a+b}=w_t^{a}+w_t^b$ we have:
\begin{equation}
    w_t^{a+b}=P_\varepsilon w_{t+1}^{a+b}+\sum_{i:s_i-2=t}\frac{a_i}{\hat p(t_i,x_i)}e_i+\sum_{i:s_i-1=t}\frac{b_i}{\hat p(t_i,x_i)}e_i
\end{equation}
Putting all this together we find that:
\begin{equation}
    -m\frac{\partial L_p}{\partial \phi_k}=\sum_{t=0}^{S-2}p_t^T\frac{\partial P_\varepsilon}{\partial \phi_k}w_t^a+\sum_{t=0}^{S-1}p_t^T\frac{\partial P_\varepsilon}{\partial \phi_k}w_t^b=p_{S-1}^T\frac{\partial P_\varepsilon}{\partial \phi_k}w_{S-1}^b+\sum_{t=0}^{S-2}p_t^T\frac{\partial P_\varepsilon}{\partial \phi_k}w_t^{a+b}
\end{equation}
So just as before we can compute the gradient in a single forward pass.

All of the above requires the analytic computation of $\frac{\partial P_\varepsilon}{\partial \phi_k}$, and this depends on the ODK $K_\varepsilon$ and the form of the normalisation of $P_\varepsilon$ in terms of $K_\varepsilon$.  The elements of the normalised transition operator $P_\varepsilon$ take the form:
\begin{equation}
    P_{ij}\coloneqq (P_\varepsilon)_{ij}=\frac{K(x_i,x_j)}{Z(x_i)q(x_j)},\qquad\text{where:}\qquad Z(x_i)=\sum_{j=1}^m\frac{K(x_i,x_j)}{q(x_j)}
\end{equation}
where we drop the $\varepsilon$ subscript of $P_\varepsilon$ for notational clarity.  As a result, the $(i,j)^\text{th}$ component of the derivative of $P_\varepsilon$ takes the form:
\begin{align}
    \frac{\partial P_{ij}}{\partial \phi_k}&=\frac{1}{Z_iq_j}\frac{\partial K_{ij}}{\partial \phi_k}-\frac{K_{ij}}{Z_i^2q_j}\sum_{n=1}^m\frac{1}{q_n}\frac{\partial K_{in}}{\partial \phi_k}=\frac{P_{ij}}{K_{ij}}\frac{\partial K_{ij}}{\partial \phi_k}-P_{ij}\sum_{n=1}^m\frac{P_{in}}{K_{in}}\frac{\partial K_{in}}{\partial \phi_k}\\
    &=P_{ij}\left[\frac{\partial \log (K_{ij})}{\partial \phi_k}-\sum_{n=1}^mP_{in}\frac{\partial \log(K_{in})}{\partial \phi_k}\right]\label{appB:dP_dphi}
\end{align}
which is in terms of the normalised transition operator $P_\varepsilon$ and the natural logarithm of the raw ODK $K_\varepsilon$.  The presence of the natural logarithm points us towards using the prototypical kernel, i.e. where $h(u)=\exp(-u/2)$, therefore we consider the ODK:
\begin{equation}\label{appB:odk}
    K_{ij}\coloneqq K_\varepsilon(x_i,x_j)=\exp\left(-\frac{D_{ij}^2}{2\varepsilon\rho_i^2}-\frac{(l_j-l_i-\varepsilon\kappa_i\tau_i)^2}{2\varepsilon\kappa_i^2\rho_i^2}+\frac{(l_j-l_i)^2}{2\varepsilon\kappa_i^2\rho_i^2}\right)
\end{equation}
where we use subscripts, $f_i$, to denote evaluation of a function $f(x)$ at the sample indexed by $i$ in $\mathcal{D}$, i.e. $f_i=f(x_i)$, and where $D_{ij}=\lVert x_i-x_j\rVert$.  Since the ODK, $K_{ij}$, given in equation \eqref{appB:odk}, has a rigid structure in terms of the input functions $\tau(x),\ \rho(x)$ and $l(x)$ we can compute the derivative of $\log(K_{ij})$ with respect to any $\phi_k$ (where $\phi=\tau,\rho,l$) analytically.  We first write $\log(K_{ij})$ as:
\begin{equation}
    \log (K_{ij})=-\frac{D_{ij}^2}{2\varepsilon\rho_i^2}+\frac{\tau_i}{\kappa_i\rho_i^2}(l_j-l_i)-\frac{\varepsilon\tau_i^2}{2\rho_i^2}
\end{equation}
then we can read of the derivatives of $\log(K_{ij})$ with respect to $\tau_k,\ \rho_k$ and $l_k$, before plugging the result back into equation \eqref{appB:dP_dphi} to compute the respective derivatives with respect to $P_{ij}$:
\begin{align}
    \frac{\partial \log(K_{ij})}{\partial \tau_k}&=\frac{(l_j-l_i-\varepsilon\kappa_i\tau_i)}{\kappa_i\rho_i^2}\delta_{ki}&\implies&\qquad\qquad\frac{\partial P_{ij}}{\partial\tau_k}=\frac{P_{ij}}{\kappa_i\rho_i^2}\left(l_j-\sum_{n=1}^mP_{in}l_n\right)\delta_{ki}\label{appB:dK_dP_dtau}\\
    \frac{\partial \log(K_{ij})}{\partial \rho_k}&=-\frac{2}{\rho_i}\log(K_{ij})\delta_{ki}&\implies&\qquad\qquad\frac{\partial P_{ij}}{\partial\rho_k}=\frac{2P_{ij}}{\rho_i}\left(\sum_{n=1}^mP_{in}\log(K_{in})-\log(K_{ij})\right)\delta_{ki}\label{appB:dK_dP_drho}\\
    \frac{\partial \log(K_{ij})}{\partial l_k}&=\frac{\tau_i}{\kappa_i\rho_i^2}(\delta_{kj}-\delta_{ki})&\implies&\qquad\qquad\frac{\partial P_{ij}}{\partial l_k}=\frac{\tau_i}{\kappa_i\rho_i^2}P_{ij}\left(\delta_{kj}-P_{ik}\right)\label{appB:dK_dP_dl}
\end{align}
When $\phi=\tau,\rho$ we notice that the only row of $\tfrac{\partial P_\varepsilon}{\partial\phi_k}$ that is not identically zero will be the $k^\text{th}$ row, and thus we can define:
\begin{equation}
    \frac{d P_\varepsilon}{d\phi}=\sum_{k=1}^m\frac{\partial P_\varepsilon}{\partial \phi_k}
\end{equation}
as a method to store all derivative information in a single matrix.  We could also re-write the full derivative of the $L_q$ and $L_p$ objectives as:
\begin{align}
    -m\frac{dL_q}{d\phi}&=\sum_{s=1}^S\bar r(s)\sum_{t=0}^{s-1}p_t^T\text{Diag}\left(\frac{dP_\varepsilon}{d\phi}\hat\omega_{s-1-t}\right)\\
    -m\frac{dL_p}{d\phi}&=\sum_{t=0}^{S-1}p_t^T\text{Diag}\left(\frac{dP_\varepsilon}{d\phi}\omega_{t}\right)
\end{align}
But in practice, we do not compute the gradient using the diagonal matrices, but we do store all gradient information in a single matrix.  This wouldn't work when $\phi=l$, but if you were training for the ordering you would naturally set $\tau(x)$ and $\kappa(x)$ to be constants, so every iteration you would just have to invert $\rho(x)^2$ to compute gradients with respect to $l$, since we notice in equation \eqref{appB:dK_dP_dl} that we do not have any costly sums over the rows of $P_\varepsilon$ as we do in equations \eqref{appB:dK_dP_dtau} \& \eqref{appB:dK_dP_drho}, and thus require us to precompute and store the derivatives at the beginning of our calculations.

\subsubsection{Laplacian Eigenfunction Decomposition}\label{appB:LE_decomp}

In an attempt to enforce smoothness without adding constraints to the differential structure of dynamic parameters in our ODK as is done when using gradient or Laplacian regularisation given in equations \eqref{eq:sde_gradreg} \& \eqref{eq:sde_lapreg}, respectively.  We can compute a basis of functions that will naturally be smooth.  In the spirit of this paper, and with reference to the work of Coifman \& Lafon in Diffusion Maps \cite{coifman2006diffusion}, we will extract the top $K$ eigenfunctions of the Laplace-Beltrami operator on $\mathcal{M}$ by using a symmetric (diffusion maps) kernel with the form:
\begin{equation}\label{appBeq:appB_Psym}
    P_\varepsilon^\text{sym}(x_i,x_j)=\frac{K_\varepsilon^\text{sym}(x_i,x_j)q(x_j)^{-1}}{\sum_{k=1}^mK_\varepsilon(x_i,x_k)q(x_k)^{-1}},\qquad\text{with:}\qquad K_\varepsilon^\text{sym}(x_i,x_j)=h\left(\frac{\lVert x_i-x_j\rVert^2}{\varepsilon}\right)
\end{equation}
Which is an ODK with no ordering and $\rho\equiv 1$, and by Theorem \ref{thm:odk_backward_expansion} (or equivalently, by Theorem 2, Pg. 15, of \cite{coifman2006diffusion}), we know that:
\begin{equation}\label{appBeq:appB_Lsym}
    \lim_{\varepsilon\to 0}\varepsilon^{-1}(P_\varepsilon^\text{sym}-I)f(x)=\frac{1}{2}\sigma_h^2\Delta f
\end{equation}
In the examples given in section \ref{subsec:numerical_results_sde}, we take $h(u)=\exp(u/2)$, which makes the right-hand side of equation \eqref{appBeq:appB_Lsym}, equal $\tfrac{1}{2}$.  We can write $P_\varepsilon^\text{sym}$ as $P_\varepsilon^\text{sym}=D_L^{-1}D_R^{-1}K_\varepsilon^\text{sym}D_R^{-1}$, where $K_\varepsilon^\text{sym}$ is the matrix of pairwise unnormalised kernel evaluations, $D_R=\diag(q)$ for the sampling density $q(x)$, and $D_L=\diag(D_R^{-1}K_\varepsilon^\text{sym}D_R^{-1}\mathbb{1})$.  We note that introducing the first $D_R$ on the left, diverges from the usual form of the symmetric transition matrix, but since this simply scales each row by $q_i^{-1}$, and $(D_L)_{ii}$ equals the sum of the $i^\text{th}$ row, this additional term immediately cancels and we are left with the standard transition matrix given in equation \eqref{appBeq:appB_Psym}.  We now define the symmetrised transition matrix as:
\begin{equation}
    \bar P_\varepsilon^\text{sym}\coloneqq D_L^{1/2}P_\varepsilon^\text{sym} D_L^{-1/2}=D_L^{-1/2}D_R^{-1}K_\varepsilon^\text{sym}D_R^{-1}D_L^{-1/2}
\end{equation}
since the kernel matrix $K_\varepsilon^\text{sym}$ is, by construction symmetric, so too will $\bar P_\varepsilon^\text{sym}$.  We now note that if $\varphi_k$ is an eigenvector of $\bar P_\varepsilon^\text{sym}$ with eigenvalue $\lambda_k$, i.e. $P_\varepsilon^\text{sym}\varphi_k=\lambda_k\varphi_k$, then $v_k\coloneqq D_L^{1/2}\varphi_k$ satisfies:
\begin{equation}
    \bar P_\varepsilon^\text{sym}v_k=D_L^{1/2}P_\varepsilon^\text{sym} D_L^{-1/2}v_k=D_L^{1/2}P_\varepsilon^\text{sym}\varphi_k=\lambda_kD_L^{1/2}\varphi_k=\lambda_kv_k
\end{equation}
the reverse is also true, that being that if $v_k$ is an eigenvector of $\bar P_\varepsilon^\text{sym}$ with eigenvalue $\lambda_k$ then $\varphi_k=D_L^{-1/2}v_k$ is an eigenvector of $P_\varepsilon^\text{sym}$ with the same eigenvalue.  It is of practical importance to formulate this in terms of the symmetrised transition matrix since it is far more computationally stable to compute the eigendecomposition of a symmetric matrix than an asymmetric matrix.  Thus, after computing an eigendecomposition of the $\bar P_\varepsilon^\text{sym}$, we can take the top $K$ eigenvectors, corresponding the $K$ largest eigenvalues of $\bar P_\varepsilon^\text{sym}$, and then form the basis $\mathcal{B}=\{\varphi_1,\dots,\varphi_K\}$.  All the eigenvalues of $\bar P_\varepsilon^\text{sym}$ are in the range $(0,1]$ (as was shown, for example, in \cite{chung1997spectral}), thus for a given value of $\varepsilon$, all the eigenvalues of $L_\varepsilon^\text{sym}=\varepsilon^{-1}(P_\varepsilon^\text{sym}-I)$ lie in the range $(-1/\varepsilon,0]\xrightarrow{\varepsilon\to 0}(-\infty,0]$.  Since $L_\varepsilon^\text{sym}\xrightarrow{\varepsilon\to 0}\Delta$ we know that we can interpret the eigenvectors of $L_\varepsilon^\text{sym}$ as being harmonics, with their corresponding eigenvalues being their frequency (where frequency controls variability over $\mathcal{M}$).  Thus, by taking the top $K$ eigenfunctions, we approximate a basis of functions that is both smooth and not overly variable (rapid variability is effectively the definition of not being smooth when one is restricted to discrete data).

It is important to note that although the eigenvectors of $\bar P_\varepsilon^\text{sym}$ are orthogonal (by the spectral theory of symmetric matrices), the transformed eigenvectors of $P_\varepsilon^\text{sym}$ instead satisfy:
\begin{equation}
    \langle\varphi_i,\varphi_j\rangle_{L^2(\mathcal{M},q)}\approx\varphi_iD_L\varphi_j=v_i^TD_L^{-1/2}D_LD_L^{-1/2}v_j=\langle v_i,v_j\rangle_2=\delta_{ij}
\end{equation}
So, unless the data is uniformly distributed, we are not guaranteed that $\varphi_k$'s will be orthogonal.  Thus, for practical reasons, we often orthogonalise $\mathcal{B}$ using the classical Gram-Schmidt algorithm \cite{leon2013gram}.

Once we have a basis of eigenfunctions, $\mathcal{B}$, we can parametrise $l(x)$, $\tau(x)$ or $\rho(x)$ as:
\begin{equation}
    \tau(x)=\sum_{k=1}^K\alpha_k^\tau\varphi_k(x),\qquad l(x)=\sum_{k=2}^K\alpha_k^l\varphi_k(x),\qquad\rho(x)=\exp\left(\sum_{k=1}^K\alpha_k^\rho\varphi_k(x)\right).
\end{equation}
where $\tau(x)$ is simply a weighted sum of all $K$ eigenfunctions, $l(x)$ is a weighted sum of all but the trivial (constant) eigenfunction (thus building into the parametrisation our knowledge that ODKs are invariant to addition by a constant), and $\rho(x)$ is the exponential of the weighted sum of all $K$ eigenfunctions.  The purpose of log-transforming the $\rho$'s is to ensure the noise is bounded away from zero; other functions could also be used to enforce positivity, such as $x^2+\rho_\text{min}$, but while testing we found the $\rho=\exp(\cdot)$ to be the most reliable.  We could also consider taking the modulus, ReLU, or exponential of our decomposition for $\tau(x)$ since if we have fixed the ordering we might want to stop the process from going backwards.

Other than acting as an almost parameter-free smoothness regularisation, taking a Laplacian eigenbasis also massively reduces the number of parameters that we have to optimise for, thus improving the stability and rate of convergence of both our $L_q$ and $L_p$ optimisation problems.  In figure \ref{fig:nl_results}, Panels B and E show the results of applying this parameterisation framework to data sampled from the NL process of section \ref{subsec:numerical_results_sde}, we can see that the estimate is still relatively noisy indicating our choice of $K=250$ was too high, but we also see that upon smoothing we have very nearly perfectly captured the true noise function.  This is noteworthy because the gradient- and Laplacian-regularised optimisation didn't capture this structure as faithfully, showing that the signal was weak enough that even a small amount of additional regularisation was enough to distort the minima; thus, using a Laplacian eigendecomposition might allow us to more faithfully reconstruct the weak signal from the diffusion in the data.

We note that, as discussed in section \ref{subsec:selftuning_odk}, we could replace the symmetric diffusion maps kernel with a symmetric variable bandwidth kernel:
\begin{equation}
    K_\varepsilon^S(x_i,x_j)=h\left(\frac{\lVert x_i-x_j\rVert^2}{\varepsilon\rho(x)\rho(y)}\right)
\end{equation}
which would still allow for us to stably take an eigendecomposition and adapts better to inhomogeneously sampled data, and indeed converges to (a drifted version of) the Laplacian at a faster rate, see \cite{berry2016variable} for more details.

\subsection{Identifiability of Drift and Diffusion Coefficients}\label{appBsec:identifiability}

Recovering the drift and diffusion from data sampled from an SDE is a classical inverse problem and is central to many applied domains.  However, it is not obvious whether the problems we pose for our optimisers to solve are well-posed; the theory of \textit{identifiability} for these problems seeks to answer the question of whether, given some idealised form of our data, it would be theoretically possible to infer the exact parameters of our underlying process; more precisely, let $\mathcal{P}=\{P_\theta:\theta\in\Theta\}$ be the set of all possible models we could learn given the parameter space $\Theta$, then --- under the assumption that our model is expressive enough to represent the observed realities --- if the mapping from $\theta\mapsto P_\theta$ is injective (i.e. $P_{\theta_1}=P_{\theta_2}\implies \theta_1=\theta_2$ for all $\theta_1,\theta_2\in\Theta$) then our model is said to be identifiable (see \cite{lehmann1998theory}, Chapter 1, Definition 5.2).  We are to imagine this as follows: under the assumption that, given some data $\mathcal{D}$, our model can learn the optimal parameterisation, then if our problem is identifiable, our learned model will be a faithful representation of the underlying system.

For the inference of parameters for SDEs, the difficulty of this identifiability problem depends on the data modalities available.  For example, if sample paths (samples of whole trajectories) are observed, then the quadratic variation along these trajectories will determine the diffusion coefficients, and the Kramers-Moyal expansion will determine the drift.  However, if we are not able to sample trajectories, and can instead only sample independent snapshots of the process --- which is the canonical setting of scRNA-seq data, where in order to measure the gene expression of a cell you must destroy it --- then the problem of parameter identifiability is tremendously important as fundamentally our problem may be ill-posed and we would not wish to expend energy on a futile problem.  It is datasets of this sort that we wished to study in this paper --- principally, we aimed to design a kernel which was theoretically well-posed to represent the underlying dynamics of data generated by a stochastic process, but could be easily tuned and reparametrised by users.

From the perspective of identifiability, we have a fundamentally under-determined problem; writting the Fokker-Planck equation in divergence form $\partial_t p_t=-\nabla \cdot J_t$, where $J_t=\mu p_t-\tfrac{1}{2}\nabla (\sigma^2p_t)$ is the probability current, we can see that observations of the marginals will only determine $J_t$ up to the addition of a divergence-free field, and hence drift and diffusion may be exchanged without altering $\{p_t\}_{t\in[0,T]}$.  To make progress with such problems, we need to make some structural restrictions; for example, restricting the drift to be generated by a potential $b=-\nabla U$ allows us to avoid the divergence-free gauge freedom since non-linear gradients cannot generate divergence-free terms.  

The theory of identifiability for gradient-flow SDEs with constant isotropic diffusions is now relatively complete, with parameters being identifiable, under various sets of constraints, from the observation of as few as three time marginals.  If one has a known, or highly constrained, diffusion, then it has been (see, for example, \cite{hashimoto2016learning}) that the potential is recoverable from time marginals; Lavenant et al. \cite{lavenant2024toward} established that the full law on paths is identified through a minimum-entropy variational principle, with a single sample per time point sufficing when the measurement times are dense in $[0,T]$; Chizat et al. \cite{chizat2022trajectory} designed a convergent mean-field Langevin algorithm which was later extended by Neklyudov et al. \cite{neklyudov2023action} to time-inhomogeneous known diffusivity.  More relevant to our own work is the result of Guan et al. \cite{guan2026gradientflow} who showed that a gradient-flow SDE with unknown constant isotropic diffusion is jointly identifiable if and only if your observations occur out of steady state, and furthermore, one only needs to observe three time marginals to guarantee uniqueness.

Our $L_q$ objective, given in equation \eqref{eq:sde_q_objective}, would likely only be identifiable in this restricted domain of gradient flows and constant isotropic diffusions, but we note that our assumptions are significantly weaker than knowing all time marginals.  We assume that we know the initial distribution $p_0(x)$, the time-sampling distribution $r(t)$, and the transient sampling distribution $q(x)=\int_0^Tr(t)\int_{X_0}p_t(x|x_0)p_0(x_0)\mathrm{d}x_0\mathrm{d}t$.  To show identifiability, we must show that given two processes have the same initial condition, they can only produce the same transient distribution if they have identical time-marginals for all $t\in[0,T]$, which then by the result of Guan et al., we know that the processes must be identical.  This might not be as infeasible as it first seems, since the gradient dynamics precludes recurrence and thus if our model is to add probability to the estimated transient distribution at a given sample $x$, it only really has one chance to do so --- i.e. it must add all the required probability exactly when the probability is flowing through this point.  Furthermore, every path evolved up to some $t\in [0,T]$ will have quadratic variation equal to $t\sigma^2$, so if this diffusion isn't properly tuned, the model's estimated variability will be incongruous with the observed data.  Thus, we conjecture that the spread of our observed data determines the magnitude of the diffusion, and the geometry determines the drift; by constraining ourselves to gradient-flow SDEs with constant isotropic diffusion we believe, as with existing results, these determining factors separate cleanly, giving us some form of identifiability.

However, if we wish to relax the condition of constant isotropic noise, the problem becomes significantly more challenging, as the signals from the drift and diffusion can begin to mix together and one's model must not only be able to identify the correct regime, but then at the finer scale, must untangle the balance between the drift and diffusion.  As we have mentioned several times in the main text, when the diffusion is state-dependent, you effectively introduce an effective advection term originating from the diffusion; considering the Fokker-Planck equation $\partial_tp_t=\mathcal{L}^*p_t$, we have:
\begin{equation}
    \partial_t p_t=-\nabla\cdot\left[\left(\mu-\tfrac{1}{2}\nabla\sigma^2\right)p_t\right]+\tfrac{1}{2}\nabla\cdot(\sigma^2\nabla p_t)
\end{equation}
Hence, the model must be able to isolate by which mechanism the probability is pushed forward.  A real experimental demonstration of this was given by Hoffmann et al. \cite{hoffmann2008noise} who showed that one could model the gene expression of a differentiating stem cell --- a process usually modeled with ODEs --- using an entirely diffusive model; thus, if your system is likely generated by a process with state-dependent noise, then understanding when you can, or cannot, reasonably identify ones parameters is critical.  We can also see how state-dependent noise produces issues with identifiability by assuming that $\mathcal{L}_1^*$ and $\mathcal{L}_2^*$ are two distinct forward operators that give rise to identical time-marginals, then it must be true that $(\mathcal{L}_1^*-\mathcal{L}_2^*)p_t\equiv 0$; when diffusion is constant Guan et al. exploit the fact that the diffusions factor outside of the Laplacian as $\tfrac{1}{2}(\sigma_1^2-\sigma_2^2)\Delta p_t$, and we may assume that $\sigma_1>\sigma_2$ as otherwise previously mentioned identifiability results for known diffusions would guarantee the potentials be equal \cite{lavenant2024toward}; now since the residual diffusion is non-degenerate and the residual drift $-\nabla(\Psi_1-\Psi_2)$ is Lipschitz and satisfies linear growth, thus we can instantiate a uniqueness result from PDE theory (namely, Theorem 4.1.6, of \cite{bogachev2022fokker}) to prove that $p_t=p_0$ for all $t\in [0,T]$ which contradicts the system not being in steady state.  However, if we cannot make this factorisation, we have one PDE, namely $0=(\mathcal{L}_1^*-\mathcal{L}_2^*)p_t$, in two unknowns, giving us an underdetermined PDE which precludes the application of existing uniqueness results for PDEs.  It was shown by Cristofol et al. \cite{cristofol2017simultaneous} that a general univariate SDE is identifiable if one observes the time evolution of the first- and second-order moments over some time interval; this result is intuitive given our previous discussion since we are now trying to derive uniqueness from two PDEs in terms of the forward operators, $(\mathcal{L}_1-\mathcal{L}_2)x_t=0$ and $(\mathcal{L}_1-\mathcal{L}_2)x_t^2=0$.

Our $L_p$ objective is posed towards this significantly more challenging problem where one has state-dependent noise; furthermore, we still assume that we have relatively little dynamical information, namely, it is assumed that we know the initial time-marginal $p_0(x)$, and we observe independent samples of trajectories and when, along these trajectories, they were sampled.  As a result, we have effectively observed the transient distribution $q(x)$, the time sampling distribution, $r(t)$, and a very noisy representation of the joint distribution $p(x,t)$.  The results we obtained from applying the $L_p$ objective to data sampled from the highly complex NL process showed empirically that we can faithfully reconstruct the first-order deterministic dynamics while simultaneously approximating the diffusion.  However, our estimates of the diffusion were highly noisy, and the functional changes observed by introducing a very small smoothness regularisation indicate that the signal from the diffusion, given this objective, is very weak, but nonetheless it is present.  In Panel F of figure \ref{fig:nl_results}, we show the ground truth and estimated noise functions plotted against the ground truth ordering function (which gives us a measure of how far along the process we are), we observe that our Laplacian Eigenbasis estimate, although noisy, very faithfully captures the noise (after some ad hoc smoothing), even features that would seem unimportant given the form of our objective, i.e. the small but sharp increase in diffusion outside of $\lVert x\rVert>2\pi$.  However, this estimate was highly noisy and required ad hoc smoothing to appreciate that it really had captured the signal from the diffusion; this indicates that the loss landscape is very flat around the minima, which we interpret as evidence of some weak identifiability.  The objective function is isolating the noise by using ODKs to numerically integrate the Fokker-Planck equation; given that we observe sampling times, if the path integral of the diffusion has been too small, then the observed uncertainty in the sampling times would be misaligned with the numerical integration.  Thus, we conjecture that the $L_p$ objective, computed by numerically integrating ODKs, is capable of untangling deterministic and stochastic contributions to the dynamics by adaptively matching the mean-field dynamics and the uncertainty in sampling times.  It is not clear if this is sufficient in general, but it may provide some evidence that, potentially with additional constraints, that gradient-flow SDEs with state-dependent isotropic diffusions are identifiable from observations of the joint distribution $p(x,t)$ and the initial distribution $p_0(x)$, under the assumption that the process is not in steady-state --- it may be required that we assume that the time-sampling distribution has a specific form, such as uniform or exponentially distribution; it also may be necessary to assume that one knows the steady-state distribution and this is reached at time $T$, i.e. we also may need to know $p_T(x)$.  The assumption that one is not in steady-state is essential, since otherwise the SDE $dx_t=-\nabla\Psi(x_t)\mathrm{d}t+\sigma\mathrm{d}w_t$, and the rescaled SDE $dy_t=-\alpha\nabla\Psi(x_t)\mathrm{d}t+\sqrt{\alpha}\sigma\mathrm{d}w_t$, would have identical time-marginals, and thus one has a counterexample for identifiability.

\end{document}